
\input miniltx
\input graphicx 


\includegraphics[height=4.4in,width=5.5in,angle=0]{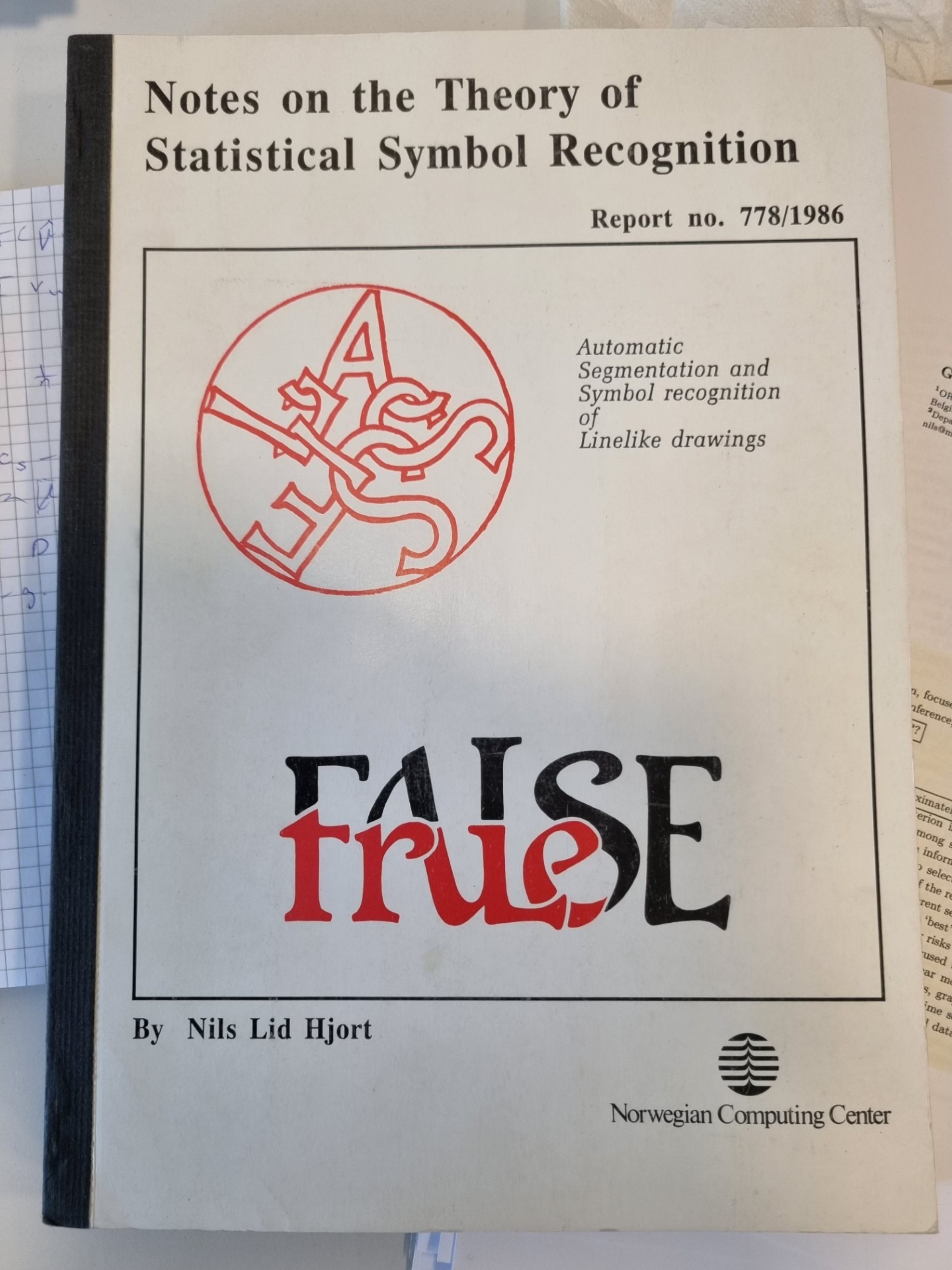}

\bigskip\bigskip\noindent 
{{\baselineskip13pt
This document is a pdf generated from old plain-TeX files of 1986,
of Nils Lid Hjort's {\it Notes on the Theory of Statistical
Symbol Recogntion}, a limited circulation 207-pages monograph
published at the Norwegian Computing Centre, as Report no.~778 / 1986.
It gives the basics of the statistical pattern recognition
theory developed to suit the needs of several industrial projects,
related to contracts with the Norwegian-German firm {\it SysScan},
the Royal Norwegian Research Council, and yet others, 
involving symbol recognition and classification analysis
from noisy images, related to maps, documents, satellite imaging, etc.  

\medskip\noindent
The methods and algorithms developed also needed to fit the
technology of that time, anno c.~1986, with machines scanning
documents, converting these to vector representation,
within computational and machine system boundaries. 
There is an accompanying and also limited circulation
booklet, {\it Statistical Symbol Recognition: Development of a System},
by Knut Br\aa ten, Erik Holb\ae k-Hanssen, and Torfinn Taxt
(Report No.~777 / 1986, Norwegian Computing Centre, Oslo),
detailing the system developments. Thus developments
took form and shape on two frontiers, in close collaboration,
Hjort's statistical methods and getting the technology
to work, with its multiple components, hardware and software.  

\medskip\noindent
About the cover: Helga Lid Hjort made the ASSEL figure
(Working out methods for recognising rotated symbols
has been one of the tasks of the project). The project
name ASSEL (Automatic Segmentation and Symbol
recognition of vectorised Linelike drawings) is the German
word for centipede. The true/false figure is taken
from Scott Kim's {\it Inversions},
BYTE Books, McGraw-Hill, 1981. \smallskip}}




\let\lt=<    \let\gt=>       \let\sla=/
\catcode`\/=0 
\catcode`\<=1 
\catcode`\>=2 
/catcode`/{=13 /let{=/ae
/catcode`/[=13 /let[=/AE
/catcode`/}=13 /let}=/aa
/catcode`/]=13 /let]=/AA
/catcode`/\=13 /let\=/O
/catcode`/|=13 /let|=/o
/escapechar`//

/def/nils#1<>

%
%
/font/bbf=cmbx10 scaled/magstep1
/font/Bbf=cmbx10 scaled/magstep3
/font/ninerm=cmr9
/font/ninesl=cmsl9
/font/ninebf=cmbx9
%
%
/mag/magstephalf
/input texini
/vsize22truecm
/hsize15truecm

/baselineskip=14pt

%
%

/newbox/xxbox
/def/nl</hfill/break>
/def/header#1</vfill/supereject/vbox to 2in<>%
/noindent</Bbf #1>/bigskip/noindent/ignorespaces>
/def/chapter#1#2</message<(Chapter #1:#2)>%
/vfill/supereject/vbox to 2in<>%
/noindent</Bbf Chapter #1/medskip/noindent #2>%
/bigskip/noindent/ignorespaces>
/newbox/seca /newbox/secb /newbox/threebox
/def/section#1#2#3</message<(#1)>/goodbreak/bigskip/nobreak/noindent%
/setbox/threebox/hbox<#3>
/setbox/xxbox/hbox</bbf #1/quad>%
/setbox/seca/hbox</bbf #1/quad #2>%
/setbox/secb/hbox</hbox to/wd/xxbox<>/bbf #3>%
/begingroup%
/bbf/hbox to/hsize</hss/vbox<%
/box/seca%
/ifdim/wd/threebox=0pt /else/box/secb/fi>/hss>%
/nobreak/smallskip/nobreak/ignorespaces%
/endgroup>
/def/subsection#1#2</message<(#1)>/goodbreak/bigskip/nobreak%
</setbox/xxbox/hbox</bf #1/quad>%
/par/parindent=0pt /hangindent=/wd/xxbox /hangafter=1%
/bf#1/quad #2/par>/nobreak/smallskip/nobreak/noindent/ignorespaces>

/def/Btheorem</medbreak/noindent</csc Theorem:>/quad/it/ >
/def/Etheorem</rm/hfil/medskip>

/def/Blemma</medskip/noindent</csc Lemma:>/quad/it>
/def/Elemma</rm/hfil/medskip>

/def/Bproof</medskip/noindent</csc Proof:>/quad/ >
/def/Eproof</qed/hfil/medskip>

/def/note#1</medskip</csc #1>>
/let/Enote=/medskip

/def/subchapt#1#2</section<#1><#2>>


/def/dl</delta>
/def/lam</lambda>
/def/Lam</Lambda>
/def/eps</varepsilon>
/def/ssg</sigma>
/def/sg</Sigma>
/def/th</theta>
/def/Th</Theta>
/def/om</Omega>
/def/som</omega>
/def/ph</phi>


/def/est<</rm est>>

/def/rx<</rm x>> /def/ry<</rm y>> /def/rz<</rm z>>
/def/rX<</rm X>> /def/rY<</rm Y>> /def/rZ<</rm Z>> 

/def/doubt<</rm DOUBT>> /def/out<</rm OUT>>

/def/pmc<</rm pmc>> /def/pcc<</rm pcc>> /def/pd<</rm pd>>
/def/apmc<</rm </overline p>mc>> /def/hpmc<</rm </wh p>mc>>
/def/op<</rm op>>

/def/cross<</rm (CROSS)>>
/def/boot<</rm (BOOT)>>
/def/jack<</rm (JACK)>>

/def/Pr<</it Pr>> /def/el<</rm EL>> /def/Tr<</rm Tr>>

/def/lr<</rm LR>> /def/ise<</rm ISE>> /def/mise<</rm MISE>>
/def/At<A_t> /def/Bt<B_t> /def/Ct<C_t> /def/Dt<D_t>

/def/exp</,</rm exp>/,> /def/log</,</rm log>/,> /def/const<</rm const.>/,>
/def/dim<</rm dim>> /def/aff</,</it aff>/,> /def/diag</,</rm diag>/,>


/def/cc<</cal C>> /def/ck<</cal K>>
/def/cl<</cal L>> /def/cx<</cal X>>


/def/rr<I/!/!R>
/def/halv<<1/over 2>/,>
/def/nil<<(0)>>
/def/nk<<n_k>>
/def/nkp<</lp/nk/rp>>
/def/upto<,/allowbreak/ldots,>
/def/qed</hskip3pt/vrule height4pt width3pt depth2pt>
/def/co<C_0>
/def/tmpa#1<</ssg/log f(A_j,Y_j)/over #1>>
/def/tmpb#1<</ssg/log f(x)/over #1>>
/def/thcurl</theta/kern-7pt/lower5pt/hbox<$/sim$>>
/def/eqtop#1</,</buildrel <#1>/over =>/,>
/def/totop#1</,</buildrel <#1>/over /longrightarrow>/,>
/def/k<<(k)>>
/def/xkj<x_j^/k>
/def/Xkj<X_j^/k>
/def/hmk</hm_k>
/def/fk<f_k>
/def/tsk</ts_k>
/def/hsk</hs_k>
/def/lb</left/lbrace>
/def/rb</right/rbrace>
/def/lp</left(>
/def/rp</right)>
/def/fkx<f_k(x)>
/def/bk<</bf k>>
/def/ca<C^/ast>
/def/hco</hc_0>
/def/hCo</hC_0>
/def/ftfk<f_1/upto f_K>
/def/pkx<P(k/mid x)>
/def/ptx<P(t/mid x)>
/def/tho</th_0>
/def/mtlek<</rm max>_<t/leq K>/,>
/def/pr</parallel>
/def/phx<P_h(x)>
/def/ndel<<1/over n>/,>
/def/icx</int_/cx>
/def/lsb</left/lbrack>
/def/rsb</right/rbrack>
/def/lfull<L_</rm full>>
/def/twice#1<#1#1>
/def/mo#1</max_<#1>>
/def/mio#1</min_<#1>>
/def/mom</mo</om>>
/def/miom</mio</om>>
/def/ktok</,k=1/upto K>
/def/ave#1</,</scriptstyle/rm ave/atop <#1>>/,>
/def/xtxn<x_1/upto x_n>
/def/XtXn<X_1/upto X_n>
/def/nil<<(0)>>
/def/tmpa#1<</partial/log f(A_j,Y_j)/over #1>>
/def/tmpb#1<</partial/log f(x)/over #1>>
/def/halv<<1/over 2>>
/def/rr<I/!/!R>
/def/kcurl<k/kern-7pt/lower5pt/hbox<$/sim$>>
/def/co<C_0>
/def/cov<</rm cov>>
/def/eqtop#1<</buildrel <#1>/over =>>
/def/totop#1<</buildrel #1/over /longrightarrow>>
/def/tsk</ts_k>
/def/prt</partial>
/def/hsk</hs_k>
/def/hvlam</vert/Lam/vert^<1/sla 2>>
/def/pietc<(2/pi)^<-d/sla 2>>
/def/mh<<-1/sla 2>>
/def/mdh<<-d/sla 2>>
/def/bd</overline/delta>
%
%
/def/wh</widehat>
/def/hph<</wh/ph>> /def/hc<</wh c>> /def/hC<</wh C>>
/def/hl<</wh/lam>> /def/hP<</wh P>>
/def/hth<</wh/th>> /def/hs<</wh/sg>> /def/ha<</wh a>> /def/hA<</wh A>>
/def/hm<</wh/mu>> /def/hp<</wh/pi>> /def/hf<</wh f>> /def/hd<</wh/dl>>
/def/heps<</wh/eps>> /def/hg<</wh/gamma>> /def/ho<</wh/omega>>
/def/hr<</wh/rho>> /def/hssg<</wh/ssg>>
%
%
/def/dhat#1</widetilde #1>
%
%
/def/wt</widetilde>
/def/tm</wt/mu> /def/tr</wt/rho> /def/ts</wt/sg> /def/tss</wt/ssg>
/def/tf</wt f>
%
%
/def/db<</bar /dl>>

/def/invert#1</vert#1/vert>
/def/vsgv</invert/sg>
/def/vhsv</invert/hs>

/def/inpar#1</parallel#1/parallel>

%
%
/def/emodel</eps_</rm model>>
/def/etrue</eps_</rm true>>
/def/eideal</eps_</rm ideal>>
/def/bj<</overline j>>
/def/uj<</underline j>>


/header<Preface>
The Norwegian Computing Centre has been involved with symbol and pattern 
recognition work since 1982, through joint projects with the Norwegian
firm SysScan, and with support from the Royal Norwegian Council for
Scientific and Industrial Research. One particular task has been to design
and build a general-purpose statistical classification system capable of
recognising handwritten and machine-printed symbols on maps. A broad overview
of the field is offered by Llewellyn, Preston, Kittler and Harris (1982);
see also Harris, Kittler, Llewellyn and Preston (1982).

These reports constituted the starting point and an important basis for the
ASSEL project (Automatic Segmentation and Symbol recognition of
vectorised Linelike drawings) at the NCC. The ASSEL project has been
carried out on a contract with SysScan since 1983. SysScan has been
responsible for the preprocessing machinery, including various scanning,
skeletonising, boundary-finding, vectorising, and segmentation procedures.
A complete system for automatic digitisation of maps, including many of
the algorithms described in the present report, is marketed
by SysScan under the name GEOREC.

The classification system that has been developed at the NCC is a
production program, but is at the same time designed as a laboratory tool,
where new methods for feature extraction and for statistical classification 
are easy to test. The system is presented by Br}ten, Holb{k-Hanssen and
Taxt (1986a, limited circulation; 1986b). It is also worth pointing out that
the general classification methodology as well as central parts of the
software system can be applied in other and totally different areas of
application.

The ASSEL project has also required reviewing and extensions of the general
theory of pattern recognition, in particular the statistical approach.
The present report provides much of the theoretical background for the
various procedures that have been implemented, and acts as a sister report
to Br}ten, Holb{k-Hanssen and Taxt (1986a).
/vfill/eject
The art of automatically recognising symbols consists of three separate parts:
a machine preprocessing system must be built to isolate and represent 
candidate symbols in some suitable form (SysScan provides such a system); one
must decide upon certain measurements to take for each candidate symbol,
constituting together a vector of feature components; and a discriminant
or classification procedure must be constructed, a rule that for each
given feature vector assigns a class label. The present report is mostly
concerned with the third and most intrinsically statistical part of the
problem: what is the best classification rule for the given feature 
extraction method? Statistical theory can also provide a basis for giving
solutions to the second part of the problem, that of extracting features
with optimal discriminatory power, but this will be reported on
elsewhere.

The report reviews some of the basic theory of statistical pattern
recognition, and goes further in some places by contributing new methods
and new results about known methods. This is commented upon in the
``concluding remarks'' that end each chapter; see also the brief 
introduction and the outline of the contents of the report that are
offered subsequently.

Although parts of the report are written in the style of a textbook and 
could be
used in a graduate level course on classification techniques and discriminant
analysis, it does not claim to be one. Other parts have the character
of a research treatise. A ``textbook version'' would have included more of the
standard methods that are sometimes mentioned only briefly here, exactly
because they are well treated in the literature; a separate chapter on
feature extraction methods and their merits; and would have toned down
some of the other aspects dwelt on here.

Details about feature extraction methods and experimental results assessing
the performance of specific classification methods are deliberately
left out in this report. Some of these explorations are summarised in Br}ten,
Holb{k-Hanssen and Taxt (1986a). The results have been remarkably good.
Further results will be recorded in forthcoming reports. The joint project
between SysScan and the NCC will continue, requiring new theoretical as well
as practical work, with the aim of producing even better results in the
future.

The project group has on an average consisted of one statistician, allocating
1/sla 4 of his time to the project, and three computer scientists, each
working about half time on the project. It continues to be a pleasant,
challenging, and inspiring experience for me as a statistician to work so 
close to the actual applications of the methodology and in close
collaboration with computer scientists.
/goodbreak

The original project group consisted of Knut Br}ten, Erik Holb{k-Hanssen,
and Bryn Llewellyn in addition to myself. Later on  Thomas Bjerch and
Torfinn Taxt have joined the project group. I wish to thank
all of them for continually posing tough questions (and demanding
answers to them) and for their sharing of knowledge. Thanks are also
due to my fellow statisticians and classification theoreticians Jon
Helgeland and Erik Mohn for many discussions. Jon Helgeland has
co-authored Chapter 2 of the present report, and Erik Mohn has worked with me
on classification methods for remotely sensed data. I am also grateful
to Helge Roald and the other SysScan-people, for paying us to do interesting
work, and for contributing to it. The NCC has kindly supported me by giving
me time to finish this report. Finally I owe special thanks to Bj|rn
Larsen, who is responsible for the beautiful type-setting of my manuscript.
/nobreak/bigskip/nobreak
/hbox</hfil/vbox</hfil Oslo, December 1986/hfil/medskip
/hfil Nils Lid Hjort/hfil>/hbox to 3truecm<>>
/vfill
/begingroup/parindent=12pc
 /advance/rightskip/parindent

/noindent About the cover: Helga Lid Hjort made the ASSEL figure. 
(Working out methods for recognising rotated symbols has been one of the
tasks of the project.)  The project name
ASSEL (Automatic Segmentation and Symbol/break recognition of vectorised
Linelike drawings) is the German word for centipede. The
true/sla false figure is taken from Scott Kim's ``Inversions'',
BYTE Books, McGraw-Hill, 1981.

/endgroup
/eject


/header<Introduction and overview>
The idealised classification problem involves $K$ classes with class
densities $f_1/upto f_K$ and so-called </it prior> or 
</it a priori>/nils<Gjetter du vil ha disse i italics, som n}r du innf|rer
dem i de senere kapitlene.> probabilities $/pi_1/upto/pi_K$. The
task is to construct a rule that for any given feature vector $x$ (which
could be the measurements taken on a new candidate symbol) assigns a class
label (and hopefully the correct one).

The $/pi_k$'s may sometimes be known in advance; for example, it is
reasonable to put $/pi_k=1/sla10$ for each of the numbers $0$,
$1/upto9$ in an application involving a large number of handwritten
numerals. The class densities are almost always unknown, however.
Chapters 1 and 2 are nevertheless concerned with the idealised situation
in which these are known. The optimal classification rules developed
there, and which depend upon $f_1/upto f_K$, provide Platonian
ideals to aim at for later approximations.

Included in Chapters 1 and 2 is discussion of the ``doubt'' decision option.
The construction of the optimal doubt region of the feature space
derived in Chapter 2 with user-defined restrictions on all individual 
error rates is new.

Training data, i.e.~a collection of vectors whose true class labels are known
(by human inspection, say), make it possible to estimate
the class densities. Chapter 3 and 4 employ parametric models for these, and
discuss relevant methods of estimation. The result is a parametric
classifier with the parameters in question estimated from
training data. Realising that a given model described with a finite
number of parameters never can be completely correct, it becomes important
to understand exactly what fitting such a model to data means. This is
discussed in Section 3.1, where arguments also are provided to show that even
parametric models that are incorrect, strictly speaking,
may give rise to very useful classification
rules. Some of the results of Sections 3.1 and 4.1 appear to be new.

In Chapter 5 several nonparametric and semiparametric methods are developed.
These rules tend to outperform parametric procedures when the training sets
are large. General theoretical background is provided in Section 5.1. The
orthogonal expansion approach discussed in Sections 5.3 and 5.4 seems
particularly valuable, in that some of these procedures rely on </it some//>
structure being present in the data, i.e.~avoid being </it too//>
nonparametric. Several of the methods of Section 5.4 are new. Also 
treated in Chapter 5 is the $k$-nearest-neighbour approach.

Occasionally the preprocessing machinery may come up with a candidate symbol
that does not belong to any of the predefined classes, i.e.~an ``outlier''.
Procedures for filtering out such wild things are developed in Chapter 6.
Those that are based on a ``simultaneous likelihood ratio approach'' 
and take estimation variability into account seem to be new.

The nonparametric procedures of Chapter 5 are better than the parametric
ones of Chapters 3 and 4, </it provided//> enough training data are available.
It is also the case that the performance of any given classification rule
becomes steadily better with more training data. Getting hold of such training
samples is often a bottleneck, however, requiring sometimes too many hours
spent on tedious label assigning and editing. On the other hand it is 
sometimes inexpensive to get hold of a moderate to large number of 
feature vectors from </it unclassified//> symbols. Chapter 7 considers
the potential for these to update class descriptions, i.e.,
the initial class density estimates $/hf/upto/hf_K$ may be replaced by more
precise versions $f_1^*/upto f_K^*$, exploiting these new vectors with
</it unknown//> class labels. Also discussed is the possibility of passing 
from one level of sophistication, say the multinormal model, to a more
sophisticated description. Methods are also proposed for updating the
$k$-nearest-neighbour rule and the general orthogonal expansion procedure.

Chapter 8 provides solutions to a commonly occurring problem in pattern
recognition, especially in situations with high feature vector dimension,
namely that of singular covariance matrices (for some of the classes).
This may actually be used to enhance discriminatory power.

While much of the theory reviewed and developed in this report is 
geared towards feature methods with continuously valued components, it is
obviously the case that discrete valued features may provide good
discrimination ability. Chapter 9 treats several useful models for
discrete and mixed discrete and continuous data, and methods for 
estimating the class probability distributions. Classification rules based
on such mixed features can be very effective.

Several important questions must be (tentatively or accurately)
answered when a new feature extraction method is explored. How well are
the classes separated? Between which pairs of classes can confusion
be expected? Can the simple-minded multinormal density description 
adequately capture the data structure, or is a more sophisticated method
necessary? What size should a secondary training set have, if needed?
How small are the error rates based on the appropriately constructed
classification rule likely to be? Will the updating methods of Chapter 7
succeed in lowering error rates?

Chapters 10 and 11 aim at giving clues to correct answers to these
admittedly ambitious questions. Chapter 10 proposes a collection of
descriptive statistics and graphical plots that could constitute
``standard output'' for each class, once a feature method is chosen.
Also discussed in Chapter 10 is a new measure of class separability,
called here the generalised Mahalanobis distance. Chapter 11 considers 
checking procedures for various popular models, in particular the
multinormal model. A new test for normality based on distance from
data points to mean vector is developed.

Finally Chapter 12 is concerned with the important topic of assessing the
quality or performance of a classification method. Methods for
estimating error rates are discussed, in particular the cross
validation or leave-one-out method. Also included is a more ambitious
and computer intensive procedure based on resampling.

The field of pattern recognition is still rapidly growing, and along an
expanding boundary. Several topics of importance for statistical 
symbol recognition are not treated in this report, or only briefly
mentioned. The project group at the Norwegian Computing Centre has gathered
experience in some of these areas, including special methods for 
recognising rotated symbols, construction of hierarchical classifiers,
feature selection methods, structural and syntactical methods, and 
cluster analysis. Future work may include problems of ``reducing 
dimensionality'' (for which a good reference is the book edited by Krishnaiah
and Kanal, 1982); applying modern computer intensive nonparametric 
methods like classification trees (see Breiman, Friedman, Olshen, and
Stone, 1984) and projection pursuit (see Huber, 1985, and the
concluding remarks ending Chapter 5); further exploration of the 
automatic updating approaches to reduce or avoid costly training
(some suggestions are already in the present Chapter 7); and
hybrid procedures involving both structural and statistical methods.


/begingroup
/parindent=0pt
/newbox/chno
/newbox/secno
/def/CH#1#2#3</setbox/chno=/hbox</bf#1>/leftskip=0pc /medskip /par%
/goodbreak/hbox to/hsize</hbox to 2pc</bbf#1./hfil>/bbf#2/hfill/rm#3>>
/def/ch#1#2</hbox to/hsize</hbox to 2pc<>/bbf#1/hfill/rm#2>>
/def/SEC#1#2#3</setbox/secno=/hbox</unhcopy/chno/bf.#1>%
/smallskip /par%
/hbox to/hsize</hbox to 2pc<>/hbox to 2pc</unhcopy/secno/hfil>/bf#2/hfill/rm#3>>
/def/sec#1#2</hbox to/hsize</hbox to 4pc<>/bf#1/hfill/rm#2>>
/def/SS#1#2< /par%
/hbox</hbox to 4pc<>/hbox to 3pc</unhcopy/secno/bf.#1/hfil>/rm#2>>
/def/ss#1</hbox</hbox to 7pc<>/rm#1>>
/header<Table of contents>
/noindent/hbox to 2pc<></bbf Preface/hfill/rm i>/medskip
/noindent/hbox to 2pc<></bbf Introduction and overview/hfill/rm iv>/medskip
/noindent/hbox to 2pc<></bbf Table of contents/hfill/rm vii>/medskip
/CH 1 <Optimal classification when class densities><>
/ch    <and prior probabilities are known><1>
    /SEC 1 <Optimal classification><1>
    /SEC 2 <Concluding remarks><4>
/CH 2 <Optimal classification with bounds on error rates><>
/ch    <(With Jon Helgeland)><5>
    /SEC 1 <Bounding error rates><5>
    /SEC 2 <Concluding remarks><8>
/CH 3 <Classification based on parametric models:><>
/ch    <The estimative approach><9>
    /SEC 1 <General considerations><9>
        /SS A <Bayes risk consistency>
        /SS B <Fitting parametric families when they are wrong>
    /SEC 2 <Examples><14>
        /SS A <The multi-exponential model>
        /SS B <The multinormal model and the best quadratic rule>
        /SS C <Common covariance matrix: the best linear rule>
        /SS D <A mixed model for discrete and continuous components>
        /SS E <Mixtures of normals>
    /SEC 3 <Concluding remarks><26>
/CH 4 <Classification based on parametric models:><>
/ch    <The predictive approach><29>
    /SEC 1 <Theory><29>
    /SEC 2 <Examples><33>
        /SS A <The exponential model>
        /SS B <Multinormal densities with known covariance matrices>
        /SS C <Unknown covariance matrices>
        /SS D <Common covariance matrix>
        /SS E <Mixed discrete and continuous data>
    /SEC 3 <Concluding remarks><42>
/CH 5 <Nonparametric classification><44>
    /SEC 1 <General considerations><45>
    /SEC 2 <Density estimation by the kernel method><48>
    /SEC 3 <Density estimation by orthogonal expansions><50>
        /SS A <A density on $/lbrack0,1/rbrack$>
        /SS B <General case>
        /SS C <Estimating a $d$-dimensional density>
    /SEC 4 <Extended theory, and some practical solutions in><>
    /sec   <high-dimensional space><60>
        /SS A <Extended theory>
        /SS B <Solution 1: The multiplicative cosine expansion estimator>
        /SS C <Solution 2: The multiplicative beta with cosine expansions>
        /SS D <Robustification>
        /SS E <Penalising higher order terms>
        /SS F <Solution 3: The third order correction to normality>
        /SS G <Taking initial estimation variability into account:>
	/ss   <A more sophisticated inclusion rule>
        /SS H <Solution 4: From orthogonal expansions to projection pursuit>
    /SEC 5 <$k$-nearest-neighbour methods><84>
        /SS A <</it k-NN//> density estimation>
        /SS B <</it k-NN//> classification>
        /SS C <Remarks on computation>
    /SEC 6 <Concluding remarks><87>
/CH 6 <Detecting outliers><89>
    /SEC 1 <Finding outliers when class densities are known><89>
        /SS A <Multinormal case, common covariance matrix>
        /SS B <General case>
        /SS C <A nonparametric outlier test>
    /SEC 2 <Finding outliers when class descriptors are unknown><94>
        /SS A <Common covariance matrix>
        /SS B <Different covariance matrices>
    /SEC 3 <Concluding remarks><100>
/CH 7 <Updating parameter estimates using><>
/ch   <unclassified objects><101>
    /SEC 1 <Estimating prior probabilities on the given map><102>
        /SS A <Maximum likelihood estimators>
        /SS B <Bayes estimators>
    /SEC 2 <Simultaneous updating of class parameters and><>
    /sec   <prior probabilities><107>
        /SS A <General program>
        /SS B <Normal case>
    /SEC 3 <Passing from parametric to semiparametric><>
    /sec   <and nonparametric models><112>
        /SS A <Updating coefficients in an expansion density estimator>
        /SS B <Updating from normality to third order corrected normality>
	/SS C <Updating the </it k-NN//> rule>
    /SEC 4 <Concluding remarks><115>
/CH 8 <Handling singular covariance matrices><120>
    /SEC 1 <Exact theory><120>
    /SEC 2 <Concluding remarks><123>
/CH 9 <Discrete feature vectors><124>
    /SEC 1 <Binary feature components><125>
        /SS A <The Bahadur-Lazarsfeld expansion>
        /SS B <Inclusion rules>
        /SS C <Chow expansions>
    /SEC 2 <General discrete feature vectors><130>
    /SEC 3 <Mixed discrete and continuous variables><131>
    /SEC 4 <Concluding remarks><132>
/CH <10> <Class descriptions and separability measures><135>
    /SEC 1 <Basic descriptive measures for classes><136>
        /SS A <Robust estimation of mean and covariance matrix>
        /SS B <Further ``standard output''>
    /SEC 2 <Error rates and separability measures><142>
        /SS A <Error rates>
        /SS B <Affinity, distance measures, and a generalised Mahalanobis distance>
        /SS C <Other models>
        /SS D <More than two classes>
    /SEC 3 <Estimating the generalised Mahalanobis distance><155>
        /SS A <Adjustments for bias>
        /SS B <Amount of noise in the distance estimates>
    /SEC 4 <Concluding remarks><163>
/CH <11> <Checking model assumptions><165>
    /SEC 1 <Introduction><165>
    /SEC 2 <Checking normality><167>
        /SS A <Distribution of distance and direction>
        /SS B <A rigorous test based on distances>
    /SEC 3 <Checks for other density estimation schemes><176>
    /SEC 4 <Concluding remarks><177>
/CH <12> <Estimating error rates><178>
    /SEC 1 <The leave-one-out method><179>
        /SS A <Apparent error rates and cross validation>
        /SS B <Leave-one-out estimates for the best quadratic rule>
        /SS C <Leave-one-out estimates for the best linear rule>
        /SS D <Leave-one-out estimates for some nonparametric methods>
    /SEC 2 <Bootstrapping and other resampling methods><185>
        /SS A <Bootstrap and jackknife estimators>
        /SS B <The <$.632$> estimator>
    /SEC 3 <Concluding remarks><191>
/bigskip
/noindent/hbox to 2pc<></bbf References/hfill/rm 193>/medskip
/endgroup


/input bldef.tex 

/chapter<1><Optimal classification when class densities/nl 
and prior probabilities are known>/pageno=1
The problem is the following: A vector $X$ is observed in a feature space 
$/cx $. $X$ is known to be of exactly one of $K$ possible types, corresponding
to classes $1/upto K$. Observations from class $k$ are distributed according to
the density $f_k(x)$, and the prior probabilities are $/pi_1/upto /pi_K$.
(The densities are w.r.t.~a common $/sigma$-finite measure $/mu$.)
$X$ is to be classified, i.e. one of $K+1$ possible decisions $1/upto K,D$ is 
to be reached on the basis
of the observed value $X=x$;  decision $k$ of course corresponds to the claim
``$X$ is from class $k$'', whereas $D$ means ``being in doubt''.
/section<1.1><Optimal classification><>
In order to choose a ``good'' (or perhaps even the ``best'') classification
procedure one has to agree on reasonable overall criteria, for example
involving the misclassification probabilities 
$$/pmc(k)=Pr/lb C^/ast(X)/not=k,C^/ast(X)/in
/lb1/upto K/rb/mid C=k/rb/eqno(1.1)$$
and the reject or doubt probabilities 
$$/pd(k)=Pr/lb C^/ast(X)=D/mid C=k/rb./eqno(1.2)$$
Here $C/in/lbrace 1/upto K/rbrace$ denotes the true class
in question, and 
$$C^/ast:/cx /to T=/lbrace 1/upto K,D/rbrace$$
is a classification procedure.

The usual way of formalising goodness criteria is by introducing a </it loss
function//>. Let $L(k,t)$ be the loss incurred by making decision $t/in T$ if the
true class is $C=k$. One should have $L(k,k)=0$, all $k$, and maybe 
$L(k,D)=c$, all $k$, whereas
the $L(k,t)$'s, $k/not=t/in/lbrace 1/upto K/rbrace$, could in principle 
be any set of positive numbers.  If
every misclassification is equally serious, then
$$L(k,t)=/cases<0&if $t=k$ (correct decision)/cr
1&if $t/not=k$ and $t/in/lbrace1/upto k/rbrace$ (wrong decision)/cr
c&if $t=D$ (being in doubt)>/eqno(1.3)$$
for
$$k/in/Omega=/lbrace1/upto K/rbrace,/eqno(1.4)$$
$$t/in T=/lbrace1/upto K,D/rbrace,/eqno(1.5)$$
is reasonable.

In what follows we employ the loss function (1.3), since it is indeed
reasonable, and since important concepts are well illustrated.  Statements
and arguments similar to those below can easily be written down for the 
general loss function.

The </it risk function//> for procedure $C^/ast$ is the expected loss when using 
it, as
a function of the unknown class $k=1/upto K$:
$$/eqalign<R(C^/ast,k)&=E/left/{L(k,C^/ast(X))/mid C=k/right/}/cr
&=/sum^K_<t=1>/,L(k,t)/,Pr/left/{C^/ast(X)=t/mid C=k/right/}/cr
&/qquad+L(k,D)/,Pr/left/{C^/ast(X)=D/mid C=k/right/}/cr
&=/pmc(k)+c/,/pd(k)./cr>/eqno(1.6)$$
The </it total risk//>, or </it Bayes risk//>, is the total expected loss, 
viewing both the class $C$ and the vector $X$ as random:
$$/eqalign<R(C^/ast)&=ER(C^/ast,C)/cr
&=/sum^K_<k=1>/pi_k/,R(C^/ast,k)/cr
&=/sum_<k=1>^K/pi_k/,/pmc(k)+c/sum_<k=1>^K/pi_k/,/pd(k),/cr>/eqno(1.7)$$
which is seen to be the overall misclassification probability plus $c$ times
the overall reject (or doubt) probability.

We may derive the optimal classification procedure, i.e.~the one
minimising the total risk, under the assumption that prior probabilities
and class densities are </it known//>.  It is seen that
$$/eqalign<R(C^/ast)&=E/left/{EL(C,C^/ast(X))/mid X/right/}/cr
&=/int_/cx E/left/{L(C,C^/ast(x))/mid X=x/right/}/,f(x)/,d/mu(x)>$$
where
$$f(x)=/sum^K_<k=1>/pi_k/,f_k(x)/eqno(1.8)$$
is the marginal density for $X$.  It suffices to minimise the integrand
$$/sum^K_<k=1>L(k,C^/ast(x))/,P(k/mid x)$$
w.r.t. $C^/ast(x)/in T$, for each $x$, where 
$$/eqalign<P(k/mid x)&=Pr/left/{C=k/mid X=x/right/}/cr
&=/pi_k/,f_k(x)/sla/sum^K_<t=1>/pi_t/,f_t(x)/cr>/eqno(1.9)$$
is the </it posterior probability//> of class $k$ given $X=x$.  The integrand 
becomes
$$1-P(1/mid x)/upto1-P(K/mid x),c$$
when $C^/ast(x)=1/upto K,D$ respectively. We are, therefore, to maximise
$$P(1/mid x)/upto P(K/mid x),1-c,$$
and end up with the </it optimal solution//>
$$/co(x)=/cases<D&if every $P(k/mid x)/leq1-c$,/cr
k&if $P(k/mid x)=<>^<max>_<t/leq K>P(t/mid x)/gt1-c$./cr>/eqno(1.10)$$

This optimal procedure is also referred to as the </it Bayes solution//>, and
the corresponding value of the total risk, $R(/co)$ is called 
</it minimum Bayes risk//>.
$R(/co)$ is the best one can achieve if the $/pi_k$'s
and $f_k$'s are known.
/note<Remark 1.> The constant $c$ acts as a safety threshold, and should
in principle be specified by the user of the resulting classifier.  The
inconveniences caused by a reject has to be judged against the consequences of 
a misclassification.  In practice one usually tries a couple of $c$-values
on a training set of vectors with known classes, and obtains estimates of
misclassification and doubt rates, cf./ Chapter 12, before
a ``final value'' is chosen.  If $c$ is near $0$ then ``doubt'' is inexpensive.
This will lead to low error rates but on few classified vectors and a high doubt 
rate.
If on the other hand $c/geq1-<1/over K>$ then decision $D$ is so expensive
that it never will be used.
/note<Remark 2.> There are no restrictions on the type of densities 
$f_1/upto f_K$
in (1.9), for example, they need not be absolutely continuous w.r.t.
Lebesgue measure. As long as $f_k=dP_k/sla d/mu$ for some $/sigma$-finite measure $/mu$
dominating the class distributions $P_1/upto P_k$ both (1.9) and (1.10) 
continue to hold. Thus some or all of the $P_k$'s may have discrete components,
they may represent normal distributions with singular covariance matrices, 
etc.  
/section<1.2><Concluding remarks><>
The small piece of theory presented here is fairly standard, although the
rigorous derivation of the optimal reject region, by means of the loss
function, is less known. The most popular special cases of the optimal rule
are the normal distribution cases with common or different covariance 
matrices, cf.~Chapters 3 and 4. The methods that correspond to these
assumptions are known as the best linear and the best quadratic rule,
respectively.

Discriminant or classification analysis historically started with Fisher's 
1936-paper, where he derived the best linear rule, but from a different
perspective than above.

It may happen that the observed $X$ vector comes from a source different
from each of the $K$ defined classes. This motivates including also
``$/out$'' as an option for the classification procedure, i.e.~$X$ is an
outlier as seen from each of the $K$ theories $f_1/upto f_K$. This topic
is treated in Chapter 6.

Of course one rarely knows the class densities and prior probabilities
exactly. The optimal rule presented above can then be taken as an ``ideal
procedure'' for approximations based on relevant data to aim at. This is
indeed what we work on in several of the following chapters.


\input norsk

/def/nils#1<>

%
%
/font/bbf=cmbx10 scaled/magstep1
/font/Bbf=cmbx10 scaled/magstep3
/font/ninerm=cmr9
/font/ninesl=cmsl9
/font/ninebf=cmbx9
%
%
/mag/magstephalf
/input texini
/vsize22truecm
/hsize15truecm

/baselineskip=14pt

%
%

/newbox/xxbox
/def/nl</hfill/break>
/def/header#1</vfill/supereject/vbox to 2in<>%
/noindent</Bbf #1>/bigskip/noindent/ignorespaces>
/def/chapter#1#2</message<(Chapter #1:#2)>%
/vfill/supereject/vbox to 2in<>%
/noindent</Bbf Chapter #1/medskip/noindent #2>%
/bigskip/noindent/ignorespaces>
/newbox/seca /newbox/secb /newbox/threebox
/def/section#1#2#3</message<(#1)>/goodbreak/bigskip/nobreak/noindent%
/setbox/threebox/hbox<#3>
/setbox/xxbox/hbox</bbf #1/quad>%
/setbox/seca/hbox</bbf #1/quad #2>%
/setbox/secb/hbox</hbox to/wd/xxbox<>/bbf #3>%
/begingroup%
/bbf/hbox to/hsize</hss/vbox<%
/box/seca%
/ifdim/wd/threebox=0pt /else/box/secb/fi>/hss>%
/nobreak/smallskip/nobreak/ignorespaces%
/endgroup>
/def/subsection#1#2</message<(#1)>/goodbreak/bigskip/nobreak%
</setbox/xxbox/hbox</bf #1/quad>%
/par/parindent=0pt /hangindent=/wd/xxbox /hangafter=1%
/bf#1/quad #2/par>/nobreak/smallskip/nobreak/noindent/ignorespaces>

/def/Btheorem</medbreak/noindent</csc Theorem:>/quad/it/ >
/def/Etheorem</rm/hfil/medskip>

/def/Blemma</medskip/noindent</csc Lemma:>/quad/it>
/def/Elemma</rm/hfil/medskip>

/def/Bproof</medskip/noindent</csc Proof:>/quad/ >
/def/Eproof</qed/hfil/medskip>

/def/note#1</medskip</csc #1>>
/let/Enote=/medskip

/def/subchapt#1#2</section<#1><#2>>


/def/dl</delta>
/def/lam</lambda>
/def/Lam</Lambda>
/def/eps</varepsilon>
/def/ssg</sigma>
/def/sg</Sigma>
/def/th</theta>
/def/Th</Theta>
/def/om</Omega>
/def/som</omega>
/def/ph</phi>


/def/est<</rm est>>

/def/rx<</rm x>> /def/ry<</rm y>> /def/rz<</rm z>>
/def/rX<</rm X>> /def/rY<</rm Y>> /def/rZ<</rm Z>> 

/def/doubt<</rm DOUBT>> /def/out<</rm OUT>>

/def/pmc<</rm pmc>> /def/pcc<</rm pcc>> /def/pd<</rm pd>>
/def/apmc<</rm </overline p>mc>> /def/hpmc<</rm </wh p>mc>>
/def/op<</rm op>>

/def/cross<</rm (CROSS)>>
/def/boot<</rm (BOOT)>>
/def/jack<</rm (JACK)>>

/def/Pr<</it Pr>> /def/el<</rm EL>> /def/Tr<</rm Tr>>

/def/lr<</rm LR>> /def/ise<</rm ISE>> /def/mise<</rm MISE>>
/def/At<A_t> /def/Bt<B_t> /def/Ct<C_t> /def/Dt<D_t>

/def/exp</,</rm exp>/,> /def/log</,</rm log>/,> /def/const<</rm const.>/,>
/def/dim<</rm dim>> /def/aff</,</it aff>/,> /def/diag</,</rm diag>/,>


/def/cc<</cal C>> /def/ck<</cal K>>
/def/cl<</cal L>> /def/cx<</cal X>>


/def/rr<I/!/!R>
/def/halv<<1/over 2>/,>
/def/nil<<(0)>>
/def/nk<<n_k>>
/def/nkp<</lp/nk/rp>>
/def/upto<,/allowbreak/ldots,>
/def/qed</hskip3pt/vrule height4pt width3pt depth2pt>
/def/co<C_0>
/def/tmpa#1<</ssg/log f(A_j,Y_j)/over #1>>
/def/tmpb#1<</ssg/log f(x)/over #1>>
/def/thcurl</theta/kern-7pt/lower5pt/hbox<$/sim$>>
/def/eqtop#1</,</buildrel <#1>/over =>/,>
/def/totop#1</,</buildrel <#1>/over /longrightarrow>/,>
/def/k<<(k)>>
/def/xkj<x_j^/k>
/def/Xkj<X_j^/k>
/def/hmk</hm_k>
/def/fk<f_k>
/def/tsk</ts_k>
/def/hsk</hs_k>
/def/lb</left/lbrace>
/def/rb</right/rbrace>
/def/lp</left(>
/def/rp</right)>
/def/fkx<f_k(x)>
/def/bk<</bf k>>
/def/ca<C^/ast>
/def/hco</hc_0>
/def/hCo</hC_0>
/def/ftfk<f_1/upto f_K>
/def/pkx<P(k/mid x)>
/def/ptx<P(t/mid x)>
/def/tho</th_0>
/def/mtlek<</rm max>_<t/leq K>/,>
/def/pr</parallel>
/def/phx<P_h(x)>
/def/ndel<<1/over n>/,>
/def/icx</int_/cx>
/def/lsb</left/lbrack>
/def/rsb</right/rbrack>
/def/lfull<L_</rm full>>
/def/twice#1<#1#1>
/def/mo#1</max_<#1>>
/def/mio#1</min_<#1>>
/def/mom</mo</om>>
/def/miom</mio</om>>
/def/ktok</,k=1/upto K>
/def/ave#1</,</scriptstyle/rm ave/atop <#1>>/,>
/def/xtxn<x_1/upto x_n>
/def/XtXn<X_1/upto X_n>
/def/nil<<(0)>>
/def/tmpa#1<</partial/log f(A_j,Y_j)/over #1>>
/def/tmpb#1<</partial/log f(x)/over #1>>
/def/halv<<1/over 2>>
/def/rr<I/!/!R>
/def/kcurl<k/kern-7pt/lower5pt/hbox<$/sim$>>
/def/co<C_0>
/def/cov<</rm cov>>
/def/eqtop#1<</buildrel <#1>/over =>>
/def/totop#1<</buildrel #1/over /longrightarrow>>
/def/tsk</ts_k>
/def/prt</partial>
/def/hsk</hs_k>
/def/hvlam</vert/Lam/vert^<1/sla 2>>
/def/pietc<(2/pi)^<-d/sla 2>>
/def/mh<<-1/sla 2>>
/def/mdh<<-d/sla 2>>
/def/bd</overline/delta>
%
%
/def/wh</widehat>
/def/hph<</wh/ph>> /def/hc<</wh c>> /def/hC<</wh C>>
/def/hl<</wh/lam>> /def/hP<</wh P>>
/def/hth<</wh/th>> /def/hs<</wh/sg>> /def/ha<</wh a>> /def/hA<</wh A>>
/def/hm<</wh/mu>> /def/hp<</wh/pi>> /def/hf<</wh f>> /def/hd<</wh/dl>>
/def/heps<</wh/eps>> /def/hg<</wh/gamma>> /def/ho<</wh/omega>>
/def/hr<</wh/rho>> /def/hssg<</wh/ssg>>
%
%
/def/dhat#1</widetilde #1>
%
%
/def/wt</widetilde>
/def/tm</wt/mu> /def/tr</wt/rho> /def/ts</wt/sg> /def/tss</wt/ssg>
/def/tf</wt f>
%
%
/def/db<</bar /dl>>

/def/invert#1</vert#1/vert>
/def/vsgv</invert/sg>
/def/vhsv</invert/hs>

/def/inpar#1</parallel#1/parallel>

%
%
/def/emodel</eps_</rm model>>
/def/etrue</eps_</rm true>>
/def/eideal</eps_</rm ideal>>
/def/bj<</overline j>>
/def/uj<</underline j>>

/chapter<2><Optimal classification with/nl bounds on error rates>
The traditional Bayes discriminant  rule picks the class with highest
posterior probability, and is (1.10) with the threshold $c$ chosen so large 
that decision $D$ will never be reached ($c/geq 1-<1/over K>$). It is clear 
that this rule may produce
large misallocation probabilities when some classes are easily confused,
i.e. when more than one of the $/pi_k/,f_k(x)$'s may be large in the same region.  
Using a
lower threshold and the doubt option is one way of reducing the
misclassification rates, and (1.10) provides the optimal solution if the
loss function (1.3) is considered to be the appropriate measure of
consequences.

Procedure (1.10) with a given $c/in/,(0,1-<1/over K>)$ offers no 
guarantee against
high error probabilities, however, it only provides the minimal cost in terms of
a combination of the </it average//> misclassification rate
$$/pmc=/sum_<k=1>^K/pi_k /pmc(k)/eqno(2.1)$$
and the average doubt rate
$$/pd=/sum_<k=1>^K/pi_k/,/pd(k),$$
where
$$/eqalign</pmc(k)&=Pr/left/{ C^/ast(X)/in /left/{1/upto k-1,k+1/upto K
/right/}/mid C=k/right/},/cr
/pd(k)&=Pr/left/{ C^/ast(X)=D/mid C=k/right/}.>/eqno(2.2)$$
/section<2.1><Bounding error rates><>
In some applications it may be important to guarantee low </it individual//> 
error rates (as opposed to low </it average//> error rate).  Specifically, one might
demand
$$/pmc(k)/leq/eps_k,/ k=1/upto K,/eqno(2.3)$$
resorting if necessary to the reject or doubt option.

The aim of the present chapter is to show that one's error probabilities
actually may be controlled in an optimal way.  (We assume
that the prior probabilities and class densities are known, as in Chapter 1.)
We want to find the procedure with highest success rate
$$/pcc=/sum_<k=1>^K/pi_k/,/pcc(k),/eqno(2.4)$$
where
$$/pcc(k)=Pr/lbrace C^/ast(X)=k/mid C=k/rbrace$$
is the probability of correct classification, among all procedures
obeying (2.3).

A (randomised) classification procedure may be defined as a function
$$/phi(x)=/left(/phi_o(x),/phi_1(x)/upto /phi_k(x)/right)$$
on the feature space $/cx $, interpreted via
$$/eqalign</phi_k(x)&=Pr/lbrace</rm decision>/,k/mid X=x/rbrace,/ k=1/upto K,/cr
/phi_o(x)&=Pr/lbrace</rm decision>/,D/mid X=x/rbrace.>$$
The $/phi_k(x)$'s are non-negative and sum to one.  Usually there are sets 
$A_0,A_1/upto A_K$
providing a partition of feature space such that $/phi_k(x)=I/lbrace x/in A_k
/rbrace,/ k=0,1/upto K$,/break
$I$ denoting the indicator function, so that the classifier is
</it non-randomised//>.  (Note that the Bayes rule found in Chapter 1 is (or 
rather: can be chosen) non-randomised, but that it is still optimal in
the larger class of all randomised rules.)

Let $/cc$ be the class of all procedures $/phi$ obeying 
$$/pmc(k)=/int/sum_<1/leq t/leq K, t/not=k>/phi_t/,f_k/leq /eps_k,
/ k=1/upto K,/eqno(2.5)$$
and let $/cc_0$ be the 
subclass having $/pmc(k)=/eps_k,/ k=1/upto K$, for a given set of numbers
$/eps_1/upto /eps_K$.
/Btheorem Assume there exist $/lam_1/upto /lam_K$ such that the 
classification procedure $/phi^/ast=(/phi_0^/ast,/phi_1^/ast/upto /phi_K^/ast)$
defined by
$$/phi_k^/ast(x)=I/lbrace H_k(x)=</rm max>_<0/leq t/leq K>/,H_t(x)/rbrace,
/ k=0,1/upto K/eqno(2.6)$$
belongs to $/cc_0$, where
$$/eqalign<H_k(x)&=(/pi_k+/lam_k)/,f_k(x),/ k=1/upto K,/cr
H_0(x)&=/sum_<k=1>^K/lam_k/,f_k(x).>$$
(If more than one $k$ achieves $H_k(x)=</rm max>_<0/leq t/leq K>
/,H_t(x)$, then $x$ may for example be assigned to the class
with lowest label.)
Then $/phi^/ast$ maximises the success probability (2.4) among all procedures 
in $/cc_0$.

Further, if the $/lam_k$'s are non-negative, then $/phi^/ast$ also maximises 
(2.4) among all procedures in $/cc$.
/Etheorem
/Bproof Let $/phi$ be an arbitrary procedure in $/cc_0$. By construction,
$$/sum_<k=0>^K /phi_k(x)/,H_k(x)/leq/sum_<k=0>^K/phi_k^/ast(x)/,H_k(x)$$
for every $x$.  Hence $/int/sum_<k=0>^K(/phi_k^/ast-/phi_k)/,H_k/geq0$.
We may subtract the function $H_0$ from the $H_k$'s to obtain
$$/int/sum_<k=1>^K(/phi_k^/ast-/phi_k)/lp/pi_k/,f_k-
/sum_<1/leq t/leq K, t/not=k>
/lam_t/,f_t/rp/geq0.$$
Accordingly
$$/eqalignno</sum_<k=1>^K/pi_k/,/int(/phi_k^/ast-/phi_k)/,f_k&/geq/sum_<k=1>^K
/sum_<1/leq t/leq K, t/not=k>/int(/phi_k^/ast-/phi_k)/lam_t/,f_t/cr
&=/sum_<t=1>^K/lam_t/int/sum_<1/leq k/leq K, t/not=k>(/phi_k^/ast-/phi_k)/,f_t/cr
&=/sum_<t=1>^K/lam_t/lbrace/pmc^/ast(t)-/pmc(t)/rbrace,&(2.7)/cr>$$
using $/pmc^/ast(t)$ and $/pmc(t)$ to denote the probability of misallocating 
an $X$
from class $t$ for the rules $/phi^/ast$ and $/phi$ respectively.  These 
are both equal to $/eps_t$ when $/phi^/ast$ and $/phi$ belong to $/cc_0$.
Hence
$$/pcc^/ast=/sum_<k=1>^K/pi_k/int/phi_k^/ast/,f_k/geq/sum_<k=1>^K/pi_k/int/phi_k/,f_k=/pcc,$$
proving the first part of the theorem. 
The second part, assuming only $/phi/in/cc$, also follows from (2.7). /Eproof

It is interesting to consider the case of $K=2$ classes. The feature space 
$/cx $ is to be divided into three sets $A_0,A_1,A_2$ such that vector $X$ is 
allocated to class $k$ if $X/in A_k,k=1,2$, whereas no final class assignment
is made if $x/in A_0$. We are to consider the functions
$$/lam_1f_1+/lam_2 f_2,(/pi_1+/lam_1)f_1,(/pi_2+/lam_2)f_2,$$
and determine the regions in which one of them dominates the other ones.
Class 1 will be assigned if
$$<f_1(x)/over f_2(x)>/gt </rm max>/lb</lam_2/over/pi_1>,
</pi_2+/lam_2/over/pi_1+/lam_1>/rb,$$
class 2 is chosen if
$$<f_1(x)/over f_2(x)>/lt</rm min>/lb</pi_2/over/lam_1>,
</pi_2+/lam_2/over/pi_1+/lam_1>/rb,$$
while decision ``doubt'' is reached when
$$</pi_2/over/lam_1>/leq<f_1(x)/over f_2(x)>/leq</lam_2/over/pi_1>.$$
(We assume for simplicity that the probability distributions are continuous,
so that the boundaries of the above regions get zero probability.)  It is
assumed here that $/lam_1,/lam_2/geq 0$.

The regions $A_0,A_1,A_2$ obtained above are recognised as those that can 
be derived 
from Neyman-Pearson's Fundamental Lemma (Lehmann 1959, p.65) when a doubt option
is used, if necessary, to control both the level $/alpha$ and the power 
$/beta$ of a test. 

Thus our theorem may be viewed as a generalisation of the Fundamental Lemma.
The generalisation is different from Lehmann's generalisation (1959, p.83),
but the methods of proof, as far as the sufficient conditions are concerned,
are similar.

The theorem is also similar to one presented by J.A.~Anderson (1969), who 
used bounds on all
$$/pmc(k,t)=/int/phi_t/,f_k,/ t/not=k.$$
Questions about existence, necessity, and uniqueness are not as easily answered.
The reader is referred to
the careful work of Anderson, which can be extended to cover the 
case under consideration here.

It would be useful to have an algorithm capable of finding the appropriate
$/lam_1/upto/lam_K$ for a given set of $/pi_k$'s, $f_k$'s, $/eps_k$'s, say 
for multinormal populations.  This is not a
trivial task, however.  Anderson's work (op.cit.) might be a starting point.
/section<2.2><Concluding remarks><> The doubt
option was considered in Chapter 1. The present chapter has provided another
way of constructing the doubt region. Bounding individual error rates is
certainly sensible, but using that as the fundamental criterion, as opposed
to trusting the loss function (1.3) for example, may be more appropriate 
in other applications than in pattern recognition problems.

In some situations some of the $K$ individual error rates are more important 
to keep low than other ones. A theorem similar to the one proved above can
be constructed for the situation where only a subset of the error rates
are to be bounded.

Constructing good numerical algorithms for finding the decision regions in 
practice is a solvable but non-trivial problem. J.A.~Anderson (1969) provides
some ideas (for a different problem) that can be used here too.

Of course class densities must be estimated in practice before the 
theorem above can be applied. Estimation uncertainty should also be taken
into account. The general theory presented in Chapters 3 and 4 applies to 
this end.

Jon Helgeland has co-authored Chapter 2.


/chapter<3><Classification based on parametric models:/hfil/break
The estimative approach>
/section<3.1><General considerations><>
In Chapter 1 the optimal classification procedure was derived under the loss
function (1.3), assuming that prior probabilities and class densities
were known.  The winning procedure was
$$/co(x)=/cases<D&if every $P(k/mid x)/leq1-c$,/cr
k&if $P(k/mid x)=<>_<t/leq K>^</rm max>P(t/mid x)/gt1-c$,/cr>$$
where
$$P(k/mid x)=</pi_k/,f_k(x)/over/sum_<t=1>^K/pi_t/,f_t(x)>,$$
and resulted in the minimum Bayes risk
$$/eqalign<R(/co)&=EL(C,/co(X))/cr
&=/sum_<k=1>^K/pi_k/lbrack/pmc_0(k)+c/,/pd_0(k)/rbrack./cr>$$
(Here $/pmc_0(k)$ and $/pd_0(k)$ denote misclassification and doubt 
probabilities for procedure $/co$.)

The unfriendly world in which statisticians live is however full of
densities and probabilities that are (at least partly) unknown, and
their task becomes one of providing good alternative procedures with
Bayes risks as close to $R(/co)$ as possible.  The usual approach has been to
utilise training data, i.e. a set of feature vectors for each class, to
construct approximations to $/co(x)$ in some sense.

It is assumed in the present chapter that the $/pi_k$'s are known and that the
class densities belong to parametric families, say
$$f_k(x)=f_k(x,/th),/eqno(3.1)$$
where $/th$ is an unknown vector of parameters.  In formulation (3.1) $/th$ 
may be thought of as a collection of component parameters, some of which
are class-dependent while others are common to all classes.  If the $K$ class
densities all belong to the </it same//> parametric family, differing only in the
values of the parameter, then
$$f_k(x)=f(x,/th_k)/eqno(3.2)$$
might be a better formalisation; we will encounter both types.

Assume that a training set
$$Z=/left/{X_j^<(k)>;/ j=1/upto n_k,/ k=1/upto K/right/}/eqno(3.3)$$
is available, where it is known that $X_1^<(k)>/upto X_<n_k>^<(k)>$ come 
from class $k$.  These data
may give rise to an estimate $/wh/th$ of $/th$ in (3.1).  A natural proposal 
is then the classification rule
$$/widehat/co(x)=/cases<D&if every $/widehat P(k/mid x)/leq1-c$,/cr
k&if $/widehat P(k/mid x)=<>_<t/leq K>^</rm max>
/widehat P(t/mid X)/gt1-c$,/cr>/eqno(3.4)$$
where $/wh/th$ is inserted in the class densities to produce approximate 
posterior probabilities
$$/widehat P/lp k/mid x/rp=</pi_k/,f_k(x,/wh/th)/over
/sum_<t=1>^K/pi_t/,f_t(x,/wh/th)>,/eqno(3.5)$$
$/wh/co$ above is called a </it plug-in rule//>, since it is constructed by just
plugging in an estimate $/wh/th$ for $/th$ in $/co$, which is known to  be 
optimal if $/th$ is the true value of the parameter. We could also write
$/hC_0=/hC_0(x;z)$ in that it depends upon the training set through the
estimate $/hth$.

It remains to decide exactly which estimator should be plugged in.  The maximum
likelihood (ML) estimator has been the most popular choice, or modifications
of it to make it unbiased.  The widespread use of the plug-in rule with the ML
estimator has been caused by the good general reputation the ML estimator
enjoys and the fact that several of the pioneers in the field, like Fisher
(1936), Mahalanobis (1936), Rao (in a series of papers starting in 1948),
and Anderson (1958), have directly or implicitly recommended it.

The use has been rather uncritical, though.  Although $/wh/th$ may be 
excellent as an estimator of $/th$ there is no guarantee that $f_k(x,/wh/th)$
is the best guess for $f_k(x,/th)$, nor is $/widehat/co(x)$ itself 
necessarily a good approximation to $/co(x)$. The performance
of plug-in rules and other procedures should really be judged by the criterion
of total risk $R(C^/ast)$ defined in (1.7), if (1.3) is still considered to be
the appropriate loss function, or by other criteria more tied to classification
accuracy than to the behaviour of $/wh/th$ as an estimator for $/th$.

These questions and related problems are returned to in Chapter 4.  Now some
general comments are offered pertaining to the use of parametric models in
discriminant analysis, before a list of examples is presented.
/subsection<3.1.A><Bayes risk inconsistency>
A reasonably simple observation that has been taken as support for the use
of plug-in parametric rules, is the following:  As the training set 
increases, i.e. $n_1/upto n_K$ go to infinity in (3.3), then under mild 
regularity $/lbrack/wh/th$ is 
consistent and the class densities are continuous in $/th/rbrack$ $/widehat/co$
in (3.4) becomes
identical to the optimal $/co$ and its Bayes risk $R(/widehat/co)$ converges to 
minimum Bayes risk $R(/co)$. Many plug-in rules, corresponding to a large 
class of possible estimators $/wh/th$, have this property.

There is an important assumption behind the argument above, however, namely
that the class densities $f_1/upto f_K$ in fact obey the parametric structure 
in question.  As even statisticians admit, their parametric models are only 
approximations to reality, implying in the present context that even when the size
of the training set increases beyond bounds, $/widehat/co$ in (3.4) will become 
close to only an approximation to $/co$ in (1.10), and the total risk 
$R(/widehat/co)$ will converge to a number greater than $R(/co)$.

It is most often possible to construct procedures that are </it Bayes risk 
consistent//> in the sense that the sequence of Bayes risks converges to 
minimum Bayes risk $R(/co)$. Unless one firmly believes in a certain parametric
model the Bayes risk consistent rules will necessarily involve 
</it nonparametric//> estimation of the densities, a topic considered in Chapter 5.

These comments are not meant to imply that parametric
models are useless; they may indeed constitute good and compact approximations
to more
complicated models.  Classifiers built on parametric assumptions may work
excellently.  Nonparametric density estimation demands for its successful
application far larger training sets than parametric alternatives.  Thus there
is a trade-off between perhaps simple, easily implementable algorithms that
work well even for moderately sized training sets, and heavier nonparametric
ones that may behave awkwardly for small to moderate training sets, but 
nevertheless will (nearly always) win if a sufficient amount of training data
is available.
/subsection<3.1.B><Fitting parametric families when they are wrong>
We will end the present section with a brief discussion of the behavior
of ML estimators when the underlying parametric model is not necessarily true.

Assume that $X_1/upto X_n$ are independent and identically distributed 
(i.i.d.) with a density $f(x)$, and that the parametric model $f_/th(x)$ is 
forced on the data, $/th$ being a $p$-dimensional parameter belonging to a set
$/Theta$. The ML estimator $/wh/th_n$ maximises the log-likelihood function
$$A_n(/th)=<1/over n>/sum_<i=1>^n/log/,f_/th(X_i)$$
w.r.t.~$/th$. By the strong law of large numbers
$$A_n(/th)/totop</rm a.s.>A(/th)=E_f/,/log/,f_/th(X_i)=/int/log/,f_/th(x)
/cdot f(x)/,d/mu./eqno(3.6)$$
Very often the function $A(/th)$ has a unique maximum for $/th=/th_0$.
$/th_0$ is not necessarily
the ``true value'' because we have not assumed that $f$ belongs to the family
of $f_/th$'s.
In a sense $/th_0$ is the value of $/th$ making $f_/th$ most close to the true 
$f$, however, in that it minimises the Kullback-Leibler distance
$$I(f,f_/th)=/int f(x)/,/log/,<f(x)/over f_/th(x)>/,d/mu./eqno(3.7)$$
One may now show that under reasonable regularity conditions
$/wh/th_n/totop</rm a.s.>/th_0$, thus 
generalising the classical consistency result for ML estimators.  
Consequently the ML plug-in rule $/widehat/co$ converges pointwise to a 
rule $C_1$ defined analogously to (1.10) but with posterior probabilities
of the form
$$P(k/mid x)=</pi_k/,f_k(x,/th_0)/over /sum_<t=1>^K/pi_t/,f_t(x,/th_0)>.$$
Also
$$R(/widehat/co)/totop</rm a.s.>R(C_1)/gt R(C_0)$$
under mild regularity.

The classical result on the limiting distribution of
$/sqrt<n>(/wh/th_n-/th_0)$
may also be generalised to the present agnostic state of affairs where 
$/lbrace f_/th;/th/in /Theta/rbrace$ not necessarily
contains the true $f$. $/wh/th_n$ solves the equations
$$U_<n,j>(/th)=</partial/over /partial/th_j>/,A_n(/th)=<1/over n>/sum_<i=1>^n
</partial/log f_/th(X_i)/over /partial/th_j>=0,/ j=1/upto p,$$
writing $/th=(/th_1/upto /th_p)'$.
Taylor expansion and analysis of remainder terms ensure
$$U_<n,j>/lp/wh/th_n/rp/doteq/ U_<n,j>(/th_0)+/sum_<l=1>^p I_<n,jl>
(/th_0)/lp/wh/th_<n,l>-/th_<0,l>/rp$$
where
$$I_<n,jl>(/th)=</partial/over/partial/th_l>/,U_<n,j>(/th)=<1/over n>
/sum_<i=1>^n</partial^2/log f_/th(X_i)/over/partial/th_j/,/partial/th_l>,/ j,l=1/upto p,$$
and $X_n/doteq/ Y_n$ means asymptotic equivalence in the sense 
$X_n-Y_n/totop<P> 0$. Hence
$$/sqrt<n>/,U_n(/th_0)/doteq/ -I_n(/th_0)/,/sqrt<n>/lp/wh/th_n-/th_0/rp,$$
where $U_n(/th)=(U_<n,1>(/th)/upto U_<n,p>(/th))'$ and $I_n(/th)$ is the 
$p/times p$ matrix of $I_<n,jl>(/th)$ entries.
Now 
$$/sqrt<n>U_n(/th_0)/totop<D>N_p(0,K(/th_0))$$
and
$$-I_n(/th_0)/totop<P>T(/th_0),$$
where
$$K(/th)=/left/lbrace E_f</partial/log f_/th(X_i)/over /partial/th_j>
</partial/log f_/th(X_i)/over /partial/th_l>/right/rbrace_<j,l/leq p>,/eqno(3.8)$$
$$T(/th)=/left/lbrace -E_f</partial^2/log f_/th(X_i)/over
/partial/th_j/partial/th_l>/right/rbrace_<j,l/leq p>./eqno(3.9)$$
Hence, by traditional arguments,
$$/sqrt<n>(/wh/th_n-/th_0)/totop<D>
N_p(0,T(/th_0)^<-1>/,K(/th_0)T(/th_0)^<-1>)./eqno(3.10)$$
Observe that if the true $f$ really belongs to $/lbrace f_/th;/th/in/Theta/rbrace$,
then $f=f_</th_0>$ with $/th_0$ as above, i.e.~$/th=/th_0$ minimises the
Kullback-Leibler distance (3.7) (and $I(f,f_</th_0>)=0$),
and $K(/th_0)=T(/th_0)$. (3.10) now simplifies to the classical result
$/sqrt<n>(/wh/th_n-/th_0)/totop<D>N_p(0,T(/th_0)^<-1>)$.

It is not difficult to estimate both $T(/th_0)$ and $K(/th_0)$ consistently.
This fact, combined with (3.10), makes it possible to obtain confidence
regions for and test hypotheses about the least false parameter value $/th_0$.
/note<Example 1.> Let $X$ have density $f(x)$ for $x$ in 
$/lbrack 0,/infty)$, and suppose $f_/th(x)=/th e^<-/th x>, x/gt 0$, is the
parametric structure. Minimising (3.7) is equivalent to maximising (3.6), 
that is
$$/int_0^/infty(/log/th-/th x)/,f(x)/,dx=/log/th-/th/mu_f$$
w.r.t.~$/th$, where $/mu_f=/int_0^/infty xf(x)/,dx$ is assumed finite.
Accordingly $/th_0=1/sla /mu_f$ is the ``pseudo-true'' or least false 
parameter, and the ML estimator $/wh/th_n=(<1/over n>/sum_<i=1>^nX_i)^<-1>$
converges to $/th_0$.
/note<Example 2.> Consider the multinormal density
$$/eqalign<f(x;/mu,/Sigma)&=N_p(/mu,/Sigma)(x)/cr
&=(2/pi)^<-p/sla 2>/,/vert/Sigma/vert^<-1/sla 2>
/,e^<-<1/over 2>(x-/mu)'/Sigma^<-1>(x-/mu)>/cr>$$
as an approximation to a given density $f$ on $/rr^p$. One may show that the 
parameter values $(/mu,/Sigma)=(/mu_0,/Sigma_0)$ that provides the best
approximation according to the Kullback-Leibler criterion are
$$/mu_0=E_f X=/int x/,f(x)/,dx,$$
$$/sg_0=</rm VAR>_f/,X=/int(x-/mu_0)(x-/mu_0)'/,f(x)/,dx.$$
Thus when the multinormal model is used to describe data from a density $f$ 
that perhaps is known </it a priori//> not to be multinormal itself, and the
ML estimators $/wh/mu_n,/widehat/Sigma_n$ nevertheless are computed, the
preceding theory shows that what they really estimate are $/mu_0,/Sigma_0$
above.
/note<Example 3.> A wide class of densities on the unit interval $(0,1)$
is provided by the </it beta distributions//>,
$$f_</alpha,/beta>(x)=</Gamma(/alpha+/beta)/over/Gamma(/alpha)/,/Gamma(/beta)>
x^</alpha-1>(1-x)^</beta-1>,$$
where $/alpha/gt0$ and $/beta/gt0$ are parameters. One has
$$EX=</alpha/over/alpha+/beta>,/ </rm Var/ >X=</alpha/beta/over(/alpha+
/beta)^2(/alpha+/beta+1)>$$
under the model, so the natural moment estimators are
$$/wt/alpha=/lb/hm(1-/hm)/sla/hssg^2-1/rb/,/hm,/ /wt/beta=/lb/hm(1-/hm)/sla
/hssg^2-1/rb/,(1-/hm),$$
where $/hm=</overline X>$ and $/hssg^2=<1/over n-1>/sum_<i=1>^n(X_i-</overline X>)^2$.
In particular these aim at certain ``best fit'' parameter values
$$/alpha_1=/lb/mu(1-/mu)/sla/ssg^2-1/rb/,/mu,/ /beta_1=/lb/mu(1-/mu)/sla/ssg^2
-1/rb(1-/mu),$$
writing $/mu$ and $/ssg$ for the underlying $f$'s </it true//> mean and
standard deviation.

The maximum likelihood estimators $/wh/alpha$, $/wh/beta$, on the other
hand, aim at other parameter values $/alpha_0$ and $/beta_0$, determined
from
$$/int_0^1/log x/,f(x)/,dx=/psi(/alpha_0)-/psi(/alpha_0+/beta_0),$$
$$/int_0^1/log (1-x)/,f(x)/,dx=/psi(/beta_0)-/psi(/alpha_0+/beta_0).$$
Only in exceptional cases (including the idealised one)
will $(/alpha_0,/beta_0)$ and $(/alpha_1,/beta_1)$ agree.
/Enote
The results for these three examples could have been obtained directly, but it is 
important
to note that the general considerations provide information about the behaviour
of ML estimators also when more complicated parametric structures
are tried on a set of data, for example the mixture models studied in 
3.2.E and in Chapter 7.
/section<3.2><Examples><>
This section considers some parametric models that may be called upon
to describe different types of data sets.
/subsection<3.2.A><The multi-exponential model>
Suppose the vector $X=(X_1/upto X_d)'$ has independent and 
exponentially distributed components, that is a density of the form
$$f(x,</underline/th>)=/th_1/,e^<-/th_1x_1>/ldots/th_d e^<-/th_d x_d>,
/ x_1/upto x_d/gt 0./eqno(3.11)$$
If a set $X_1/upto X_n$ of observed vectors is available whose density $f(x)$
is modelled by (3.11), then
$$/wh/th_1=/left(<1/over n>/sum_<i=1>^n X_<i,1>/right)^<-1>/upto /wh/th_d=
/left(<1/over n>/sum_<i=1>^n X_<i,d>/right)^<-1>/eqno(3.12)$$
are the ML estimators, writing $X_i=(X_<i,1>/upto X_<i,d>)'$.
(Improvements are possible if $d/geq 3$, see Berger (1980b).)

If $K$ populations have densities of this type, say with parameters
$$</underline /th>_k=(/th_<k,1>/upto /th_<k,d>)',/ k=1/upto K,$$
then the ML plug-in rule would be $/widehat C_0(x_1/upto x_d)$
defined as in (3.4), (3.5), with
$$f_k/lp x,/wh/th/rp=/prod_<i=1>^d /wh/th_<k,i>/,
/exp/lp-/wh/th_<k,i>/,x_i/rp./eqno(3.13)$$
$/wh/th_<k,1>/upto /wh/th_<k,d>$ are as in (3.12), based on a training set 
for class $k$. A candidate vector $x=(x_1/upto x_d)$ would be assigned to the class with the highest value of
$$/log/,/pi_k+/sum_<i=1>^d/lp/log /wh/th_<k,i>-/wh/th_<k,i>/,x_i/rp,$$
unless this highest value is below
$$/log (1-c)+/log/left/lbrace/sum_<t=1>^K /pi_t/,f_t/lp x,/wh/th/rp/right/rbrace,$$
in which case the doubt option should be used.
/subsection<3.2.B><The multinormal model and the best quadratic rule>
Let $X$ be multinormally
distributed in $d$-space, i.e. with density 
$$f(x;/mu,/Sigma)=(2/pi)^<-d/sla 2>/vert/Sigma/vert^<-1/sla 2></rm exp>/lb-
<1/over 2>(x-/mu)'/Sigma^<-1>(x-/mu)/rb/eqno(3.14)$$                       
for a vector $/mu$ and a positive definite symmetric matrix $/Sigma$. The ML 
estimators of $/mu$ and $/Sigma$ based on a random sample
$X_1/upto X_n$ from a distribution modelled by (3.14) are well known to be
$$/wh/mu=<1/over n>/sum_<j=1>^nX_j,/eqno(3.15)$$
$$/widehat/Sigma=<1/over n>/sum_<j=1>^n(X_j-/wh/mu)(X_j-/wh/mu)',/eqno(3.16)$$
see for example Anderson (1958, p.~47).

Assume now that $K$ populations are described in this manner, that is 
$$X/sim N_d(/mu_k,/Sigma_k)</rm/ if/ >X</rm/ is/ from/ class/ >k./eqno(3.17)$$
Assume further that a training set is available for each class, as in (3.3), so
that ML estimates for the individual class descriptors $/mu_k$, $/Sigma_k$
may be obtained:
$$/wh/mu_k=<1/over n_k>/sum_<j=1>^<n_k>X_j^<(k)>,/eqno(3.18)$$
$$/wh/Sigma_k=<1/over n_k>/sum_<j=1>^<n_k>(X_j^<(k)>-/wh/mu_k)
(X_j^<(k)>-/wh/mu_k)'./eqno(3.19)$$ 
The ML plug-in discriminant rule becomes as in (3.4), with estimated
posterior probabilities being proportional to 
$$/pi_k(2/pi)^<-d/sla 2>/vert/wh/Sigma_k/vert^<-1/sla 2>/exp/left/lbrace-
<1/over 2>(x-/wh/mu_k)'/wh/Sigma_k^<-1>(x,/wh/mu_k)/right/rbrace.$$
This is known as the </it best quadratic rule//> of discriminant 
analysis.
/subsection<3.2.C><Common covariance matrix: the best linear rule>
An important variation on (3.17) is
$$X/sim N_d(/mu_k,/Sigma)</rm/ if/ >X</rm/ is/ from/ class/ >k,/eqno(3.20)$$
i.e.~the $K$ populations are assumed to have a common covariance matrix 
$/Sigma$, so that the classes differ only in mean vectors.
The ML estimators $/wh/mu_1/upto/wh/mu_K,/wh/Sigma$ are needed on the basis
of the training set (3.3). One may show that the total likelihood
$$/prod_<k=1>^K/prod_<j=1>^<n_k>(2/pi)^<-d/sla 2>/vert/Sigma/vert^<-1/sla 2>
/,</rm exp>/,/lb-<1/over 2>/lp X_j^<(k)>-/mu_k/rp'
/Sigma^<-1>/lp X_j^<(k)>-/mu_k/rp/rb$$
is maximised by $/wh/mu_k$ as in (3.18), $k=1/upto K$, and
$$/eqalign</wh/Sigma&=/sum_<k=1>^K(n_k/sla n)/wh/Sigma_k/cr
&=<1/over N>/sum_<k=1>^K/sum_<j=1>^<n_k>
/lp X_j^<(k)>-/wh/mu_k/rp/lp X_j^<(k)>-/wh/mu_k/rp',/cr>/eqno(3.21)$$
where $/wh/Sigma_k$ is as in (3.19) and $N=/sum_<k=1>^K n_k$ is the
total size of the training set.

The plug-in rule becomes: for a new $X$, find the nearest mean $/hm_k$ in 
terms of the Mahalanobis distance $(X-/hm_k)'/hs^<-1>(X-/hm_k)$. This is 
equivalent to finding the maximum of 
$2/hm_k'/hs^<-1>X-/hm_k'/hs^<-1>/hm_k$. The rule is known as the </it 
best linear rule//> of discriminant analysis.

Whereas $/wh/mu_k$ in (3.18) is unbiased, $/wh/Sigma_k$ and $/wh/Sigma$ are
not. Unbiased covariance matrix estimators can easily be obtained, however:
$$/eqalign<
/ts_k&=<1/over n_k-1>/sum_<j=1>^<n_k>/lp 
X_j^<(k)>-/hm_k/rp/lp/xkj-/hm_k/rp',/cr
/ts&=/sum_<k=1>^K<n_k-1/over N-K>/,/ts_k/cr
&=<1/over N-K>/sum_<k=1>^K/sum_<j=1>^<n_k>
/lp/xkj-/hm_k/rp/lp/xkj-/hm_k/rp'/cr>$$
possess this property.  Often plug-in discriminant procedures are used
with $/tsk$ and $/ts$ replacing $/hsk$ and $/hs$.

There are no compelling reasons for choosing say $/tsk$ instead of $/hsk$,
however.  The unbiasedness property of $/tsk$ as an estimator of $/sg_k$ is 
rather irrelevant in the classification context.  It would perhaps be 
better to find an unbiased estimator $/wh f_k(x)$ of the density $f_k(x)$ 
itself. One can prove that
$$/wh f_k(x)=/left(<1/over (n-1)/pi>/right)^<d/sla 2>
</Gamma/left(<n-1/over 2>/right)/over /Gamma/left(<n-d-1/over 2>/right)>
/vert/hs/vert^<-1/sla 2>
/left/{ 1-<1/over n-1>(x-/hm)'/,/hs^<-1>(x-/hm)/right/}_<+>^<(n-d-3)/sla 2>
/eqno(3.22)$$
has this property.  Thus several versions of the classification
procedure (3.4) necessarily exist once a particular parametric
model is chosen.  It will usually be difficult to tell in advance
which of the procedures will achieve the lowest error rates etc.,
though the procedures that use unbiased estimators of the densities
themselves seem to have a slight edge on the direct plug-in rules
in that the latter ones do not reflect sample variation in 
estimators, whereas the former ones in a sense do,  cf. (3.22).

Chapter 4 takes up this problem, and presents some interesting competitors
to plug-in rules.
/subsection<3.2.D><A mixed model for discrete and continuous components>
Consider a vector $X=(A,Y_1/upto Y_d)$,
where $A$ is a $0-1$ variable whereas $Y=(Y_1/upto Y_d)'$ is continuously
distributed.  Feature vectors of this type have proven to be
valuable in symbol recognition applications, where $A$ may indicate
the presence or absence of a certain characteristic.  A possible
statistical description of the distribution of $X$ is
the following:
$$Y/mid(A=a)/sim N_d(/mu_a,/sg_a),a=0,1,/eqno(3.23)$$
$$Pr/left/{A=0/right/}=1-p, Pr/left/{ A=1/right/}=p./eqno(3.24)$$
This may be expressed as
$$/eqalign<
f(x)&=f(a,y)=/cases<(1-p)N_d(/mu_0,/sg_0)(y)&if $a=0$/cr
p/,N_d(/mu_1,/sg_1)(y)&if $a=1$/cr>/cr
&=(1-p)^<1-a>p^a N_d(/mu_a,/sg_a)(y);/cr>/eqno(3.25)$$
$f$ is the density w.r.t. the measure $/lam/times m_d$ on 
$/lbrace 0,1/rbrace/times/rr^d$, where $/lam$ is
counting measure and $m_d$ Lebesgue measure on $/rr^d$.

Let $X_1/upto X_n$ be a random sample from a distribution modelled as
above.  The ML estimators $/wh p,/hm_0,/hs_0,/hm_1,/hs_1$ are needed.  
We have by (3.25)
$$/log f(a,y)=a/,/log p+(1-a)/log(1-p)
+<d/over 2>/log(2/pi)+<1/over 2>/log/vert/Lambda_a/vert-/halv(y-/mu_a)'/Lambda_a(y-/mu_a),$$
where $/Lambda_a=/sg_a^<-1>$. The ML estimators will be found as the solution
to the equations
$$/eqalign<
/sum_<j=1>^n</partial/log f(A_j,Y_j)/over /partial p>&=0/cr
/sum_<j=1>^n</partial/log f(A_j,Y_j)/over /partial/mu_<a,i>>&=0,/ i=1/upto d,/,a=0,1,/cr
/sum_<j=1>^n</partial/log f(A_j,Y_j)/over /partial/lam_<a,il>>&=
0,/ i,l=1/upto d,/,a=0,1,/cr>$$
writing $X_i=(A_i,Y_i),/ /mu_a=(/mu_<a,1>/upto/mu_<a,d>)'$,
$/Lambda_a=/lbrace/lam_<a,il>/rbrace,/,/sg_a=/lbrace/ssg_<a,il>/rbrace$.

Generally we have 
$$/eqalignno<</partial/over /partial/mu_i>(x-/mu)'/Lambda(x-/mu)&=-2/Lambda_<(i)>(x-/mu),&(3.26)/cr
</partial/over/partial/lam_<ii>>(x-/mu)'/Lambda(x-/mu)&=(x_i-/mu_i)^2,/cr
</partial/over/partial/lam_<il>>(x-/mu)'/Lambda(x-/mu)&=
2(x_i-/mu_i)(x_l-/mu_l),l/not=i,&(3.27)/cr
</partial/over/partial/lam_<ii>>/log/vert/Lambda/vert&=(/Lambda^<-1>)_<ii>=
/ssg_<ii>,/cr
</partial/over/partial/lam_<il>>/log/vert/Lambda/vert&=2(/Lambda^<-1>)_<il>=
2/ssg_<il>,l/not=i,&(3.28)/cr>$$
where $/Lambda_<(i)>=(/lam_<i1>/upto /lam_<id>)$ is the $i$'th
row of $/Lambda$ and $/sg=/Lambda^<-1>$. Using these results, we get
$$/eqalign<
/tmpa</partial p>&=<A_j-p/over p(1-p)>,/cr
/tmpa</partial/mu_<a,i>>&=/Lambda_<a,(i)>(Y_j-/mu_a) I/lbrace A_j=a/rbrace,/cr
/tmpa</partial/lam_<a,ii>>&=/left/lbrack/halv/ssg_<a,ii>-
/halv(Y_<j,i>-/mu_<a,i>)^2/right/rbrack I/lbrace A_j=a/rbrace,/cr
/tmpa</partial/lam_<a,il>>&=/left/lbrack/ssg_<a,il>-(Y_<j,i>-/mu_<a,i>)(Y_<j,l>-/mu_<a,l>)
/right/rbrack I/lbrace A_j=a/rbrace./cr>$$
The ML-defining equations simplify to give
$$/wh p=<1/over n>/sum_<j=1>^n A_j=<N_1/over n>,/eqno(3.29)$$
$$/eqalign<
/Lambda_a/sum_<j:A_j=a>(Y_j-/mu_a)&=0,/,a=0,1,/cr
</rm i.e.>/qquad/hm_a=<1/over N_0>/sum_<j:A_j=0>Y_j,/ /hm_1&=<1/over N_1>
/sum_<j:A_j=1> Y_j,/cr>/eqno(3.30)$$
$$/eqalign<
/sum_<j=1>^n &/left/lbrace/sg_a-(Y_j-/mu_a)(Y_j-/mu_a)'/right/rbrace/,
I/lbrace A_j=a/rbrace=0,/cr
</rm i.e.>/qquad/hs_a&=<1/over N_a>/sum_<j:A_j=a>(Y_j-/hm_a)(Y_j-/hm_a)',
/,a=0,1./cr>/eqno(3.31)$$
Here $N_1=/sum_<j=1>^n A_j$ is the observed number of 
feature vectors with $A_j=1$ and $N_0=n-N_1$ is the number of vectors
with $A_j=0$. One may show that $/wh p,/hm_0,/hs_0,/hm,/hs$ above really 
maximise the likelihood.

Suppose now that $K$ populations give rise to $K$ different densities
of the type above, i.e.
$$f_k(x)=f(a,y;p_k,/mu_<k,0>,/sg_<k,0>,/mu_<k,1>,/sg_<k,1>)./eqno(3.32)$$
If a training set of the form (3.3) is available a ML
plug in classification procedure of the general type (3.4),(3.5) may
be constructed, with estimates
$$/wh p_k,/hm_<k,0>,/hs_<k,0>,/hm_<k,1>,/hs_<k,1>$$
obtained for each class as in (3.29)--(3.31) and then plugged in.

The approach outlined above involves a lot of parameters,
indeed $1+2(d+/halv d(d+1))=1+3d+d^2$ for each class.  Relatively large 
training sets are required to estimate all of them with sufficient precision.
(Just how large depends on the separability between classes
and also the dimension $d$.)  Some variations on the full model (3.32) are
possible, reducing the number of parameters.
/note<Reduction 1> puts $/sg_<k,0>=/sg_<k,1>$ in (3.32):
$$f_k(x)=f(a,y;p_k,/mu_<k,0>,/sg_k,/mu_<k,1>,/sg_k)./eqno(3.33)$$
Accordingly ML estimators $/wh p,/hm_0,/hm_1,/hs$ are needed in the model 
(3.23)--(3.24) when $/sg_0=/sg_1=/sg$
is common.  Going through arguments and algebraic manipulations
similar to those that led to (3.29)--(3.31) we end up
with $/wh p,/hm_0,/hm_1$ as before and 
$$/hs=<N_0/over n>/hs_0+<N_1/over n>/hs_1./eqno(3.34)$$
The number of parameters is now $1+<5/over 2>d+/halv d^2$ for each class.
/note<Reduction 2> assumes that the $K/;/sg_<k,0>$ matrices may be taken to be 
the same $/sg_0$ and similarly puts $/sg_<1,1>=/cdots=/sg_<k,1>$ equal to 
a common $/sg_1$:
$$f_k(x)=f(a,y;p_k,/mu_<k,0>,/sg_0,/mu_<k,1>,/sg_1)./eqno(3.35)$$
The total number of parameters is now $K(1+2d)+d(d+1)$, which may be
substantially less than the $K(1+<5/over 2>d+/halv d^2)$ involved in 
Reduction 1.

ML estimators $/wh p_k,/hm_<k,0>,/hm_<k,1>,/hs_0,/hs_1$ are needed, on the 
basis of a training set of the form (3.3). Using the general identity
$$/sum_<j=1>^n(X_j-/mu)'/Lambda(X_j-/mu)=Tr(/Lambda B)+n(/mu-/hm)'/Lambda(/mu-/hm),$$
where $B=/sum_<j=1>^n (X_j-/hm)(X_j-/hm)'$ and $/hm=<1/over n>/sum_<j=1>^n 
X_j$, we get by (3.25) the following expressions for the likelihood of the data
from class $k$, $X_j^<(k)>=(A_j^<(k)>,Y_j^<(k)>),/ j=1/upto n_k:$
$$/prod_<j=1>^<n_k>/lsb p_k^<A_j^<(k)>>(1-p_k)^<1-A_j^<(k)>>/vert
   /sg_<A_j^<(k)>>/vert^<-1/sla 2> /exp/lb-/halv/lp Y_j^<(k)>-/mu_<k,A_j^<(k)>>/rp'
   /sg_<A_j^<(k)>>^<-1>/lp Y_j^<(k)>-/mu_<k,A_j^<(k)>>/rp/rb/rsb$$
$$=p_k^<N_1^<(k)>>(1-p_k)^<N_0^<(k)>>/vert/Lambda_0/vert^<N_0^<(k)>/sla 2>
  /exp -/halv/lsb/Tr/lp/Lambda_0 B_0^<(k)>/rp+N_0^<(k)>(/mu_<k,0>-/hm_<k,0>)'
  /Lambda_0(/mu_<k,0>-/hm_<k,0>)/rsb$$
$$/qquad/vert/Lambda_1/vert^<N_1^<(k)>/sla 2>/exp-/halv/lsb/Tr
   /lp/Lambda_1 B_1^<(k)>/rp+N_1^<(k)>(/mu_<k,1>-/hm_<k,1>)'
   /Lambda_1(/mu_<k,1>-/hm_<k,1>)/rsb,$$
where $N_1^<(k)>=/sum_<j=1>^<n_k>A_j^<(k)>,N_0^<(k)>=n_k-N_1^<(k)>,/Lambda_0=/sg_0^<-1>, /Lambda_1=/sg_1^<-1>$,
$$B_a^<(k)>=/sum_<j=1>^<n_k>(Y_j^<(k)>-/hm_<k,a>)(Y_j^<(k)>-/hm_<k,a>)'=N_a^<(k)>/hs_<k,a>,a=0,1,/eqno(3.36)$$
$$/hm_<k,0>=<1/over N_0^<(k)>>/sum_<j:A_j^<(k)>=0>Y_j^<(k)>,/hm_<k,1>=<1/over N_1^<(k)>>/sum_<j:A_j^<(k)>=1>Y_j^<(k)>./eqno(3.37)$$
The total likelihood for $Z$ of (3.3) under model (3.35) becomes, using 
$$N_a=/sum_<k=1>^K N_a^<(k)>, B_a=/sum_<k=1>^K B_a^<(k)>=/sum_<k=1>^K
N_a^<(k)> /hs_<k,a>,/ a=0,1,$$
$$/displaylines<
/prod_<k=1>^K/left/lbrace p_k^<N_1^<(k)>>(1-p_k)^<N_0^<(k)>>/right/rbrace/vert/Lambda_0/vert^<N_0/sla 2>/vert/Lambda_1/vert^<N_1/sla 2>/!
/exp -/halv/lbrace Tr(/Lambda_0 B_0)+Tr(/Lambda_1 B_1)/rbrace/cr
/qquad/exp/left/lbrace-/halv/sum_<k=1>^K/sum_<a=0>^1 N_a^<(k)>
(/mu_<k,a> -/hm_<k,a>)'/Lambda_a(/mu_<k,a>-/hm_<k,a>)/right/rbrace./cr>$$  
This entails that $/hm_<k,0>,/hm_<k,1>$ of (3.37) indeed are ML estimators
also in the present model, and similary for
$$/wh p_k=N_1^<(k)>/sla n_k,/ k=1/upto K./eqno(3.38)$$
Taking derivatives of $/halv N_a/log/vert/Lambda_a/vert
-/halv Tr(/Lambda_a B_a)$ using (3.28), or appealing to
Lemma 3.2.2 in Anderson (1958), we arrive at
$$/hs_a=<1/over N_a> B_a=/sum_<k=1>^K (N_a^<(k)>/sla N_a)/hs_<k,a>,
/ a=0,1,/eqno(3.39)$$
i.e.~$/hs_0,/hs_1$ are obtained by pooling across the $K$ classes.
/note<Reduction 3> combines the two previous reductions, and pools
both across classes and across outcomes $a=0,1$:
$$f_k(x)=f(a,y;p_k,/mu_<k,0>,/sg,/mu_<k,1>,/sg)./eqno(3.40)$$
Applying this parametric model, which has only $K(1+2d)+/halv d(d+1)$ 
parameters, may be successfull even if the $2K$ covariance
matrices involved really differ.  Often $p_k=Pr/lbrace A=1/mid C=k/rbrace$ and the mean vectors
$/mu_<k,0>,/mu_<k,1>$ provide the most important discriminating information.

By arguments similar to those presented in Reductions
1 and 2 one may show that the ML estimators $/wh p_k,/hm_<k,0>,/hm_<k,1>$
are as before,
see (3.37)--(3.38), whereas the ML estimator $/hs$ of the
common $/sg$ may be expressed in various forms, illustrating
that it is obtained through pooling across both classes
and A-outcomes:
$$/eqalign</hs&=<1/over N>/sum_<a=0>^1 /sum_<k=1>^K /sum_<j=1>^<n_k>
/twice</lp Y_j^<(k)>-/hm_<k,a>/rp>'/cr
&=<1/over N>/lp N_0/hs_0+N_1/hs_1/rp/cr
&=<1/over N>/sum_<k=1>^K /lp N_0^<(k)>/hs_<k,0>+N_1^<(k)>/hs_<k,1>/rp./cr>/eqno(3.41)$$

The resulting ML plug-in discriminant procedures are in each of the cases
considered above of the general form (3.4)--(3.5). For example, suppose that 
a candidate vector $x=(a,y)$ is observed and that the model corresponding 
to (3.40) is used. The the procedure works in he following way:
if $a=1$, find the class $k$ giving maximum value of 
$$/log/pi_k+/log/wh p_k-/halv(y-/hm_<k,1>)'/hs^<-1>(y-/hm_<k,1>),$$
and if $a=0$, find the class maximising 
$$/log/pi_k+/log(1-/wh p_k)-/halv(y-/hm_<k,0>)'/hs^<-1>(y-/hm_<k,0>).$$
In both cases the candidate vector is assigned to the class giving maximum
value, unless
$$</pi_k f(a,y;/wh p_k,/hm_<k,0>,/hs,/hm_<k,1>,/hs)/over
/sum_<t=1>^K /pi_t f(a,y;/wh p_t,/hm_<t,0>,/hs,/hm_<t,1>,/hs)>/leq 1-c$$
for every $k$, in which case the doubt reject option is used.
/subsection<3.2.E><Mixtures of normals> We have experienced several instances of important feature
calculation methods in symbol recognition that have resulted in class densities
with a clear bi-modal structure. Consequently, histograms
for even well chosen transformations of data from such  distributions are not
well described by fitting a multinormal density.
The bi-modality suggests studying mixtures of two multinormal 
distributions as a means of describing class densities.

Let $X_1/upto X_n$ be a random sample from a density which we want
to describe parametrically by 
$$/eqalign<f(x)&=f(x;p,/mu_1,/sg_1,/mu_2,/sg_2)/cr
&=(1-p)N_d(/mu_1,/sg_1)(x)+pN_d (/mu_2,/sg_2)(x)./cr>/eqno(3.42)$$
The important problem of fitting data to this class of densities, 
i.e.~obtaining good parameter estimates 
$/wh p,/hm_1,/hm_2,/hs_1,/hs_2$ on the basis of $X_1/upto X_n$,
is a difficult one and is not yet satisfactorily solved in
the literature.

The model is not properly defined until a restriction of the parameter
set is made to avoid problems of identifiability:
we may exchange $(/mu_1,/sg_1)$ and $(/mu_2,/sg_2)$ and rename $p$
as $1-p$ to get two representations of the same density.
The model is identifiable if one demands $p/in/lbrack 0,/halv/rbrack$ or 
$/mu_<1,1>/leq/mu_<2,1>$ for example, writing 
$/mu_a=(/mu_<a,1>/upto/mu_<a,d>)',a=1,2$.
One may check by drawing graphs in the one-dimensional case,
however, that curves with rather different sets of parameters may still
come close to each other, making estimation of parameters
a more confusing and difficult task than usual. The curve
(3.42) is not necessarily bi-modal, cf. Eisenberger (1964).

The maximum likelihood program does not work as smoothly as in the previous 
examples. We will nevertheless review this approach.

A maximum likelihood estimator does not exist in the usual sense;
it may be seen that the likelihood 
$$L_n=/prod_<j=1>^n f(X_j;p,/mu_1,/sg_1,/mu_2,/sg_2)$$
is unbounded, it has in fact several singularities. For
example $L_n/to/infty$ as $/mu_1=X_1$ and $/sg_1/to 0$, corresponding to the
``explanation'' $p=1-<1/over n>,X_1/sim N(/mu_1,0),X_2/upto X_n$ i.i.d.
$N(/mu_2,/sg_2)$. Clearly this is not the solution we want.
$L_n$ will however tend to have several local maxima,
and one of these corresponds to the $n$-th element in a sequence of stationary points that converge a.s.~to the true parameter values.

Write (3.42) in the form
$$f(x)=(1-p)f_1(x)+p f_2(x),$$
for convenience, and introduce
$$Q(1/mid x)=<(1-p)f_1(x)/over f(x)>,/;Q(2/mid x)=<p f_2(x)/over
f(x)>,/eqno(3.43)$$
which are the posterior probabilities for a given $x$ to belong to type
1 and 2 respectively, imagining the population
under study as a mixture of type 1 and 2.
To find local maxima, we derive/message<(Before)>
$$/eqalign<
/tmpb</partial p>&=<1/over f(x)>(-f_1(x)+f_2(x))/cr
&=-<1/over 1-p>Q(1/mid x)+<1/over p>Q(2/mid x),/cr>$$/message<Break here?>
$$/eqalign<
/tmpb</partial/mu_<1,i>>&=<1/over f(x)>(1-p)f_1(x)
        </sigma/log f_1(x)/over /sigma/mu_<1,i>>/cr
&=Q(1/mid x)/Lambda_<1,(i)>(x-/mu_1),/cr
/tmpb</partial/mu_<2,i>>&=Q(2/mid x)/Lambda_<2,(i)>(x-/mu_2),/cr
/tmpb</partial/lam_<1,il>>&=<1/over f(x)>(1-p)/,f_1(x)/,/tmpb</partial/lam_<1,il>>/cr
&=/lp1-/halv/dl_<il>/rp Q(1/mid x)/lbrace/sigma_<1,il>-(x_i-/mu_<1,i>)(x_l-/mu_<1,l>)/rbrace,/cr
/tmpb</partial/lam_<2,il>>&=/lp1-/halv/dl_<il>/rp/,Q(2/mid x)/lbrace/sigma_<2,il>-(x_i-/mu_<2,i>)(x_l-/mu_<2,l>)/rbrace,/cr
>$$/message<(After)>
writing $/Lambda_1=/sg_1^<-1>,/Lambda_2=/sg_2^<-1>,/sg_a=/lbrace/sigma_<a,il>/rbrace,
/Lambda_a=/lbrace/lam_<a,il>/rbrace,/ a=1,2$. Setting derivatives of the likelihood
$L_n$ equal to zero leads to the following equivalent set of equations:
$$/eqalign<
p&=<1/over n>/sum_<j=1>^n Q(2/mid X_j),/cr
1-p&=<1/over n>/sum_<j=1>^n Q(1/mid X_j),/cr
/mu_1&=</sum_<j=1>^n Q(1/mid X_j)X_j/over
/sum_<j=1>^n Q(1/mid X_j)>=<1/over n(1-p)>/sum_<j=1>^n Q(1/mid X_j)X_j,/cr
/mu_2&=</sum_<j=1>^n Q(2/mid X_j)X_j/over /sum_<j=1>^n Q(2/mid X_j)>=<1/over np>/sum_<j=1>^n Q(2/mid X_j)X_j,/cr
/sg_1&=<1/over n(1-p)>/sum_<j=1>^n Q(1/mid X_j)(X_j-/mu_1)(X_j-/mu_1)',/cr
/sg_2&=<1/over np>/sum_<j=1>^n Q(2/mid X_j)(X_j-/mu_2)(X_j-/mu_2)'./cr>$$
Note that the parameters $p,/mu_1,/sg_1,/mu_2,/sg_2$ are
involved in $Q(1/mid X_j), Q(2/mid X_j)$ on the right hand sides,
making these equations rather intricate.

A solution $(/wh p,/hm_1,/hs_1,/hm_2,/hs_2)$ of these equations
may be arrived at via an iterative computational procedure. Start out with an initial guess $(p^/nil,/mu_1^/nil,/sg_1^/nil,/mu_2^/nil,/sg_2^/nil)$, and
compute sequentially
$$p^<(t+1)>=<1/over n>/sum_<j=1>^n Q^<(t)>(2/mid X_j),/eqno(3.44)$$
$$/mu_a^<(t+1)>=</sum_<j=1>^n Q^<(t)>(a/mid X_j)X_j/over
/sum_<j=1>^n Q^<(t)>(a/mid X_j)>,/ a=1,2,/eqno(3.45)$$
$$/sg_a^<(t+1)>=</sum_<j=1>^n Q^<(t)>(a/mid X_j)/lp X_j-/mu_a^<(t+1)>/rp
/lp X_j-/mu_a^<(t+1)>/rp'/over /sum_<j=1>^n Q^<(t)>(a/mid X_j)>,/ a=1,2,/eqno(3.46)$$
$t=0,1,2,/ldots$ until specified convergence criteria are met,
where $Q^<(t)>(a/mid x), a=1,2$ are as in (3.43), with 
$(p^<(t)>,/mu_1^<(t)>,/sg_1^<(t)>,/mu_2^<(t)>,/sg_2^<(t)>)$.

The procedure above, with a given starting point, is easy to implement,
and will always lead to a solution of the preceding equations, i.e.~a
stationary point for the likelihood function, but not
necessarily a local maximum. The procedure may be viewed as an application
of the so-called EM-algorithm (Dempster, Laird, and Rubin, 1977).

It may be shown that the estimators $(/wh p,/hm_1,/hs_1,/hm_2,/hs_2)$ defined
above, or rather the corresponding sequence as $n$ grows, are consistent
and enjoy further good asymptotic properties, </it if//> the starting point 
$(p^/nil/upto /sg_2^/nil)$ is consistent itself. The interpretation of this
asymtotic statement is that the starting point has to be chosen reasonably 
close to the true values $(p/upto /sg_2)$. Several procedures providing 
perhaps crude, but perhaps good enough starting points may be proposed,
for example using one of the available clustering routines, asking for the 
division of $X_1/upto X_n$ in $d$-space into two clusters.

It has proven difficult to obtain a completely automatic and always reliable
algorithm solving the above problem.
An approach based on results by Ranneby (1984) appears to solve the problem in the one-dimensional case, but seems difficult to generalise to the 
multi-dimensional case. Furthermore, Hathaway (1985) showed that 
maximisation of the log-likelihood under certain constraints led to consistent
solutions.
One may have to resort to the semi-automatic approach, where
several starting points are found by inspection of histograms etc., and then
tried in the EM machine (3.44)--(3.46), checking also for the
local maximum property, and so on.
/note<Reduction 1.> The estimation problem considered above
for the general model (3.42) was seen to be a difficult one, because many 
parameters  $(1+2d+d(d+1)=1+3d+d^2)$ are involved with the frustrating property
that fairly well separated parameter points may still produce densities pretty 
close to each other.

A reduction, in terms of number of parameters and likelihood difficulties,
is provided by the restriction $/sg_1=/sg_2$, i.e.
$$f(x)=(1-p)/,N_d(/mu_1,/sg)(x)+pN_d(/mu_2,/sg)(x)./eqno(3.47)$$
The likelihood of a random sample $X_1/upto X_n$ is now free
of singularities, but may still be plagued with several local maxima.

One may now mimic the elaborations presented above in the 
$(/sg_1,/sg_2)$ case. One needs 
$$/eqalign</tmpb</sigma/lam_<il>>&=/lp1-/halv/dl_<il>/rp
/biggl/lbrace Q(1/mid x)/left(/sigma_<il>-(x_i-/mu_<1,i>)(x_l-/mu_<1,l>)/right)/cr
&/qquad+Q(2/mid x)/left(/sigma_<il>-(x_i-/mu_<2,i>)(x_l-/mu_<2,l>)/right)/biggr/rbrace,/cr>$$
writing $/Lambda=/sg^<-1>$, and ends up with equations 
$$p=<1/over n>/sum_<j=1>^n Q(2/mid X_j),/eqno(3.48)$$
$$/mu_a=</sum_<j=1>^n Q(a/mid X_j)X_j/over /sum_<j=1>^n Q(a/mid X_j)>,/,a=1,2,/eqno(3.49)$$
$$/sg=<1/over n>/sum_<j=1>^n/left/lbrace Q(1/mid X_j)(X_j-/mu_1)(X_j-/mu_1)'+Q(2/mid X_j)(X_j-/mu_2)(X_j-/mu_2)'/right/rbrace,/eqno(3.50)$$
that any local maximum has to obey. $Q(1/mid x), Q(2/mid x)$
are as in (3.43) but with $/sg_1=/sg_2$.

An EM-type iterative computational procedure with a starting point
$(p^/nil,/allowbreak/mu_1^/nil,/mu_2^/nil,$/break$/sg^/nil)$ may now be defined
analogously to (3.44)--(3.46), using this time (3.48)--(3.50) as
``iterative equations''.
/note<Reduction 2.> We have not forgotten the discriminant analysis context,
but have silently assumed that some or all of $K$ class densities
$f_1/upto f_K$ were to be parametrised by the mixture of two normals.
Model (3.42) for each class would lead to 
$$f_k(x)=(1-p_k)/,N_d(/mu_<k,1>,/sg_<k,1>)(x)+p_k/,N_d(/mu_<k,2>,/sg_<k,2>)(x)/eqno(3.51)$$
whereas the simpler model (3.47) leads to
$$f_k(x)=(1-p_k)/,N_d(/mu_<k,1>,/sg_k)(x)+p_k/,N_d(/mu_<k,2>,/sg_k)(x)./eqno(3.52)$$

A further simplification is now the restriction $/sg_1=/cdots=/sg_K=/sg$, i.e.
$$f_k(x)=(1-p_k)/,N_d(/mu_<k,1>,/sg)(x)+p_k/,N_d(/mu_<k,2>,/sg)(x)./eqno(3.53)$$
One may consider $(p_k,/mu_<k,1>,/mu_<k,2>)$ to be the most important 
information about a class, i.e.~the mixing probability and the position of the
two centres, and interpret the common covariance matrix $/sg$ as securing 
against random noise from the $(p_k,/mu_<k,1>,/mu_<k,2>)$ structure.

This firm restriction on the $2K$ covariance matrices $/sg_<k,1>,/sg_<k,2>$, 
namely that they should all agree, may actually help in the consistency 
problems encountered above.

The total likelihood for a training set of the form (3.3) is
$$L=/prod_<k=1>^K f_k/lp x_1^<(k)>/rp/cdots f_k/lp x_<n_k>^<(k)>/rp.$$
Since $p_k,/mu_<k,1>,/mu_<k,2>$ are involved in the $k$-th factor only, it is
seen that a local maximum $(/wh p_1,/hm_<1,1>,/hm_<2,1>/upto 
/wh p_K,/hm_<K,1>,/hm_<K,2>,/hs)$ of $L$ must satisfy the equations
$$p_k=<1/over n_k>/sum_<j=1>^<n_k> Q_k/lp2/mid X_j^<(k)>/rp,/eqno(3.54)$$
$$/mu_<k,1>=/sum_<j=1>^<n_k> Q_k/lp1/mid X_j^<(k)>/rp X_j^<(k)>/sla /sum_<j=1>^<n_k> Q_k/lp1/mid X_j^<(k)>/rp,/eqno(3.55)$$
$$/mu_<k,2>=/sum_<j=1>^<n_k> Q_k/lp2/mid X_j^<(k)>/rp X_j^<(k)>/sla/sum_<j=1>^<n_k> Q_k/lp2/mid X_j^<(k)>/rp,/eqno(3.56)$$
where
$$Q_k(1/mid x)=<(1-p_k)/,f_<k,1>(x)/over f_k(x)>,/;Q_k(2/mid x)= <p_k/,f_<k,2>(x)/over f_k(x)>./eqno(3.57)$$
It remains to obtain an equation involving $/sg$.
Write $/Lambda=/sg^<-1>=/lbrace/lam_<il>/rbrace$, and utilise previous efforts to get 
$$/eqalign<</partial/log L/over /partial/lam_<il>>&=
/sum_<k=1>^K /sum_<j=1>^<n_k>/lp1-/halv/dl_<il>/rp/biggl/lbrace Q_k/lp1/mid X_j^<(k)>/rp/left/lbrack/sigma_<il>-/lp X_<j,i>^<(k)>-/mu_<k,1,i>/rp/lp X_<j,l>^<(k)>-/mu_<k,1,l>/rp/right/rbrack/cr
&/qquad+Q_k/lp2/mid X_j^<(k)>/rp/left/lbrack /sigma_<il>-/lp X_<j,i>^<(k)>-/mu_<k,2,i>/rp/lp X_<j,l>^<(k)>-/mu_<k,2,l>/rp/right/rbrack/biggr/rbrace./cr>$$
Hence
$$/eqalign</sg&=<1/over N>/sum_<k=1>^K/sum_<j=1>^<n_k>/Bigl/lbrace Q_k/lp1/mid X_j^<(k)>/rp/lp X_j^<(k)>-/mu_<k,1>/rp/lp X_j^<(k)>-/mu_<k,1>/rp'/cr
&/qquad+Q_k/lp2/mid X_j^<(k)>/rp/lp X_j^<(k)>-/mu_<k,2>/rp/lp X_j^<(k)>-/mu_<k,2>/rp'/Bigr/rbrace.
/cr>/eqno(3.58)$$

Once more an iterative computational procedure may be given that from a starting point $p_k^/nil,/mu_<k,1>^/nil,/mu_<k,2>^/nil;/ k=1/upto K,/sg^/nil$
arrives at a solution to (3.54)--(3.56), (3.58).
/section<3.3><Concluding remarks><>
Parametric plug-in versions of the optimal Bayes rule, using maximum
likelihood estimates or unbiased modifications, are perhaps the most
widely used classification procedures. That one sometimes can do better,
by taking estimation variability into account, is the topic of the
following chapter.

It is important to realise the fundamental shortcoming of
parametric procedures pointed out in Section 3.1: the parametric models
are never completely correct, and the least false version of the 
optimal classifier that the parametric methods aim at are almost
certainly not as good as the truly optimal one. Only nonparametric
methods have the favourable Bayes risk consistency property. But again,
the nonparametric methods need large and sometimes very large training
sets to be effective. This means that procedures with inconsistent 
Bayes risks, say parametric ones, very well can be best for the given
sample sizes.

It is relevant here to point out that the updating methods
presented in Chapter 7 can be utilised to make estimates more precise,
thereby helping the parametric classifier under consideration to
perform better, in situations where additional, unclassified data 
vectors are available. It is even possible to pass from an initial
level of sophistication (say the normal distribution description) to
a higher level (semi- or nonparametric models) using such methods.

The best linear and the best quadratic rule were seen to be
special cases of the general Bayes procedure of Chapter 1, and are both
optimal under appropriate conditions. Theoretical investigations and
simulation experiments show the best linear rule to be surprisingly
competitive even in situations with different covariance matrices, and
indeed often wins over the best quadratic rule, unless the training sets
are fairly large.

This fact points again to the dangers of using models with too
many parameters: estimation variability takes over. It also suggests
that methods taken from the small world of interesting methods lying
between the linear and the quadratic ones can perform well. A simple
variation of these methods could take proportional (but otherwise
unknown) covariance matrices as the modelling starting point. This
greatly reduces the number of unknown parameters. One can also
construct empirical Bayes type methods that let the data themselves 
decide whether the linear or the quadratic method is more appropriate.
As a final example of such modified quadratic methods we mention the
recent paper by Kimura, Takashina, Tsuruoka, and Miyake (1987). They
apply their method to automatic recognition of Chinese characters.

These remarks also point to the need for </it model choice criteria//>.
One can check each model considered with a goodness of fit test, cf.~Chapter
11, or use criteria like those of Akaike (1974) and Schwarz (1978); see
also the brief discussion in Section 7.4.

The mixture models considered in 3.2.D are important, since 
such distributions often arise for naturally chosen feature extraction
methods. It is inherently difficult to estimate the parameters of the
mixture, and it seems to be almost impossible to construct a completely
automatic, fool-proof method. Often some human interaction is necessary
in order to get proper starting values for the iterative schemes 
discussed in 3.2.D, for example. Even methods with guaranteed
consistency, as the method searching for the maximum log-likelihood
over a partly compactified parameter set (see e.g./ Hathaway (1985)),
have their difficulties, and the difficulties quickly increase with
increasing dimension. The book by Titterington, Smith and Makov (1985)
contains several illuminating examples.

There is a natural connection from mixture models to cluster
analysis. Anderberg (1973) is a good source of such methods. If such
analysis suggests that a certain class really consists of two disjoint 
subclasses, then the set of possible decisions for the classifier might
as well be extended to accomodate this. In cases where a class consists
of somewhat overlapping clusters a mixture model can be fitted, with
starting values suggested by the clustering.

We have used the so-called $k$-means clustering method (cf./ %
Anderberg (1973)) with reasonable success, see Br}ten, Holb{k-Hanssen and
Taxt (1986a, Chapter 7 and Section 9.7).

We mention finally a group of methods for fitting parametric 
families that are somewhat different from those discussed previously
in this chapter. These are based on finding suitable initial data
transformations that perhaps aim at symmetrising the data, or making the
data ``more normal''. Thus the popular Box-Cox method finds, for
each coordinate, a power transformation that makes the data come closer
to normality. In this way one-dimensional data can be fitted to a 
three-parametric family; one power parameter defining the transformation
and two to describe the normal density. While this method has much to commend
it, is not clear how well the several possible $d$-dimensional
generalisations of it perform.

There are links from Chapter~3 to other chapters. Chapter~4
discusses methods of taking  estimation variability into consideration for
parametric classification procedures. Chapter~7 provides updating
methods for making parameter estimates more precise for a given
application, and that also can be used, in some cases, to pass from an
initial parametric model into a more ambitious semi- or nonparametric
method, and such models are discussed in Chapter 5. Section 3.2.D
treated the simplest one among a large group of models for mixed
discrete and continuous data (see also Section 4.2.E). Chapter 9 expands
on this theme. It is sometimes important to use estimators that are
non-sensitive to extreme data points, i.e./ robust estimators.
Some such for the mean vector and the covariance matrix are presented in
Chapter~10, along with proposals for what could constitute ``standard output''
after choosing a feature extraction method. It appears important, for
example in view of the general considerations of Section 3.1, to have
checking methods available for at least the most popular parametric
families, like the normal distribution, see Chapter~11 for discussion
and some goodness-of-fit tests. Finally Chapter~12 is concerned with
estimating error rates for a given classification procedure, and in
particular quick and precise methods are provided for some of the rules
considered in this chapter.


/chapter<4><Classification based on parametric models:/hfil/break
The predictive approach>
/section<4.1><Theory><>
A feature vector $X$ is considered in its space $/cx$.
$X$ is drawn from  the density $/fkx$ if it originates from class
$k;/ktok$. These class densities are modelled parametrically as
$$/fkx=f_k(x,/th),/ktok,/eqno(4.1)$$
where $/th$ belongs to a region $/Th$ in $/rr^p$.
$/th$ may be thought of as a collection of component parameters, some of
which are common to all classes and some of which are class-specific.
$X$ is to be classified, and 
%
%
the measure of loss is
$$L/lsb(k,/th),t/rsb=/cases<
0&if $t=k$ (correct decision)/cr
1&if $t/not=k$ and $t/in/lb1/upto K/rb$ (wrong decision)/cr
c&if $t=D$ (being in doubt)/cr>/eqno(4.2)$$
for $k/in/Omega=/lb1/upto K/rb,/,/th/in/Th,/,t/in T=/lb1/upto K,D/rb$.
The prior probabilities are $/pi_k=/Pr/lb C=k/rb,/ktok$.

This was the problem considered in the previous chapter, where it was dealt 
with in a somewhat </it ad hoc//> manner, by simply plugging in estimates
for unknown parameters.
The loss function implicitly involved in Chapter 3 was indeed (4.2), which 
differs from (1.3) only by the presence of the unknown 
$/th$. The methods of that chapter aimed however directly
at approximating the rules that were known to be optimal for the known 
parameters case, and in a 
%
%
way the original measure of consequences (4.2) was not kept in sight.

It is clear that there is small room for improvement over the estimative
(or plug-in) approach when large training sets, compared with the complexity of the statistical model, are available. Then the estimates being inserted are 
very close to the true (least false) values, and the corresponding classification rules
closely approximate their ideals. If the training sets are only small or 
moderate, however, and the number of parameters to be estimated is high, then 
the approach described below may very well give a significant
%
%
improvement.

Assume that a training set of the form $Z$ in (3.3) is available.
Consider a general classification rule $/ca=/ca(x;/,z):/cx/times/cx^N
/to/lb1/upto K,D/rb,/,N=/sum_<k=1>^Kn_k$ being the total size of the
training set. Its </it risk function//> is
$$/eqalignno<
R(/ca,k,/th)&=E_<k,/th>/,L/lsb(k,/th),/ca(X;/,Z)/rsb/cr
&=P_<k,/th>/lb/ca(X;/,Z)/in/lb1/upto k-1,k+1/upto K/rb/rb+
	cP_<k,/th>/lb/ca(X;/,Z)=D/rb/cr
&=/pmc(k,/th)+c/,/pd(k,/th),/,k/in/Omega,/,/th/in/Th,&(4.3)/cr>$$
realising that the misclassification and doubt probabilities depend upon
the unknown $/th$, too.
$P_<k,/th>$ and $E_<k,/th>$ are used to signalise probabilities
%
%
and expectations when $X$ has class $C=k$ and the parameter value is $/th$.
$k$ may be averaged over as in (1.7), giving risk (expected loss) as a 
function of $/th$ only:
$$/eqalignno<
R(/ca,/th)&=E_/th L/lsb(C,/th),/ca(X;/,Z)/rsb/cr
&=E_/th/lb E_/th L/lsb(C,/th),/ca(X;/,Z)/rsb/mid C/rb/cr
&=/sum_<k=1>^K/pi_k R(/ca,k,/th)/cr
&=/ave<</it k>>/pmc(k,/th)+c/ave<</it k>>/pd(k,/th).&(4.4)/cr>$$

The theory of Chapter 1 provides for each given $/th$ a rule that minimises
this expression. Since (4.4) cannot be minimised
%
%
uniformly in $/th$ an overall criterion is needed.
One possibility is to average over $/th$, using a non-negative weight 
function $/lam(/th)$:
$$/eqalignno<
R(/ca,/lam)&=/int_/Th R(/ca,/th)/,/lam(/th)/,d/th/cr
&=/sum_<k=1>^K/pi_k/,/lsb/int_/Th/pmc(k,/th)/lam(/th)d/th+
	c/int_/Th /pd(k,/th)/lam(/th)d/th/rsb/cr
&=/ave</th>/ave<</it k>>/pmc(k,/th)+c/!/ave</th>/ave<</it k>>/pd(k,/th).&(4.5)/cr>$$
$/lam(/th)$ could for example be constant, or could be chosen large in regions
where it is important to ensure low values of the risk function
$R(/ca,/th)$.

If $/int/lam(/th)d/th/lt/infty$ we may as well assume that the weight function
is scaled so as to give $/int/lam(/th)d/th=1$. In such cases $/lam$
is often referred to as the </it prior distribution//>
%
%
or </it prior density//> for $/th$. This is the Bayesian formulation
of the problem, where $/lam$ is interpreted as a quantification of the
statistician's uncertanity regarding the unknown $/th$. However, $/lam$ may 
simply be
taken as a weight function, and  the minimisation
of (4.5) will be a reasonable criterion, regardless of faith.

It is indeed often possible to construct a unique, non-subjective 
$/lam$-function appropriate for a given parametric model, namely the
</it canonical non-informative prior density//>, see
Box and Tiao (1973) and Berger (1980a) for full discussions of this concept.

The following theorem provides the
%
%
Bayes solution to the classification problem, i.e.~the
rule that minimises (4.5), when $/lam$ is a proper density. Note that (4.5) 
then may be interpreted as ``total expected loss'',
$$R(/ca,/lam)=</it EL//>/lsb(C,/th),/ca(X;/,Z)/rsb,$$
where the four elements $/th,Z,C,X$ are all stochastic. Observe also that 
$/pmc(k,/th)$, $/pd(k,/th)$ in (4.3), (4.5) are evaluated w.r.t.~the
simultaneous distribution of $(Z,X)$, and that they of course depend upon 
the particular procedure $/ca$.
%
%

/Btheorem<> Assume that $/lam$ in (4.5) is a proper density.
Let $/lam(/th/mid z)$ be the posterior density for $/th$ given
$Z=z$ and define the </rm predictive class densities>
$$g_k(x,z)=/int_/Th f_k(x,/th)/lam(/th/mid z)d/th,/ktok./eqno(4.6)$$
Then the rule
$$/wh C(x;/,z)=/cases<
D&if every $P(k/mid x,z)/leq 1-c$/cr
k&if $P(k/mid x,z)=/mo<t/leq K>P(t/mid x,z)/gt 1-c$/cr>/eqno(4.7)$$
minimises the total risk (4.5), where
$$P(k/mid x,z)=/pi_kg_k(x,z)/sla/sum_<t=1>^K/pi_tg_t(x,z)./eqno(4.8)$$  
/Etheorem
/Bproof
$(/th,Z,C,X)$ have simultaneous density
$$/lam(/th)h(z,/th)/pi_kf_k(x,/th)$$
w.r.t.~the appropriate dominating product measure
in $/Th/times/cx^N/times/Omega/times/cx$, where
%
%
$$h(z,/th)=/prod_<k=1>^K/prod_<j=1>^<n_k>f_k/lp x_j^/k,/th/rp/eqno(4.9)$$
is the density of $Z$ given $/th$. The total risk
$$/eqalign<R(C^/ast,/lam)&=</it EL//>/lsb(C,/th),C^/ast(X;/,Z)/rsb/cr
&=/sum_<k=1>^K/pi_k/int_/Th/int_/cx/int_</cx^N>L/lsb(k,/th),C^/ast
	(x;/,z)/rsb h(z,/th)/,f_k(x,/th)/lam(/th)/,dz/,dx/,d/th/cr>$$
may now be evaluated by conditioning on the observable data $(X,Z)$.

Write
$$/lam(/th)/,h(z,/th)=h(z)/,/lam(/th/mid z),$$
i.e.~$h(z)$ is the marginal density for $Z$. The unconditional density of 
$(X,Z)$ becomes
$$/sum_<k=1>^K/int_/Th/lam(/th)/,h(z,/th)/,/pi_k/,f_k(x,/th)/,d/th=
	/sum_<k=1>^K/pi_k/,g_k(x,z)/,h(z),$$
so that $(C,/th)$ given $(x,z)$ has distribution 
$$</pi_k/,f_k(x,/th)/,/lam(/th/mid z)/over
	/sum_<t=1>^K/pi_t/,g_t(x,z)>.$$
Thus we may write
%
%
$$/eqalign<
R(C^/ast,/lam)&=E/lb EL/lsb(C,/th),C^/ast(X;/,Z)/rsb/mid X,Z/rb/cr
&=/int_/cx/int_</cx^N>/lb/sum_<k=1>^K/int_/Th L/lsb(k,/th),
	C^/ast(x;/,z)/rsb</pi_k/,f_k(x,/th)/,/lam(/th/mid z)/over
	/sum_<t=1>^K/pi_t/,g_t(x,z)>d/th/rb/cr
&/qquad/qquad/qquad/sum_<t=1>^K/pi_t g_t(x,z)/,h(z)/,dz/,dx,/cr>$$
and it suffices to choose $C^/ast(x;/,z)$, for each given
$(x,z)$, so as to minimise the inner integral. This conditional expected
loss may be simplified to 
$$/sum_<k=1>^K L/lsb(k,/th),/,C^/ast(x;/,z)/rsb/,P(k/mid x,z)$$
since the loss function (4.2) is independent of $/th$.

That $/wh C(x;/,z)$ is the best choice follows now by arguments similar to
those needed in the simpler situation studied in Chapter 1./Eproof
/note<Remark 1.> In expression (4.5) 
%
%
$/pmc(k,/th)$ and $/pd(k,/th)$ are as defined in (4.3), i.e.~w.r.t.~the
simultaneous distribution for $(X,Z)$. One may also consider a criterion 
conceptually different from (4.5), namely the ``expected loss given earlier
experience'', i.e.
$$/eqalign<R_z(C^/ast,/lam)&=E/bigl/lbrace L/lsb(C,/th),/,C^/ast(X;/,Z)/rsb/mid Z=z/bigr/rbrace/cr
&=/sum_<k=1>^K/pi_k/int_/Th/lsb/pmc_z(k,/th)+c/,/pd_z(k,/th)/rsb
	/lam(/th/mid z)/,d/th,/cr>$$ 
which can be expressed as (4.5), but with $/lam(/th/mid z)$ being
the weight function for $/th$, and with 
$$/eqalign</pmc_z(k,/th)&=P_<k,/th>/bigl/lbrace C^/ast(X;/,z)/in
	/lb1/upto k-1,k+1/upto K/rb/bigr/rbrace,/cr 
/pd_z(k,/th)&=P_<k,/th>/bigl/lbrace C^/ast(X;/,z)=D/bigr/rbrace./cr>$$
%
%
It may appear more satisfactory to minimise $R_z(C^/ast,/lam)$,
which is
$$/ave<</it k>>/pmc_z(k,/th)+c/ave<</it k>>/pd_z(k,/th)$$
averaged according to ``the current knowledge'' $/lam(/th/mid z)/,d/th$ about
the unknown $/th$.

It is reassuring, then, to see that $/wh C$ given in the theorem actually
minimises $R_z(C^/ast,/lam)$ too, for each given $z$. We omit
the proof. $/lam(/th/mid z)$ is often a surface with a well defined
peak close to some estimate of $/th$, and with small values outside the
most probable region for $/th$.
/note<Remark 2.> An obvious criticism against $/wh C$ given in Theorem 4.1 is 
its dependency
%
%
on the initial prior density $/lam(/th)$. One may show, however, that 
different choices for $/lam$ give almost the same posterior density 
$/lam(/th/mid z)$, and hence approximately identical classification rules,
when the training sets are large.
/note<Remark 3.> The theorem assumed $/lam(/th)$
to be proper. The criticism mentioned in the previous remark leads
to the use of objective, canonical non-informative prior
densities for $/th$. These are often improper, however, i.e.~$/int/lam(/th)
/,d/th=/infty$. Nevertheless it is usually possible to write
$$/lam(/th)/,h(z,/th)=h(z)/,/lam(/th/mid z)/eqno(4.10)$$
%
%
where both $h(z)=/int h(z,/th)/,/lam(/th)/,d/th$ and the resulting
$/lam(/th/mid z)$ have finite integrals, i.e.~are proper.
Then the program of the theorem can be pursued, and usually
results in a classification procedure identical to the one
that emerges as the limit of totally proper ones, say corresponding to
a sequence of priors $/lb/lam_n(/th)/rb$ with $/int/lam_n(/th)/,d/th=1$
and $c_n/lam_n(/th)/to/lam(/th)$ for some positive constants $/lb c_n/rb$.
/Enote
The term ``predictive procedures'' has its origins in Bayesian analysis,
cf. Geisser (1982). While elements of the theory presented above have appeared
earlier in the literature,
%
%
it seems not to have been explicitly applied to discrimination (between
multinormal populations only) until Aitchison, Habbema and Kay (1977), 
who  praised it. Further studies by for example Moran and Murphy (1979) are more
cautious, but seem to indicate that the predictive approach in general can be
expected to perform better than the estimative (plug-in) approach.
The main advantage of the former is that sampling variability
of estimates are taken into account, as will be seen in the examples
below. These examples also point out, however, that the two approaches lead
to very similar rules when the training sets are large.
%
%
/section<4.2><Examples><>
The predictive rules involve the predictive class densities (4.6). By (4.10) 
these may also be found as
$$g_k(x,z)=/int_/Th f_k(x,/th)/,h(z,/th)/,/lam(/th)/,d/th/sla h(z).
	/eqno(4.11)$$
It is seen from (4.7) and (4.8) that the marginal density
$h(z)$ is unimportant; it ensures only that $/int g_k(x,z)/,dx=1$
above. If deemed convenient we are free to use
$$/bar g_k(x,z)=/int_/Th f_k(x,/th)/,h(z,/th)/,/lam(/th)/,d/th/eqno(4.12)$$
in (4.7), (4.8); that is, we do not explicitly need the posterior
density $/lam(/th/mid z)$.
/subsection<4.2.A><The exponential model> Let $X$'s from class $k$ be exponentially
distributed with parameter $/th_k;
%
%
/ktok$. A training set $Z$ of the form (3.3) is available.

The plug-in rule uses (3.4)--(3.5) with 
$$f_k/lp x,/hth/rp=/hth_k/exp/lp-/hth_kx/rp$$
and $/hth_k=(<1/over n_k>/sum_<j=1>^<n_k>X_j^/k)^<-1>$. The predictive
approach uses (4.7)--(4.8), with
$$/eqalign<
g_k(x,z)&=/int_0^/infty/th_ke^<-/th_kx>/lam(/th_k/mid z)/,d/th_k/cr
	&=E/lb f_k(x,/th)/mid z/rb./cr>$$

The canonical non-informative prior has $/th_1/upto/th_K$ independent
with the same improper density $/th^<-1>,/,/th/gt0$,
cf. Berger (1980a, Ch.~6). Hence
$$/eqalign</int_</lbrack0,/infty)^K> h(z,/th)/,/lam(/th)/,
	d/th&=/prod_<k=1>^K/int_0^/infty/th_k^<n_k>e^<-/th_ks_k>/th_k^<-1>/,
	d/th_k/cr
&=/prod_<k=1>^K/Gamma(n_k)/sla s_k^<n_k>,/cr>$$
where $s_k=/sum_<j=1>^<n_k>x_j^/k=n_k/sla/hth_k$. It follows that
%
%
$$/bar g_k(x,z)=/prod_<t/not=k>(/Gamma(n_t)/sla s_t^<n_t>)/,
	/Gamma(n_k+1)/sla(s_k+x)^<n_k+1>.$$
Cancelling out the common factor $/prod_<t=1>^K(/Gamma(n_t)/sla s_t^<n_t>)$ we
arrive at
$$/eqalignno<
g_k(x,z)&=/const n_k/,s_k^<n_k>/sla(s_k+x)^<n_k+1>/cr
	&=/const<n_k/over s_k>/lp s_k/over s_k+x/rp^<n_k+1>/cr
	&=/const/hth_k/lp1+<1/over n_k>/hth_kx/rp^<-(n_k+1)>.&(4.13)/cr>$$
Note that $g_k(x,z)/approx/const/hth_k/exp(-/hth_kx)$ when $n_k$ is
large enough. (The const.~above is in fact $1$.)

One may also start out with Gamma densities, say 
$/th_k/sim</rm/ Gamma>(/alpha_k,/beta_k)$, independently,
$/ktok$, that is 
$$/lam_k(/th_k)=(/beta_k^</alpha_k>/sla/Gamma(/alpha_k))/th_k^</alpha_k-1>
	e^<-/beta_k/th_k>,/,/th_k/gt0./eqno(4.14)$$
Then it is easily seen that $/lam_k(/th_k/mid z)$, in 
%
%
obvious notation, is Gamma$(/alpha_k+n_k,/,/beta_k+s_k)$.
Furthermore,
$$/eqalignno<
g_k(x,z)&=<(/beta_k+s_k)^</alpha_k+n_k>/over/Gamma(/alpha_k+n_k)>/,
		</Gamma(/alpha_k+n_k+1)/over
		(/beta_k+s_k+x)^</alpha_k+n_k+1>>/cr
	&=</alpha_k+n_k/over/beta_k+s_k>/,/lp/beta_k+s_k/over
		/beta_k+s_k+x/rp^</alpha_k+n_k+1>&(4.15)/cr>$$
is arrived at, using
$$EY^ae^<-bY>=</beta^/alpha/over/Gamma(/alpha)>/,</Gamma(/alpha+a)/over
	(/beta+b)^</alpha+a>>,/,/alpha+a/gt0,/,/beta+b/gt0,/eqno(4.16)$$
which holds when $Y/sim</rm Gamma>(/alpha,/beta)$.

Again it can be seen that $g_k(x,z)/approx f_k(x,/hth)$ for sufficiently large $n_k$ 
(and then $/hth_k/approx/th_k$), illustrating Remark 2.

The Gamma densities (4.14) are proper if $/alpha_k/gt0,/,/beta_k/gt0$, whereas the non-informative prior used above corresponds to $/alpha_k=/beta_k=0$. 
Remark 3 is illustrated by noting that 
%
%
the predictive class densities in (4.13), which  are the ones many 
statisticians will prefer in practice on the grounds of objectivism, can be
obtained as the limits of those in (4.15) as $/alpha_k/to0,
/beta_k/to0$. The same type of phenomenon will be observed
if we let 
$$/th_k/sim/lp/log<R/over/eps>/rp^<-1><1/over/th_k>,/,/eps/leq/th_k/leq R,$$
and afterwards let $/eps/to0,/,R/to/infty$.

These operations (thinking properly first and taking limits afterwards)
are mainly cosmetic, however. If one wants to follow the predictive program 
with a non-informative improper prior, then one might as well use
(4.10) -- (4.12) directly.
%
%
/subsection<4.2.B><Multinormal densities with known covariance matrices> Assume
$$/fkx=N_d(/mu_k,/sg_k)(x),/ktok,/eqno(4.17)$$
where the $/mu_k$'s are unknown, whereas the $/sg_k$'s
are known. The estimative aproach uses
$$/eqalignno<
f_k/lp x,/hth/rp&=N_d(/hm_k,/sg_k)(x)/cr
&=(2/pi)^<-d/sla 2>/vert/sg_k/vert^<-1/sla 2>/exp/lb-/halv(x-/hm_k)'/sg_k^<-1>
	(x-/hm_k)/rb&(4.18)/cr>$$
with $/hm_k=(1/sla n_k)/sum_<j=1>^<n_k>x_j^/k$ obtained from the
training set $z$.

To construct a predictive rule, aiming at taking variability of $/hm_k$ into 
account, we may start out with $/mu_1/upto/mu_K$ independent
and $/mu_k/sim N_d(/mu_<k,0>,A_k)$. Then it may be shown that
%
%
$$/mu_k/mid z/sim N_d/lp/hm_k-<1/over n_k>/sg_k/lp<1/over n_k>/sg_k
+A_k/rp^<-1>
/lp/hm_k-/mu_<k,0>/rp, /lp n_k/sg_k^<-1>+A_k^<-1>/rp^<-1>/rp./eqno(4.19)$$
The canonical non-informative prior distribution is the uniform
and improper one,
$$/mu_k/sim1/,d/mu_<k,1>/cdots d/mu_<k,d></rm/ in/ >/rr^d,$$
corresponding to $A_k/to/infty$ (in the sense that the eigenvalues of
$A_k$ tend to infinity).
The posterior densities are proper:
$$/mu_k/mid z/sim N_d/lp/hm_k,<1/over n_k>/sg_k/rp,/ktok./eqno(4.20)$$

The predictive class density
$$g_k(x,z)=E_z(2/pi)^<-d/sla 2>/vert/sg_k/vert^<-1/sla 2>/exp/lb-/halv
	(/mu_k-x)'/sg_k^<-1>(/mu_k-x)/rb$$
is needed, where $E_z$ denotes expectation w.r.t.~the
posterior distribution (4.20). Under
%
%
this distribution 
$$(/mu_k-x)'/lp<1/over n_k>/sg_k/rp^<-1>(/mu_k-x)/sim/chi_d^2/lp n_k/dl_k(x)^2/rp,$$
i.e.~a non-central $/chi_d^2$ distribution with eccentricity parameter
$$n_k/dl_k(x)^2=n_k(x-/hm_k)'/sg_k^<-1>(x-/hm_k)./eqno(4.21)$$
Hence
$$/eqalignno<
g_k(x,z)&=(2/pi)^<-d/sla 2>/vert/sg_k/vert^<-1/sla 2>E/,/exp
		/lb-/halv/,<1/over n_k>/chi_d^2(n_k/dl_k(x)^2)/rb/cr
	&=(2/pi)^<-d/sla 2>/vert/sg_k/vert^<-1/sla 2>
		/lp n_k/over n_k+1/rp^<d/sla 2>/exp/lb-/halv/,<n_k/over n_k+1>
		/,/dl_k(x)^2/rb,&(4.22)/cr>$$
compare (4.18).

We used the following identity in order to arrive at (4.22):
$$E/,/exp/lp-s/chi_d^2/lp/dl^2/rp/rp=/lp1+2s/rp^<-d/sla 2>/exp/lp-<s/dl^2/over 1+2s>/rp,
	/eqno(4.23)$$
which can be proved utilising the fact that if $J/sim</rm Poisson>(/halv/dl^2)$ and 
$W/mid(J=j)/sim/chi_<d+2j>^2$, then $W/sim/chi_d^2(/dl^2)$.
%
%
/subsection<4.2.C><Unknown covariance matrices> We reconsider the model (4.17), but assume now that also the
$/sg_k$'s are unknown. Then $/hm_k, /hs_k$ as in (3.18), (3.19) are the
estimates to plug in if the ML-estimative approach is taken.

In order to pursue the predictive program in this situation
$$g_k(x,z)=/int/int f_k(x,/mu_k,/sg_k)/,/lam(/mu_k,/sg_k/mid z)/,d/mu_k/,d/sg_k$$
is needed, where $/lam(/mu_k,/sg_k/mid z)$ is the posterior distribution of 
$/mu_k,/sg_k$ given the training set. Above $/mu_k$ is in $/rr^d$ while $/sg_k$
runs through the set $M^d$ of all $d/times d$ symmetric positive definite 
matrices.

Since $(/mu_k,/sg_k)$ is relevant only for
%
%
$/fkx$, only the $x_j^/k$'s from class $k$ in the total training
set $z$ have any influence in $/lam(/mu_k,/sg_k/mid z)$, and we
might as well solve the following one-class problem:

Let $X_1/upto X_n$ be i.i.d.~$N_d(/mu,/sg)$, and let
$(/mu,/sg)$ have a prior distribution $/lam(/mu,/sg)/,d/mu/,d/sg$ in 
$/rr^d/times M^d$. Find the posterior distribution $/lam(/mu,/sg/mid/xtxn)$
and then the predictive density
$$g(x,z)=E_z(2/pi)^<-d/sla 2>/vert/sg/vert^<-1/sla 2>/exp
	/lb-/halv(/mu-x)'/sg^<-1>(/mu-x)/rb,$$
where $z=(/xtxn)$ and $E_z$ again denotes expectation w.r.t.~the
posterior distribution.
%
%

There is no general agreement on what the ``canonical'' non-informative prior
distribution for $(/mu,/sg)$ should be in this multi-parameter situation.
The left and right invariant Haar measures on the appropriate group
of transformations are both of the form
$$/lam_0(/mu,/sg)=/vert/sg/vert^<-a_0/sla 2>/eqno(4.24)$$
but with different $a_0$'s. Berger (1980a, Ch.~6) and Box and Tiao (1973, Ch.~8)
seem to prefer respectively the right and left invariant Haar measure.
There have also been proposals of other values for $a_0$ in the literature, so
it will be useful to do the calculations for this
%
%
general constant $a_0$ (depending upon the dimension $d$).

The manipulations needed become easier when $/Lam=/sg^<-1>$  is used as 
parameter. $(/mu,/Lam)$ has prior density
$$/eqalignno<
/lam(/mu,/Lam)&=/lam_0(/mu,/sg(/Lam))
		/vert</partial/sg/over/partial/Lam>/vert/cr
	      &=/vert/Lam/vert^<a_0/sla 2>/,/vert/Lam/vert^<-(d+1)>/cr
	      &=/vert/Lam/vert^<-a/sla 2>&(4.25)/cr>$$
for $a=2(d+1)-a_0$. We have used an expression for 
$/vert/partial/sg/sla/partial/Lam/vert$ available in Box and Tiao (1973, A8).

The density for $(X_1/upto X_n)$ given $/mu,/Lam$ may be written
$$(2/pi)^<-nd/sla 2>/vert/Lam/vert^<n/sla 2>/exp/lb-/halv/lsb Tr(/Lam A)+n(/mu-/hm)'/Lam(/mu-/hm)/rsb/rb,$$
where $A=/sum_<j=1>^n/twice<(X_j-/hm)>'=n/hs$. Hence
%
%
$(X_1/upto X_n,/mu,/Lam)$ has density 
$$(2/pi)^<-nd/sla 2>/vert/Lam/vert^<(n-a)/sla 2>/exp/lb-/halv/lsb
	Tr(/Lam A)+n(/mu-/hm)'/Lam(/mu-/hm)/rsb/rb.$$
Successive integrations give that $(X_1/upto X_n,/Lam)/sim$
$$(2/pi)^<-(n-1)d/sla 2>/vert/Lam/vert^<(n-a-1)/sla 2>n^<-d/sla 2>e^<-/halv 
	Tr(/Lam A)>$$
and that $(X_1/upto X_n)/sim$
$$(2/pi)^<-(n-1)d/sla 2>n^<-d/sla 2>/,/vert A/vert^<-/halv(n-a+d)>/,2^</halv 
	d(n-a+d)>/Gamma_d/lp n-a+d/over 2/rp,/eqno(4.26)$$
where $/Gamma_d(x)=/pi^<d(d-1)/sla 4>/,/prod_<j=1>^d
/Gamma/bigl(x-/halv(d-j)/bigr)$, cf. for example Box and Tiao (1973, 8.2.21).
(4.26) is valid when $n/gt a+2d-1$.

It follows that $(/mu,/Lam)$ has posterior distribution proportional to
$$(2/pi)^<-d/sla 2>/vert n/Lam/vert^<1/sla 2>/exp/lb-/halv(/mu-/hm)'(n/Lam)
(/mu-/hm)/rb/vert/Lam/vert^<(n-a-1)/sla 2>e^<-/halv Tr(/Lam A)>,$$
%
%
which shows that
$$/eqalignno<
/Lam/mid z&/sim</rm/ Wishart>_d(n+d-a,A^<-1>),&(4.27)/cr
/mu/mid(/Lam,z)&/sim N_d/lp/hm,/ndel/sg/rp,&(4.28)/cr>$$
cf. Mardia, Kent and Bibby (1979, p.~85). Accordingly,
$$/eqalign<
g(x,z)&=/pietc E_z/hvlam/exp/lb-/halv(/mu-x)'/Lam(/mu-x)/rb/cr
&=/pietc E_z/lsb/hvlam E_z/exp/lb-/halv(/mu-x)'/Lam(/mu-x)/rb/mid/Lam/rsb/cr
&=/pietc E_z/hvlam E/,/exp/lb-/halv/,/ndel/chi_d^2
	/bigl(n(x-/hm)'/Lam(x-/hm)/bigr)/rb/cr
&=/pietc E_z/hvlam/lp1+/ndel/rp^<-d/sla 
	2>/exp/lb-/halv<n/over n+1>(x-/hm)'/Lam(x-/hm)/rb/cr>$$
by (4.23). We need further the following 
/Blemma Let $/Lam/sim</rm/ Wishart>_d(m,B)$, i.e.~the
density is
$$k_d(m,B)/vert/Lam/vert^<(m+1-d)/sla 2-1>/exp/lb-/halv/,Tr(/Lam B^<-1>)/rb,$$
where
$$k_d(m,B)=/lsb 2^<md/sla 2>/pi^<d(d-1)/sla 4>/vert B/vert^<m/sla 2>
	/prod_<i=1>^d/Gamma/lp/halv(m+1-i)/rp/rsb^<-1>.$$
%
%
Then
$$/eqalignno<
&E/vert/Lam/vert^<b/sla 2>e^<-/halv/,cy'/Lam y>/cr
&/qquad=2^<bd/sla 2>/prod_<i=1>^d</Gamma/bigl(/halv(m+b-i+1)/bigr)
	/over/Gamma/bigl(/halv(m-i+1)/bigr)>/,/vert B/vert^<b/sla 2>
	(1+cy'By)^<-(m+b)/sla 2>,/cr>$$
for $b/gt -(m+1-d),/,c/geq0,/,y/in/rr^d$.
/Elemma
/Bproof We use the identities 
$$cy'/Lam y=Tr(cyy'/Lam),$$
$$/vert I+cByy'/vert=1+cy'By,$$
and the fact that the cited Wishart density integrates to one. /Eproof
%
%
Now $g(x,z)$ can be found explicitly. Some manipulations give
$$/eqalignno<
 g(x,z)&=/pi^<-d/sla 2>(n+1)^<-d/sla 2>
	</Gamma/bigl(/halv(n+d-a+1)/bigr)/over
	/Gamma/bigl(/halv(n-a+1)/bigr)>/cr
&/qquad/vert/hs/vert^<-1/sla 2>/lsb1+<1/over n+1>(x-/hm)'
	/hs^<-1>(x-/hm)/rsb^<-(n+d-a+1)/sla 2>&(4.29)/cr>$$
which is a multivariate $t$-distribution, more precisely,
$$g(x,z)=t_d/lp n-a+1,/hm,<n+1/over n+1-a>/,/hs/rp(x)$$
according to the definition in Berger (1980a, p.~395).

The picture we should have in mind is one with two smooth surfaces
$N_d(/hm,/hs)(x)$ and $g(x,z)$
in $/rr^d$, both having a clear peak at $x=/hm$ and both having
the same system of ellipsoids for $/lb x/mid</rm density>=</rm constant>/rb$,
but where $g(x,z)$ is slightly more spread out than the 
estimative version $N_d(/hm,/hs)(x)$. In fact, the covariance matrix
for a random vector drawn from $g(x,z)$ can be shown to be
$<n+1/over n-a-1>/hs$.  Note that the difference between the densities
is small when $n$ is large.

The resulting predictive classification rule is as in
Theorem 4.1 with
$$/eqalignno<
 g_k(x,z)&=/pi^<-d/sla 2>(n_k+1)^<-d/sla 2></Gamma/bigl((n_k+d-a+1)
	/sla 2/bigr)/over/Gamma/bigl((n_k-a+1)/sla 2/bigr)>/cr
&/qquad/vert/hs/vert^<-1/sla 2>/lsb1+<1/over n_k+1>(x-/hm_k)'/hs_k^<-1>
	(x-/hm_k)/rsb^<-(n_k+d-a+1)/sla 2>&(4.30)/cr>$$
Aitchison, Habbema and Kay (1977) propose this method with $a=d+1$
in (4.25), corresponding to $a_0=d+1$ in (4.24), which gives
the vague prior distribution advocated by Box and Tiao (1973),
who use Jeffreys' (1961) information principle.  Berger's
(1980a) approach seems to lead to  $a_0=2d$ and $a=2$ in the 
formulae above.
/subsection<4.2.D><Common covariance matrix>  It is interesting to see the effect of the 
predictive approach in the model 
$$/fkx=N_d(/mu_k,/sg),/ktok,/eqno(4.31)$$
where both the $/mu_k$'s and the common $/sg$ are unknown.

Let us use the vague prior
$$(/mu_1/upto/mu_k,/Lam)/sim/vert/Lam/vert^<-a/sla 2>d/mu_1/cdots d/mu_k/,d/Lam
	/eqno(4.32)$$
where $/Lam=/sg^<-1>$. The density for $Z$ given $/mu_k$'s and $/Lam$ can be
written
$$/prod_<k=1>^K/lsb(2/pi)^<-d/sla 2>/vert n_k/Lam/vert^<1/sla 2>/exp
	/lb-/halv(/mu_k-/hm_k)'(n_k/Lam)(/mu_k-/hm_k)/rb/rsb$$
where $A_k=n_k/hs_k$ and $A=/sum_<k=1>^KA_k=N/hs$.
When the $/mu_k$'s are integrated out in the simultaneous distribution
for $(Z,/mu_1/upto/mu_k,/Lam)$ it is seen that
$$(Z,/Lam)/sim(2/pi)^<-(N-K)d/sla 2>/prod_<k=1>^Kn_k^<-d/sla 2>/vert/Lam/vert
	^<(N-K-a)/sla 2>e^<-/halv Tr(/Lam A)>$$
and that
$$Z/sim(2/pi)^<-(N-K)d/sla 2>/prod_<k=1>^Kn_k^<-d/sla 2>/lsb
	k_d(N-K-a+d+1,A^<-1>)/rsb^<-1>,$$
where $k_d(m,B)$ was defined in the Lemma in 4.2.C. It follows that
$$/displaylines<
/Lam/mid z/sim</rm Wishart>_d(N-K-a+d+1,A^<-1>),/cr
/qquad/mu_1/upto/mu_K/mid(/Lam,z)/sim/prod_<k=1>^KN_d/lp/hm_k,
	<1/over n_k>/sg/rp./cr>$$
These results are valid when $N/geq K+a-1$.
%
%

Now the predictive densities can be found explicitly, reasoning as in 4.2.C:
$$/eqalign<
g_k(x,z)&=E_z(2/pi)^<-d/sla2>/hvlam/,/exp/lb-/halv(x-/mu_k)'/Lam(x-/mu_k)/rb/cr
&=(2/pi)^<-d/sla2>E_z/lsb/hvlam E_z/,/exp/lb-/halv(/mu_k-x)'
	/Lam(/mu_k-x)/rb/mid/Lam/rsb/cr
&=(2/pi)^<-d/sla2>E_z/hvlam E/,/exp/lb-/halv/,<1/over n_k>/chi_d^2
	/bigl((x-/hm_k)'n_k/Lam(x-/hm_k)/bigr)/rb/cr
&=(2/pi)^<-d/sla2>E_z/hvlam/lp1+<1/over n_k>/rp^<-d/sla2>/exp/lb-/halv<n_k/over
	n_k+1>(x-/hm_k)'/Lam(x-/hm_k)/rb/cr
&=/pi^<-d/sla2></Gamma/bigl((N-K-a+d+2)/sla2/bigr)/over
	/Gamma/bigl((N-K-a+2)/sla2/bigr)>N^<-d/sla2>/cr
&/qquad/lp n_k/over n_k+1/rp/vert/hs/vert^<1/sla2>/lsb1+<1/over N>/,
	<n_k/over n_k+1>(x-/hm_k)'/hs^<-1>(x-/hm_k)/rsb
	^<(N-K-a+d+2)/sla2>./cr>$$
Note that the predictive and ML-estimative procedures are equivalent (apart 
from decisions about the Doubt option) if (and only if) the
$n_k$'s are equal.

As in 4.2.C $a=d+1$ may be used, for example, to please Jeffreys (1961),
or $a=2$, giving a certain right invariant Haar measure.
%
%
/subsection<4.2.E><Mixed discrete and continuous data> Let us include also a non-standard example, namely the model for
vectors $X=(A,Y)$ studied in 3.2.D, where $A/in/lb0,1/rb$ and $Y/in/rr^d$.
More general models along the lines of Chapter 9 can be handled similarly.

Let us choose the ``non-reduced'' model
$$/fkx=f(a,y;/,p_k,/mu_<k,0>,/sg_<k,0>,/mu_<k,1>,/sg_<k,1>)$$
of (3.32) for concreteness. The cases corresponding to Reductions
1, 2, 3 in 3.2.D can also be taken care of by similar means.

The predictive approach needs a prior density for the full family of 
parameters $p_k,/mu_<k,0>,/sg_<k,0>,/mu_<k,1>,/sg_<k,1>;/ktok$.
If we choose a prior distribution making these $K$ 
%
%
sets of parameters statistically independent, then the
arguments presented in 4.2.C can be repeated, and are seen to imply that 
$$g_k(x,z)=E(/fkx/mid z_k),$$
say, where $z_k$ is the set of $x_j^/k$'s from class $k$
only. Thus $k$ may be dropped in the notation.

Consider the prior density
$$/lam(p,/mu_0,/Lam_0,/mu_1,/Lam_1)=/vert/Lam_0/vert^<-b/sla2>
	/vert/Lam_1/vert^<-b/sla2>,/eqno(4.33)$$
where $/Lam_0=/sg_0^<-1>$ and $/Lam_1=/sg_1^<-1>$ are preferred
as parameters. It takes $p,/mu_0,/sg_0,/mu_1,/sg_1$ to be independent,
$p$ is uniform on $/lsb 0,1/rsb$, $/mu_0,/mu_1$ are uniform on $/rr^d$, whereas
$/sg_0$ and $/sg_1$ are given the same vague prior. Only $p$ has a
proper prior distribution.

The density for a sample $X_1/upto X_n$, where
$X_j=(A_j,Y_j)$, may be written 
$$/displaylines<
(1-p)^<N_0>(2/pi)^<-N_0d/sla 2>/vert/Lam_0/vert^<N_0/sla 2>e^<-/halv 
	Tr(/Lam_0B_0)>/exp/lb-/halv/,N_0(/mu_0-/hm_0)'/Lam_0(/mu_0-/hm_0)/rb/cr
p^<N_1>(2/pi)^<-N_1d/sla2>/vert/Lam_1/vert^<N_1/sla2>e^<-/halv Tr(/Lam_1B_1)>
	/exp/lb-/halv/,N_1(/mu_1-/hm_1)'/Lam_1(/mu_1-/hm_1)/rb,/cr>$$
where
$$/eqalign<
N_a&=/sum_<j=1>^nI/lb A_j=a/rb,/cr
/hm_a&=<1/over N_a>/sum_<j:A_j=a>Y_j,/cr
B_a&=/sum_<j:A_j=a>/twice<(Y_j-/hm_a)>'=N_a/hs_a,/cr>$$
$a=0,1$. Hence $(X_1/upto X_n,p,/mu_0,/Lam_0,/mu_1,/Lam_1)$ has
simultaneous density
$$/displaylines<
(2/pi)^<-nd/sla2>(1-p)^<N_0>p^<N_1>/prod_<a=0,1>/Bigg/lbrace/vert/Lam_a/vert^<(N_a-b-1)/sla2>e^<-/halv Tr(/Lam_aB_a)>/cr
/vert/Lam_a/vert^<1/sla2>/exp/lsb-/halv(/mu_a-/hm_a)'N_a/Lam_a(/mu_a-/hm_a)/rsb/Bigg/rbrace./cr>$$

It may be shown after the necessary
%
%
amount of algebraic manipulations that $p, (/mu_0,/Lam_0),/allowbreak (/mu_1,/Lam_1)$ are 
three independent sets of parameters given the data, and that
$$/eqalignno<
p/mid z&/sim</rm Beta>(N_0+1,N_1+1),/cr
/Lam_a/mid z&/sim</rm Wishart>_d(N_a+d-b,A_a^<-1>),/cr
/mu_a/mid(/Lam_a,z)&/sim N_d/lp/hm_a,<1/over N_a>/sg_a/rp,&(4.34)/cr>$$
$a=0,1$. $z$ is used as short hand notation for the data set $X_1/upto X_n$.

This makes it possible to determine the predictive density, which will exist
provided the Wishart and normal densities above are proper,
i.e.~$N_a/geq b,/,a=0,1$. (This happens with a probability tending to one.)

We need
%
%
$$g(x,z)=g(a,y,z)=E/lb f(a,y)/mid z/rb$$
for $a=0,1$ and $y/in/rr^d$. By (4.34) and the work already carried out in
4.2.C we get
$$/displaylines<
g(a,y,z)=<N_a+1/over n+2>/pi^<-d/sla 2>(N_a+1)^<-d/sla 2>
	</Gamma/bigl((N_a+d-b)/sla2/bigr)/over
	/Gamma/bigl((N_a-b)/sla 2/bigr)>/cr
/vert/hs_a/vert^<-1/sla2>/lsb1+<1/over N_a+1>(y-/hm_a)'/hs_a^<-1>(y-/hm_a)/rsb
	^<-(N_a+d-b+1)/sla2>./cr>$$
/section<4.3><Concluding remarks><>
Chapters~3 and~4 have been concerned with the same problems, but the
present chapter has sought to take the sampling variability that every
estimator is plagued with into account. The investigations led us
naturally to the class of predictive methods. These have connections
to  Bayesian methods, see Geisser (1982), but do not necessarily demand
as input subjective knowledge in the form of </it a priori> distributions
for the parameters.

We have seen that the predictive versions tend to be more
complicated than the plug-in versions, but that they are reasonably 
similar when the training sample sizes are moderate or large. Only
when the training sets are small or the number of classes large will
there be significant differences between estimative and predictive
rules. See Aichison, Habbema and Kay (1977) for an example where the 
difference is dramatic.

The important question of which model to employ is of course always present,
see the comments made at the end of Chapter 3. There is an interesting method 
of </it smoothing over several models//> that fits within the predictive
framework: The resulting $/hf_k(x)$ is obtained as a linear combination of
predictive densities for the respective models, with weights equal to the 
likelihood of the models, as dictated by the data.

There are other ways of incorporating estimation variability.
One could employ unbiased estimators of the densities themselves,
see (3.22) for the normal case. An interesting facet of (3.22) is that
it is equal to zero outside a certain finite region. Incoming $X$-vectors
falling in this zero region would then be taken as outliers, see Chapter~6
for alternative but similar strategies.

A related approach that may appear somewhat more natural in the
classification context, where decisions hinge on how large the ratios
$f_k(x)/sla f_l(x)$ are, finds unbiased estimators for the </it logarithm//>
of the densities, i.e./ uses $/wt f_k(x)$ to estimate $f_k(x)$, where
$/log/wt f_k(x)$ is unbiased for $/log f_k(x)$. It is easy to construct 
such procedures in the normal case, using results available in Section 10.1.B.

Later on in this work several principles are illustrated using,
essentially, plug-in estimative methods. But predictive versions along
the lines of the present chapter could also be constructed and used in
their place. This remark applies to the updating methods of Chapter~7
and some of the estimation methods for mixed discrete and continuous
data in Chapter~9.


/chapter<5><Nonparametric classification>
Consider once more the optimal discriminant rule derived in
Chapter 1 when prior probabilities and class densities
were assumed known:
$$C_0(x)=/cases<D&if every $/pkx/leq1-c$,/cr
                k&if $/pkx=/max_<t/leq K>/,P(t/mid x)/gt1-c$/cr>/eqno(5.1)$$
where $/pkx=/pi_k/fkx/sla/sum_<t=1>^K/pi_t f_t(x)$.
In Chapters 3 and 4 the approach was, when encountering unknown
class densities, to model them according to a parametric structure,
and then estimate the parameters involved from training
sets. A warning was given in Section 3.1, however, observing that
the resulting procedures did not actually approximate the optimal $C_0$,
as they might seem to pretend to, but rather a possibly worse $C_1$,
having $f_k(x,/tho)$ instead of the genuine $/fkx$,
where $/tho$ was the least false parameter value according to
the Kullback-Leibler information distance.

The present chapter is concerned with the possibility of estimating the
densities nonparametrically, that is, with few or none assumptions 
imposed on them. Such an approach will usually demand large training sets
to be successful (but see Chapter 7).

Section 5.1 reviews some results that suggest that nonparametric classifiers
will tend to surpass the parametric ones if a sufficient amount of
training data is available. The popular kernel density estimation
method is considered in 5.2. These estimators cannot be
condensed to compact formulae that use only a small or moderate number of
characteristics for the classes, and in general the full set of training
data must be traversed in order to classify a new feature vector. This
rather unsuitable quality is not shared by an interesting class of density
estimators studied in Sections 5.3 and 5.4, the orthogonal expansion
estimators. These fit more nicely into the ``usual design'' where 
$/wh f_k(x)$ is given in terms of a moderate number of class description
parameters. 5.4 in particular presents methods that are more semiparametric
than fully nonparametric in nature, in that an initial parametric description
is taken as a starting point for further ``nonparametric corrections''.
Section 5.5 is ``fully nonparametric'' again, and treats
$k$-nearest-neighbour methods.
/section<5.1><General considerations><>
Consider a general classification procedure 
$/ca=/ca(x;z):/cx/times/cx^N/to/lb 1/upto K,D/rb$
based on a training set in the form (3.3),
$N=/sum_<k=1>^K n_k$ being the training set's total size.
Let $A_k=/lb x/mid/ca(x;z)=k/rb,/ k=1/upto K$, and
$A_0=/lb x/mid/ca(x;z)=D/rb$ partition the feature space
$/cx$, and write $/phi_k(x)=I/lb x/in A_k/rb$ for
the indicator functions. The conditional risk of $/ca$,
given $z$, becomes
$$/eqalignno<
R(/ca)&=R(/ca,/ftfk)/cr
&=EL(C,/ca(X;z))/cr
&=/int_/cx E/lb L(C,/ca(x;z))/mid X=x/rb/,f(x)/,d/mu/cr
&=/int_/cx/sum_<k=1>^K L(k,/ca(x;z))/,/pkx/,f(x)/,d/mu/cr
&=/int_<A_0> cf(x)/,d/mu+/sum_<k=1>^K /int_<A_k>(1-/pkx)/,f(x)/,d/mu/cr
&=/int_/cx /lb c/phi_0(x)+/sum_<k=1>^K/phi_k(x)(1-/pkx)/rb/,f(x)/,d/mu/cr
&=1-/int_/cx/lb(1-c)/,/phi_0(x)+/sum_<k=1>^K/phi_k(x)/,/pkx/rb/,f(x)/,d/mu,&
(5.2)>$$
still employing the loss function (1.3), and where
$f=/sum_<k=1>^K/pi_k/,f_k$ is the marginal density of
$X$. It may be seen that the above expression is retained also for
randomised procedures, with the interpretation
$$/phi_k(x)=Pr/lb/ca(x;z)=k/mid X=x/rb.$$

The optimal rule (5.1) has the minimum possible risk,
$$/eqalignno<
R_0&=R_0(/ftfk)/cr
   &=R(C_0,/ftfk),&(5.3)>$$
i.e.~minimum Bayes risk is expression (5.2) with
$$/eqalign<
/phi_k(x)&=I/lb/pkx=/mtlek/ptx/gt1-c/rb,/ k=1/upto K,/cr
/phi_0(x)&=I/lb/mtlek/ptx/leq1-c/rb./cr>$$
(Ties may be broken in an arbitrary way without affecting
$R(C_0,/ftfk)$.)

Assume now that density estimates $/wh f_1/upto/wh f_K$
have been obtained on the basis of $z$. The natural analogue
of (5.1) is
$$/hCo(x;z)=/cases<D&if every $/widehat/pkx/leq1-c$,/cr
                   k&if $/widehat/pkx=/mtlek/widehat/ptx/gt1-c$,/cr>/eqno(5.4)$$
where
$$/widehat/pkx=</pi_k/,/wh/fkx/over/sum_<t=1>^K/pi_t/wh f_t(x)>,/ k=1/upto K./eqno(5.5)$$
Its risk, conditional on $z$, is
$$/displaylines<
/qquad R(/hCo,/ftfk)/cr
/qquad/qquad=1-/int_/cx/lb(1-c)/wh/phi_0(x)+/sum_<k=1>^K/hph_k(x)/,/pkx/rb/,f(x)/,d/mu/quad/cr>$$
where
$$/eqalign<
/hph_k(x)&=I/lb/wh/pkx=/mtlek/wh/ptx/gt1-c/rb,/cr
/hph_0(x)&=I/lb/mtlek/wh/ptx/leq1-c/rb.>$$
Thus
$$/eqalign<
0&/leq R(/hCo,/ftfk)-R_0(/ftfk)/cr
&=/int_/cx/lb(1-c)/lp/ph_0(x)-/hph_0(x)/rp+/sum_<k=1>^K/lp/ph_k(x)-/hph_k(x)/rp
/,/pkx/rb/,f(x)/,d/mu./cr>/eqno(5.6)$$
We should also take interest in the </it global risk//> of
$/hCo(x;z)$, defined by taking expectation w.r.t.~the random training set $Z$:
$$</rm GR>(/hCo,/ftfk)=ER(/hCo,/ftfk)$$
$$=1-/int_/cx/lb(1-c)/,E/hph_0(x)+/sum_<k=1>^K E/hph_k(x)/pkx/rb f(x)/,d/mu.$$

It is now clear that both the conditional and the global risk of $/hCo$ will 
be close to minimum Bayes risk (5.3) if the density estimates
$/wh f_1/upto /wh f_K$ based on $z$ are so good that
$/hph_k(x)$ and $E/hph_k(x)$ are close to the optimal 
$/ph_k(x)$.

The classification procedure $/hCo=/hCo(x;z)$ is termed
</it Bayes risk consistent//> if
$$R(/hCo,/ftfk)/totop<P>R_0(/ftfk)/eqno(5.7)$$
when $n_1/upto n_K$ in (3.3) go to infinity.
Then an immediate consequence is </it global  
risk consistency//>,
$$</rm GR>(/hCo,/ftfk)/longrightarrow R_0(/ftfk),/eqno(5.8)$$
using Lebesgue's dominated convergence theorem.
(Sometimes (5.8) may hold and not (5.7).)

Sufficient for Bayes risk consistency are the 
conditions
$$/displaylines<
/icx/lp/ph_0(x)-/hph_0(x)/rp/,f(x)/,d/mu/totop<P>0,/cr
/icx/lp/ph_k(x)-/hph_k(x)/rp/,/fkx/,d/mu/totop<P>0,/ k=1/upto K,/cr>$$
using (5.6). Writing $/ph_k(x)=I/lb x/in A_k/rb$,
$/hph_k(x)=I/lb x/in/wh A_k/rb$, we get
$$/icx/lp/ph_k-/hph_k/rp f_k/,d/mu=/int_<A_k>/lp1-/hph_k/rp f_k/,d/mu-/sum_<t/not=k>/int_<A_t>/hph_k/,f_k/,d/mu.$$
Even weak assumptions will imply that
$$/displaylines<
/mu/lb x/in A_k/mid/hph_k(x)/not=1/rb/to0,/cr
/mu/lb x/in A_t/mid/hph_k(x)/not=0/rb/to0,/cr>$$
ensuring Bayes risk consistency.
/Btheorem<> Suppose
$$/wh/fkx/totop<P>/fkx,</rm a.e.>/,(/mu),/ k=1/upto K./eqno(5.9)$$
Then $/hCo$ defined in (5.4) is Bayes risk consistent.
/Etheorem/Bproof
 Using the bounded convergence theorem for random functions 
(Glick, 1974),
this slightly extended version of Greblicki's (1978) result is obtained. /Eproof

It seems to have been Van Ryzin (1966) who introduced
the Bayes risk consistency term. His sufficient conditions were
$$/icx E/mid/wh/fkx-/fkx/mid d/mu/to0,/ k=1/upto K./eqno(5.10)$$

There are large classes of nonparametric density estimators that satisfy both 
(5.10) and the condition in the theorem. On the
other hand, a classification  procedure based on a parametric model,
as in (3.4)--(3.5), will most often have a conditional risk converging
in probability to a $R(C_1,/ftfk)/gt R_0(/ftfk)$,
cf. Section 3.1. The important consequence is that nonparametric classifiers
will tend to outperform any given parametric one, if 
sufficiently large training sets are used.

As is usual when man thinks he has obtained a piece of truth,
some warnings need to be given. Firstly, the statement above is really
concerned with a sequence (as the training set increases) of
nonparametric procedures versus a sequence of procedures derived under a 
</it fixed//> parametric model. It is possible to increase the complexity of a 
parametric model as more and more training data become available,
and then the sequence of parametric procedures might win after all.
Secondly, at least some of the popular nonparametric density estimation 
methods are heavy and may need much computing effort to reach
say $/wh/fkx$, and maybe the full set of observations
$X_1^<(k)>/upto X_<n_k>^<(k)>$ needs to be traversed in the process 
(as opposed to using a compact formula involving a small number of class
descriptors (parameter estimates)). While this may be perfectly acceptable when only a small number of curves are to be produced, say,
it may be impractical and perhaps impossible in an automatic pattern
recognition context.

/section<5.2><Density estimation by the kernel method><>
Let $X_1/upto X_n$ be i.i.d.~with a (preferably continuous) density
$f(x),/,x/in/rr^d$. Introduce a </it kernel function//> $K$ on $/rr^d$,
assumed to fulfil
$$/eqalignno<
</rm sup>_<y/in/rr^d>/vert K(y)/vert&/lt/infty&(5.11</rm a>)/cr
/int/vert K(y)/vert d/mu&/lt/infty,&(5.11</rm b>)/cr
/pr y/pr^d K(y)&/to0/,</rm as>/,/pr y/pr/to/infty,&(5.11</rm c>)/cr
/int K(y)/,dy&=1.&(5.11</rm d>)>$$
The kernel method estimator of $f(x)$ is
$$/eqalign<f_n(x)&=n^<-1>h(n)^<-d>/sum_<j=1>^n K((x-X_j)/sla h(n))/cr
&=/int h(n)^<-d> K/lp<x-y/over h(n)>/rp/,dF_n(y),/cr>/eqno(5.12)$$
where $F_n=<1/over n>/sum_<j=1>^n/dl(X_j)$ is the empirical
distribution. $f_n(x)$ may be viewed as the density of $X+Y$,
where $X/sim F_n$ (i.e.~$Pr/lb X=X_j/rb=<1/over n>,/ j=1/upto n$), and
$$Y/sim h(n)^<-d> K(y/sla h(n))$$
is independent of $X$.

Usual additional requirements imposed on $K(y)$ is
$$/eqalignno<
K(y)&/geq0,/,y/in/rr^d,&(5.13</rm a>)/cr
K(y)&=K(-y),/,y/in/rr^d.&(5.13</rm b>)/cr>$$
It is not unusual to put $K(y)=N_d(0,/wh/sg)(y)$,
where $/wh/sg$ is the estimated covariance matrix for the sample.

Some results on the behaviour of the kernel estimator may be listed 
before we go on to discuss its relevance in the symbol recognition
context.

Assume (5.11 a -- d) and that $h(n)/to0$. Then $Ef_n(x)/to f(x)$ as $n/to/infty$.
Suppose further that (5.13 a -- b) holds, and that
$$h(n)/to0,/,n/,h(n)^d/to/infty./eqno(5.14)$$
Then
$$E/lb f_n(x)-f(x)/rb^2/to0/eqno(5.15)$$
whenever $f$ is continuous in $x$, in particular, $f_n(x)$ is
consistent for $f(x)$. The consistency is uniform if
$$/eqalignno<
h(n)&/to0,/,n/,h(n)^<2d>/to/infty,&(5.16)/cr
/widetilde<K>(u)&=/int_</rr^d>e^<iu'y>/,K(y)/,dy&(5.17)/cr>$$
is absolutely integrable, and if (5.11 a -- d), (5.13 a -- b)
continue to hold. More precisely,
$$/sup_<x/in/rr^d>/,/vert f_n(x)-f(x)/vert/totop<P>0./eqno(5.18)$$

Now consider the possibility of using a particular kernel method
estimator in a pattern recognition context, say to obtain a 
classification rule of the form (5.4) when devising an automatic
symbol recogniser.

First, there is the problem of choosing the kernel function $K(y)$ and the 
``window size'' $h(n)$. It turns out the latter problem is both more
important and more difficult than the first one.
Accordingly, we might as well try
$$K(y)=N_d(0,/wh/sg)(y),/eqno(5.19)$$
with $/wh/sg$ defined in (3.19), and then employ one of a variety
of methods proposed in the literature to choose a suitable $h(n)$
based on the data, see for example Silverman~(1976, 1978),
Rudemo~(1982), Duin~(1976), and Scott, Tapia, and Thompson~(1980).
It is important to realise that every method is subjective to
some extent.

Having chosen $K$ and $h(n)$ a more practical problem remains.
$f_n(x)$ in the form (5.12) is very awkward in that a lot of computation is 
needed to produce the estimate, and in that the full data set
$X_1/upto X_n$ has to be traversed. This may be acceptable if a moderate
number of feature candidate vectors are to be classified, in fact
Hermans~and Habbema~(1976) report success with the use of a kernel method
discriminant analysis program package intended for use in medical
statistics, and Hand (1982) devotes a full book to kernel
discriminant analysis, unbothered by computational aspects.
In the technological pattern recognition world, however, where a
large number of candidate vectors are to be classified 
automatically, $f_n(x)$ in the form (5.12) is highly impractical.

Several clever approximations and shortcuts are however available, and 
can make some of the kernel method estimators much more feasible
computationally.

It was pointed out that $f_n(x)$ may be viewed as the density
of a convolution $X+Y$, where
$$/eqalignno<
Ee^<is'X>&=<1/over n>/sum_<j=1>^n e^<is'X_j>,&(5.20)/cr
Ee^<is'Y>&=/int e^<is'y>/,h(n)^<-d>/,K(y/sla h(n))/,dy&(5.21)/cr>$$
are their characteristic functions. In particular,
if (5.19) is chosen, then
$$/int e^<is'x>/,f_n(x)/,dx=<1/over n>/sum_<j=1>^n e^<is'X_j>/exp/lp
-/halv h(n)^2 s'/wh/sg s/rp,$$
which equivalently may be stated
$$/widetilde<f>_n(s)=u(s)/exp/lp-/halv h(n)^2 s'/wh/sg s/rp,/eqno(5.22)$$
where $/widetilde<f>_n$ is the Fourier transform of $f_n$ and
$u(s)=(2/pi)^<-1/sla 2>/,<1/over n>/sum_<j=1>^n/exp(is'X_j)$
is the Fourier transform of the data.

The kernel method estimator $f_n(x)$ may now be
approximated by fast inverse Fourier algorithms. Some
technical suggestions related to this are in Silverman (1982) and
Jones and Lotwick (1984).

If a kernel with bounded support is used only the nearest neighbours
to $x$ among $X_1/upto X_n$ enter into $f_n(x)$. References to smart
algorithms for finding nearest neighbours are given in 5.5.C.

There are several important variants of the kernel approach discussed 
above. One might smooth differently in different directions, i.e./ use
$$/wh f_n(x)=(nh_1/cdots h_d)^<-1>/sum_<j=1>^n K/lp(x_1-X_<j,1>)/sla
h_1/upto(x_d-X_<j,d>)/sla h_d/rp$$
instead of (5.12). And one might use </it locally varying//> smoothing 
parameters; see Breiman, Meisel, and Purcell (1977) for a paper on 
variable kernel estimates, and see Section 5.5 for some related ideas.

</it Histograms//> are less sophisticated but useful alternatives
to kernel estimators. The cells should not necessarily be of equal volume.
Other binning-the-data procedures that reduce computations are mentioned
in Hjort~(1986).

Before we leave the kernel method estimator it should be pointed out that
(5.15) and (5.18) already guarantee, in view of (5.9), that
the resulting nonparametric classification procedures are Bayes
risk consistent, provided only that the class densities $/ftfk$ are
continuous.

We conclude this subsection by mentioning another approach to 
density estimation that is related to the kernel type estimators, but
are computationally much simpler and faster. Scott (1985a, 1985b)
recently studied two variations on the ordinary histogram, namely the
frequency polygon and the average shifted histogram, and found that they
are able to compete with (for example) kernel type estimators while
retaining the advantage of being conceptually and computationally
simple. Hjort (1986) improved Scott's results, and the ``generalised
frequency polygon of the average shifted histogram'' described there is
an attractive procedure in the present context.
/section<5.3><Density estimation by orthogonal expansions><>
The method of orthogonal expansions provide us with
a new, large class of nonparametric density estimators, and
a fair number of special choices may give rise to valuable classification 
algorithms.  Although the approach will be seen to be rather different from the
kernel method several results in the literature indicate that the two types of 
estimators may produce reasonably similar results on the whole.
Each of them will have its own small peculiar dis- and advantages, however.

The approach is nicely illustrated in the following simple framework involving
a (univariate) density on $/lsb0,1/rsb$. The method will be generalised
afterwards.
/subsection<5.3.A><A density on $/lsb0,1/rsb$> Let $X_1/upto X_n$ be i.i.d.~with a
continuous density $f(x)$ on $/lsb0,1/rsb$. It is possible to extend
$f$ to an even function on $/lsb-1,1/rsb$, i.e.~let $f(x)=f(-x)$. 
Then $f$ may be expanded in a Fourier series
$$f(x)=a(0)+/sum_<j=1>^/infty a(j)/cos(j/pi x)./eqno(5.23)$$

Consider the finite-sum approximation 
$$f_m(x)=a(0)+/sum_<j=1>^m a(j)/cos(j/pi x)/eqno(5.24)$$
to $f(x)$. If the fit is measured by integrated square error
$$/ise(f,f_m)=/int_0^1(f-f_m)^2/,dx,$$
then it is easily seen that the values to choose are
$$/eqalign<
a(0)&=1,/cr
a(j)&=2/int_0^1/cos(j/pi x)/,f(x)/,dx=2E/cos(j/pi X),/,j/geq1./cr>/eqno(5.25)$$
A natural estimator of $f(x)$ is then
$$/eqalignno<
/wh f_m(x)&=1+/sum_<j=1>^m /wh a(j)/cos(j/pi x);&(5.26)/cr
/wh a(j)&=2/ndel/sum_<t=1>^n/cos(j/pi X_t),/,j/geq1.&(5.27)>$$

The problem left is to decide when to stop, i.e.~which $m$ should be used. 
Consider 
$$/eqalign<
/ise(f,/wh f_m)&=/int_0^1/lb/sum_<j=m+1>^/infty a(j)/cos(j/pi x)
      +/sum_<j=1>^m(a(j)-/wh a(j))/cos(j/pi x)/rb^2/,dx/cr
&=/halv/sum_<j=m+1>^/infty a(j)^2+/halv/sum_<j=1>^m(/wh a(j)-a(j))^2,/cr>$$
which gives a mean integrated squared error
$$/eqalign<</rm MISE>(f,/wh f_m)&=E/lb /ise(f,/wh f_m)/rb/cr
&=/halv/lb<1/over n>/sum_<j=1>^m /sigma(j)^2+/sum_<j=m+1>^/infty a(j)^2/rb,/cr>/eqno(5.28)$$
where
$$/sigma(j)^2=</rm Var>/lb2/cos(j/pi X_t)/rb,/,j/geq0./eqno(5.29)$$

We should like to choose the $m$ that makes $/mise(f,/wh f_m)$ as small
as possible, thus meeting the trade-off problem with the two types
of error: $/sum_<j=m+1>^/infty a(j)^2$ is
the error made by stopping $f_m$ too early
($/wh f_m$ is seen to approximate $f_m$ rather than $f$), whereas
$<1/over n>/sum_<j=1>^m/sigma(j)^2$ accounts for the variability of
the estimates of $a(j), j/leq m$.

Now it may be seen that 
$$/eqalign<
D(j)&=/mise(j)-/mise(j-1)/cr
    &=/halv/lb<1/over n>/sigma(j)^2-a(j)^2/rb./cr>/eqno(5.30)$$
For many smooth densities $/mise(m)=/mise(f,/wh f_m)$ decreases
initially but then starts climbing, so that we might expect $D(j)/lt0$ for the 
first $j$-values and then $D(j)/gt0$, say for $j/geq m_0+1$.
In such a case $/wh f_<m_0>$ should be the best choice.

This procedure can only be mimicked in that the $D(j)$'s
are unobservable. They may be estimated from the data in an unbiased way, 
however. One has
$$/eqalign<
/sigma(j)^2&=4/int_0^1/cos^2(j/pi x)/,f(x)dx-a(j)^2/cr
           &=2+a(2j)-a(j)^2,/,j/geq1,/cr>/eqno(5.31)$$
so that
$$D(j)=/halv/lb/ndel(2+a(2j))-<n+1/over n>/,a(j)^2/rb.$$
Since
$$/eqalign<
E/wh a(j)^2&=a(j)^2+/ndel/sigma(j)^2/cr
            &=<n-1/over n>/,a(j)^2+/ndel(2+a(2j)),/cr>$$
$<n/over n-1>/,/wh a(j)^2-<1/over n-1>(2+/wh a(2j))$ is unbiased for
$a(j)^2$. Hence
$$/eqalignno<
/wh D(j)&=/halv/,/ndel/lb2+/wh a(2j)-<n(n+1)/over n-1>/,/wh a(j)^2
                +<n+1/over n-1>(2+/wh a(2j))/rb/cr
        &=<1/over n-1>/lb2+/wh a(2j)-/halv(n+1)/,/wh a(j)^2/rb&(5.32)/cr>$$
is an unbiased estimator of $D(j);/,j/geq1$.

Accordingly, an easily implementable nonparametric density
estimation procedure is the following: Calculate for example the first
$/lsb n^<1/sla 2>/rsb$ estimates $/wh D(j)$. Let $m_0+1$ be
the first $j$ having $/wh D(j)$ positive. Then use
$$/wh f_<m_0>(x)=1+/sum_<j=1>^<m_0>/wh a(j)/cos(j/pi x),/,0/leq x/leq1.
/eqno(5.33)$$

This rule was (essentially) proposed in the classical expansion estimator
paper by Kron/-mal and Tarter (1968). It has certain weaknesses, though,
and safer rules are developed in Section 5.4.
/subsection<5.3.B><General case> Let $X_1/upto X_n$ be i.i.d.~with
an unknown, preferably continuous density $f(x)$ on a set $/cx$. $/cx$ will be 
a subset of $/rr^d$ in the applications we have in mind. Assume that 
$f=dP/sla d/mu/in/cl^2(/mu)$, and assume further that an explicit,
complete orthonormal basis $/lb/psi_j;/,j/geq0/rb$ for $/cl^2(/mu)$ is 
available, i.e.
$$/int/psi_k(x)/,/psi_j(x)/,d/mu=/dl_<k,j>=I/lb k=j/rb./eqno(5.34)$$
(Integrals are over the full space $/cx$ unless otherwise noted.)

$f$ may be expanded in a series
$$f(x)=/sum_<j=0>^/infty a(j)/psi_j(x),/,x/in/cx./eqno(5.35)$$
The best finite-sum approximation 
$$f_m(x)=/sum_<j=0>^m a(j)/psi_j(x)$$
in terms of integrated square error results when
$$a(j)=/int/psi_j(x)f(x)/,d/mu=E_f/psi_j(X_t),/,j/geq0./eqno(5.36)$$
The traditional orthogonal series estimator of $f$, based on the chosen
set of orthonormal functions, is
$$/wh f_m(x)=/sum_<j=0>^m /wh a(j)/psi_j(x)/eqno(5.37)$$
where
$$/wh a(j)=/ndel/sum_<t=1>^n/psi_j(X_t),/,j/geq 0./eqno(5.38)$$
$/wh a(j)$ is the natural estimator of $a(j)$.

Above $m$ should be thought of as dependent upon sample size $n$,
i.e.~$m=m(n)$. One may show that the essential condition
under which
$$/mise/lp/wh f_<m(n)>,f/rp=E/,/int/lp/wh f_<m(n)>-f/rp^2/,d/mu/to0/eqno(5.39)$$
is guaranteed, is
$$m(n)/to/infty,/,m(n)/sla n/to0./eqno(5.40)$$

The estimation procedure is not fully specified until a rule defining 
$m=m(n)$ for a given sample $X_1/upto X_n$ is defined. (5.40) is only a crude
guideline. One may actually consider (at least) two
aspects of the problem. Firstly, one may concentrate
on the asymptotic framework, and for example obtain conditions under which 
Bayes risk consistency of the resulting classification procedures is ensured,
cf. (5.9) and (5.10). Secondly, the finite-sample
problem should be dealt with, devising an explicit stopping rule for $m$, 
perhaps taking the asymptotic results into account.

We may quote some results from the literature pertaining to the asymptotic 
aspect. If the functions constituting the complete orthonormal
basis are uniformly bounded, i.e.~$/vert/psi_j(x)/vert/leq M$
for all $j$ and $x$, then the condition
$$m(n)/to/infty,/,m(n)^2/sla n/to0/eqno(5.41)$$
will imply
$$/sup_<x/in/cx>/,E/vert/wh f_<m(n)>(x)-f(x)/vert^2/to0,/eqno(5.42)$$
for every continuous $f$, and hence the nonparametric classification procedures
built as in (5.4), (5.5) with orthogonal series estimators for class densities 
will be Bayes risk consistent. This follows from Theorem 5.1.
More generally, if
$$M_n=/sup_<0/leq j/leq m(n)>/,/sup_<x/in/cx>/vert /psi_j(x)/vert
/lt/infty,/,n/geq1,$$
then $M_n^4 m(n)^2/sla n/to0$ together with $m(n)/to/infty$
will secure (5.42), cf. Prakasa Rao (1983, Section 2.2).

It is also possible under some regularity conditions
to determine the sequence of $m(n)$'s that will lead to the best estimators
in terms of $/mise(/wh f_<m(n)>,f)$ for a given $f$.
Systematic investigations by Bosq (1970) indicate that $m(n)$
should satisfy
$$/vert a(m(n))/vert/approx</rm const.>/sla(n+1)^<1/sla 2>,$$
where const. depends on the values of $/int/psi_j^2/,f/,d/mu$.

A reasonable stopping rule for $m$ based on the actual sample
$X_1/upto X_n$ may be constructed in a way similar to the rule proposed in 
5.3.A leading to (5.33). Consider
$$/eqalignno<
/mise(m)&=E/int/lp/wh f_m-f/rp^2 d/mu/cr
        &=E/int/lb/sum_<j=m+1>^/infty a(j)/,/psi_j(x)+/sum_<j=0>^m
          /lb a(j)-/wh a(j)/rb/psi_j(x)/rb^2/,d/mu/cr
        &=/sum_<j=m+1>^/infty a(j)^2+/ndel/sum_<j=0>^m/sigma(j)^2,/cr
/noalign</rm where>
/sigma(j)^2&=</rm Var>_f/psi_j(X_t)=/int/lb/psi_j(x)-a(j)/rb^2/,f(x)/,
d/mu,/,j/geq0.&(5.43)/cr>$$
Then 
$$/eqalignno<D(j)&=/mise(j)-/mise(j-1)/cr
		&=/ndel/sigma(j)^2-a(j)^2,/,j/geq1.&(5.44)/cr>$$ 
Now $S(j)^2$ is unbiased for $/sigma(j)^2$, where
$$S(j)^2=<1/over n-1>/sum_<t=1>^n /lb/psi_j(X_t)-/wh a(j)/rb^2,/,j/geq0.
/eqno(5.45)$$
Combining this with $E/wh a(j)^2=a(j)^2+/ndel/sigma(j)^2$, an unbiased 
estimator of the increase in MISE by including also the term 
$/wh a(j)/,/psi_j(x)$, i.e.
$D(j)$, may be constructed:
$$/wh D(j)=<2/over n>S(j)^2-/wh a(j)^2,/,j/geq1./eqno(5.46)$$
A simple stopping rule is: calculate the first $/wh D(j)$'s from the sample.
If $m_0+1$ is the smallest $j$ for which $/wh D(j)/gt0$, then use 
$/wh f_<m_0>(x)$.

The rule proposed above may not always be appropriate.
It depends strongly upon the chosen complete orthonormal
set $/lb/psi_j(x);/,j/geq0/rb$, and in particular also on the ordering within
the set.  In the situation of 5.3.A a natural ordering
exists, and the proposed rule seems reasonable.  Kron/-mal~and
Tarter~(1968) studied such trigonometric density
estimators and found them to be convenient and reasonably
accurate.  They were also happy to report that for many
smooth estimands  $f$ only few terms in the expansion were needed,
i.e. $m_0=m_0(n)$ above tended to be moderate in size even for
large values of $n$.
In other situations the orthonormal set may become more
complex, and perhaps no particular ordering seems natural.
Then variations of the above rule could be devised, as
exemplified in 5.3.C and 5.4 below.

Another approach to the determination of the terms $/wh a(j)/,/psi_j(x)$
which should be included in a final estimator of the 
form
$$/wh f_/ck(x)=/sum_<j/in/ck>/wh a(j)/,f_j(x),/eqno(5.47)$$
$/ck$ a finite subset of indexes, can be based on the simple result
$$/sqrt<n>/lb/ha(j)-a(j)/rb/totop<D>N/lb0,/sigma(j)^2/rb,/,j/geq0./eqno(5.48)$$ 
Under the hypothesis that $a(j)=0$, $/sqrt<n>/ha(j)/sla/wh/sigma(j)$ is 
approximately standard normal, where $/wh/sigma(j)$ is either $S(j)$ of
(5.45) or $/lbrace/ndel/sum_<t=1>^n/psi_j(X_t)^2/rbrace^<1/sla 2>$.
Hence one could evaluate the first $/lbrack n^<1/sla 2>/rbrack$ or so
$/ha(j)$'s and include in (5.47) only the ``significant'' ones,
for example those having $/vert/sqrt<n>/ha(j)/sla/wh/sigma(j)/vert/gt$ 1.96.
/note<Remark 1.>  It is known that orthogonal expansion estimators
may be bothered by the so-called Gibb's phenomenon, that
is, if many terms are included in the density estimator,
it may provide good fit to the underlying true density
in large portions of the space $/cx$, but may produce 
perhaps wild values near the boundary. The reason for this possible shortcoming
of the orthogonal expansion estimator can be traced to
the use of integrated squared error as loss function.

A couple of simple remedies may be devised, however,
in the case  where it is possible to give a fair
picture of the density in advance.

Suppose that $f_0(x)$ is such a ``prior guess'' or ``initial estimate'' of the 
density $f(x)$. Then write
$$f(x)=f_0(x)+g(x),/,x/in/cx,/eqno(5.49)$$
so that intuitively $g(x)$ is a curve exhibiting smaller
variation (around $0$) than does $f(x)$. Consequently $g(x)$ might
be more successfully estimated by the method of
orthogonal expansions, and might hopefully need fewer
terms.  An estimator of $g(x)$ is
$$/wh g_m(x)=/sum_<j=0>^m/wh b(j)/,/psi_j(x),$$
where $/wh b(j)$ estimates the best coefficient, which is seen
to be
$$/displaylines<
b(j)=/int/psi_j(x)/,g(x)/,d/mu=a(j)-a_0(j),/cr
a(j)=/int/psi_j(x)/,f(x)/,d/mu,/,a_0(j)=/int/psi_j(x)/,f_0(x)/,d/mu,>$$
in terms of $/int(g-/wh g_m)^2/,d/mu$. Hence 
$/wh b(j)=/wh a(j)-a_0(j)$ is the natural choice.  The
resulting estimator of $f(x)$ becomes
$$/eqalignno<
/dhat<f>_m(x)&=f_0(x)+/sum_<j=0>^m/lp/ha(j)-a_0(j)/rp/,/psi_j(x)/cr 
	     &=/wh f_m(x)+f_0(x)-/sum_<j=0>^m a_0(j)/,/psi_j(x).&(5.50)>$$	
It may sometimes be reasonable to use $f_0(x)=N_d(/hm,/hs)(x)$, with $/hm,/hs$ defined in (3.15)--(3.16), as the initial guess of $f(x)$.

Another possibility is to write
$$f(x)=f_0(x)/,h(x),/,x/in/cx,/eqno(5.52)$$
and estimate $h(x)$. The following procedure may be successful
when the quotient $h(x)=f(x)/sla f_0(x)$ varies less (around 1) than $f(x)$.
If $f_0(x)$ is a good estimate already, then perhaps only
few terms are needed in an orthogonal expansion estimator
of $h(x)$.  Writing $h(x)=/sum_<j=0>^/infty c(j)/,/psi_j(x)$, where
$$c(j)=/int/psi_j(x)/,h(x)/,d/mu=E_f/lb/psi_j(X_t)/sla f_0(X_t)/rb$$
is the best coefficient in terms of $/int/lbrace h-/sum_<j=0>^m c(j)
/psi_j/rbrace^2/,
d/mu$, an estimator of the form
$$/eqalign</wh h_m(x)&=/sum_<j=0>^m /wh c(j)/psi_j(x),/cr
           /wh c(j)&=/ndel/sum_<t=1>^n /psi_j(X_t)/sla f_0(X_t),/cr>$$
suggests itself.  Rules for the determination of $m$ based
on the observations $X_1/upto X_n$ may be developed in a way similar to those
discussed for the original $/wh f_m(x)$.

In some ways (5.52) is a more natural starting point than (5.49). The
optimal Bayes rule directly involves ratios of densities, for example.
Some of the practical procedures developed in 5.4 below expand on these
ideas.
/note<Remark 2.>  We have limited ourselves to the discussion
of estimators of the type $/wh f_m(x)=/sum_<j=0>^m /wh a(j)/psi_j(x)$. A 
more general estimator is
$$/widetilde<f>(x)=/sum_<j=0>^/infty/lam_j/,/wh a(j)/,/psi_j(x)$$ 
where $/lb/lam_j;/,j/geq0/rb$ is a sequence of weights.  Some decrease in 
$/mise$ can be achieved.  The best theoretical values for
the weights may be shown to be
$$/lam_j=a(j)^2/sla/lb a(j)^2+/ndel/sigma(j)^2/rb,$$
and these may again be estimated.  It is seen that little
can be gained if $n$ is moderate or large, and having
automatic pattern recognition in mind we stick to the
classical ones, noting the computational difficulties
that would result if $/widetilde<f>(x)$ above were to be used.

Similar remarks apply to another type of generalisation,
sometimes referred to as matrix density estimation, cf.
Crain (1976).
/subsection<5.3.C><Estimating a $d$-dimensional density>
The general theory of 5.3.B of course applies to the 
estimation of $d$-dimensional densities.  We will now look
into the practical difficulties encountered if the general
approch is used and the dimension $d$ is high.  It will 
be seen that the most general solutions may become 
impractical, suggesting that alternative methods, based
on approximations or assumptions about the unknown density,
should be worked out.

Let $X_1/upto X_n$ be i.i.d.~with unknown density $f(x)=f(/rx_1/upto/rx_d)$
w.r.t.~Lebesgue measure $/mu$ in a region $/cx$ in Euclidian $d$-space $/rr^d$.
The first practical difficulty one encounters is the 
choice of an orthonormal complete set $/lb/psi_j(x);/,j/geq0,/,x/in/cx/rb$ in 
$/cl^2(/cx,/mu)$. While such a set always exists, it may be impossible to 
obtain one in an explicit manner, if the region $/cx$ is anything
but well-structured.  The easy cases seem to be the 
class of examples where
$$/cx=/prod_<i=1>^d /cx_i$$
and each $/cx_i$ is an interval, bounded or unbounded.  If
$/lb/psi_<i,j>(/rx_i);/,j/geq0,/,/rx_i/in/cx_i/rb$ is a complete orthonormal 
set for $/cx_i,/ i=1/upto d$, then
$$/lb/psi_<j_1/upto j_d>(/rx_1/upto/rx_d);/,j_i/geq0,/,/rx_i/in/cx_i,/,
i=1/upto d/rb$$  
may be seen to constitute a complete orthonormal set for
$/prod_<i=1>^d/cx_i$, where
$$/psi_<j_1/upto j_d>(/rx_1/upto /rx_d)=/prod_<i=1>^d/psi_<i,j_i>(/rx_i),
/,x/in/cx./eqno(5.52)$$

Very often there is a </it bounded//> region within which all $X$'s are 
guaranteed  to fall, for all practical purposes. In such cases intervals
$/lbrack a_i,b_i/rbrack$ can be found such that 
$$a_i/leq/rX_i/leq b_i,/ i=1/upto d/eqno(5.53)$$
always, writing $X=(/rX_1/upto /rX_d)'$. Hence it is easy to
put up a set of orthonormal functions on $/prod_<i=1>^d/lbrack a_i,b_i/rbrack$,
which however may be larger than the feature space $/cx$. The approach
described below may be used to produce a curve $/wh f(x)$ for $x/in/prod_<i=1>^d/lbrack a_i,b_i/rbrack$, after which a more appropriate density estimator 
is obtained by defining
$$/dhat f(x)=/cases</wh f(x)&if $x/in/cx$,/cr
	            0&if $x/not/in/cx$./cr>$$

Vectors $X$ fulfilling (5.53) may be linearly transformed
to $/lbrack0,1/rbrack^d$, viz. $/rX_i'=(/rX_i-a_i)/sla(b_i-a_i),/ i=1/upto d$.
We assume such a transformation has been performed, and consider
training data $X_1/upto X_n$ in $/lbrack0,1/rbrack^d$.
Let $/lb/psi_0(x),/psi_1(x),/ldots/rb$ be a complete orthonormal set in
$/lbrack0,1/rbrack$ with $/psi_0(x)=1$. Important examples are 
constituted by $/lbrace 1,/sqrt<2>/cos(/pi s)/allowbreak,/sqrt<2>/cos$/break$(2/pi s),/ldots/rbrace$,
the normalised Legendre polynominals, and the normalised  Chebyshev
polynomials (Abramowitz~and Stegun, 1964, Chapter 22).
Further trigonometric constructions are in Kron/-mal~and Tarter (1968).

An orthogonal expansion of the unknown density $f$ takes the form
$$/eqalignno<
f(x)&=/sum_<j_1=0>^/infty/cdots/sum_<j_d=0>^/infty a(j_1/upto j_d)/,
		/psi_<j_1>(/rx_1)/ldots/psi_<j_d>(x_d),&(5.54)/cr
a(j_1/upto j_d)&=/int_0^1/cdots/int_0^1
		/psi_<j_1>(/rx_1)/cdots/psi_<j_d>(/rx_d)/,f(x)/,dx.&(5.55)/cr
>$$
The orthogonal expansion estimator is 
$$/wh f_<m_1/upto m_d>(x)=/sum_<j_1=0>^<m_1>/cdots
/sum_<j_d=0>^<m_d> /ha(j_1/upto j_d)/,/psi_<j_1>(/rx_1)/cdots/psi_<j_d>(/rx_d)/eqno(5.56)$$
where
$$/wh a(j_1/upto j_d)=/ndel/sum_<t=1>^n
	/psi_<j_1>(/rX_<t,i>)/cdots/psi_<j_d>(/rX_<t,d>),/eqno(5.57)$$ 
writing $X_t=(/rX_<t,i>/upto/rX_<t,d>),/ t=1/upto n$.

It remains to provide a rule for the determination of $m_1/upto m_d$. This 
can be done in various ways, the following proposal is but one
possibility.

Introduce $/uj=(j_1/upto j_d)$ for ease of notation. In a sense all
terms $/wh a/lp/uj/rp /,/psi_<j_1>(/rx_1)/cdots/allowbreak/psi_<j_d>(/rx_d)$
having the same value of $/sum_<i=1>^d j_i=j$ are of the same order. Let
$$/cl_j=/lb/uj:/sum_<i=1>^d j_i=j/rb,/,j/geq0,/eqno(5.58)$$
so that $f(x)=/sum_<j=0>^/infty/sum_</uj/in/cl_j> a/lp/uj/rp /psi_<j_1>(/rx_1)/cdots/psi_<j_d>(/rx_d)$.
A natural estimator is now 
$$/wh f_m(x)=1+/sum_<j=1>^m /sum_</uj/in/cl_j>/wh a/lp/uj/rp/psi_<j_1>(/rx_1)/cdots/psi_<j_d>(/rx_d),/eqno(5.59)$$
where $/psi_0(x)$ is taken to be the constant function $1$, for
convenience.
This estimate takes terms of order $/leq m$ into account.

A rule for  choosing $m$ can be given along the lines of 5.3.B. Let
$$/eqalign</mise(m)&=E/int_0^1/cdots/int_0^1(/wh f_m(x)-f)^2/,dx/cr
&=/sum_</uj/not/in/cl_m>a/lp/uj/rp^2+/ndel/sum_</uj/in/cl_m>/sigma/lp/uj/rp^2./cr>$$
Then the increase in expected loss, as measured by the integrated
squared error, when passing from $/wh f_<m-1>$ to $/wh f_m$, is
$$/eqalign<D(m)&=/mise(m)-/mise(m-1)/cr
&=/sum_</uj/in/cl_m>/lb/ndel/sigma/lp/uj/rp^2-a/lp/uj/rp^2/rb,/cr>$$
which may be estimated in an unbiased way by
$$/wh D(m)=/sum_</uj/in/cl_m>/lb<2/over n>S/lp/uj/rp ^2-/wh a/lp/uj/rp ^2/rb.
/eqno(5.60)$$
Here $/wh a/lp/uj/rp $ is defined in (5.57) whereas
$$/displaylines<
/sigma/lp/uj/rp ^2=/int_0^1/cdots/int_0^1 /lb/psi_<j_1>(/rx_1)/cdots/psi_<j_d>(/rx_d)-/wh a/lp/uj/rp /rb^2/,f(x)/,dx,/cr
S/lp/uj/rp ^2=<1/over n-1>/sum_<t=1>^n/lb/psi_<j_1>(/rX_<t,1>)/cdots/psi_<j_d>(/rX_<t,d>)-/wh a/lp/uj/rp /rb^2./cr>$$
The rule proposes $/wh f_<m_0>(x)$ as estimator, where
$m_0+1$ is the first $m$ having $/wh D(m)$ positive.

As an illustration, consider
$$/eqalign<
/wh f_3(x)&=1+/sum_<i=1>^d /wh a(0/upto1/upto0)/psi_1(/rx_i)/cr
&/qquad+/sum_<i=1>^d/wh a(0/upto2/upto0)/psi_2(/rx_i)^2+/sum_<i/lt j>/wh a(0/upto1/upto1/upto0)/psi_1(/rx_i)/psi_1(/rx_j)/cr
&/qquad+/sum_<i=1>^d/wh a(0/upto3/upto0)/psi_3(/rx_i)^2+/sum_<i/not=j>/wh a(0/upto2/upto1/upto0)/psi_2(/rx_i)^2/,/psi_1(/rx_j)/cr
&/qquad+/sum_<i/lt j/lt k>/wh a(0/upto1/upto1/upto1/upto0)/psi_1(/rx_i)/psi_1(/rx_j)/psi_1(/rx_k)./cr>$$
The number of coefficients $/wh a(j_1/upto j_d)$ that needs to be stored for 
the repeated computation of $/wh f_3(x)$ is quite large, for example, this 
number is 83 when the dimension is 6. A further drawback is that at most three
components in $x=(/rx_1/upto /rx_d)$ are involved in the 
individual terms; we have to use at least $/wh f_d(x)$ to see the term
$/psi_1(/rx_1)/cdots/psi_1(/rx_d)$. This may indicate that the density will 
not be very well estimated unless a very high number of terms
are used, in which case an enormous amount of training data is needed. Adding
this to  problems related  to a possibly difficult implementation and
long computing time, we see that more practical, but perhaps not totally
nonparametric methods are called for, when high-dimensional feature
vectors are used.
/note<Remark 3.> The assumption pertaining to the boundedness of the feature 
space $/cx$ was  made partly for concreteness and mainly because it is met so often
in practice. The arguments are easily seen to go through,
however, for a general $/cx=/prod_<i=1>^d /cx_i$ and 
orthonormal functions as in (5.52). In particular if $/cx=/rr^d$, then the set 
of normalised Hermitian polynomials, multiplied by $/exp(-/halv x^2)$, is 
the most usual choice of an orthonormal basis, see 5.4.F below.
/section<5.4><Extended theory, and some practical><solutions in high-dimensional space>
The general orthogonal expansion density estimators (5.56), (5.59) were seen to
be awkward and perhaps impractical in the case of a high dimension $d$. 
Whether they are useful at all depends of course upon the situation at 
hand, the densities that are to be estimated, the amount of available 
training data, whether the shape of a histogram invites a preliminary 
rough parametric guess $f_0$ to be made so that the suggestions of Remark 
1 may be utilised, etc. In any case some more convenient and easily 
implementable procedures should be devised. 

A general point that could be made is that many of the nonparametric 
density estimation methods are in a sense ``</it too//> nonparametric'',
in that they do not dare assume anything at all about the $f$ at hand, 
whereas most data sets met in practice, even ``dirty'' ones, do display 
</it some//> reasonable structure. The price paid by being too pessimistic, 
as opposed to taking a modest amount of cooperation from Nature for
granted, may be too high, in spite of guaranteed optimality for (possibly 
enormous) large training sets, if the dimension is higher than say four.

One possible way of toning down this curse of dimensionality is to find 
some transformation to approximately independent sets of components. 
Assume $n$ realisations of a six-dimensional feature vector 
$X=(</rm X>_1/upto </rm X>_6)$ are observed, and imagine a ``Transformation 
Pursuit'' algorithm is put to work and produces a smooth transformation 
to $Y=(</rm Y>_1/upto </rm Y>_6)$, where $(</rm Y>_1,</rm Y>_2,</rm Y>_3)$ 
and $(</rm Y>_4, </rm Y>_5,</rm Y>_6)$ become 
practically independent. Then the densities for 
$(</rm Y>_1,</rm Y>_2,</rm Y>_3)$ and 
$(</rm Y>_4,</rm Y>_5,</rm Y>_6)$ could be estimated using some
of the general methods 
described above, in a computationally and statistically efficient way. 
The estimated density for $Y$ is the product of these, and finally the 
density $f$ for the original $X$ is obtained by inverse transformation.

Some of the solutions invented below are ``nonparametric but optimistic'' 
in that they cautiously rely on some rough initial description being 
reasonable or use devices like the transformation trick above.
/subsection<5.4.A><General theory>
Allow us to expand somewhat on the general theory of 5.3 before we 
describe concrete applications.

We wish to apply the technique of orthogonal expansions, and therefore 
in some cases start out transforming $X$'s into a convenient product 
space, and that will be $/lsb0,1/rsb^d$, $/lsb-1,1/rsb^d$, or $/rr^d$ in 
our applications. The methods are intended for reasonably large sample 
sizes, and we sometimes (formally) ignore the sampling variability that 
is present in such an initial transformation $X/to Y$; see 5.4.G below, 
however.

Forget, therefore, the $X$'s for a moment, and study $Y_1/upto Y_n$. Let 
$g_0(y)$ be some initial estimate of the density, e.g.~the normal density 
with estimated mean and covariance matrix, and write the correct density 
as
$$/eqalignno<
g(y)&=g_0(y)/,h(y)/cr
&=g_0(y)/lb c(0)/psi_0(y)+c(1)/psi_1(y)+/cdots/rb,&(5.60)/cr>$$
where 
$$/int/psi_j/psi_/ell/omega/,dy=/delta_<j,/ell>,$$
i.e.~$/psi_0,/psi_1,/ldots$ are orthonormal w.r.t.~the weigth function 
$w(y)$. Minimising 
$$/int/lb/sum_0^m c(j)/psi_j-h/rb^2 w/,dy$$
gives
$$/eqalignno<
c(j)&=/int h/psi_j w/,dy/cr
&=/int g</psi_j w/over g_0>/,dy/cr
&=E/psi_j(Y)/,w(Y)/sla g_0(Y),&(5.61)/cr>$$
a natural estimate for which is
$$/wh c(j)=<1/over n>/sum_<t=1>^n /psi_j(Y_t)/,w(Y_t)/sla 
g_0(Y_t)./eqno(5.62)$$
Accordingly
$$/wh g_m(y)=g_0(y)/sum_<j=0>^m/wh c(j)/psi_j(y)/eqno(5.63)$$
is proposed.

Assume $h$ is within reach of the $/psi_j$-system, and consider 
$$/eqalign<
</rm ISE>(m)&=/int(/wh h_m-h)^2 w/,dy/cr
&=/int(/wh g_m-g)^2<w/over g_0^2>/,dy/cr
&=/int/lsb/sum_0^m/lb/wh c(j)-c(j)/rb/psi_j-/sum_<m+1>^/infty 
c(j)/psi_j/rsb^2 w/,dy/cr
&=/sum_0^m/lb/wh c(j)-c(j)/rb^2+/sum_<m+1>^/infty c(j)^2/cr>$$
and the expected integrated squared error
$$/eqalignno<
</rm MISE>(m)&=/sum_0^m<1/over n>/tau(j)^2+/sum_<m+1>^/infty c(j)^2/cr
&=/sum_0^/infty c(j)^2-/sum_0^m/lb c(j)^2-<1/over 
n>/tau(j)^2/rb.&(5.64)/cr>$$
The $m$ we should like to use (among $0,1/upto m_<up>$, say) is the one 
maximising
$$A(m)=/sum_0^m/lb c(j)^2-<1/over n>/tau(j)^2/rb./eqno(5.65)$$
Here $/tau(j)^2=</rm Var>/lb/tau_j(Y)/,w(Y)/sla g_0(Y)/rb$, an unbiased
estimator for which is provided by
$$/wh/tau(j)^2=<1/over n-1>/sum_<t=1>^n/lb/psi_j(Y_t)<w(Y_t)/over g_0(Y_t)>- /wh 
c(j)/rb^2./eqno(5.66)$$
One finds that
$$/wh A(m)=/sum_0^m/lb/wh c(j)^2-<2/over n>/wh/tau(j)^2/rb/eqno(5.67)$$
is an unbiased estimator of $A(m)$. Accordingly we propose to use (5.64) 
with $m$ chosen to maximise $/wh A(m)$, among $0,1/upto m_<up>$. (There 
is seldom need to go beyond $m_<up>=10$ for one-dimensional cases, for 
example.)
/subsection<5.4.B><Solution 1: The multiplicative cosine expansion method>
Matters simplify dramatically if one dares assume independence between the 
components of the feature vector $X=(X_1/upto X_d)$. Then
$$f(x_1/upto x_d)=f_1(x_1)/,/cdots/,f_d(x_d),$$
say, and the individual univariate component densities may be estimated 
with good precision by one of the methods of 5.2 or 5.3. The consequences 
of falsely assuming independence in such a direct manner may be grave, 
however, and the following method will tend to produce better results.

Compute mean vector $/hm$ and empirical covariance matrix $/hs$ from the 
data set $X_1/upto X_n$, as in (3.15)---(3.16). Find the spectral 
decomposition
$$/wh P/hs/wh P'=/wh D=</rm diag>/lb/hl_1/upto /hl_d/rb,/eqno(5.68)$$
where $/hl_1/le/cdots/le/hl_d$ are the eigenvalues and $/wh P'$ has the 
normalised eigenvectors. $/wh P$ is orthonormal (unitary), i.e.~$/wh P/wh 
P=I=/wh P'/wh P$.  Define new variables 
$$Z_t=/wh P(X_t-/hm),/ t=1/upto n/eqno(5.69)$$
These have mean vector $0$ and covariance matrix
$$/hs_z=<1/over n-1>/sum_<t=1>^n(Z_t-/bar Z)(Z_t-/bar Z)'=/wh D,$$
by construction. Thus $Z_1/upto Z_n$ are pseudo-data with uncorrelated 
components. Admitting that there is a long way from absence of 
correlation to independence, and also that $Z_1/upto Z_n$ are not 
independent among themselves (since $/wh P$ depends upon the full set 
$X_1/upto X_n$), we may nevertheless approximate the density $q(z)$ for a 
$Z_t$ as follows:
$$q(z_1/upto z_n)=q_1(z_1)/cdots q_d(z_d).$$
Next estimate the component densities using the method presented in 
5.3.A, after a linear transformation has taken place to produce data in 
$/lsb 0,1/rsb$. The result is
$$/wh q_i(z_i)=/wh g_i/lp<z_i-a_i/over b_i-a_i>/rp<1/over 
b_i-a_i>,/,z_i/in/lsb a_i, b_i/rsb,/eqno(5.70)$$
where for $i=1/upto d$
$$/wh g_i(y_i)=1+/sum_<j=1>^<m_i>/wh c_i(j)/cos(j/pi 
y_i),/,y_i/in/lsb0,1/rsb,/eqno(5.71)$$
$$/wh c_i(j)=<2/over n>/sum_<t=1>^n/cos/lp j/pi<Z_<t,i>-a_i/over b_i - 
a_i>/rp, j/le1./eqno(5.72)$$
It is assumed here that
$$a_i/le Z_<t,i>=/wh P_<(i)>(X_t-/hm)/le b_i, </rm all/ > t,/eqno(5.73)$$
writing $/wh P_<(i)>$ for the $i$'th row of the matrix $/wh P$.

The final estimate for $f(x)$ becomes
$$/eqalignno<
/wh f(x)&=/wh f(x_1/upto x_d)/cr
&=/wh q/lp/wh P(x-/hm)/rp/inpar</wh P>/cr
&=/wh q_1(z_1)/cdots/wh q_d(z_d)&(5.74)/cr>$$
for $z=/wh P(x-/hm)$ in the box $/prod_<i=1>^d/lsb a_i,b_i/rsb$, and zero 
outside of it. 

$a_i$ and $b_i$ of (5.73) should be chosen slightly to the left and right 
of the observed range of $Z_<1,i>/upto Z_<n,i>$. Optimal unbiased 
estimates for the underlying population parameters, under a simple 
uniform model, are 
$$/eqalignno<
a_i&=/min_t Z_<t,i>-<2/over n-1>/lp/max_t Z_<t,i>-/min_t Z_<t,i>/rp,/cr
b_i&=/max_t Z_<t,i>+<2/over n-1>/lp/max_t Z_<t,i>-/min_t 
Z_<t,i>/rp,&(5.75)/cr>$$
and there is no particular need to be more sophisticated. The stopping 
rule for the number of terms in (5.71) is to choose the $m_i$ that 
maximises
$$/hA_i(m)=1+/sum_<j=1>^m/lb/wh c_i(j)^2-<2/over n>/wh/tau_i(j)^2/rb,$$
$m=0,1/upto10$ (say), where
$$/wh/tau_i(j)^2=<1/over n-1>/sum_<t=1>^n/lb2/cos/lp j/pi<<Z_<t,i>-a_i/over 
b_i-a_i>>/rp-/wh c_i(j)/rb^2,$$
and where $/hA_i(0)=1$.

Since the coefficient estimates that are involved are based on a random 
transformation it is clear that a large sample size may be necessary in 
order for the method to be effective, i.e.~$/wh P$ ought to be stable, 
which may be difficult when the underlying eigenvalues are not well 
separated.

Let us descripe how this density estimation method works in the context 
of a classification problem. Below $k=1/upto K$ is index for classes, 
$i=1/upto d$ is coordinate (component) number, $j=0,1/upto m(k,i)$ is 
index for the expansion coefficient, and $t=1/upto n_k$ indexes the 
replicates in the $k$'th training set.

Based on $X_1^<(k)>/upto X_<n_k>^<(k)>$, compute $/hm_k$, $/hs_k$,
and $/wh P_k$ as in $/wh P_k/hs_k/wh P_k'=/wh D_k$. Next let
$Z_t^<(k)>=/wh P_k(X_t^<(k)>-/hm_k)$ and find
$$</rm Box>(k)=/prod_<i=1>^d/lsb a_<k,i>,b_<k,i>/rsb,/eqno(5.76)$$
inside which all $Z_t^<(k)>$ vectors fall, using (5.75).
Then evaluate $/wh c_<k,i>(j)$ and $/wh/tau_<k,i>(j)$ as mean and standard 
deviation for $2/cos(j/pi(Z_<t,i>^<(k)>-a_i)/sla
(b_<k,i>-a_<k,i>)),$ $t=1/upto n_k$, and compute the stopping 
criterion statistic
$$/wh B_<k,i>(m)=/sum_<j=1>^m/lb/wh c_<k,i>(j)^2-<2/over 
n_k>/wh/tau_<k,i>(j)^2/rb$$
along the way, $m=1/upto 10$ (and $/wh B_<k,i>(0)=0$).
Choose $m(k,i)$ to maximise it. The class density $f_k$ is estimated as
$$/wh f_k(x)=/prod_<i=1>^d/lb1+/sum_<j=1>^<m(k,i)>/wh c_<k,i>(j)/cos/lp
j/pi<Z_<k,i>-a_<k,i>/over b_<k,i>-a_<k,i>>/rp/rb<1/over 
b_<k,i>-a_<k,i>>$$
for $z_k=(z_<k,1>/upto z_<k,d>)'=/wh P_k(x-/hm_k)$ inside Box$(k)$ and is 
zero outside it.
/subsection<5.4.C><Solution 2: The multiplicative beta with cosine 
expansions> Above the initial $X_t$ vectors were transformed into 
$Y_t$'s in $/lsb0,1/rsb^d$, via an initial orthogonal transformation and 
$Y_<t,i>=(Z_<t,i>-a_i)/sla(b_i-a_i)$. The density for $Y_<t,i>$ (or 
rather, for the $i$-th component of a $Y$ obtained from a future $X$ 
using the transformation established from the training set) was then estimated
using cosine term corrections to the constant 1. This procedure works 
perhaps best if the corrections to be made are few and small, i.e.~if
the real density is not too far from being uniform. In many cases it </it 
is//> far from the uniform, however, having for example some clearly 
defined peak or perhaps a $J$ shape. While a cosine expansion can cope 
with such shapes by adding enough terms, this may lead to too much 
``noise'' in the estimates of coefficients or may require a too high 
sample size to be effective.

It should in many cases be more efficient, then, to start an expansion 
with a better initial ``guess'' of the shape than the uniform one. Such a 
procedure is described now.

The beta distribution with parameters $(/alpha,/beta)$ on $/lsb0,1/rsb$ 
has density
$$B(/alpha,/beta)y^</alpha-1>(1-y)^</beta-1>,$$
where 
$B(/alpha,/beta)=/Gamma(/alpha+/beta)/sla/Gamma(/alpha)/Gamma(/beta)$ and 
$/Gamma$ is the gamma function. It has mean $</alpha/over/alpha+/beta>$ 
and variance $</alpha/beta/over(/alpha+/beta)^2(/alpha+/beta+1)>$, and a 
natural way of matching a beta density to an observed set $Y_<1,i>/upto 
Y_<n,i>$ is to match moments, i.e.~to solve the equations
$$/eqalign<
</alpha/over/alpha+/beta>&=<1/over n>/sum_<t=1>^nY_<t,i>=/wh/xi_i,/cr
</alpha/beta/over(/alpha+/beta)^2(/alpha+/beta+1)>&=<1/over 
n-1>/sum_<t=1>^n/lp Y_<t,i>-/wh/xi_i/rp^2=/wh/tau_i^2./cr>$$
This leads to
$$/wh g_<0,i>(y)=B/lp/wh/alpha_i,/wh/beta_i/rp y^</wh/alpha_i-1> 
(1-y)^</wh/beta_i-1>/eqno(5.77)$$
as the initial estimate of the density for the $i$-th component of $Y$, 
where
$$/eqalign<
/wh/alpha_i&=/lb/wh/xi_i/lp1-/wh/xi_i/rp/sla/wh/tau_i^2-1/rb/wh/xi_i,/cr
/wh/beta_i&=/lb/wh/xi_i/lp1-/wh/xi_i/rp/sla/wh/tau_i^2-1/rb 
/lp1-/wh/xi_i/rp./cr>/eqno(5.78)$$
The maximum likelihood program could equally well be pursued. It consists 
in this case of matching expectations of $/log Y_<t,i>$ and 
$/log(1-Y_<t,i>)$ with their observed counterparts, i.e.~solving
$$/eqalign<
<1/over n>/sum_<t=1>^n/log Y_<t,i>&=/psi(/alpha)-/psi(/alpha+/beta),/cr
<1/over n>/sum_<t=1>^n/log(1-Y_<t,i>)&=/psi(/beta)-/psi(/alpha+/beta),/cr>$$
where $/psi=/Gamma'/sla/Gamma$, by some numerical technique starting from 
$/wh/alpha_i$, $/wh/beta_i$ above. While the maximum likelihood estimators 
have a superior performance </it if//> the beta model reasonably fits the data
it is not clear that they have when the beta model is wrong, i.e.~must be 
corrected on, as is our intention here, and we will usually be
content with (5.78). (See also Example 3 of Section 3.1.B.)

The class of beta densities is in some ways more flexible then the normal 
one, in that they can be non-symmetric and can have a $J$ shape and even 
$U$ shape for suitable choices of parameters.

Next let us follow the general theory outlined before Solution 1, 
choosing again $1$, $/sqrt<2>/cos(/pi y)$, $/sqrt<2>/cos(2/pi y)$, 
$/ldots$ as basis functions, orthogonal w.r.t.~$w(y)=1$ on $/lsb0,1/rsb$. 
The result is a generalisation of (5.71)---(5.72):
$$/wh g_i(y_i)=/wh g_<0,i>(y_i)/lb/wh c_i(0)+/sum_<j=1>^<m_i>/wh 
c_i(j)/sqrt<2>/cos(j/pi y_i)/rb,/eqno(5.79)$$
where
$$/eqalign<
/wh c_i(0)&=<1/over n>/sum_<t=1>^n 1/sla/wh g_<0,i>(Y_<t,i>),/cr
/wh c_i(j)&=<1/over n>/sum_<t=1>^n/sqrt<2>/cos(j/pi Y_<t,i>)/sla/wh 
g_<0,i>(Y_<t,i>)./cr>$$
The number of recommended terms $m_i$ is arrived at by inspecting
$$/hA_i(m)=/sum_<j=0>^m/lb/wh c_i(j)^2-<2/over n>/wh/tau_i(j)^2/rb$$
for $m=0,1/upto10$ (say) and choosing the maximising one, where
$$/eqalign<
/wh/tau_i(0)^2&=<1/over n-1>/sum_<t=1>^n/lb1/sla/wh g_<0,i>(Y_<t,i>)-/wh c_i(0)/rb^2,/cr
/wh/tau_i(j)^2&=<1/over n-1>/sum_<t=1>^n/lb/sqrt<2>/cos(j/pi Y_<t,i>)/sla/wh 
g_<0,i>(Y_<t,i>)-/wh c_i(j)/rb^2./cr>$$

Solution 1 is the special case where $/wh/alpha_i$ and $/wh/beta_i$ are 
set to $1$ instead of being allowed to fit the data in an initial 
description (5.77). The present generalisation should in most cases lead 
to lower values of $m_i$, and apart from aiming at greater predictive 
efficiency it therefore requires less computer time.
/subsection<5.4.D><Robustification>
There is a potential problem with the method above, if it is applied directly.
$/wh c_i(j)$ and $/wh/tau_i(j)$ are the mean and standard deviation for the set
$/sqrt<2>/cos(j/pi Y_<t,i>)/sla/wh g_<0,i>(Y_<t,i>),/ t=1/upto n$, and
the denominator $/wh g_<0,i>(Y_<t,i>)$ can become devastatingly small for
some but very few $Y_<t,i>$'s, so that $/wh c_i(j)$ and $/wh/tau_i(j)$ are in
danger of being blown up. Thus the estimator (5.79) has been seen to 
perform ridiculously in cases where $/wh/alpha_i$ and $/wh/beta_i$ have 
been large.

Thus a robustification is absolutely necessary. Several ways are possible, 
one could for example use robust estimates of 
$c_i(j)=/int_0^1/sqrt<2>/cos(j/pi y)<w(y)/over/wh g_<0,i>(y)> g_i(y) dy$ 
and $/tau_i(j)^2=</rm Var>/lb/sqrt<2>/cos(j/pi Y_<t,i>) w(Y_<t,i>)/sla/wh 
g_<0,i>(Y_<t,i>)/rb$ instead of $/wh c_i(y)$ and $/wh/tau_i(j)^2$. We recommend 
avoiding the very smallest ones of the observed $/wh g_<0,i>(Y_<t,i>)$ by 
adjusting them somewhat upwards, for example ``recoding'' the $r$ 
smallest ones to the $(r+1)$-st smallest value, where $r$ equal to the 
nearest integer to $0.03n$ has been an effective choice. (This could mean 
``Winsorising'' with 1.5 percent in each tail.) Another way of describing 
this robustification is that it amounts to using a slightly adjusted 
version $/overline g_<0,i>(y)$ of (5.77), pulling this initial density 
estimate up from zero in the tails, if either $/wh/alpha_i$ of 
$/wh/beta_i$ is greater than one.

If we reconsider the general theory outlined before Solution 1 it becomes 
clear that the potential trouble can be traced to the use of
$$/int/lp/wh h-h/rp^2 w dy=/int/lp/wh g-g/rp^2<w/over g_0^2> dy$$
as the goodness criterion, which punishes $/wh g$ to harshly for not 
being close to $g$ in the (tail) areas where $g_0$ is very small, and 
where very few observations are likely to fall. The problem can be 
overcome by a correspondingly chosen weight function $w(y)$, making 
$w/sla g_0^2$ less severe in the tails. This is done for Solution 3 
below, but is not always feasible, since natural orthogonal basis 
functions may be hard to construct. Thus pulling one's initial $g_0$ 
slightly up frome zero in extreme areas, as we did above, seems to be the 
generally best solution.
/subsection<5.4.E><Penalising higher order terms>
The stopping criterion we have used to determine the number of terms in 
an expansion is 
$$/hA(m)=/sum_<j=0>^m/lb/wh c(j)^2-<2/over n>/wh/tau(j)^2/rb,$$
in the general context of (5.60)---(5.67), and which is an unbiased 
estimator for
$$A(m)=/sum_<j=0>^m/lb c(j)^2-<1/over n>/tau(j)^2/rb.$$
The following phenomenon is often observed. The criterion function 
$/hA(m)$ that we want to maximise climbs sharply in the beginning up to 
some $m_1$, say, after which it continues to climb, but only very 
moderately, or maybe fluctuates somewhat until it reaches its maximum for 
$m_0$, say. If $/hA(m_0)$ is only slightly larger than $/hA(m_1)$, and 
$m_0$ is several steps beyond $m_1$, then perhaps $m_1$ should be used 
after all. Too many terms in an expansion means much ``estimation 
noise'', a less smooth curve, and also an increased likelihood for seeing 
negative values of $/wh g_m(y)=g_0(y)/wh h_m(y)$.

This motivates adjustments to the criteria $A(m)$ and $/hA(m)$. One has
$$/hA(m)=/sum_<j=0>^m<1/over n>/wh/tau(j)^2/lb<n/hc(j)^2/over 
/wh/tau(j)^2>-2/rb,$$
and a theoretical reason for the empirically observed phenomenon 
described above can be seen: if $c(j)$ really is zero, so that one's 
criterion shouldn't encourage including term no.~$j$, then $/hA$ 
increases for $j$ nevertheless, with probability
$$/Pr/lb<n/hc(j)^2/over /wh/tau(j)^2>/gt2/rb/doteq Pr/lb/chi_1^2/gt2/rb=0.157
/eqno(5.80)$$
This means that $/hA$ often is tempted to pick too many terms.

A version that is less encouraging for higher order terms is
$$/eqalignno<
/hA_/eps(m)&=/sum_<j=0>^m<1/over n>/wh/tau(j)^2/lb<n/hc(j)^2/over 
/wh/tau(j)^2>-(2+/eps j)/rb/cr
&=/sum_<j=0>^m/lb/wh c(j)^2-<2+/eps j/over n>/wh/tau(j)^2/rb,&(5.81)/cr>$$
where $/eps$ is some positive number. Its expectation is 
$$/eqalignno<
A_/eps(m)&=/sum_<j=0>^m/lb c(j)^2-<1+/eps j/over n>/tau(j)^2/rb/cr
&=/sum_<j=0>^m<1/over n>/tau(j)^2/lb<nc(j)^2/over/tau(j)^2>-(1+/eps 
j)/rb,&(5.82)/cr>$$
and is related to the expected integrated and weighted squared error
$$/mise(m)=E/int/lp/wh g_m-g/rp^2<w/over g_0^2>dy=E/int/lp/wh h_m-h/rp^2 w dy$$
by
$$A_/eps(m)=/sum_<j=0>^/infty c(j)^2-/lb/mise(m)+/sum_<j=0>^m<1/over 
n>/eps j/tau(j)^2/rb./eqno(5.83)$$
Thus choosing $m=m_0$ to maximise $/hA_/eps(m)$ aims at minimising the 
``penalised'' MISE criterion
$$/eqalignno<
/mise_/eps(m)&=/mise(m)+/eps/sum_<j=0>^m<1/over n>j/tau(j)^2/cr
&=E/lsb/int/lp/wh g_m-g/rp^2<w/over g_0^2>dy+/eps/sum_<j=0>^m j/lb/wh c(j)- 
c(j)/rb^2/rsb.&(5.84)/cr>$$

One could start from scratch, employing the penalised MISE criterion 
(5.84), and devise techniques that found the smoothing parameters $/eps$ 
and $m$ from data, perhaps using cross validation techniques. A 
satisfactory solution is however one that settles on some small, fixed 
$/eps$ and then finds $m=m_0$ to maximise $/hA_/eps(m)$. Using 
$/eps=0.20$, for example, leads to a ``false encouragement'' for term $j$ 
with probability
$$/Pr/lb<n/,/hc(j)^2/over/wh/tau(j)^2>/gt2+/eps j/vert c(j)=0/rb/doteq 
/Pr/lb/chi^2_1/gt2+<j/over5>/rb$$
which is $4.6/

Of course another slowly increasing penalty function than $/eps j$ could 
have been chosen, and it should not be allowed to increase as fast as 
linearly after the first ten or so terms. The penalty could for example 
be brought to a halt after a certain number of terms.

/pageinsert/vbox to5truecm<>
/note<Example.></it/ The figure above shows the (robustified) beta times cosine
expansion estimator at work in a one-dimensional example. A certain
feature component is computed for each of $400$ hand-written ``7''
symbols. These are transformed to data $Y_1/upto Y_<400>$ in the unit
interval $/lsb0,1/rsb$.

The initial beta density estimate has $/wh/alpha =0.780$ and $/wh/beta =1.357$.
The table below displays coefficients $/wh c(j)$ in the cosine expansion
for the correction factor; standard deviation estimates $/wh/tau(j)$;
the test ratios $z(j)=/sqrt<n>/wh c(j)/sla/wh/tau(j)$; and finally
$/hA(m)$ and $/hA_<.20>(m)$ (see 5.4.E). The ``Solution 2'' method
chooses the ninth order expansion as the best density estimator, and
so does in this case also the modified version of 5.4.D. The $z(j)$
column suggests that terms of order $0$, $2$, $3$, $5$, $6$, $8$ are
most important.

Shown in the figure are the initial beta density estimate; the third
order estimator; the winning ninth order estimator; and a histogram of
the data.>
/smallskip
/noindent/hbox to /hsize</hfil/vbox</halign<&/quad/hfil#/cr
order/hfil&$/wh c(j)$/hfil&$/wh/tau(j)$/hfil&$z(j)$/hfil&$/hA(m)$/hfil%
&$/hA_<.20>(m)$/hfil/cr
/noalign</smallskip>
 0& 0.903702&0.40522&44.6029&0.81586&0.81586/cr
 1& 0.003185&0.80667& 0.0790&0.81261&0.81229/cr
 2&-0.343390&0.83960&-8.1799&0.92701&0.92597/cr
 3& 0.270378&1.02863& 5.2570&0.99482&0.99220/cr
 4&-0.051568&1.03711&-0.9944&0.99210&0.98733/cr
 5& 0.194279&0.94870& 4.0957&1.02534&1.01833/cr
 6& 0.246293&0.97465& 5.0540&1.08125&1.07139/cr
 7&-0.077885&0.97732&-1.5938&1.08254&1.06933/cr
 8& 0.147567&1.01424& 2.9099&1.09918&1.08185/cr
 9& 0.107503&0.99048& 2.1707&1.10583&1.08409/cr
10& 0.023575&0.99334& 0.4747&1.10145&1.07478/cr>>/hfil>
/vfill
/endinsert
/subsection<5.4.F><Solution 3: The third order correction to normality> 
Below a method is described that takes the usual fitted normal density 
$N_d(/hm,/hs)(x)$ as a starting point and then performs a ``third order 
correction''.

Compute once more $/hm$, $/hs$, and the orthonormal transformation matrix 
$/wh P$ satisfying $/wh P/hs/wh P'=/wh D=/diag/lbrace/wh/lambda_1 
/upto/wh/lambda_d/rbrace$. $/hs^<-1/sla2>=/wh P'/wh D^<-1/sla2>/wh P$ is the 
unique symmetric inverse square root if $/hs$. Transform this time to 
pseudo-data 
$$Y_t=/hs^<-1/sla2>(X_t-/hm),/ t=1/upto n;/eqno(5.85)$$
this is sometimes referred to as ``sphering the data''. These have mean 
zero and the identify $I$ for empirical covariance matrix. Using the 
normal as one's initial description or estimate of the density $f$ for 
$X$ amounts to using
$$/eqalign<
g_0(y)&=N_d(0,I)(y)/cr
&=/phi(y_1)/cdots/phi(y_d)/cr
&=(2/pi)^<-d/sla2>/exp/lp-/halv/sum_<i=1>^d y_i^2/rp/cr>$$
for $Y$.

We want again to employ the general machinery presented before Solution 
1. A convenient choice for the weight function $w$ is $g_0$ itself. 
Products
$$/psi_<j_1/upto j_d>(y_1/upto y_d)=H_<j_1>(y_1)/cdots H_<j_d>(y_d)$$
of Hermitian polynomials become orthogonal w.r.t.~$w(y)=g_0(y)$. We use 
the normalised versions right away, i.e.
$$/eqalign<
H_0(y)&=1,/cr
H_1(y)&=y,/cr
H_2(y)&=<1/over/sqrt<2>>/lp y^2-1/rp,/cr
H_3(y)&=<1/over/sqrt<6>>/lp y^3-3y/rp,/cr
H_4(y)&=<1/over/sqrt<24>>/lp y^4-6y^2+3/rp,/cr
H_5(y)&=<1/over/sqrt<120>>/lp y^5-10y^3+15y/rp,/ldots/cr>/eqno(5.86)$$
Further Hermitian polynomials can be found e.g.~in Abramowitz and Stegun
(1964, Ch.~22).

Write in what follows $e_i=(0/upto 1/upto 0)$ for the $i$-th basis 
vector, so that we are allowed to write $e_i+e_j+e_k$ instead of $(0/upto 
0,1,0/upto0,1,0/upto0,1,0/upto0)$, etc. A correction expansion to order 
$m$ terms to the normal initial estimate $g_0$ is 
$$/wh g_m(y)=g_0(y)/wh h_m(y),$$
where $/wh h_m(y)$ is the estimated version of 
$$/eqalign<
h_m(y)&=c(0/upto0)+/sum_<i=1>^d c(e_i)H_1(y_i)/cr
&/quad+/sum_<i=1>^d c(2e_i)H_2(y_i)+/sum_<i/lt j>c(e_i+e_j)H_1(y_i) 
H_1(y_j)/cr
&/quad+/sum_</scriptstyle i,l,s;/atop </rm order> j_1+j_2+j_3=3>
c(e_i+e_l+e_s)H_<j_1>(y_i)H_<j_2>(y_l)H_<j_3>(y_s)/cr
&/quad+/cdots+/sum_</scriptstyle i_l/upto i_m;/atop </rm order> =m>
c(e_<i_1>+/cdots+e_<i_m>)H_<j_1>(y_<i_1>)/cdots H_<j_d>(y_<i_d>)./cr>$$
But from (5.62), $c(0/upto0)=1$, and 
$$/eqalign<
c(e_i)&=/int H_1(y_i)g(y)dy=/int y_i g(y)dy/cr
&/wh =/ /wh c(e_i)=<1/over n>/sum_<t=1>^n Y_<t,i>=0,/cr
c(2e_i)&=/int H_2(y_i)g(y)dy=/int <1/over/sqrt<2>>(y_i^2-1)g(y)dy/cr
&/wh=/ /wh c(2e_i)=<1/over/sqrt<2>><1/over n>/sum_<t=1>^n/lb/lp 
Y_<t,i>/rp^2-1/rb=0,/cr
c(e_i+e_j)&=/int H_1(y_i) H_1(y_j)g(y)dy=/int y_i y_j  g(y) dy/cr
&/wh=/ /wh c(e_i+e_j)=<1/over n>/sum_<t=1>^n Y_<t,i>Y_<t,j>=0,/cr>$$
i.e.~the first and second order terms vanished.

The first non-trivial term is the third. One has
$$/eqalignno<
c(3e_i)&=/int<1/over/sqrt<6>>(y_i^3-3y_i) g(y)dy/cr
&/wh=/ <1/over/sqrt<6>>/wh/gamma_i,/ /wh/gamma_i=<1/over 
n>/sum_<t=1>^n(Y_<t,i>)^3,/cr
c(2e_i+e_j)&=/int<1/over/sqrt<2>>(y_i^2-1) y_j g(y)dy/cr
&/wh=/ <1/over/sqrt<2>>/wh/delta_<i,j>,/ /wh/delta_<i,j>=<1/over 
n>/sum_<t=1>^n(Y_<t,i>)^2Y_<t,j>,/cr
c(e_i+e_j+e_l)&=/int y_i y_j y_l g(y) dy/cr
&/wh=/ /wh/eps_<i,j,l>=<1/over n>/sum_<t=1>^n Y_<t,i>Y_<t,j>Y_<t,l>.&(5.87)/cr>$$
Thus
$$/eqalign<
/wh h_3(y)&=1+/sum_i 
<1/over/sqrt<6>>/wh/gamma_i<1/over/sqrt<6>>(y_i^3-3y_i)/cr
&/quad+/sum_<i/not=j><1/over/sqrt<2>>/wh/delta_<i,j><1/over/sqrt<2>>(y_i^2-1) 
y_j+/sum_<i/lt j/lt l>/wh/eps_<i,j,l>/, y_i y_j y_l,/cr>$$
and an interesting semi-parametric density estimator evolves:
$$/eqalignno<
/wh g_3(y)&=g_0(y)/wh h_3(y),/cr
/wh f_3(x)&=/wh g_3/lp/hs^<-1/sla2>(x-/hm)/rp/invert</hs^<-1/sla2>>/cr
&=N_d/lp/hm,/hs/rp(x)/bigg/lbrace1+/sum_i<1/over6>/wh/gamma_i(y_i^3-3y_i)/cr
&/quad+/sum_<i/not=j>/halv/wh/delta_<i,j>(y_i^2-1)y_j+/sum_<i/lt j/lt 
l>/wh/eps_<i,j,l>/, y_i y_j y_l/bigg/rbrace,&(5.88)/cr>$$
in which $y=y(x)=/hs^<-1/sla2>(x-/hm)$ is computed first.

Is it worthwhile to include the third order correction factor? Of course 
the added variability introduced by the coefficients $/wh/gamma_i$, 
$/wh/delta_<i,j>$, $/wh/eps_<i,j,l>$ may lead to a less precise estimate 
than the rough normal description, and we need a diagnostic tool to 
decide for us. The loss incurred when using the proposed (5.88) is
$$/eqalign<
/ise(3)&=/int/lb/wh h_3(y)-h_3(y)/rb^2 w(y) dy/cr
&=/int/Bigg/lbrack/sum_i<1/over6>/lp/wh/gamma_i-/gamma_i/rp/lp y_i^3-3y_i/rp/cr
&/quad+/sum_<i/not=j>/halv/lp/wh/delta_<i,j>-/delta_<i,j>/rp(y_i^2-1)y_j 
+/sum_<i/lt j/lt l>(/wh/eps_<i,j,l>-/eps_<i,j,l>) y_i y_j y_l/cr
&/quad-/sum_</scriptstyle </rm order> /ge4>c(j_1/upto j_d)H_<j_1>(y_1)
/cdots H_<j_d>(y_d)/Bigg/rbrack^2 g_0(y) dy/cr
&=/sum_i <1/over6>(/wh/gamma_i-/gamma_i)^2+/sum_<i/not=j> 
/halv(/wh/delta_<i,j>-/delta_<i,j>)^2/cr
&/quad+/sum_<i/lt j/lt l>(/wh/eps_<i,j,l>-/eps_<i,j,l>)^2+/sum_<
/scriptstyle </rm order>/ge4>c(j_1/upto j_d)^2,/cr>$$
where $/gamma_i=/int y_i^3 g(y) dy$, $/delta_<i,j>=/int y_i^2 y_j g(y) 
dy$, and $/eps_<i,j,l>=/int y_i y_j y_l g(y) dy$. On the other hand, the 
loss by being non-sophisticated and sticking to the simple $g_0(y)$ for
$Y$ is
$$/eqalign<
/ise(0)&=/int/lb1-h(y)/rb^2 w(y) dy/cr
&=/sum_</scriptstyle </rm order>/ge3>c(j_i/upto j_d)^2./cr>$$
Denoting the expected integrated squared error loss by MISE,  one has
$$/mise(3)/lt/mise(0)$$
if and only if
$$/eqalign<
&/sum_i<1/over6>Var/nils<Burde disse kanskje v{re i "roman"?>/wh/gamma_i+/sum_<i/not=j>/halv Var/wh/delta_<i,j>+
/sum_<i/lt j/lt l>Var/wh/eps_<i,j,l>/cr
&/quad/lt/sum_</scriptstyle </rm order>=3>c(j_1/upto j_d)^2=/sum_i 
<1/over6>/gamma_i^2+/sum_<i/not=j>/halv/delta_<i,j>^2+/sum_<i/lt j/lt 
l>/eps_<i,j,l>^2./cr>$$
A natural criterion is therefore the following: Compute
$$/eqalignno<
T&=/sum_i<1/over6>/lb/wh/gamma_i^2-<2/over n>est. var (Y_i^3)/rb/cr
&/quad+/sum_<i/not=j>/halv/lb/wh/delta_<i,j>^2-<2/over n>est. var(Y_i^2 
Y_j)/rb/cr
&/quad+/sum_<i/lt j/lt l>/lb/wh/eps_<i,j,l>^2-<2/over n>est. var(Y_i 
Y_j Y_l)/rb,&(5.89)/cr>$$
where
$$/eqalign<
est. var(Y_i^3)&=<1/over n-1>/sum_<t=1>^n/lp Y_<t,i>^3-/wh/gamma_i/rp^2, 
/cr
est. var(Y_i^2 Y_j)&=<1/over n-1>/sum_<t=1>^n/lp Y_<t,i>^2 
Y_<t,j>-/wh/delta_<i,j>/rp^2,/cr
est. var(Y_i Y_j Y_l)&=<1/over n-1>/sum_<t=1>^n/lp Y_<t,i>Y_<t,j>Y_<t,l>-
/wh/eps_<i,j,l>/rp^2./cr>$$
An observed positive $T$ indicates that including the third order 
correction factor is worth the trouble.

It is also possible to include a fourth order correction to normality, 
but the high number of coefficients that must be estimated based on the 
single set of observations disencourage this, if the dimension is higher 
than, say, three. But in dimensions one, two, and three higher order 
correction terms could lead to improved performance.

This also invites hybrids to be constructed. One can imagine a 
high-dimensional feature vector $X$ being transformed to some $Y$, where 
specific subsets are found to be mutually independent, and with subset 
size less than or equal to three. Each subset of variates now gets its 
low-dimensional density estimated using (5.88) or some higher order 
version. A fifth order corrected normal density estimate looks for 
example as follows, in the univariate case:
$$/eqalignno<
/wh 
f_5(x)&=N(/hm,/wh/ssg^2)(x)/Bigg/lbrace1+<1/over6>/wh/gamma_3(y^3-3y)/cr
&/quad+<1/over24>(/wh/gamma_4-3)(y^4-6y^2+3)/cr
&/quad+<1/over120>(/wh/gamma_5-10/wh/gamma_3)(y^5-10y^3+15y)/Bigg/rbrace,
&(5.90)/cr>$$
where $y=(x-/hm)/sla/wh/ssg$ and $/wh/gamma_j=<1/over 
n>/sum_<t=1>^n(Y_t)^j$.
/subsection<5.4.G><Taking initial estimation variability into account:/nl
A more sophisticated inclusion rule> $X_1/upto X_n$ come from an unknown
density $f$. 
The general orthogonal expansion procedure for estimating $f$ involves the 
following steps (see 5.4.A): First perform an initial variable
transformation $X/to Y=R(X)$; then decide on a benchmark density
$g_0(y)$ for $Y$; find a convenient weight function $w(y)$
w.r.t.~which functions $/lbrace/psi_j/rbrace$ become orthonormal and
complete; and finally use
$$/wh g_m(y)=g_0(y)/sum_<j=0>^m/hc(j)/,/psi_j(y)/eqno(5.91)$$
as estimator for the density $g$ of $Y$. Accompanying (5.91) is the
estimator for the original $f$:
$$/wh f_m(x)=/wh g_m(y(x))/,/invert</partial y(x)/sla/partial x>./eqno(5.92)$$
Let us also recapitulate how the order $m$ is determined: The
integrated, squared, weighted error is
$$/eqalign<
/ise(m)&=/int/lp/wh g_m-g/rp^2<w/over g_0^2>/,dy/cr
&=/sum_<j=0>^/infty c(j)^2-/Delta(m),/cr>$$
where
$$/Delta(m)=/sum_<j=0>^m/lsb c(j)^2-/lb/hc(j)-c(j)/rb^2/rsb
/eqno(5.93)$$ 
ought to be as large as possible. Its expectation is
$$A(m)=E/Delta(m)=/sum_<j=0>^m/lb c(j)^2-/ndel/tau(j)^2/rb,
/eqno(5.94)$$ 
for which an unbiased estimator is
$$/hA(m)=/sum_<j=0>^m/lb/hc(j)^2-<2/over n>/wh/tau(j)^2/rb.
/eqno(5.95)$$ 
$m$ is often chosen as the one that maximises $/hA(m)$ (but see
5.4.E). 

The orthonormal expansion estimator is most effective when only few
terms are needed. This is achieved if the benchmark density $g_0$ is
already close to $g$, cf.~(5.91). Solution 2 above was motivated by
this, and effectively used $/wh g_0(y)=g_0(y;/,/wh/alpha_1,
/wh/beta_1/upto/wh/alpha_d,/wh/beta_d)$, say, fitting a product beta
density to the (transformed) data. The trouble with this, however, is
that the </it estimation variability//> present in $/wh g_0(y)$
destroys parts of the mathematical derivations summarised above; in
particular $A(m)$ of (5.94) is not equal to $E/Delta(m)$, and $/hA(m)$
of (5.95) is not unbiased for $E/Delta(m)$ (nor for $A(m)$). Also
disturbing is noise in the initial data transformation $y=R(x)$, say
$y=/hs^<-1/sla2>(x-/hm)$ as in 5.4.F. (That even $w(y)$ sometimes is
chosen as a data-based $/wh w(y)$ is not as important.)

This subsection looks into the consequences of using a data-based $/wh
g_0(y)$ as benchmark density for $Y$. In particular another criterion
than (5.95) will be found to be (approximately) unbiased for
$E/Delta(m)$.

The analysis that follows does not solve each of the
noise-in-transformation/sla noise-in-benchmark-density/sla
noise-in-weight-function problems, but the mathematical methods we use
are applicable in some generality.

Assume
$$/wh g_0(y)=g_0(y,/wh/th)/eqno(5.96)$$
is the initial estimator of $g(y)$, i.e.~a parametric family of
densities is fitted to $Y_1/upto Y_n$. Assume further, for
concreteness, that $/wh/th$ is the maximum likelihood estimator,
solving $/sum_<t=1>^n/partial/log g_0(Y_t,/wh/th)/sla/partial/th_l
=0$, $l=1/upto p$, denoting by $p$ the number of components in $/th$.
According to Section 3.1 $/wh/th$ is consistent for a certain
parameter value $/th_0$ that is least false (most fitting) as measured
by the Kullback-Leibler distance, and
$$/sqrt<n>/lp/wh/th-/th_0/rp- T^<-1>/sqrt<n>/,</overline L>/totop<P>0,$$
where
$$</overline L>=/ndel/sum_<t=1>^n L(Y_t)=/ndel/sum_<t=1>^n</partial
/log g_0/lp Y_t,/th_0/rp/over/partial/th>/eqno(5.97)$$
and
$$T=-/int</prt^2/log g_0(y,/th_0)/over/prt/th./,/prt/th.>/ g(y)/,dy.
/eqno(5.98)$$ 
$T$ is a $p/times p$ matrix, and $L(y)=/prt/log g_0(y,/th_0)/sla
/prt/th$ is a $p$-vector.
$$/wh T=-/ndel/sum_<t=1>^n</prt^2/log g_0/lp Y_t,/wh/th/rp/over
/prt/th./,/prt/th.>/eqno(5.99)$$ 
is consistent for $T$. It was also shown in 3.1 that
$$/sqrt<n>/lp/wh/th-/th_0/rp/totop<D>/ T^<-1>N_p(0,K)=N_p/lp0,
T^<-1>KT^<-1>/rp,/eqno(5.100)$$ 
writing
$$K=/int L(y)/,L(y)'/,g(y)/,dy.$$

That $/wh g_0(y)$ is used instead of (the unavailable)
$g_0(y)=g_0(y,/th_0)$ necessitates some notational adjustments. One
has $g=/wh g_0h$, where $h$ has coefficients
$$/eqalign<
c(j)_/est &=/int/psi_j wh/,dy=/int/psi_j<w/over/wh g_0>/,g/,dy/cr
&=E/psi_j(Y_0)<w(Y_0)/over/wh g_0(Y_0)>,/cr>/eqno(5.101)$$
the ``estimated version'' of the ideal
$$c(j)=/int/psi_j<w/over g_0>/,g/,dy./eqno(5.102)$$
$Y_0$ denotes a fresh, new vector, independent of $Y_1/upto Y_n$.
$c(j)_/est $ could be estimated in a cross validation manner,
using 
$$c(j)_/est ^/ast =/ndel/sum_<t=1>^n/psi_j/lp Y_t/rp/,w/lp
Y_t/rp/sla/wh g_<0,(t)>/lp Y_t/rp,$$
where $/wh g_<0,(t)>=g_0(y,/wh/th_<(t)>)$ is as before, but skipping
$Y_t$ in the evaluation of $/wh/th$. We use the simpler
$$/wh c(j)_/est =/ndel/sum_<t=1>^n/psi_j(Y_t)/,<w(Y_t)/over /wh
g_0(Y_t)>,/eqno(5.103)$$ however, which is the estimated version of
$$/hc(j)=/ndel/sum_<t=1>^n/psi_j(Y_t)/,<w(Y_t)/over g_0(Y_t)>.
/eqno(5.104)$$ 
The analogue of (5.93) is easily derived:
$$/eqalignno<
/ise(m)&=/int/lp/wh g_m-g/rp^2/ <w/over/wh g_0^2>/,dy/cr
&=/int/lp/wh h_m-h/rp^2/,w/,dy/cr
&=/sum_<j=0>^/infty c(j)_/est ^2-/Delta(m),/cr
/Delta(m)&=/sum_<j=0>^m/lsb c(j)_/est ^2-/lb /hc(j)_/est -
c(j)_/est /rb^2/rsb.&(5.105)/cr>$$

The task ahead of us consists of ($i$) obtaining a usable expression
for $E/Delta(m)$, and ($ii$) constructing an (almost unbiased)
estimator for $E/Delta(m)$.

The starting point for solving these problems is to Taylor expand
$1/sla/wh g_0(y)$:
$$/eqalignno<
<1/over/wh g_0(y)>&/doteq<1/over g_0(y)>+/sum_<l=1>^p/lb</prt/over
/prt/th_l>/ <1/over g_0(y,/th_0)>/rb/lp/wh/th_l-/th_<0,l>/rp/cr
&=<1/over g_0(y)>/lb1-L(y)'/lp/wh/th-/th_0/rp/rb.&(5.106)/cr>$$
It follows that 
$$/eqalignno<
c(j)_/est&/doteq/int/psi_j/,<w/over g_0>/lb1-L'/lp/wh/th-/th_0
/rp/rb/,g/,dy/cr 
&=c(j)-b(j)'/lp/wh/th-/th_0/rp&(5.107)/cr>$$
and that
$$/eqalignno<
/hc(j)_/est &/doteq/ndel/sum_<t=1>^n/psi_j(Y_t)/,<w(Y_t)/over
g_0(Y_t)>/lb1-L(Y_t)'/lp/wh/th-/th_0/rp/rb/cr
&=/hc(j)-B(j)'/lp/wh/th-/th_0/rp.&(5.108)/cr>$$
Here
$$/eqalignno<
B(j)&=/ndel/sum_<t=1>^n/psi_j(Y_t)/,<w(Y_t)/over g_0(Y_t)>/,L(Y_t)/cr
&/totop<P>/ b(j)=/int/psi_j/,<w/over g_o>/,Lg/,dy.&(5.109)/cr>$$
These expressions, together with
$$/wh/th-/th_0/doteq T^<-1>/,</overline L>/eqno(5.110)$$
(see (5.97)), make it possible to find an expression for
$E/Delta(m)$.

First, (5.107) gives
$$/eqalign<
E c(j)_/est^2 &/doteq c(j)^2+E/lb b(j)'/lp/wh/th-/th_0/rp/rb^2/cr
&/doteq c(j)^2+/ndel/,b(j)'/,T^<-1>KT^<-1> b(j),/cr>$$
see (5.100). Second, (5.107) and (5.108) in concert imply
$$/hc(j)_/est -c(j)_/est /doteq/hc(j)-c(j)-/lb
B(j)-b(j)/rb'/lp/wh/th-/th_0/rp,$$ 
and all three terms $/hc(j)-c(j)$, $B(j)-b(j)$, $/wh/th-/th_0/doteq
T^<-1> </overline L>$ are of the form $/ndel/sum_<t=1>^nA_t$,
$EA_t=0$. By appealing to the lemma below one can show that
$$/eqalign<
E/lb/hc(j)_/est -c(j)_/est /rb^2&= E/lb/hc(j)-c(j)/rb^2 +
O/lp<1/over n^2>/rp/cr
&=/ndel/ </rm Var>/lb/psi_j(Y_0)/,<w(Y_0)/over
g_0(Y_0)>/rb+O/lp<1/over n^2>/rp=/ndel/tau(j)^2+O/lp<1/over
n^2>/rp./cr>$$
Combining these facts, task ($i$) set up for us above is done:
$$E/Delta(m)=/sum_<j=0>^m/lb c(j)^2-/ndel/tau(j)^2+/ndel b(j)' T^<-1>
K T^<-1> b(j)/rb+O/lp<1/over n^2>/rp./eqno(5.111)$$
Notice the generalisation from (5.94).

The following lemma conveniently summarises facts about expectations
of products of averages. It was used above and will be appealed to
again. Its proof is left as an exercise.
/goodbreak
/Blemma Let $(A_t,B_t,C_t,D_t)$ be independent quadruples, $t=1/upto
n$, and let $EA_t=EB_t=EC_t=ED_t=0$. Then
%
%
/begingroup
/def/At<A_t> /def/Bt<B_t> /def/Ct<C_t> /def/Dt<D_t>
/def/E(#1#2)<E/lp #1_t#2_t/rp/,>
/def/ol#1<</overline #1>>
$$/eqalign<
E/ol A_n/ol B_n&=/ndel E/At/Bt;/cr
E/ol A_n/ol B_n/ol C_n&=<1/over n^2>E/At/Bt/Ct;/cr
E/ol A_n/ol B_n/ol C_n/ol D_n&=<1/over n^3>E/At/Bt/Ct/Dt/cr
+<n-1/over n^3>&/Bigg/lbrace /E(AB)/E(CD)+/cr
&/qquad/E(AC)/E(BD)+/E(AD)/E(BC)/Bigg/rbrace/cr
&/doteq<1/over n^2>/lbrace/E(AB)/E(CD)/cr
&/qquad/qquad+/E(AC)/E(BD)+/E(AD)/E(BC)/rbrace;/cr
E/ndel/sum_<t=1>^n/At/Bt/ndel/sum_<t=1>^n/Ct&=/ndel E/At/Bt/Ct;/cr
E/ndel/sum_<t=1>^n/At/Bt/ /ndel/sum_<t=1>^n/Ct/ /ndel/sum_<t=1>^n/Dt&=
<1/over n^2>E(/At/Bt/Ct/Dt)+<n-1/over n^2>/E(AB)/E(CD)/cr
&/doteq/ndel/E(AB)/E(CD)./cr>$$ /endgroup /Elemma
/note<Remark 4.> It is interesting to see the amount of extra burden the
estimation variability in $/wh g_0(y)=g_0(y,/wh/th)$ causes the
expected $/ise$. From
$$/ise(m)=/sum_<j=0>^m/lb/hc(j)_/est -c(j)_/est /rb^2+
/sum_<j/gt m> c(j)_/est ^2$$
and the analysis above, we find
$$/mise(m)=/sum_<j=0>^m/ndel/tau(j)^2+/sum_<j/gt m>/lb c(j)^2+/ndel
b(j)' T^<-1>KT^<-1>b(j)/rb+O/lp<1/over n^2>/rp,$$
which can be compared to $/mise(m)$ for the ideal case:
$/sum_<j=0>^m/ndel/tau(j)^2+/sum_<j/gt m>c(j)^2$.
/Enote
Task ($ii$) remains: $E/Delta(m)$ of (5.111) must be estimated from
the data. Start out considering 
$$/hc(j)_/est /doteq c(j)+/hc(j)-c(j)-b(j)'
/lp/wh/th-/th_0/rp-(B(j)-b(j))'/lp/wh/th-/th_0/rp.$$
If one uses $/wh/th-/th_0/doteq T^<-1></overline L>$ again, and
results
from the lemma, then one can show that 
$$/eqalignno<
E/hc(j)_/est ^2&/doteq c(j)^2+/ndel/tau(j)^2+/ndel b(j)'
T^<-1>KT^<-1>b(j)/cr 
&/quad-2/ndel b(j)' T^<-1>b(j)-2/ndel/Tr(T^<-1>U)+O/lp<1/over n^2>/rp,
&(5.112)/cr>$$ 
where $U$ is the matrix with elements
$$U_<l,m>=/cov/lp/psi_j<w/over g_0>L_l,L_m/rp=/int/psi_j<w/over
g_0>L_l L_m/,g/,dy./eqno(5.113)$$
A consistent estimator is
$$/wh U=/ndel/sum_<t=1>^n/psi_j(Y_t)<w(Y_t)/over/wh g_0(Y_t)>
</prt/log g_0(Y_t,/wh/th)/over/prt/th.>/lp< /prt/log g_0(Y_t,/wh/th)
/over /prt/th.>/rp.$$

Let
$$/wh/tau(j)_/est ^2=<1/over
n-1>/sum_<t=1>^n/lb/psi_j(Y_t)<w(Y_t)/over/wh g_o(Y_t)>-/hc(j)_</rm
est>/rb^2$$ 
and
$$/wh B(j)_/est =/ndel/sum_<t=1>^n/psi_j(Y_t)<w(Y_t)/over /wh
g_0(Y_t)> </prt/log g_0(Y_t,/wh/th)/over/prt/th>.$$
Judicious calculations then show that
$$/hc(j)_/est ^2-<2/over n>/wh/tau(j)_/est ^2+<2/over n>/wh
B(j)_/est '/wh T^<-1>/wh B(j)_/est+<2/over n>/hc(j)_/est/Tr/lp/wh T^<-1>
/wh U/rp$$
has expectation
$$c(j)^2-/ndel/tau(j)^2+/ndel b(j)'T^<-1>KT^<-1>b(j)$$
up to this $/ndel$-order, i.e.,
$$/eqalignno<
A^/ast (m)&=/sum_<j=0>^m/Big/lbrace/hc(j)_/est^2-<2/over n>/wh/tau(j)_/est^2/cr
&/quad+<2/over n>/wh B(j)_/est'/wh T^<-1>/wh B(j)_/est+<2/over
n>/hc(j)_/est /Tr/lp/wh T^<-1>/wh U/rp/Big/rbrace&(5.114)/cr>$$
is unbiased, up to $O(<1/over n^2>)$, for $E/Delta(m)$.

In conclusion, one should use
$$/wh g_m(y)=/wh g_0(y)/sum_<j=0>^m/hc(j)_/est/psi_j(y)$$
to estimate $g(y)$, and this is the same estimator as (5.63),
advocated long ago! But the order $m$ should be arrived at in a more
sophisticated manner: choose $m$ to maximise $A^/ast (m)$ of (5.114)
instead of simply 
$$/hA(m)=/sum_<j=0>^m/lb/hc(j)_/est^2-<2/over
n>/wh/tau(j)_/est^2/rb,$$
which is the criterion derived when the estimation variability in
$g_0$ was disregarded, cf.~(5.67).

/subsection<5.4.H><Solution 4: From orthogonal expansions to projection 
pursuit> The previously discussed solutions have started from some (more 
or less explicit) ``initial guess'' $/wh f_0$ for the unknown density at 
hand, and then tried to provide a good correction factor to it, using a 
single or a product of expansions. It is natural, then, to recycle the 
good idea, and correct the obtained density estimate in the same way, 
i.e.~provide a second correction factor to the first one.

A general iterative scheme is the following. The starting point is an 
initial estimate $/wh g_<(0)>(y)$ for the unknown density $g(y)$, and 
orthonormal basis functions $/psi_0,/psi_1,/ldots$ w.r.t.~the weight 
function $w$. Writing
$$g(y)=/wh g_<(0)>(y) h_<(1)>(y)/doteq/wh g_<(0)>(y)/sum_<j=0>^<m(1)>
c_<(1)>(j)/psi_j(y)$$
and applying the
$$/int/lb/wh h_<(1)>-h_<(1)>/rb^2w dy=/int/lb/wh g_<(1)>-g/rb^2<w/over/wh 
g_<(0)>^2>dy$$
criterion, one arrives at
$$/wh g_<(1)>(y)=/wh g_<(0)>(y)/sum_<j=0>^<m(1)>/wh 
c_<(1)>(j)/psi_j(y)/eqno(5.115)$$
as before. Here
$$/wh c_<(1)>(j)=<1/over n>/sum_<t=1>^n/psi_j(Y_t)<w(Y_t)/over /wh 
g_<(0)>(Y_t)>,$$
and $m(1)$ is decided upon studying 
$$/hA_<(1)>(m)=/sum_<j=0>^m/lb/wh c_<(1)>(j)^2-<2/over 
n>/wh/tau_<(1)>(j)^2/rb$$
or a version penalising higher $j$'s, writing
$$/wh/tau_<(1)>(j)^2=<1/over n-1>/sum_<t=1>^n/lb/psi_j(Y_t)<w(Y_t)/over/wh 
g_0(Y_t)>-/wh c_<(1)>(j)/rb^2.$$
(A more sophisticated rule, as in (5.114), could also be used.)

The next step writes
$$g_<(2)>(y)=/wh g_<(1)>(y)h_<(2)>(y)/doteq/wh 
g_<(1)>(y)/sum_<j=0>^<m(2)>c_<(2)>(j)/psi_j(y),$$
and the second set of coefficients can be determined using the 
$/int/lbrace h_<(2)>-/sum_<j=0>^<m(2)>c_<(2)>(j)/psi_j/rbrace^2/,/allowbreak  w/, dy$
criterion. The  
result is
$$c_<(2)>(j)=/int/psi_j h_<(2)> wdy=/int/psi_j<w/over/wh g_<(1)>> gdy,$$
which is the expected value of $/psi_j(Y_0)w(Y_0)/sla/wh g_<(1)>(Y_0)$ for a 
new $Y_0$, independent of $Y_1/upto Y_n$. We could use a cross-validated 
estimate
$$<1/over n>/sum_<t=1>^n/psi_j(Y_t)<w(Y_t)/over/wh 
g_<(1),t>(Y_t)>/eqno(5.116)$$
for this expectation, where $/wh g_<(1),t>(y)$ is the first stage density 
estimate based on all observations except $Y_t$. This becomes quite 
intricate, however, and we settle on
$$/wh c_<(2)>(j)=<1/over n>/sum_<t=1>^n/psi_j(Y_t)<w(Y_t)/over/wh g_<(1)> 
(Y_t)>/eqno(5.117)$$
which continues to be a good estimate in spite of the subtle difficulty 
alluded to above, since $/wh g_<(1)>(y)/doteq/wh g_<(1),t>(y)$. Thus 
$$/wh g_<(2)>(y)=/wh g_<(0)>(y)/sum_<j=0>^<m(1)>/wh c_<(1)>(j)/psi_j(y) 
/sum_<j=0>^<m(2)>/wh c_<(2)>(j)/psi_j(y)/eqno(5.118)$$
is the second stage density estimator, with $m(2)$ found from 
maximisation of
$$/hA_<(2)>(m)=/sum_<j=0>^m/lb/wh c_<(2)>(j)^2-<2/over 
n>/wh/tau_<(2)>(j)^2/rb,$$
or, possibly, from a related ``cross-validated'' criterion.

At step $M$ of this iterative process a product type density
$$/eqalignno<
/wh g_<(M)>(y)&=/wh g_0(y)/wh h_<(1)>(y)/cdots/wh h_<(M)>(y)/cr
&=/wh g_0(y)/prod_<s=1>^M/sum_<j=0>^<m(s)>/wh 
c_<(s)>(j)/psi_j(y)&(5.119)/cr>$$
is produced, with $/wh c_<(s)>(j)$ and $/wh/tau_<(s)>(j)$ determined from
$$/psi_j(Y_1)<w(Y_1)/over/wh g_<(M-1)>(Y_1)>/upto/psi_j(Y_n)<w(Y_n)/over 
/wh g_<(M-1)>(Y_n)>,$$
and $m(s)$ determined from $/hA_<(s)>(m)$.

We should be rather careful about allowing ``false terms'' here, since
including a $/wh c_<(2)>(j)$, say, if the underlying $c_<(2)>(j)$ is 
close to zero, could hurt the following iteration steps. Thus we 
recommend a ``stingy'' version of the $/hA_<(s)>(m)$ criterion, for 
example using the penalising
$$/hA_<(s),/eps>(m)=/sum_<j=0>^m/lb/wh c_<(s)>(j)^2-<2+/eps j/over 
n>/wh/tau_<(s)>(j)^2/rb,/eqno(5.120)$$
cf.~the discussion ending Solution 2, and also putting all 
``non-significant'' $/wh c_<(s)>(j)$ equal to zero, $j/lt m(s)$, in order 
for the final version (5.119) to be parsimonious. 

It remains to decide upon a single one of the candidates (5.119). Consider 
the distance criterion
$$/int(/wh g-g)^2/rho dy=/int/wh g^2/rho dy-2/int/wh g g/rho dy+/int 
g^2/rho dy,$$
where $/rho$ is some weight function, e.g.~the identity. The last term 
does not involve the data at all, so if we can predict with some 
precision which of the candidates $/wh g_<(0)>, /wh g_<(1)>,/ldots$ has 
largest
$$Q=2/int/wh g g/rho dy-/int/wh g^2/rho dy,$$
then we have a reasonable selection procedure. But
$$/wh Q=<2/over n>/sum_<t=1>^n/wh g(Y_t)/rho(Y_t)-/int/wh g(y)^2 /rho(y) 
dy/eqno(5.121)$$ 
can be computed, and estimates $Q$. A more careful version observes that 
$/int/wh g g/rho dy$ really is $E_/ast /wh g(Y_0)/rho(Y_0)$ for a fresh $Y_0$, 
independent of $Y_1/upto Y_n$. Hence we recommend computing the 
cross-validated 
$$/wh Q_<M,/ast >=<2/over n>/sum_<t=1>^n/wh g_<(M),t>(Y_t)/rho(Y_t)-/int /wh 
g(y)^2/rho(y) dy/eqno(5.122)$$
for our candidates $M=0,1/upto M_<UP>$ (say), where $/wh g_<(M),t>$ is as 
in (5.119), but computed with observation $Y_t$ deleted in $/wh g_0$ and 
$/wh c_<(s)>(j)$. The final estimator is $/wh g_<(M)>$ for the $M$ that 
has the largest $/wh Q_<M,/ast >$.

Another possibility uses the Kullback-Leibler information distance
$$I(g,/wh g)=/int g/log<g/over/wh g>dy=/int g/log g dy-/int g/log/wh g 
dy,$$
cf.~Section 3.1. It suffices to estimate $/int g/log/wh g dy$, which is
$E_/ast /log/wh g(Y_0)$ for a fresh $Y_0$. Hence
$$I_<M,/ast >=<1/over n>/sum_<t=1>^n/log/wh g_<(M),t>(Y_t)$$
could be evaluated for each $M$, and the largest value determines the 
choice $/wh g_<(M)>$.

Our estimator (5.119) is built up in a manner similar to the </it
projection pursuit//> density estimator, see Friedman, Stuetzle, and
Schroeder  (1984) and Huber (1985). There one searches, at each step,
for the direction where a correction would be most advantageous, and
constructs such a correction in a ``univariate manner'' after having
found the best direction. $/wh g_<(M)>$ is similar, and in effect
finds both a ``direction'' and an explicit correction at each step.
/section<5.5><$k$-nearest-neighbour classification><> A very
simple-minded classification scheme is the following: for a given
feature vector $x$, find its nearest neighbour among those of the full
training set, and assign the class of this nearest neighbour to $x$.
This is (of course) called the </it nearest-neighbour//> ($NN$) rule.
The present section briefly discusses the $NN$ rule and its
generalisation, the </it k-NN//> method. The </it k-NN//> rule finds first the
$k$ nearest neighbours to $x$, and then assigns $x$ to a class label among
those present, by a majority vote.

It is not immediately seen that such a rule is related to the main,
nonparametric theme developed in this chapter; estimate $f_1/upto f_K$
and find the maximum among $/pi_1/hf_1/upto/pi_K/hf_K$. Such a
relation exists, however, and will be pointed out below.
/subsection<5.5.A><</it k-NN//> density estimation> Forget the
classification problem for a moment, and consider a sample $X_1/upto
X_n$ from a continuous density $f$ on $/rr^d$. For each $x$, let
$S(x,v)$ be the sphere around $x$ with volume $v$. Let $A_k(x)$ be the
smallest ones of these that encloses exactly $k$ of the $n$ data
points. The </it k-NN//> density estimator for $f(x)$ is
$$f_<k-NN>(x)=<k/sla n/over V_k(x)>,/eqno(5.123)$$
where $V_k(x)$ is the volume of $A_k(x)$.

The notion of a ``sphere'' is tied to a notion of distance. This could
be ordinary Euclidian distance, or, for reasons of scaling,
$$d(y,z)=/lb(y-z)'/sg^<-1>(y-x)/rb^<1/sla2>./eqno(5.124)$$
Here $/sg$ could be some crude or more serious estimate of the covariance
matrix for $f$. With this choice,
$$S(x,v)=/lb y:/ (y-x)'/sg^<-1>(y-x)/le/rho(v)^2/rb,$$
with the radius of the ellipsoid chosen to give 
$$</rm volume/ >/lb S(x,v)/rb=/invert</sg>^<1/sla2></pi^<d/sla2>/over
/Gamma/lp<d/over2>+1/rp>/ /rho(v)^d=v,/eqno(5.125)$$
i.e.~$/rho(v)$ proportional to $v^<1/sla d>$. (The simplest case is
the one-dimensional one, in which $S(x,v)=/lbrack x-/halv v,x+/halv
v/rbrack$.)

It is interesting to contrast $f_<k-NN>$ above to histograms and
kernel density estimators. In (5.123) $k$ is fixed but $V_k(x)$
allowed to vary with $x$. If instead $V_k(x)=V$ is the fixed volume of
a sphere $S(x,v)$ around $x$, and $k/sla n$ is the random fraction of
data points falling within it, then
$$<k/sla n/over V>=<1/over nV>/sum_<i=1>^n I/lb x-X_i/in
S(0,V)/rb/eqno(5.126)$$ 
takes the form of a kernel density estimator, with a Parzen-type
either/sla or kernel. Since $k/sla n$ estimates
$/int_<S(x,v)>f(y)/,dy$, which is close to $f(x)V$ when $V$ is small,
(5.126) is a reasonable estimate. This also provides, in a somewhat
indirect manner, motivation for $f_<k-NN>$ of (5.123).

If space is divided into fixed cells with volumes $V_1$, $V_2/upto$
and $k_i$ is the number of $X$'s falling into cell $i$, then $<k_i/sla
n/over V_i>$ is the histogram estimator, and is, again, of the same
form as (5.123).

A disadvantage of (5.126), and the more general kernel estimator, is
that the same degree of smoothing is applied all over the space.
Presumably (5.126) would perform better if $V$  is smaller in regions
in which $X$'s are dense than in sparse regions. One remedy is the
</it variable bandwidth//> kernel estimator, see Breiman, Meisel, and
Purcell (1977); another is to fix $k$ and let $V$ be determined by the
density of data around $x$, which is the $f_<k-NN>$ approach.

Let us briefly study some of the properties of $f_<k-NN>$ before we
connect it to the classification framework. The essential requirement
for $k$, as a function of sample size $n$, is that 
$$k/to/infty,/ k/sla n/to0./eqno(5.127)$$
This is a necessary and sufficient condition for pointwise convergence
in probability. One can show that 
$$/sup_<x>/invert<f_<k-NN>(x)-f(x)>/to 0$$
with probability one, provided $k/sla/log n/to/infty$, if $f$ is
uniformly continuous.

That $V_k(x)$ in (5.123) is greater than $v$ is equivalent to the
statement that $S(x,v)$ captures at most $k-1$ data points. Hence
$$/eqalign<
/Pr/lb V_k(x)/gt v/rb&=/sum_<i=0>^<k-1><n/choose i>/lp/int_<S(x,v)>
f/,dy/rp^i/lp1-/int_<S(x,v)> f/,dy/rp^<n-i>/cr
&/doteq/sum_<i=0>^<k-1> e^<-/lam_n(x,v)></lam_n(x,v)^i/over
i!>,/cr>$$
using the Poisson approximation, which is good for $n$ large and 
$$/lam_n(x,v)=n/int_<S(x,v)>f/,dy=n/,f(y_x)v$$
not too large. Here $y_x$ is a point somewhere in $S(x,v)$, tending to
$x$ when $V$ tends to zero. But
$$/sum_<i=0>^<k-1>e^<-/lam></lam^i/over i!>=1-/Pr/lb/chi_<2k>^2
/le 2/lam/rb.$$
It follows that
$$/Pr/lb V_k(x)/le v/rb/doteq/Pr/lb/chi_<2k>^2/le 2n/,f(y_x)v/rb,$$
at least for moderate values of $nv$, i.e.
$$V_k(x)/doteq</chi_<2k>^2/over2n/,f(y_x)>./eqno(5.128)$$
To the extent that this approximation can be trusted, therefore, 
$$E/,f_<k-NN>(x)/doteq<k/over n>/,E/,<2n/,f(y_x)/over/chi_<2k>^2>
/doteq <k/over k-1>/,f(x).$$
An approximation to its variance can also be written down. Requiring
asymptotic unbiasedness and limiting variance zero amounts exactly to
(5.127).

This suggests that $f_<k-NN>$ can be expected to overshoot its aim,
and indeed, this corresponds to a somewhat too small $A_k(x)$ sphere
in (5.123). $A_k(x)$ is the very smallest one containing $k$ data
points; a ``continuity correction'' is
$$/wt f_<k-NN>(x)=<k/sla n/over/halv/lb V_k(x)+V_<k+1>(x)/rb>.$$
/subsection<5.5.B><</it k-NN//> classification> Suppose now that class
densities $f_1/upto f_K$ are estimated by the </it k-NN//> method described
in 5.5.A, producing $f_<i, k-NN>(x)$, $i=1/upto K$. (Since the ``k'' in
</it k-NN//> is the widely accepted prescript we use $i=1/upto K$ as class
index in this particular section.) Then the usual construction leads
to the following method: maximise
$$/pi_if_<i, k-NN>(x)=/pi_i<k/sla n_i/over V_<i,k>(x)>/eqno(5.129)$$
w.r.t.~$i$, and assign this $i$ as class label to $x$. If $/pi_i$
equals $n_i/sla N$, $N=/sum_<i=1>^K n_i$, or is estimated in this way,
then the scheme becomes that of minimising $V_<i,k>(x)$ w.r.t.~$i$:
Which class will first have $k$ members within an increasing sphere
$S(x,v)$? 

Another version of the </it k-NN//> theme has become more popular, however.
Pool the separate training sets $/lbrace X_1^<(i)>/upto
X_<n_i>^<(i)>/rbrace$ into one, say $Z$, and let $W_k(x)$ be the
volume of the smallest sphere containing exactly $k$ points. Then use
$$/wh f_<i, k-NN>(x)=<k_i/sla n_i/over W_k(x)>,/eqno(5.130)$$
where $k_i/sla n_i$ is the fraction of class $i$ points that fall
within the smallest sphere. (One usually employs the $/sg$-weighted
distance function (5.124), where $/sg$ is an estimate of the pooled
covariance matrix.)

Notice that $k=k_1+/cdots+k_N$. $k$ would usually be chosen larger
than the ``individual $k$'' appearing in 5.5.A.

If $/pi_i$ is estimated as $n_i/sla N$ again, then
$$/wh/pi_i/hf_<i, k-NN>(x)=<k_i/sla N/over W_k(x)>$$
and
$$/wh P_<k-NN>(i/mid x)=<k_i/over k>,/ i=1/upto K/eqno(5.131)$$
are the remarkably simple estimates for the posterior probabilities.

If $/pi_i$ is known, or is estimated differently from $n_i/sla N$
(because of a sampling distortion, for example), then one should
maximise over $/pi_i/,k_i/sla n_i$, $i=1/upto K$. The doubt option is
used if every
$$/wh P_<k-NN>(i/mid X)=</pi_i/over n_i/sla N>/,k_i/sla/sum_<j=1>^K
</pi_j/over n_j/sla N>/,k_j$$
is less than a threshold $(1-c)$ for all classes.

The problem of choosing a value of $k$ is akin to that of choosing a
bandwidth smoothing parameter, and is a difficult one. One might
consider plots of the (5.130) functions and compare to histograms, or
one could carry out a cross-validation plan. The best way seems to be
trial and error, however, estimating error rates for a small number of
$k$ values. 	
/subsection<5.5.C><Remarks on computation> It is easy to program one
or another of the </it k-NN//> methods, by merely computing every distance
needed. This quickly becomes a heavy burden, however, because of the
number of computations and because every training sample point must be
stored and used for each new vector. Fortunately several tricks are
available for reducing the computational task.

Hart (1968) proposed a </it condensed-nearest-neighbour//> method,
Gates (1972) devised a related </it reduced-nearest-neighbour//> rule,
and Hand and Batchelor (1978) invented a combination called the </it
edited-nearest-neighbour//> rule. Fukunaga and Narendra (1975) start out
with a $k$-means clustering method and go on using a combinatorial
optimisation technique called a branch and bound algorithm to quickly
find the $k$ nearest neighbours. Their technique has recently been
improved upon by Kamgar-Parsi and Kanal (1985). Other impressive
</it k-NN//> algorithms are steadily being published, see for example Ruiz
(1986).

Observe that these algorithms can be extremely useful for computing
kernel estimators for cases where the kernel has bounded support.
/section<5.6><Concluding remarks><>
Of course there are more nonparametric classification methods than
those chosen for discussion in Chapter~5. In particular, each
published density estimation method gives rise to a classification
procedure by inserting estimates in the Bayes rule of Chapter~1. Among
the nonparametric estimation procedures not considered here are the
methods of splines, stochastic approximation, delta sequences, and
penalty function smoothing. A convenient source of theoretical 
information for these methods is Prakasa Rao (1983).

We mentioned projection pursuit density estimation methods in
Section 5.4.H, and saw a connection to these from orthogonal
expansions. Basic references for such estimation procedures are
Friedman, Stuetzle, and Schroeder~(1984) and Huber~(1985). There are
also projection pursuit methods more directly geared towards
classification, without explicitly computing density estimates. One of
these is presented in Friedman (1986), and uses a SMART approach
(smooth multiple additive regression technique.) Another one is offered
by Henry~(1983), generalising the usual logistic discriminant analysis
scheme by allowing ratios of posterior probabilities to be a sum of
smooth functions of linear combinations of $x$-components.

A general, nonparametric approach to classification is discussed
in Breiman, Friedman, Olshen, and Stone's~(1984) CART-book. They give
methods for automatic construction of classification and regression
trees. A computer intensive search is carried out at each decision point of 
the tree to find the optimal split.

It is easy to make up classification methods based on principal
components type description of the individual classes. S.~Wold~(1976)
and Sj/"ostr/"om and S.~Wold~(1980) take this approach, and show how
high-dimensional problems sometimes can be handled successfully that way.
Examples can be found, however, where most of the discriminatory
information lies in directions with low eigenvalues, in which case the
mentioned methods perform poorly.

There are more simulation studies published that consider
the kernel method than the orthogonal expansion method. They are in
many situations expected to perform similarly. The fact that the 
orthogonal expansion approach fitted better into the general system design
developed at the Norwegian Computing Centre, since only a small number
of class descriptions parameters need to be stored for this method,
prompted the research reported on in Section 5.4. These methods also
provide a convenient machinery for handling discrete and mixed data,
see Chapter~9.

In many applications methods that rely on </it some//> structure
being present in the data, but avoids being (too) parametric, will be
best. Semiparametric methods like those considered in Section~5.4 are
promising, and our admittedly limited experience so far has showed them
to perform well.

Methods developed in Chapter~7 have the important ability of
updating and even refining a class description from an initial one,
based perhaps on a moderate training set, to more sophisticated ones,
if unclassified data vectors are available. In particular a normal
distribution description can be refined by incorporating third order
or higher order corrections to it, again relying on the orthogonal
expansion approach of the present chapter. Also the </it k-NN//> rule and 
the kernel method can be updated using incoming, unconfirmed vectors.


/chapter<6><Detecting outliers>
The task of an automatic pattern recogniser is to assign class labels,
according to a chosen procedure, to each in a sequence
of ``new'' feature vectors. These feature vectors have usually been
extracted from objects belonging to one of a number $K$ of classes.
Occasionally the preprocessing machinery may have included ``alien objects'',
however, or the feature vector may have been incorrectly evaluated, etc.
Thus there is a need for a procedure committed to the detection of these
``incredible vectors'', thereby avoiding incorrect forced classification. Such a
%
%
procedure would necessarily make outcasts out of a small percentage of
perfectly legal objects. However, this could very well improve classification
accuracy, too, in that these most freakish among normal feature vectors would
be difficult to label correctly anyway.
/section<6.1><Finding outliers when class densities are known><>
/subsection<6.1.A><Multinormal case, equal covariance matrices>  A possible
procedure would be to declare as outlier any feature vector all of whose
components deviate with more than say $2.5$ standard deviations from all
corresponding components of the $K$ class means. To arrive at a less heuristic
outlier definition, assume that  the class densities are known to be
%
%
$$/fkx=N_d(/mu_k,/sg)(x),/ k=1/upto K./eqno(6.1)$$
In practice the descriptors $/mu_k,/sg$ will have been estimated
from training sets, but with good precision, i.e.~the sampling variability 
in these estimates will be close to being negligible compared
to the sampling variance in a single feature vector. (Nevertheless, see
6.2 below.)

One way in which to formalise the problem is to test the hypothesis
$$H_0:f/in/lb N_d(/mu_1,/sg)/upto N_d(/mu_k,/sg)/rb,/eqno(6.2)$$
where $f$ denotes the density from which the observed candidate vector $X$ is 
drawn. A natural </it a priori//> assumption in view of (6.1), (6.2)
%
%
is that
$$f/in/lb N_d(/mu,/sg);/,/mu/in/rr^d/rb./eqno(6.3)$$
Such an </it a priori//> assumption is generally needed to contrast the null
hypothesis, and in particular, it is needed in order to pursue the likelihood 
ratio program below. (6.3) certainly defines a rich enough class
in that any outcome of $X$ is perfectly explicable by one of its theories.
It should also be pointed out that several other choices each give
exactly the same resulting procedure.

The </it likelihood ratio//> test considers
$$/lr(x)=</mo<H_0>f(x)/over/mom f(x)>/eqno(6.4)$$
and rejects $H_0$ if $/lr(x)/leq a$, where the constant
%
%
$a$ is chosen to achieve 
$$Pr/lb/lr(X)/leq a/mid H_0</rm/ is/ true>/rb/leq/eps./eqno(6.5)$$
$/eps$ is the chosen size or significance level of the test, and $/om$
denotes the set of </it a priori//> accepted theories.

It is easy to simplify $/lr(x)$ in the present case:
$$/eqalign<
/lr(x)&=</max/lb(2/pi)^<-d/sla 2>/vert/sg/vert^<-1/sla 2>e^<-/halv(x-/mu)'
	/sg^<-1>(x-/mu)>;/,/mu=/mu_1/upto/mu_K/rb/over
	/max/lb(2/pi)^<-d/sla 2>/vert/sg/vert^<-1/sla 2>
	e^<-/halv(x-/mu)'/sg^<-1>(x-/mu)>;/,/mu/in/rr^d/rb>/cr
      &=/max_<k/leq K>/exp/lb-/halv(x-/mu_k)'/sg^<-1>(x-/mu_k)/rb/cr
      &=e^<-/halv W>,/cr>$$
where
$$W=/min_<k/leq K>/,(x-/mu_k)'/sg^<-1>(x-/mu_k)./eqno(6.6)$$
Thus $H_0$ is rejected if $W/geq b=-2/log a$, i.e.~an observed
%
%
feature vector $X$ is taken to be an outlier if and only if
all $K$ Mahalanobis distances $(X-/mu_k)'/sg^<-1>(X-/mu_k)/geq b$.

It remains to relate $b$ to the chosen size of the test.
If $/eps=0.005$ is chosen, for example, (6.5) guarantees
that in the long run only $/halv$/
will be declared outliers. (6.5) is equivalent to
$$P_k/lb/lr(X)/leq a/rb/leq/eps,/ k=1/upto K,/eqno(6.7)$$
where $P_k$ is the class distribution for class $k$.
In the present case $b$ must be chosen such that 
$$P_k/lb W/geq b/rb/leq/eps,/ k=1/upto K.$$
The exact distribution of $W$ under $P_k$ is very
%
%
complicated, being a minimum of $K$ dependent non-central 
$/chi^2_d$ distributed variables. The upper bound
$$P_k/lb W/geq b/rb/leq P_k/lb(X-/mu_k)'/sg^<-1>(X-/mu_k)/geq b/rb=Pr/lb/chi^2_d/geq b/rb$$
is however easily derived, and provides at the same time
also an approximation, in that the $k$'th variable in (6.6) will
tend to be the smallest one under $P_k$, unless some classes are easily
confused. Consequently $b=/gamma_<d;1-/eps>$ does the job, where
$Pr/lb/chi^2_d/geq/gamma_<d;1-/eps>/rb=/eps$.  

The outlier detector may also be  described as follows:
if
$$f_k(X)=N_d(/mu_k,/sg)(X)/leq(2/pi)^<-d/sla 2>/vert/sg/vert^<-1/sla 2>
	/exp/,/lp-/halv/gamma_<d;1-/eps>/rp,/ k=1/upto K,/eqno(6.8)$$
then $X$ is rejected as being incredible.
%
%
/subsection<6.1.B><General case> Let $K$ class densities $/ftfk$ be given, and
let $/om$ denote the set of </it a priori//> acceptable statistical models 
for the feature vector $X$ about to be tested. Then
$$/lr(x)=</mo<H_0>f(x)/over/mom f(x)>=</mo<k/leq K>/fkx/over M(x)>,$$
assuming $M(x)=/mom f(x)$ to be finite, and by (6.7) a constant $a$
needs to be found fulfilling 
$$P_k/lb</mo<m/leq K>>f_m(X)/leq a/,M(X)/rb/leq/eps,/ k=1/upto K./eqno(6.9)$$
The outlier criterion becomes
$$f_k(X)/leq a/,M(X),/ k=1/upto K./eqno(6.10)$$
%
%

The best possible value of $a$ will often be very difficult to find,
and besides, it would not be practical in  that its value would
depend  in an intricate way on the simultaneous aspects of
$f_1(X)/upto f_K(X)$, as opposed to only the separate characteristics
of the $K$ class densities. A conservative  and for practical purposes quite 
satisfactory approximation can be derived, as in the special case
6.1.A:
$$P_k/lb</mo<m/leq K>>f_m(X)/leq a/,M(X)/rb/leq P_k/lb f_k(X)/leq a/,M(X)/rb./eqno(6.11)$$
Hence $a$ in (6.10) will be chosen such that 
$$P_k/lb f_k(X)/leq a/,M(X)/rb/leq/eps,/ k=1/upto K./eqno(6.12)$$
/note<Remark.> One may also consider a related, but different approach to
the outlier problem. Test separately $K$ individual hypoteses, say
$H_k:f=f_k$ versus $f/in /om-/lb f_k/rb$, and declare $X$ an outlier if
it fails each of these tests. If these tests are constructed by a
likelihood ratio method then the resulting outlier criterion resembles (6.12):
$f_k(X)/leq a_k/,M(X),/ k=1/upto K$, where $P_k/lb f_k(X)/leq a_k/,M(X)/rb/leq/eps$
defines $a_k$. This ``$K$ separate tests'' approach may sometimes
be ``sharper'' than the ``simultaneous approach'', which however
concerns itself more directly with the union hypotesis $H_0:H_1/cup/cdots/cup H_K$,
and perhaps in a philosophically more natural way.

--- This second approach is illustrated at the end of the present chapter in a problem
where unknown parameters also are present.
/bigskip
Now let us turn to the important case where
$$/fkx=N_d(/mu_k,/sg_k)(x),/ k=1/upto K./eqno(6.13)$$
A natural choice for the </it a priori//> set $/om$ is
the set of all $N_d(/mu,/sg)$ distributions with $/mu/in/rr^d$ and
$/sg$ ranging over a suitable subset of positive definite symmetric
matrices related to $/sg_1/upto/sg_K$ in some way, e.g.
$$/lb/sg_1/upto/sg_K/rb,/lb/sum_<k=1>^Ka_k/sg_k;/,a_k/geq0,/ k=1/upto K,/sum_<k=1>^Ka_k=1/rb.$$
%
%
The arguments below show that the same test results for all these
choices, as long as $/miom/vert/sg/vert$ is finite.

It is seen that
$$M(x)=/mom f(x)=(2/pi)^<-d/sla 2>/mom/vert/sg/vert^<-1/sla 2>$$
is a constant, and that $/fkx/leq a/,M(x)$ is equivalent to
$$(x-/mu_k)'/sg_k^<-1>(x-/mu_k)+/log/vert/sg_k/vert/geq b$$
for some $b$. Hence
$$P_k/lb f_k(X)/leq a/,M(X)/rb=Pr/lb/chi_d^2+/log/vert/sg_k/vert/geq b/rb,$$
and the best value for $b$ obeying (6.12) is
%
%
$$b=/gamma_<d;1-/eps>+/log/vert/sg_<k_0>/vert;/eqno(6.14)$$
$$/vert/sg_<k_0>/vert=</mo<k/leq K>>/vert/sg_k/vert./eqno(6.15)$$
The outlier criterion may be written
$$f_k(X)=N_d(/mu_k,/sg_k)(X)/leq(2/pi)^<-d/sla 2>/vert/sg_<k_0>/vert^<-1/sla 2>
	e^<-/halv/gamma_<d;1-/eps>>,/ k=1/upto K./eqno(6.16)$$
The probability of falsely making an outcast out of a legal feature vector 
is less than, and approximately equal to, the chosen level $/eps$.

As a second application of the general considerations that led to 
(6.10) and (6.12), suppose
$$/fkx=/prod_<i=1>^d/th_<k,i>e^<-/th_<k,i>/rx_i>,/,
	/rx_1/upto/rx_d/gt0,/ k=1/upto K,$$
giving the model in 3.2.A. The parameter vectors $(/th_<k,1>/upto/th_<k,d>)$
are known, and our quest is to construct a criterion for outliers
based on the likelihood ratio test. We choose the </it a priori//> set $/om$ as
the full family of multi-exponential densities
$f(x)=/prod_<i=1>^d/th_i/exp(-/th_ix_i);/,/th_1/upto/th_d/gt0$. Then
it is easily seen that
$$M(x)=/mom f(x)=/prod_<i=1>^d<1/over x_i>e^<-1>.$$
(6.10) becomes equivalent to
$$/prod_<i=1>^d/th_<k,i>x_ie^<-/th_<k,i>x_i>/leq a/,e^<-d>,/ k=1/upto K./eqno(6.17)$$
Under $P_k$ $/th_<k,i>X_i=Y_i$ is a unit exponentional.
Consequently the constraint on $a$ becomes
$$/eqalign<
P_k/lb f_k(X)/leq a/,M(X)/rb&=Pr/lb/prod_<i=1>^dY_ie^<-Y_i>/leq a/,e^<-d>/rb/cr
	&=H_d(a/,e^<-d>)=/eps./cr>$$
$a$ may be found from this by numerical methods.
%
%
/subsection<6.1.C><A nonparametric outlier test>
It is also possible to construct outlier tests without any
parametric assumptions on the class densities. The well known
Chebyshev inequality bounds the probability of observing
a variable differing with more than  a given multiple of standard
deviations from the mean. If the mean vector $/mu_k$ and covariance 
matrix $/sg_k$ exist  for class $k$, Markov's inequality
can be used to obtain a multivariate generalisation 
of Chebyshev's:
$$/eqalignno<
&P_k/lb(X-/mu_k)'/sg_k^<-1>(X-/mu_k)/geq b/rb/hfill/cr
&/qquad/leq<1/over b>E(X-/mu_k)'/sg_k^<-1>(X-/mu_k)=<d/over b>.&(6.18)/cr>$$
One may show that this provides the best 
%
%
generalisation in the sense that the confidence ellipsoids
that can be derived from it get the smallest possible volumes.

One may now define $X$ to be an outlier if all $K$
Mahalanobis distances are sufficiently large,
$$(X-/mu_k)'/sg_k^<-1>(X-/mu_k)/geq d/sla /eps,/ k=1/upto K./eqno(6.19)$$
(Obviously estimates for $/mu_k$ and $/sg_k$ must be available.)
This and related procedures yield very mild bounds, however, compared to
previously presented methods based on multinormality.
%
%
/section<6.2><Finding outliers when class descriptors are unknown><>
In Section 6.1 the procedures (6.8), (6.16) were derived for the purpose of
detecting outliers not belonging to any of $K$ known multinormal distributions,
with a probability of ``a false detection'' being at most $/eps$.
We argued that the sampling variability of the needed estimates
$/hm_k,/hs_k,/hs$ was negligible. Nevertheless, it is of interest to derive
procedures that explicitly take the training sets and sampling variability
of estimates into account. The difference between the following
outlier rules and those in 6.1 can be
%
%
appreciable only when training
sets are small, however.
/subsection<6.2.A><Common covariance matrix> Assume that $/fkx=N_d(/mu_k,/sg)(x),
/ k=1/upto K$, and that a training set $Z$ of the form (3.3) is available from 
which to learn estimates of $/mu_k,/sg$. In addition a new vector $X$ is 
observed, originating from a $N_d(/mu,/sg)$ distribution
with an unspecified $/mu/in/rr^d$ </it a priori//>. The hypothesis
to be tested is 
$$H_0:/mu/in/lb/mu_1/upto/mu_K/rb./eqno(6.20)$$

The total likelihood for the observed data is
$$/eqalign<
L(z,x)&=L(z,x/mid/mu_1/upto /mu_k,/sg;/mu)/cr
&=/prod_<k=1>^K/prod_<j=1>^<n_k>N_d(/mu_k,/sg)/lp x_j^/k/rp/cdot N_d(/mu,/sg)(x)/cr
%
%
&=/prod_<k=1>^K(2/pi)^<-n_kd/sla 2>/vert/sg/vert^<-n_k/sla 2>/exp/lb-/halv
	/lsb Tr(/sg^<-1>A_k)+n_k(/hm_k-/mu_k)'/sg^<-1>(/hm_k-/mu_k)/rsb/rb/cr
&/qquad(2/pi)^<-d/sla 2>/vert/sg/vert^<-1/sla 2>/exp/lb-/halv(x-/mu)'/sg^<-1>(x-/mu)/rb/cr
&=(2/pi)^<-(N+1)d/sla 2>/vert/sg/vert^<-(N+1)/sla 2>/exp/lb-/halv Tr(/sg^<-1>/sum_<k=1>^KA_k)/rb/cr
&/qquad/exp/lb-/halv/lsb/sum_<k=1>^Kn_k(/hm_k-/mu_k)'/sg^<-1>(/hm_k-/mu_k)
	+(x-/mu)'/sg^<-1>(x-/mu)/rsb/rb,/cr>$$
where $/hm_k,/hs_k=A_k/sla n_k$, $/hs$ are as in (3.18), (3.19), (3.21).

We will pursue the likelihood ratio program, and start with
$$</mo<H_0>>L(z,x)=</mo<k/leq K>/mo</mu=/mu_k>>L(z,x).$$
It is seen that the maximum likelihood obtained when
considering only the $n_1+/cdots+n_K$ observations in $z$ becomes
$$(2/pi)^<-Nd/sla 2>/vert/hs/vert^<-N/sla 2>/exp/lb-/halv Tr(/hs^<-1>N/hs)/rb=
	(2/pi e)^<-Nd/sla 2>/vert/hs/vert^<-N/sla 2>./eqno(6.21)$$
Consequently
%
%
$$</mo</mu=/mu_k>>L(z,x)=(2/pi e)^<-(N+1)d/sla 2>
/vert/wt/sg_<(k)>/vert^<-(N+1)/sla 2>,$$
where $/wt/sg_/k$ is the adjusted maximum likelihood
estimator of $/sg$ obtained by including $x$ as an observation from class
$k$,/ i.e.
$$/wt/sg_/k=<1/over N+1>(A_1+/cdots+/wt A_k+/cdots+A_K),$$
where
$$/wt A_k=/sum_<j=1>^<n_k>/twice</lp x_j^/k-/wt/mu_k/rp>'+/twice<(x-/wt/mu_k)>',$$
$$/wt/mu_k=<1/over n_k+1>/lp/sum_<j=1>^<n_k>x_j^/k+x/rp =/hm_k+<1/over n_k+1>(x-/mu_k).$$
Some manipulations give
$$/wt A_k=A_k+<n_k/over n_k+1>/twice<(x-/hm_k)>',/eqno(6.22)$$
implying
$$/wt/sg_/k=<N/over N+1>/lsb/hs+<1/over N>/,<n_k/over n_k+1>
	/twice<(x-/hm_k)>'/rsb.$$
By the identity
$$/eqalignno<
/vert B+cyy'/vert&=/vert B(I+B^<-1>cyy')/vert/cr
&=/vert B/vert(1+cy'B^<-1>y),&(6.23)/cr>$$
%
%
which can be seen to follow from (8.4.11) in Box and Tiao
(1973),
$$/eqalign<</mo</mu=/mu_k>>L(z,x)&=(2/pi e)^<-(N+1)d/sla 2>
/lp N+1/over N/rp^<(N+1)d/sla 2>/vert/hs/vert^<-(N+1)/sla 2>/cr
&/qquad/lb1+<1/over N>/,<n_k/over n_k+1>(x-/hm_k)'/hs^<-1>(x-/hm_k)/rb^<-(N+1)/sla 2>,/cr>$$
so that
$$/eqalignno<</mo<H_0>>L(z,x)&=(2/pi e)^<-(N+1)d/sla 2>
/lp N+1/over N/rp^<(N+1)d/sla 2>/vert/hs/vert^<-(N+1)/sla 2>/cr
&/qquad</mo<k/leq K>>/lb 1+<1/over N>/,<n_k/over n_k+1>(x-/hm_k)'
/hs^<-1>(x-/hm_k)/rb^<-(N+1)/sla 2>.&(6.24)/cr>$$

Next consider $/mo<a/ priori>L(z,x)$. $/mu$ is now
free to vary in $/rr^d$ and makes $L(z,x)$ maximal when it is set equal
to $x$. Hence
$$/mo<a/ priori>L(z,x)=/mo</sg>(2/pi)^<-(N+1)d/sla 2>
	/vert/sg/vert^<-(N+1)/sla 2>e^<-/halv Tr(/sg^<-1>/sum_<k=1>^KA_k)>.$$
The maximum is attained when $/sg=<1/over N+1>/sum_<k=1>^KA_k=<N/over N+1>/hs$,
using for example Lemma 3.2.2 in Anderson (1958). Hence
%
%
$$/mo<a/ priori>L(z,x)=(2/pi e)^<-(N+1)d/sla 2>
	/lp N+1/over N/rp^<(N+1)d/sla 2>/vert/hs/vert^<-(N+1)/sla 2>.$$

By this result and (6.24) it follows that
$$/eqalignno<
/lr(z,x)&=</mo<H_0>L(z,x)/over/mo<a/ priori>L(z,x)>/cr
&=/lb1+/mio<k/leq K><1/over N>/,<n_k/over n_k+1>(x-/hm_k)'/hs^<-1>(x-/hm_k)/rb^
<-(N+1)/sla2>.&(6.25)/cr>$$
Using arguments similar to those used in Section 6.1,
cf.~(6.5), (6.7), (6.9), (6.12), it is seen that a conservative
approximation to the (exact) likelihood ratio test is:
declare $X$ an outlier provided
$$<n_k/over n_k+1>(X-/hm_k)'/hs^<-1>(X-/hm_k)/geq b,/ k=1/upto K,/eqno(6.26)$$
where $b$ is chosen to fulfill
$$P_k/lb<n_k/over n_k+1>(X-/hm_k)'/hs^<-1>(X-/hm_k)/geq b/rb/leq/eps,/ k=1/upto K.$$

It remains to relate $b$ to $/eps$. Under $P_k$,
%
%
$X-/hm_k$ is normal with zero mean and covariance matrix $<n_k+1/over n_k>/sg$,
and independent of $A=N/hs$, which is Wishart$<>_d(N-K,/sg)$.
Hence
$$/eqalign<T^2(d,N-K)&=(N-K)<n_k/over n_k+1>(X-/hm_k)'A^<-1>(X-/hm_k)/cr
	&=<N-K/over N><n_k/over n_k+1>(X-/hm_k)'/hs^<-1>(X-/hm_k)/cr>$$
is Hotelling $(d,N-K)$ distributed, cf.~Mardia, Kent and Bibby (1979, p.~73). 
Furthermore,
$$T^2(d,N-K)=<d(N-K)/over N-K-d+1>F(d,N-K-d+1)$$
relates Hotelling to Fisher (op.~cit., p.~74).
Accordingly,
$$<n_k/over n_k+1>(X-/hm_k)'/hs^<-1>(X-/hm_k)=<Nd/over N-K-d+1>/,F(d,N-K-d+1),/eqno(6.27)$$
which implies that
$$b=<Nd/over N-K-d+1>/,F^<-1>(1-/eps;/,d,N-K-d+1)/eqno(6.28)$$
is the best choice in (6.26), $F^<-1>(1-/eps, g, h)$
denoting
%
%
the $(1-/eps)$-fractile in the $F$-distribution with $g$ and
$h$ degrees of freedom.

The outlier criterion may also be written
$$/eqalignno<
/wh f_k(X)&=N_d(/hm_k,/hs)(X)/cr
&/leq(2/pi)^<-d/sla 2>/vert/hs/vert^<-1/sla 2>/exp/lp-/halv/,<n_k+1/over n_k>b/rp,/ k=1/upto K,&(6.29)/cr>$$
with $b$ given by (6.28), depending upon $/eps,d,K,N$.
It may be shown that the criterion above converges almost surely
to the ``known density case'' (6.8) as the $n_k$'s go to infinity.
%
%
/subsection<6.2.B><Different covariance matrices>
Let $/fkx=N_d(/mu_k,/sg_k)(x),/ k=1/upto K$, let $Z$ of (3.3) be available,
and face a new vector $X$ from some $N_d(/mu,/sg)$ distribution, where
$/mu/in/rr^d$ and $/sg/in/lb/sg_1/upto/sg_K/rb$ </it a priori//>.
Is $X$ from one of the $K$ classes, or is it an outlier?

We must test
$$H_0:(/mu,/sg)/in/lb(/mu_1,/sg_1)/upto(/mu_K,/sg_K)/rb,/eqno(6.30)$$
and look at the likelihood ratio statistic.
The observed data has total likelihood
$$/eqalign<
L(z,x)&=L(z,x/mid/mu_1,/sg_1/upto /mu_k,/sg_k;/mu,/sg)/cr
&=(2/pi)^<-(N+1)d/sla 2>/prod_<k=1>^K/mid/sg_k/mid^<-n_k/sla 2>/mid/sg/mid
                      ^<-1/sla 2>/cr
&/qquad/exp/Biggl/lbrace-/halv/Biggl/lbrack/sum_<k=1>^K/lp Tr(/sg_k^<-1>A_k)+n_k(/hm_k-/mu_k)'
                      /sg_k^<-1>(/hm_k-/mu_k)/rp/cr
&/qquad/qquad+(x-/mu)'/sg^<-1>(x-/mu)/Biggr/rbrack/Biggr/rbrace./cr>$$
Using arguments and techniques similar to those used in 6.2.A we get
$$/mo<(/mu,/sg)=(/mu_k,/sg_k)>L(z,x)=(2/pi e)^<-(N+1)d/sla 2>
                      /prod_<l=1>^K/mid/hs_l/mid^<n_l/sla 2>$$
$$/lp<n_k+1/over n_k>/rp^<(n_k+1)d/sla 2>/mid/hs_k/mid^<-1/sla 2>W_k,$$
where
$$W_k=/lb 1+<1/over n_k+1>(x-/hm_k)'/hs_k^<-1>(x-/hm_k)/rb^<-(n_k+1)/sla 2>./eqno(6.31)$$
We also need
$$/mo</mu/in/rr^d,/sg=/sg_k>L(z,x)=(2/pi e)^<-(N+1)d/sla 2>/prod_<l=1>^K/mid/hs_l/mid
                      ^<-n_l/sla 2>$$
$$/lp<n_k+1/over n_k>/rp^<(n_k+1)d/sla 2>/mid/hs_k/mid^<-1/sla 2>,$$
which is arrived at by similar means again.
Hence
$$/lr(z,x)=</mo<k/leq K>(1+<1/over n_k>)^<(n_k+1)d/sla 2>/mid/hs_k/mid^<-1/sla 2>
    W_k/over/mo<k/leq K>(1+<1/over n_k>)^<(n_k+1)d/sla 2>/mid/hs_k/mid^<-1/sla 2>>
    /eqno(6.32)$$
is an exact expression for the likelihood ratio.

Fisher's likelihood ratio program urges us to find a constant $a$, depending
on $/eps, d, n_1/upto n_K$, such that $P_k/lb /lr(Z,X)/leq a/rb/leq/eps,/ k=1/upto K$.
A conservative approximation results if $a$ is chosen to satisfy
$$P_k/lb/lp 1+<1/over n_k>/rp^<(n_k+1)d/sla 2>/mid/hs_k/mid^<-1/sla 2>
    W_k/leq a/,S/rb/leq/eps,/ k=1/upto K,/eqno(6.33)$$
where
$$S=/mo<k/leq K>/lp 1+<1/over n_k>/rp^<(n_k+1)d/sla 2>/mid/hs_k/mid^<-1/sla 2>./eqno(6.34)$$
$P_k$ denotes probability distribution of $(Z,X)$ when $X$ is from population $k$.
$a$ could in principle be found from this requirement, but it would be extremely
difficult. An approximation can be worked out, writing (6.33) as
$$P_k/lb W_k/leq aU_k/rb/leq/eps,/ k=1/upto K,$$
and treating
$$U_k=/lp 1+<1/over n_k>/rp^<-(n_k+1)d/sla 2>/mid/hs_k/mid^<1/sla 2>S$$
as non-random for a moment. After all most of the randomness here belongs to
$W_k$; $U_k$ is close to $/mid/sg_k/mid^<1/sla 2>/sla/mid/sg_<k_1>/mid^<1/sla 2>$,
where $/mid/sg_<k_1>/mid=/mio<k/leq K>/mid/sg_k/mid$,whereas
$$W_k/eqtop<D>/lb 1+<d/over n_k-d>F(d,n_k-d)/rb^<-(n_k+1)/sla 2>,/eqno(6.35)$$
$F(g,h)$ again denotes a $F$-distribution with $g$ and $h$ degrees of freedom.
This latter result is shown as in (6.27).

Using this device as our requirement for $a$ is that
$$Pr/lb F(d,n_k-d)/geq<n_k-d/over d>/lsb(a/,U_k)^<-2/sla(n_k+1)>-1/rsb/rb/leq/eps,$$
or
$$<n_k-d/over d>/lsb(a/,U_k)^<-2/sla(n_k+1)>-1/rsb/geq F^<-1>(1-/eps,d,n_k-d),$$
$k=1/upto K$. On these ground we choose
$$/eqalign</ha_0&=/mio<k/leq K> U_k^<-1>/lsb1+<d/over n_k-d>F^<-1>
      (1-/eps,d,n_k-d)/rsb^<-(n_k+1)/sla 2>/cr
&=/lb/mo<k/leq K> U_k/lsb 1+<d/over n_k-d> F^<-1>(1-/eps,d,n_k-d)/rsb^<(n_k+1)
      /sla 2>/rb^<-1>.>/eqno(6.36)$$
$X$ is considered an outlier if $W_k/leq/ha_0/,U_k$ for each $k=1/upto K$. This
criterion can also be written
$$/widetilde/fk(X)=(2/pi)^<-d/sla 2>/mid/hs_k/mid^<-1/sla 2>W_k/leq(2/pi)
      ^<-d/sla 2>/mid/hs_k/mid^<-1/sla 2>/ha_0/,U_k,/ktok./eqno(6.37)$$

This may be considered the appropriate fine-tuning of the coarser criterion
(6.16), taking sampling variability of estimates into account. Notice from
(6.31) that $/wt/fk(X)$ above is close to $/wh/fk(X)=N_d(/hmk,/hsk)(X)$ for
large $n_k$. Indeed, (6.37) can also be given in the form
$$/wh/fk(X)/leq(2/pi)^<-d/sla 2>/mid/hsk/mid^<-1/sla 2>/exp/lb-/halv(n_k+1)
      /lsb(/ha_0 U_k)^<-2/sla(n_k+1)>-1/rsb/rb,/ k=1/upto K./eqno(6.38)$$
It is also worth pointing out that $/mid/hsk/mid^<-1/sla 2>/ha_0 U_k$
becomes close to $/mid/sg_<k_0>/mid^<-1/sla 2>e^<-/halv/gamma_<d;1-/eps>>$ when
the $n_k$'s are large, i.e.~the r.h.s.~of (6.37), and that of (6.38), converge
to the r.h.s.~of (6.16), almost surely. To see this, start with the
approximation $S=e^<d/sla 2>$$(/mio<l/leq K>/break/mid/hs_l/mid)^<-1/sla 2>$
from (6.34) and proceed to
$$U_k/eqtop<.>e^<-d/sla 2>/mid/hsk/mid^<1/sla 2> S/eqtop<.>/mid/hsk/mid^<1/sla 2>
      /lp/mio<l/leq K>/mid/hs_l/mid/rp^<-1/sla 2>.$$
It can be shown that
$$/lsb 1 + <d /over n_k-d> F^<-1>(1-/eps,d,n_k-d)/rsb^<-(n_k+1)/sla 2>
      /longrightarrow /exp/lp<-/halv/gamma_<d;1-/eps>>/rp/eqno(6.39)$$
when $n_k$ grows to infinity. Hence
$$/ha_0/eqtop<.>/mio<k/leq K>U_k^<-1>e^<-/halv/gamma_<d;1-/eps>>/eqtop<.>
 /lp/mo<l/leq k>/mid/hs_l/mid/rp^<-1/sla 2>/lp/mio<l/leq K>/mid/hs_l/mid
      /rp^<1/sla 2> /exp/lp-/halv/gamma_<d;1-/eps>/rp.$$
Putting things together the result is
$$/mid/hsk/mid^<-1/sla 2>/ha_0 U_k/eqtop<.>/lp/mo<l/leq K>/mid/hs_l/mid/rp
      ^<-1/sla 2>/exp/lp<-/halv/gamma_<d;1-/eps>>/rp,$$
which proves the claim made above.

Another, and related, approach to the problem is as follows. Test each of the
hypoteses $H_k:(/mu,/sg)=(/mu_k,/sg_k)$ versus $/mu/in/rr^d-/lb/mu_k/rb,
/sg=/sg_k$, using the likelihood ratio criterion. Call $X$ an outlier if
all of these $K$ hypoteses are rejected. Using results of earlier efforts
one obtains the following likelihood ratio statistic for $H_k$:
$$/eqalign</lr_k(z,x)&=</mo<(/mu,/sg)=(/mu_k,/sg_k)> L(z,x)
                /over   /mo</mu/in/rr^d,/sg=/sg_k> L(z,x)>/cr
                     &=/lb 1+ <1/over n_k+1>/lp x-/hm_k/rp'/hs_k^<-1>
                /lp x-/hm_k/rp/rb^<-(n_k+1)/sla 2>/cr
                     &= W_k./cr>$$
$H_k$ is rejected at level $/eps$ provided
$$W_k/leq/lsb 1+<d/over n_k-d> F^<-1>/lp 1-/eps,d,n_k-d/rp/rsb^<-(n_k+1)/sla 2>
     =/sigma_k(/eps),/eqno(6.40)$$
 using (6.35).

The outlier criterion resulting from this second approach can be written in
several ways, also, a variety of approximations are available. We may take
$$/wt/fk(X)=(2/pi)^<-d/sla 2>/mid/hsk/mid^<-1/sla 2> W_k/leq(2/pi)^<-d/sla 2>
     /mid/hsk/mid^<-1/sla 2>/sigma_k(/eps),/ktok/eqno(6.41)$$
as a natural description, or
$$/eqalign<
(X-/hmk)'/hsk^<-1>(X-/hmk)&/geq(n_k+1)/lsb/sigma_k(/eps)^<-2/sla(n_k+1)>-1/rsb/cr
                          &=(n_k+1)<d/over n_k-d> F^<-1>(1-/eps,d,n_k-d),
                          /ktok.>$$
An approximation for the r.h.s.~here is just $/gamma_<d;1-/eps>$.

It is instructive to compare this latter proposed approach to the first one.
Using (6.40) and (6.39) we see that (6.41) essentially says $/wt/fk(X)/leq
(2/pi)^<-d/sla 2>/mid/hsk/mid^<-1/sla 2> e^<-/halv/gamma_<d;1-/eps>>, /ktok$,
whereas (6.37) essentially amounts to $/wt/fk(X)/leq(2/pi)^<-d/sla 2>
(/mo<l/leq K>/mid/hs_l/mid)$
$e^<-/halv/gamma_<d;1-/eps>>$, $/ktok.$ Thus, an outlier according to the
``simultaneous'' approach usually is an outlier according to the ``K separate
tests'' approach, but not conversely. This points to the slight disagreement
commented upon in the remark following (6.12).
/section<6.3><Concluding remarks><> The outlier problem commonly 
treated in the literature is that of finding outliers among a collection
of data points, in contrast to the problem considered above, where class
descriptions are established initially and new incoming vectors are to
be tested against these. But the first type of considerations are important
when these  initial  class  descriptions  are determined---if a couple of 
extreme data points  are allowed to artificially inflate the $/hs$-matrix 
for a class then  perhaps  no  incoming  vectors  will  ever  qualify  as 
outliers.  Methods for obtaining  outlier-resistant  or  robust estimates 
for means and covariance matrices are discussed in Section 10.1.A.

The methods here that use the simultaneous likelihood ratio approach,
and include estimator variability for the  class descriptions, seem to be 
new.

We  found  it  convenient  and important to  have  a  reject  option  for 
``incredible''  symbols,  chiefly  to  reduce error rates.   It  is  also 
important  to  filter  out  these outliers when the updating  methods  of 
Chapter 7 are used.  Parameter estimates could again be misleading if too 
strange feature vectors are allowed to enter the computations.

One  needs  outlier  criteria for classification procedures based on other 
models than  the normal. But in many cases the criteria developed here are 
satisfactory even outside  the  normal  models.    An  exception is cases 
involving discrete feature components.    We have developed outlier rules 
for the discrete-times-normal model studied in Chapter 9, but details are 
not given here.


/chapter<7><Updating parameter estimates/hfil/break
using unclassified vectors>
We have studied classification procedures based on parametric models in Chapters
3 and 4. The methods of Chapter 3 simply replace unknown 
parameters in the underlying optimal Bayes classifiers with
estimates based on training sets. The methods of Chapter 4,
although different in character, amount to essentially the same, but
take sampling variability of estimates into account. In any case
only information from the training sets is used to construct the 
classification algorithms.

Besides being instrumental in the classification of future objects
the parameter estimates of course have some interest in their
own right in that they may provide insight in the mechanisms
generating the data.

The purpose of the present chapter is to show that parameter
estimates actually may be </it updated using unclassified vectors//>.
In the symbol recognition context this could mean adjusting one's
training set based parameter estimates especially for a given ``map''
of unclassified symbols, and then performing the classification with the 
adjusted discriminant rules. Another possibility could be to update
one's class descriptions </it on the run//>, i.e.~adjusting estimates
after each new vector, in spite of the fact that these new
vectors have unconfirmed class memberships. This sequential idea
is fascinating, but the former application will usually be more
practical and important.

These updated estimates will be more reliable
than the initial ones in that much more information is drawn upon in their
construction, and should also result in better classification accuracy.
Cases where even a tiny (initial) training set is sufficient for the
updating machinery to work provide the most interesting and dramatic
examples of the usefulness of the approach developed here. Using 
these methods the sometimes costly and tedious training phase
can, in many cases, be dramatically simplified. Let the 
computer work for us!

What makes the updating approach possible is the fact that a feature vector
$X$ evaluated from an object randomly chosen from the set of unclassified
objects (i.e.~not belonging to any training set) has density 
$$f(x)=/pi_1/,f_1(x)+/cdots+/pi_K/,f_K(x),/eqno(7.1)$$
i.e.~a mixture of the $K$ class densities with the prior probabilities as 
mixture proportions. Thus it is not a contradiction in terms when 7.1 below 
provides methods giving estimates of </it a priori//> probabilities
for a given map of unclassified objects, based on known (or fixed
estimates of) class densities.
/note<Remark 1.> Occasionally our pattern recogniser may encounter alien vectors
not originating from any of the $K$ predefined classes at all. This will 
happen if a non-numeral symbol falsely is taken as a numeral by the 
preprocessing machinery, for example. It is assumed in (7.1) that
the object in question really belongs to one of the classes.
In the rest of this chapter we assume that all outliers, i.e.~vectors from 
``incredible symbols'', have been detected separately and removed, for example by one of the
methods presented in Chapter 6.
%
%
/medskip
Section 7.2 presents updating procedures for some of the parametric models 
described in Chapters 3 and 4. The resulting parameter estimates
combine information from both training samples and unclassified
feature vectors. In particular, the traditional Gaussian cases are considered.
/section<7.1><Estimating a priori probabilities on the given map><>
Assume that the class densities $/ftfk$ are </it known//>, or that
reliable estimates have been obtained, for example on the basis
of training sets and/sla or past experience.
%
%
$/fkx$ could typically be $N_d(/hm_k,/hs_k)(x)$ with
$/hm_k,/hs_k$ obtained from training sets as in (3.18), (3.19).
We wish to point out, however, that the $f_k$'s may be
completely arbitrary as far as the present section is concerned,
they might for example have been obtained via nonparametric density
estimation.

The vectors from objects outside the training set follow the mixture
density (7.1) (but see the remark following that equation).
Consider $M$ new independent vectors $X_1/upto X_M$. These have
simultaneous likelihood
$$/prod_<i=1>^M f(X_i)=/prod_<i=1>^M/lb/sum_<k=1>^K/pi_k /,f_k(X_i)/rb./eqno(7.2)$$
%
%

Now several procedures are possible for the estimation of
$/pi_1 /upto/pi_K$ on the basis of $X_1/upto X_M$. We review first the maximum
likelihood method, which has been applied to the mixed
model problem before in statistical literature, and then present
a new method.
/subsection<7.1.A><Maximum likelihood estimators> The maximum likelihood 
principle urges us to find the parameter set 
$/hp_k/upto/hp_K$, $/sum_<k=1>^K$$/hp_k/allowbreak=1$,
that maximises (7.2). Consider its logarithm
$$H(/pi_1 /upto/pi_<K-1>)=/sum_<i=1>^M/log/lb/sum_<k=1>^K/pi_k 
/,f_k(X_i)/rb,/eqno(7.3)$$
where $/pi_K =1-/sum_<k=1>^<K-1>/pi_k $. Then
$$/eqalign<
</partial H/over/partial/pi_k >&=/sum_<i=1>^M/lsb f_k(X_i)-f_K(X_i)/rsb/sla f(X_i)/cr
&=<1/over/pi_k >/sum_<i=1>^M P(k/mid X_i)-<1/over/pi_K >/sum_<i=1>^M P(K/mid X_i),/cr>$$
where
$$/pkx=/pi_k /,/fkx/sla/sum_<m=1>^K/pi_m /,f_m(x)/eqno(7.4)$$
%
%
is the posterior probability of class $k$ given that vector $X=x$ is observed.
It follows that a maximum point for $H$ must fulfil the equations
$$/pi_k =<1/over M>/sum_<i=1>^M P(k/mid X_i),/ k=1/upto K./eqno(7.5)$$
These are indeed natural equations since $/sum_<i=1>^M P(k/mid X_i)$ is the 
expected number of vectors from class $k$ among the $M$, given
the observations $X_1/upto X_M$.

The log-likelihood function $H$ is fortunately concave, so that it has a unique
maximum point, i.e.~the equations (7.5) have only one solution. $H$ is concave 
since the matrix with elements
%
%
$$-</partial^2 H/over/partial/pi_k /partial/pi_m >=
/sum_<i=1>^M/bigl/lbrack f_k(X_i)-f_K(X_i)/bigr/rbrack/bigl/lbrack
 f_m(X_i)-f_K(X_i)/bigr/rbrack/sla f(X_i)^2$$
is positive definite, as is easily shown.

The $/pi_k $'s enter into the r.h.s.~of (7.5), and the 
equations cannot be solved explicitly or recursively. The
solutions $/hp_1/upto/hp_K$ have to be arrived at by an iterative
computational procedure. The Newton-Raphson technique may be applied, but the 
easiest and completely satisfactory method is to start out with initial
values, say $/pi^<(0)>_k=<1/over K>,/ k=1/upto K$, and then use
(7.5) as iteration equations, i.e.
$$/pi^<(t+1)>_k=<1/over M>/sum_<i=1>^M P^<(t)>(k/mid X_i),/ k=1/upto K,
/eqno(7.6)$$
%
%
$t=1,2/upto$ where $P^<(t)>(k/mid x)=/pi^<(t)>_k/,f_k(x)/sla/sum_<m=1>^K/pi^<(t)>_m/,f_m(x)$.
This may be seen to be an application of the so-called EM algorithm,
cf.~Dempster, Laird, and Rubin~(1977). The general properties of
this algorithm (cf.~Wu,~1983) ensure that the sequence 
$/pi^<(t)>_k,/,t=1,2,/ldots$ will converge to the maximum likelihood
estimator $/hp_k$ regardless of the chosen initial values; $k=1/upto K$.
/subsection<7.1.B><Bayes estimators> The Bayes estimation principle requires in the
present problem the specification of an </it a priori//> distribution for 
$(/pi_1 /upto/pi_K )$. Assume 
%
%
that this distribution is given by a density
$$g(/pi)=g(/pi_1 /upto/pi_<K-1> )/eqno(7.7)$$
for $/pi=(/pi_1 /upto/pi_<K-1>)$ in the simplex
$$B=/lb/pi/mid/pi_<k>/geq0,/ k=1/upto K-1,/sum_<k=1>^<K-1>
/pi_k/leq1/rb./eqno(7.8)$$
Then $/pi$ given the data $X_1/upto X_M$ has posterior
density
$$g(/pi)/,L_M(/pi)/,d/pi/sla/int_B g(/pi)/,L_M(/pi)/,d/pi,$$
where
$$/eqalignno<
L_M(/pi)&=L_M(/pi,X_1/upto X_M)/cr
        &=/prod_<i=1>^M/lb/sum_<k=1>^K/pi_k/,f_k(X_i)/rb&(7.9)/cr>$$
is the likelihood of the observed data, given the unknown parameters. The
Bayes estimators are the posterior expectations
%
%
$$/eqalignno<
/wt</pi>_k&=E/lb/pi_k/mid X_1/upto X_M/rb/cr
           &=N_M/k/sla D_M/k,&(7.10)/cr>$$
where
$$/eqalignno<
N_M/k&=/int_B/pi_k/,g(/pi)/,L_M(/pi)/,d/pi,&(7.11)/cr
D_M/k&=/int_B g(/pi)/,L_M(/pi)/,d/pi,&(7.12)/cr>$$
$k=1/upto K$.

It is difficult to compute the Bayes estimates in practice.
The likelihood $L_M(/pi)$ may be expanded in a large sum of terms
of the type
$$/pi_<j_1>/cdots/pi_<j_M>/,f_<j_1>(X_1)/cdots f_<j_M>(X_M),/eqno(7.13)$$
where $1/leq j_i/leq K,/ i=1/upto M$. Consequently, both
the numerator and denumerator of $/wt</pi>_k$ may in theory
be written as a sum 
%
%
of $K^M$ terms, each of which is explicitly computable from moments
$$/eqalignno<
A(a_1/upto a_K)&=E/pi_1^<a_1>/cdots/pi_K^<a_K>/cr
               &=/int_B/prod_<k=1>^K/pi_k^<a_k>g(/pi)/,d/pi&(7.14)/cr>$$
(writing $/pi_K=1-/sum_<k=1>^<K-1>/pi_k$) of the prior
distribution. Such a procedure might well take a computer several
years unless $M$ is small, however, so that alternative procedures
have to be devised.

A possible method consists of obtaining $N_M/k$ and 
$D_M/k$ by numerical integration. This would also pose 
implementational difficulties since typically $L_M(/pi)$ would
be a product of a large number of possibly
%
%
small terms.

We will instead outline another method. When expanding $L_M(/pi)$
the terms (7.13) may be sorted and collected together in a way
that gives
$$L_M(/pi)=/sum_<C(M)>/pi_1^<a_1>/cdots/pi_K^<a_K>/,S_M(a_1/upto a_K)/eqno(7.15)$$
where $C(M)$ is the set of $(a_1/upto a_K)$ having 
$a_k/geq0,/ k=1/upto K$, and $/sum_<k=1>^Ka_k=M$.
$S_M(a_1/upto a_K)$ is the sum of all terms $f_<j_1>(X_1)/cdots f_<j_M>(X_M)$
that have exactly $a_k$ $f_k$-factors, $k=1/upto K$.
$S_M(a_1/upto a_K)$ has $M!/sla/lb/prod_<k=1>^Ka_k!/rb$ terms, whereas
the number of terms in the expansion (7.15) may be seen to be
%
%
$$/eqalignno<
/#C(M)&=/sum_<a_1=0>^M/sum_<a_2=0>^<M-a_1>/,/cdots/,
	/sum_<a_<K-1>=0>^<M-a_1-/cdots-a_<K-2>> 1/cr
      &=<1/over(K-1)!>(M+1)/cdots(M+K),&(7.16)/cr>$$
utilising the identity
$$/sum_<j=1>^n j(j+1)/cdots(j+b-1)=<1/over b+1>/,n(n+1)/cdots(n+b),$$
which can be proved by induction.

By (7.11), (7.12), (7.14), and (7.15)
$$/eqalign<
N_M/k&=/sum_<C(M)>A(a_1/upto a_k+1/upto a_K)/,S_M(a_1/upto a_K),/cr
D_M/k&=/sum_<C(M)>A(a_1/upto a_k/upto a_K)/,S_M(a_1/upto a_K),/cr>$$
giving us a method to obtain  $/wt</pi>_k$, provided
the moments $A(a_1/upto a_K)$ are easy to compute
(as they are w.r.t.~a Dirichlet prior density, for
%
%
example), provided $/#C(M)$ in (7.16) is a manageable number (notice
that it is dramatically less than $K^M$ in the typical pattern recognition
application), and provided the terms $S_M(a_1/upto a_K)$ can be easily
computed.

It can now be seen that 
$$S_M(a_1/upto a_K)=/sum_<k=1>^Kf_k(X_M)/,S_<M-1>(a_1/upto a_k-1/upto a_K),/eqno(7.17)$$
which holds even if some $a_k$'s are zero, if $S_M(a_1/upto a_K)$
is defined as zero when some $a_k$ equals $-1$. This makes it possible to
obtain the $S_M$'s in a recursive way. At intermediate step $m$ the full set 
of $(m+1)/cdots(m+K)/sla(K-1)!/quad$ $S_m$'s will need to be stored.
%
%

Suppose for example that the uniform prior $g_0(/pi)/equiv(K-1)!,/,/pi/in B$ of
(7.8), is chosen, which is one way of formalising ``total
ignorance </it a priori//>'' about $/pi_1/upto/pi_K$. Then
$$E/prod_<k=1>^K/pi_k^<a_k>=/lb/prod_<k=1>^Ka_k!/rb/sla K/cdots(K+M-1),$$
from known properties of Dirichlet distributions.
Hence in this case the Bayes estimates can be computed as 
$$/wt</pi>_k=
<1/over K+M>
/,</sum_<C(M)>a_1!/cdots a_K!/,(a_k+1)/,S_M(a_1/upto a_K)/over
 /sum_<C(M)>a_1!/cdots a_K!/,S_M(a_1/upto a_K)>,/eqno(7.18)$$
relying on (7.17) to find the $S_M$'s.
%
%
/note<Discussion.> The Bayes estimates (7.10) are generally harder to compute
than the maximum likelihood estimates (7.5) -- (7.6).
They may have some statistical advantages for moderate sample sizes,
however, as prior knowledge often exists in a form that could be
approximated with an appropriate Dirichlet density. The computations
needed in for example (7.18) increase considerably with the
number of classes, whereas the maximum likelihood procedure is essentially
unbothered.

One may prove that Bayes estimators and maximum likelihood estimators are
asymptotically equivalent, i.e.
%
%
$$/sqrt<M>(/wh/pi_k-/wt/pi_k)/totop<P>0$$
as $M/to/infty,/ k=1/upto K$, regardless of prior distribution (7.7).
This follows from a multidimensional extension of Theorem~6.7.2 in Lehmann
(1983) in the case where the assumed model (7.2) is absolutely correct, 
and may be shown to be true even in the more realistic case where the model is
considered to be an approximation only. It can also be shown that both 
$/wh/pi_k$ and $/wt/pi_k$ are asymptotically optimal.

Thus we may recommend the maximum likelihood procedure on the grounds
of computational simplicity and statistical objectivity.
The Bayes estimators were included in this chapter partly
%
%
because of general theoretical interest. Bayes estimators may constitute
good alternatives to maximum likelihood estimators in the more general
mixture distribution problem where the class densities have unknown parameters,
in which case there are severe problems with the latter approach.

On the whole, however, Section 7.1.B may be considered a digression from 
the main theme of the present chapter.
%
%
/vfill/eject
/section<7.2><Simultaneous updating of class><parameters and prior probabilities>
/subsection<7.2.A><General program> Section 7.1 was concerned with the model 
(7.1) in the ideal case where
the class densities $/ftfk$ were known. Assume now that there are unknown 
parameters in them, say
$$/fkx=f(x,/th_k),/eqno(7.19)$$
where $/th_k$ is a $p$-dimensional parameter characterising class
$k;/ k=1/upto K$. Hence vectors from objects whose class labels are unknown 
follow the mixture distribution
$$f(x)=/sum_<k=1>^K/pi_k/,f(x,/th_k)./eqno(7.20)$$

Let as in Section 7.1 $X_1/upto X_M$ be
%
%
feature vectors from a set of $M$ unclassified  objects, so that
their simultaneous likelihood becomes 
$$L_</rm rest>=/prod_<i=1>^M/lb/sum_<k=1>^K/pi_k/,f(x,/th_k)/rb./eqno(7.21)$$

It is now possible to construct estimators for $/pi_k,/th_k,/ k=1/upto K$ 
</it simultaneously//>, based on $X_1/upto X_M$ only. Indeed, the maximum
likelihood approach may be invoked, and results in natural likelihood equations
for the  estimates. These equations may have several solutions, however,
and another difficulty is that the likelihood (7.21) may have unpleasant 
singularities at edges of the parameter space. Nevertheless, the described
%
%
approach will be successful if consistent  estimates $/wh/th_<k,0>$ from
the training sets are available and if the needed iterative computational
procedure have $/pi_k=<1/over K>,/,/th_k=/wh/th_<k,0>,/allowbreak
/ k=1/upto K$, as initial values.

Better yet is to utilise the training data in combination with the ``new''
vectors $X_1/upto X_M$. Assume that a training set of the form
$$Z=/lb X_j^/k;/ j=1/upto n_k,/ k=1/upto K/rb/eqno(7.22)$$
is available, where $X_1^/k/upto X_<n_k>^/k$ are from class $k$.
%
%
$(Z,X_1/upto X_M)$ has total likelihood 
$$/eqalignno<
/lfull&=L_</rm training>/,L_</rm rest>/cr
      &=/prod_<k=1>^K/prod_<j=1>^<n_k>f/lp X_j^/k,/th_k/rp/cdot
	/prod_<i=1>^M/lb/sum_<k=1>^K/pi_k/,f_k(X_i)/rb.&(7.23)/cr>$$
The likelihood equations become
$$/eqalign<
</partial/log/lfull/over/partial/pi_k>&=
	<1/over/pi_k>/sum_<i=1>^MP(k/mid X_i)-<1/over/pi_K>
	/sum_<i=1>^MP(K/mid X_i),/ k=1/upto K,/cr
</partial/log/lfull/over/partial/th_<k,l>>&=
	/sum_<j=1>^<n_k></partial/log f/lp X_j^/k,/th_k/rp
	/over/partial/th_<k,l>>+
	/sum_<i=1>^MP(k/mid X_i)</partial/log 
        f(X_i,/th_k)/over/partial/th_<k,l>>,/cr>$$
$l=1/upto p,/ k=1/upto K$, where
$$/pkx=</pi_k/,f(x,/th_k)/over/sum_<m=1>^K/pi_m/,f(x,/th_m)>./eqno(7.24)$$
It follows that any local maximum $(/pi_1^/ast/upto
/pi_K^/ast,/th_1^/ast/upto/th_K^/ast)$ of $/lfull$ must satisfy
$$/eqalignno<
/pi_k&=<1/over M>/sum_<i=1>^MP(k/mid X_i),/ k=1/upto K,&(7.25)/cr
</partial/log/lfull/over/partial/th_<k,l>>&=0,/ l=1/upto p,/ k=1/upto K.&(7.26)/cr>$$
In the presence of $L_</rm training>$ the full likelihood
will most often have a global maximum 
%
%
(whereas $L_</rm rest>$ alone may be unbounded).
One may prove that the sequence of estimators obtained by iterating
equations (7.25)--(7.26) in an EM-manner with starting values
$/pi_k=<1/over K>,/,/th_k=/wh/th_<0,k>,/ k=1/upto K$, where
$/wh/th_<0,k>$ based on $X_1^/k/upto X_<n_k>^/k$, is consistent.
The asymptotic framework implicitly referred to here is one where both the 
$n_k$'s and $M$ grow towards infinity. The practical interpretation,
as far as a typical pattern recognition application is concerned, is that the 
start values above will suffice to find the correct bump on the likelihood
surface, even for 
%
%
moderately sized training sets. 
/note<Remark 2.> It is silently assumed in (7.23) that the training class
sizes $n_1/upto n_K$ have been obtained by a procedure having nothing to do 
with the $/pi_k$'s; they could for example have been decided upon in advance.
In the case where vectors are recruited to $Z$ in (7.22) independently of each
other, and with probability $/pi_k$ of coming from class $k$, then 
$L_</rm rest>$ and $L_</rm full>$ have an extra factor $/pi_1^<n_1>/cdots
/pi_K^<n_K>$. The program derived from (7.23) can be 
easily altered to give updating formulae also in this second case.
/subsection<7.2.B><Normal case>
Let us apply the program above to the important case
where 
$$/eqalignno<
/fkx&=N_d(/mu_k,/sg_k)(x)/cr
    &=(2/pi)^<-d/sla2>/vert/sg_k/vert^<-1/sla2>/exp/lb-/halv(x-/mu_k)'
	/sg_k^<-1>(x-/mu_k)/rb,&(7.27)/cr>$$
$k=1/upto K$. Then some manipulation entails
$$/displaylines<
</partial/log/fkx/over/partial/mu_<k,a>>=
	/Lambda_<k,(a)>(x-/mu_k),/cr
</partial/log/fkx/over/partial/lam_<k,ab>>=
	/lp1-/halv/dl_<ab>/rp/lsb/sigma_<k,ab>-(x_a-/mu_<k,a>)
	(x_b-/mu_<k,b>)/rsb,/cr>$$
$a,b=1/upto d$, where $/Lambda_k=/sg_k^<-1>=/lb/lam_<k,ab>/rb
=/lb/sigma_<k,ab>/rb^<-1>$ and $/Lambda_<k,(a)>$ is the $a$'th row
in $/Lambda_k$. Hence the estimates we seek must obey
$$/Lambda_<k,(a)>/lb/sum_<j=1>^<n_k>/lp X_j^/k-/mu_k/rp+/sum_<i=1>^M
	P/lp k/mid X_i/rp/lp X_i-/mu_k/rp/rb=0,$$ 
$$/displaylines</qquad/sum_<j=1>^<n_k>/lb/sigma_<k,ab>-/lp 
X_<j,a>^/k-/mu_<k,a>/rp/lp X_<j,b>^/k-/mu_<k,b>/rp/rb/hfill/cr
/hfill+/sum_<i=1>^MP(k/mid X_i)/lb/sigma_<k,ab>-/lp X_<i,a>-/mu_<k,a>/rp
	/lp X_<i,b>-/mu_<k,b>/rp/rb=0,/qquad/cr>$$
%
%
for $k=1/upto K$ and $a,b=1/upto d$, where
$$/pkx=</pi_k/,N_d(/mu_k,/sg_k)(x)/over/sum_<m=1>^K/pi_m/,N_d(/mu_m,/sg_m)(x)>./eqno(7.28)$$
The equations simplify to
$$/mu_k=</sum_<j=1>^<n_k>X_j^/k+/sum_<i=1>^MP(k/mid X_i)X_i/over
	n_k+/sum_<i=1>^MP(k/mid X_i)>,/eqno(7.29)$$
$$/sg_k=</sum_<j=1>^<n_k>/lp X_j^/k-/mu_k/rp/lp X_j^/k-/mu_k/rp'+
	/sum_<i=1>^MP(k/mid X_i)(X_i-/mu_k)(X_i-/mu_k)'/over
	n_k+/sum_<i=1>^MP(k/mid X_i)>,/eqno(7.30)$$
$k=1/upto K$. These equations may be considered as natural generalisations of known ones 
for the case where data from only the mixture distribution are available, 
cf.~Duda~and Hart~(1973, p.~200), to the case where also training data
are present.
%
%

The maximum likelihood estimators from the training sets are
$$/eqalignno<
/hm_k&=<1/over n_k>/sum_<j=1>^<n_k>X_j^/k,&(7.31)/cr
/hs_k&=<1/over n_k>/sum_<j=1>^<n_k>/lp X_j^/k-/hm_k/rp/lp X_j^/k-/hm_k/rp',
&(7.32)/cr>$$
$k=1/upto K$, cf. 3.2.B. We are
now in a position to define the updated estimates 
$/pi_k^/ast ,/mu_k^/ast ,/sg_k^/ast ,/ k=1/upto K$. Compute successive
estimates $/pi^<(t)>_k$, $/mu_k^<(t)>$, $/sg_k^<(t)>$, $k=1/upto K$,
using the iteration equations
$$/pi_k^<(t+1)>=<1/over M>/sum_<i=1>^MP^<(t)>(k/mid X_i),/eqno(7.33)$$
$$/mu_k^<(t+1)>=<n_k/hm_k+/sum_<i=1>^MP^<(t)>(k/mid X_i)X_i/over
	n_k+/sum_<i=1>^MP^<(t)>(k/mid X_i)>,/eqno(7.34)$$
/goodbreak
$$/eqalignno<
/sg_<k>^<(t+1)>&=/Biggl/lbrace n_k/hs_k+n_k/lp/mu_<k>^<(t+1)>-/hm_k/rp
		/lp/mu_<k>^<(t+1)>-/hm_k/rp'/cr
&/quad+/sum_<i=1>^MP^<(t)>/lp k/mid X_i/rp/lp X_i-/mu_<k>^<(t+1)>/rp
/lp X_i-/mu_<k>^<(t+1)>/rp'
/Biggr/rbrace/sla/lb n_k+/sum_<i=1>^MP^<(t)>(k/mid X_i)/rb;/cr
&&(7.35)/cr>$$
%
%
$t=0,1,2,/ldots$, with initial values
$$/pi^<(0)>_/k=<1/over K>,/,/mu_k^<(0)>=/hm_k,/,/sg_k^<(0)>=/hs_k,
	/ k=1/upto K./eqno(7.36)$$
Here $P^<(t)>(k/mid X_i)$ is as in (7.28), but with $t$-th generation
parameter values. The iteratively defined sequence of points $/lb/pi^<(t)>_/k,
/mu_k^<(t)>,/sg_k^<(t)>;/ k=1/upto K/rb$ in 
$K/bigl(1+d+/halv d(d+1)/bigr)$-dimensional space are ensured convergence
as $t/to/infty$ by the general theory of EM algorithms, cf.~Wu (1983).
The limiting values are our updated estimates $/pi_k^/ast,
/mu_k^/ast,/sg_k^/ast;/ k=1/upto K$. (They are also, hopefully, 
maximising the full likelihood (7.23), but there is no guarantee for 
this. The limiting values could by chance correspond to a local, 
non-global maximum, for example.)

Note that all $K/bigl(1+d+/halv d(d+1)/bigr)$ points in the $t$-th
generation need to be computed before the $(t+1)$-st generation
points can be obtained. Observe further that fewer iterations will be
needed to obtain the maximum likelihood estimates, probably, if 
$/pi_k^<(0)>=/hp_k$ is used in (7.36), where $/hp_k$ is found as in Section 
7.1.
%
%

Next consider the case where the covariance matrices are taken to be
equal, i.e.
$$/fkx=N_d(/mu_k,/sg)(x),/ k=1/upto K./eqno(7.37)$$
Then the program above may be applied with a few technical changes.
The essential equations for the maximum likelihood estimators in the previous
case were (7.25), (7.29), (7.30).
After elaborations similar to those presented there one arrives in the present
case at (7.25), (7.29) as before, whereas the $K$ equations in (7.30) are 
replaced by
$$/sg=</sum_<k=1>^K/sum_<j=1>^<n_k>/twice</lp X_j^/k-/mu_k/rp>'+
	/sum_<i=1>^M/sum_<k=1>^KP(k/mid X_i)/twice<(X_i-/mu_k)>'/over
	N+/sum_<i=1>^M/sum_<k=1>^KP(k/mid X_i)>,/eqno(7.38)$$
where $N=/sum_<k=1>^Kn_k$. An iterative computational procedure
can now be put up, resulting in a sequence of parameter points
$/lb/pi^<(t)>_k,/mu_k^<(t)>,/sg^<(t)>/rb,/,t=0,1,2,/ldots$, converging
to estimates $/pi_k^</ast/ast>,/mu_k^</ast/ast>,
/sg^</ast/ast>$. Start out with
$$/pi^<(0)>_k=<1/over K>,/,/mu_k^<(0)>=/hm_k,/,/sg^<(0)>=
	/sum_<k=1>^K<n_k/over N>/hs_k=/hs;/eqno(7.39)$$
variations exist, for example $/pi^<(0)>_/k=/hp_k$
obtained from $X_1/upto X_M$ as described in Section 7.1. Then compute 
successively
$$/pi_k^<(t+1)>=<1/over M>/sum_<i=1>^MP^<(t)>(k/mid X_i),/eqno(7.40)$$
$$/mu_k^<(t+1)>=<n_k/hm_k+/sum_<i=1>^MP^<(t)>(k/mid X_i)X_i/over
	n_k+/sum_<i=1>^MP^<(t)>(k/mid X_i)>,/eqno(7.41)$$
$$/eqalignno<
/sg^<(t+1)>&=/Biggl/lbrack N/hs+/sum_<k=1>^Kn_k/twice</lp/mu_k^<(t+1)>-
	/hm_k/rp>'/cr
  &/qquad+/sum_<i=1>^M/sum_<k=1>^KP^<(t)>(k/mid X_i)/twice</lp X_i-
	/mu_k^<(t+1)>/rp>'
	/Biggr/rbrack/cr
  &/qquad/qquad/qquad/sla/lsb N+/sum_<i=1>^M/sum_<k=1>^KP^<(t)>(k/mid X_i)/rsb,
	&(7.42)>$$
%
%
$t=0,1,2,/ldots$, until a convergence criterion is met. Here
$P^<(t)>(k/mid X_i)$ is evaluated according to $t$-th generation
parameter values.
/note<Remark 3.> We have not been specific about criteria for convergence.
One can base one criterion on the amount of change in parameter values
from iteration to iteration, and another on the increase in the 
log-likelihood function. These criteria can also be combined. It
is sometimes useful to look at a ``log'' of what happened through the
iterations in cases where convergence was slow or the final values
were surprising.
/note<Remark 4.> In our development we have stressed likelihood-based
equations and iterative computational schemes to solve them. Let us also 
point out that another option is available, namely that of optimising
the long, high-dimensional log-likelihood function directly, using powerful,
modern optimisation algorithms, and that are available in some 
mathematical software program packages.
/note<Remark 5.> This section has provided methods to obtain updated class
description parameters for a given ``map'' of new symbols,
along with estimates of prior probabilities. The resulting estimators combine 
information from training sets (the supervised part) with
information from the rest of the map (the unsupervised part) in an optimal way.
In applications one may witness the gradual adjustment
of the original class descriptions $/hm_k,/hs_k$ to the 
full-scene estimates $/mu_k^/ast,/sg_k^/ast$.
%
%
This will give the experimenter information on how representative his training
set really was. Better yet in this respect might be to compute maximum 
likelihood estimates $/wt/mu_k,/wt/sg_k$ for the unsupervised part alone, which
is possible by methods similar  to those already presented, and compare them
with the maximum likelihood estimates $/hm_k,/hs_k$
for the supervised part alone.
/section<7.3><Passing from parametric to semi-><parametric and nonparametric models>
The previous sections have developed methods for updating of parametric 
class descriptions, initially based on (perhaps small) training sets,
using (perhaps a large number of) unclassified feature vectors. This
program can also, and perhaps advantageously, be used in steps: instead of
trying to get an ambitious $N_d(/mu_k,/sg_k)$ decription for each class,
after having started with only $5$ (say) vectors from each, one could for
example update only to $N_d(/mu_k,/sg)$ with a common $/sg$ in a first
step, and then update further afterwards.

This section briefly considers methods for taking the updating business 
even further: from parametric to semiparametric and nonparametric models.
We are deliberately 
vague in places, since the success of the methods may depend upon several
circumstances with the application at hand, and since we have as 
yet not gathered any practical experience with the methods.
/subsection<7.3.A><Updating coefficients in an orthogonal expansion density 
estimator> Consider, as an introductory example, a simple univariate 
situation in which class densities on $(0,1)$ are involved, and are
expanded in cosine series
$$/fkx=1+/sum_<j=1>^<m/k>c_k(j)/cos(j/pi x),/ 0/le x/le1,$$
as in 5.3.A. Based on initial data $X_1^/k/upto X_/nk^/k$ an estimate
is obtained with
$$/hc_k(j)=<2/over/nk>/sum_<t=1>^/nk/cos/lp j/pi X_t^/k/rp,/ j=1,2/upto
/eqno(7.43)$$
and an $m(k)$ is settled on as well. Suppose now that new vectors 
$X_1/upto X_M$ from the mixture distribution $/pi_1f_1+/cdots+/pi_Kf_K$
come in. How can they be exploited to give better estimates of $c_k(j)$
than (7.43)?

The previous sections employed the maximum likelihood principle, which however
appears unfeasible now. But the parameter
$$c_k(j)=/int_0^1 2/cos(j/pi x)/fk(x)/,dx=E_f2/cos/lp j/pi X_t^/k/rp$$
is the mean of a specific random variable, and earlier
efforts gave updating formulae for such means. Following (7.29),
for example, we propose
$$c_k(j)=</sum_<t=1>^/nk 2/cos/lp j/pi X_t^/k/rp+/sum_<i=1>^MP/lp k/mid X_i/rp
2/cos/lp j/pi X_i/rp/over/nk+/sum_<i=1>^MP/lp k/mid X_i/rp>/eqno(7.44)$$
as the equations to be solved. This must be done iteratively again; $/pkx$
above is the usual ratio $/pi_k/fkx/sla/sum_<m=1>^K/pi_kf_k(x)$, which 
involves the $c_k(j)$'s.

This is a scheme that can be generalised. Consider the general situation of
Section 5.4.A, where (possibly transformed) data $Y_1^/k/upto Y_/nk^/k$,
if they stem from class $k$, have density
$$/eqalign<g_k(y)&=g_<k,0>(y)/,h_k(y)/cr
&=g_<k,0>(y)/sum_j c_k(j)/,/psi_j(y)./cr>$$
Here the $/psi_j$'s are orthonormal w.r.t.~some $w(y)$, and 
$c_k(j)=/int g_k(y)/lbrace/psi_j(y)/,w(y)/sla g_<0,k>(y)/rbrace/,dy$
is estimated by
$$/hc_k(j)=<1/over/nk>/sum_<t=1>^/nk/psi_j/lp Y_t^/k/rp/,w/lp Y_t^/k/rp
/sla g_<0,k>/lp Y_t^/k/rp./eqno(7.45)$$

If new, unclassified vectors $X_1/upto X_M$ become available, then estimates
of the posterior class probabilities $P(k/mid X_i)$ can be obtained, using
current estimates
$$/hf_k(x)=/wh g_k/lp R_k(x)/rp/invert</partial R_k(x)/sla/partial x>.$$
Here $x/to y=R_k(x)$ is the transformation used for class $k$ (see Section 
5.4 about such initial transformations). One can accordingly define
$$c_k(j)=</nk/hc_k(j)+/sum_<i=1>^M P(k/mid X_i)/,/psi_j/lp Y_<k,i>/rp/,
w/lp Y_<k,i>/rp/sla g_<0,k>/lp Y_<k,i>/rp/over/nk+/sum_<i=1>^M P(k/mid X_i)>
,/eqno(7.46)$$
where $Y_<k,i>=R_k(X_i)$. Iterating all these equations sufficiently will
hopefully lead to solutions, say $c_k(j)^/ast $.

If the initial $/wh g_k(y)$ involved $m(k)$ terms in the expansion, then
also $m(k)$ could be pushed somewhat if the new data $X_1/upto X_M$ suggest
so. One can try out $m(k)^/ast =m(k)+3$, say, and use the equations above,
only starting with $c_k(j)=0$ for $j=m(k)+1$, $m(k)+2$, $m(k)+3$ as initial
values.

The method proposed here is not without difficulties, theoretically
and practically, but ``looks reasonable'' and ought to be tried out.
Some modifications may be necessary in cases where the initial transformations
$x/to R_k(x)=y$ involve estimation variability.
/subsection<7.3.B><Updating from normality to third order corrected normality>
If $50$ vectors from each class were available initially, then the 
$N_d(/mu_k,/sg_k)$ description might happen to be the most practical one,
even if plots etc.~suggest that strict normality cannot be trusted. Can
one pass to a more sophisticated description of $f_1/upto f_K$ if a couple
of thousand unclassified vectors became available?

Let us write down the formulae that the considerations of 7.3.A suggest.

The current description can be written
$$/eqalign<
/hf_k(x)&=N_d/lp/hm_k,/hs_k/rp/lp x/rp/Bigg/lbrace1+/sum_i <1/over6>/,
/hg_<k;/,i>/lp y_<k,i>^3-3y_<k,i>/rp/cr
&/quad+/sum_<i/not=j>/halv/,/hd_<k;/,i,j>/lp y_<k,i>^2-1/rp/,y_<k,j>+
/sum_<i/lt j/lt l>/heps_<k;/,i,j,l>/,y_<k,i>y_<k,j>y_<k,l>/Bigg/rbrace,/cr>$$
but where $/hg_<k,i>$ and the other third order coefficients are zero.
Here $y_k=/hs_k^<-1/sla2>(x-/hm_k)$.

Natural formulae for revising are
$$/gamma_<k;/,i>=</sum_<t=1>^/nk/lp Y_<t,i>^<(k)>/rp^3+/sum_<s=1>^M
P(k/mid X_s)/lp Y_<k;/,s,i>/rp^3/over/nk+/sum_<s=1>^MP(k/mid 
X_s)>,/nils<Skj|nte ikke helt korrekturen din p} denne...>$$
$i=1/upto d$, $k=1/upto K$, and similarly for $/delta_<k;/,i,j>$ and
$/eps_<k;/,i,j,l>$, following 5.4.F. Here $Y_<k;s>=/hs_k^<-1/sla2>
(X_s-/hm_k)$.

It might be best to update to the full, third-order corrected 
normal in several manageable steps. First one could update only 
$/hmk$ and $/hsk$, forgetting $/gamma$'s and $/delta$'s and $/eps$'s.
After having obtained trustworthy full-sample estimates $/mu_k^/ast $,
$/sg_k^/ast $ could one try to get at the third order correction factor.
/subsection<7.3.C><Updating the </it k-NN//> rule>
As our final, specifically described updating procedure, consider the 
(very non-parametric) </it k-NN//> rules of Section 5.5. The most popular one
is based on
$$/hf_<i, k-NN>=<k_i/sla n_i/over W_k(x)>,/eqno(7.47)$$
where $k_i$ is the number of class $i$ vectors present in the smallest
sphere $S(x,w)$ around $x$ that encloses exactly $k$ data points from
the combined training set, and $W_k(x)$ is the volume of this smallest
sphere. In particular, (7.47) is based on a combined, but possibly initial
training set with $N=n_1+/cdots+n_K$ vectors whose classes are known.

Suppose again that new vectors $X_1/upto X_M$ from the mixture distribution
$/sg/pi_k/fk$ become available. How can we modify (7.47), and the </it k-NN//>
rule that follows from it (namely that of having $/wh P(i/mid x)=
</rm const./ ></pi_i/over n_i/sla N>/,k_i/sla k,/,k=/sum_<i=1>^K k_i$)?

A natural proposal is to replace $k_i/sla n_i$ with
$$<k_i+/sum_<s=1>^M P/lp i/mid X_s/rp/,I/lb X_s/in B_k(x)/rb/over
n_i+/sum_<s=1>^M P/lp i/mid X_s/rp>=<k_i+k_i^/ast /over n_i+n_i^/ast >,
/eqno(7.48)$$
where $B_k(x)$ is the smallest sphere containing $k$ points from the 
</it full, extended//> set $Z/cup/lbrace X_1/upto X_M/rbrace$.
$k_i$ is still the number of the original $X_t^<(i)>$ data that are present
in $B_k(x)$; now they do not sum to $k$, however, since some of the 
$k$ points within $B_k(x)$ may be new ones. (One would presumably use a larger
$k$ in the presence of extra data $X_1/upto X_M$ than before.) We
have instead $k=/sum_<i=1>^K k_i+/sum_<i=1>^K k_i^/ast =k^0+k^/ast $, say.

If $/pi_i$'s are unknown and the sampling procedure is such that $/pi_i$ is 
the expected fraction for class $i$ both in the initial set and in
the extra set, then
$$/pi_i^/ast =<n_i+n_i^/ast /over N+M>/eqno(7.49)$$
is a natural estimate (notice that $/sum_<i=1>^K n_i^/ast =M$). The
</it updated  k-NN rule//> then becomes: allocate a (completely
new) $x$ to class $i$ if $k_i+k_i^/ast $ is larger than the 
other $k_j+k_j^/ast $. In fact the estimated posterior class probability is
$$P(i/mid x)^/ast =<k_i+k_i^/ast /over k^0+k^/ast >=<k_i(x)+k_i^/ast (x)/over k>,
/eqno(7.50)$$
and a doubt decision rule can use this.

The rule (7.50) and the estimator (7.49) look like very reasonable 
generalisations of the classical ones. What remains to be specified,
however, is $P(i/mid X_s)$ in
$$k_i^/ast =/sum_<s=1>^M P/lp i/mid X_s/rp/,I/lb X_s/in B_k(x)/rb./eqno(7.51)$$
Unless one is specifically interested in $/pi_i^/ast $, only $P(i/mid X_s)$'s for
new vectors $X_s$ that actually fall in $B_k(x)$, i.e.~the vicinity
of $x$, are needed. It will not be possible, in general, to find exact solutions of
the form (7.48) and (7.49). But several approximations can be put forward.
A simple-minded method is to take some other, rough estimate of $P(i/mid X_s)$,
say based on Gaussian descriptions. One could also stick to
$$P(i/mid X_s)=k_i(X_s)/sla k.$$

Serious study remains before the performance of this ``dynamically revised''
</it k-NN//> rule can be assessed.
/section<7.4><Concluding remarks><>
The classical statistical symbol recognition system needs an initial stage of
``learning the symbols'' before it can be applied. Optimal performance in 
terms of lowest possible error rates requires efficient semi- and
nonparametric estimation of the underlying class densities, which in
turn, with the classical approach, needs large training sets. This
training stage is often a bottleneck, however, usually demanding many
hours of tedious labelling and editing. Thus it can be practically
unfeasible for a map firm to produce 300 different samples of each
symbol, for example.

The admirable aim of this chapter's updating methods is to reduce this
training phase to a minimum, exploiting instead the potential of
unclassified vectors to reveal the separate class structures. They
aim, in other words, at ``unsupervised learning'' after a perhaps very
brief introductory session of ``supervised learning''.

A typical application could involve the following steps: First an
initial training set is found, consisting perhaps of only five of each
symbol class. Then a file with a perhaps large number of unclassified
vectors is established, and class descriptions are updated using
these. Finally the unclassified vectors are classified, using a
discriminant rule based on updated class descriptions.

Successfully obtaining updated and precise estimates for the class
distributions in this partially unsupervised way is a delicate problem
for a number of reasons, however. The success of a particular solution
to it hinges on several circumstances, and should not be taken for
granted.

A naive proposal is to </it classify//> a new, incoming vector
according to the current classification rule, and then treat it as
really coming from that class. This way, trusting the classification
procedure, class parameters can be established using new data, even
sequentially. Procedures like this are called </it decision
directed//>, but they are generally biased and not recommended. Among
remedies that have been proposed in the literature is the </it
probabilistic teacher//>, which randomly assigns class labels to
incoming vectors according to the computed posterior probabilities,
and then update parameter estimates, again on the assumption that the
labels are correctly assigned. This method also has shortcomings,
however, and other methods, including the </it quasi-Bayes//> and the
</it probabilistic editor//> approaches, seem to have greater
potential. See Titterington, Smith and Makov (1985, Ch. 6) for
comments about these in some particular situations and for clues to
generalisations. These latter methods are particularly valuable for
</it rapid//> and </it sequential//> updating. These aspects are often
of little importance for symbol recognition applications, however.

The approach used in the present chapter seems to be the statistically
soundest one: natural models for the class densities entail a
statistical model also for the unclassified vectors, and parameters
can be estimated using well understood principles, like that of
maximum likelihood, utilising all the available data. There are still non-trivial matters to tend to,
however, both practical and theoretical, before the outlined approach
can crystallise into a completely automatic and foolproof machine for
unsupervised learning.

Remarks 3 and 4 of Section 7.2 touched upon matters related to the practical
implementation of the maximum likelihood solution. The EM algorithm
generally works satisfactorily, in spite of its sometimes slow
convergence, provided the starting values for the iterations are
acceptable. Estimates based on the initial labels-known training set
are the obvious choices. But sometimes the model that is updated is
too ambitious for the modest initial set, i.e.~the starting values may
be too crude and could lead to EM convergence to a wrong stationary
point of the log likelihood. A remedy is to start with a simpler
sub-model to update, and then use the updated parameters of this
sub-model as the starting point for updating the more ambitious model.

To illustrate, suppose the task is to arrive at full, multinormal
class descriptions $N_d(/mu_k,/sg_k)$. The number of unknown
parameters is fairly large, $Kd+Kd(d+1)/sla2$. One might now start
with a very much simplified sub-model with only $Kd+1$ parameters,
specifying that $/sg_k=/ssg^2I$, i.e./ all feature components are
independent with the same standard deviation $/ssg$. After having
updated this sub-model, ending in parameter estimates $/mu_k^<(0)>$,
$/sg_k^<(0)>=/wh/ssg^2I$, say, these can be used to start updating
iterations for the more flexible model that says $/sg_k=D=</rm diag>
(/ssg_1^2/upto /ssg_d^2)$ (with $Kd+d$ parameters). In the next stage
one could go to descriptions $N_d(/mu_k,D_k)$ with $Kd+Kd$ parameters,
and the resulting estimates could be initial values for the final step
EM analysis, passing from $N_d(/mu_k,D_k)$ to the full model.

This conjures up a picture of a system that carefully ``surfs'' from
model to model, and at each step takes care to make sure that the
starting values are acceptable and that the model can be supported by
the data. In addition to playing safe at each level surfing also leads
to quicker EM convergence, because the log-likelihood surface is more
well-behaved for simpler models, and because safer initial values are
used. This remark also applies to other approaches to the numerical
optimisation of the log-likelihood function.

The statistical information that each unclassified vector carries
about the underlying parameters is less than that of a label-known
vector. Accordingly it may happen that a model is too ambitious for
updating, i.e.~has too many parameters to succeed in that ever-present
statistical game of balancing modelling bias against estimation
variability. One interesting type of procedure is to smooth over
several models choices with weights coming from Bayes/sla empirical
Bayes considerations and computations. Another approach is to
construct data-based criteria for model choice, which can be done in
several ways. One possibility is to compute goodness-to-fit test
statistics for each model. A simpler method that appears effective in
practice is the following: Assume models $M(1)$, $M(2)/upto M(m)$ are
those considered, with model $M(j)$ having $s(j)$ parameters. Compute
$/log L_</rm max>(j)$, the maximised log-likelihood based on the $j$'th
model. Of course this number tends to be highest for models with many 
parameters, so it won't do to simply pick the model with highest
$L_</rm max>(j)$.
Choose instead the model with highest value of
the penalised criterion
$$/log L_</rm max>(j)-c/ s(j).$$
Akaike's information criterion </it AIC//> is of this type, with
$c=1$. </it AIC//> is known to have a tendency towards overfitting,
however, and we favour a related and more stingy procedure due to
Schwarz (1978) and others, where $c=/halv/log n$, $n$ being the
number of observations. (Schwarz' criterion was derived under
conditions somewhat more stringent than those present in the mixture
problem with both classified and unclassified data, but his 
arguments can be generalised to cover these cases too.)

This criterion can for example be used to find out whether the full
normal description $N_d(/mu_k,/sg_k)$ is too ambitious for the amount
of training data at hand, compared to the simpler $N_d(/mu_k,D_k)$
model with diagonal covariance matrices. Similarly the usefulness of
the full normal model can be compared to that of versions of the third
order corrected normal model considered in 7.3.B, by computing maximal
log-likelihood and subtracting the appropriate complexity term,
counting the number of parameters.

Section 3.1 provided insights into what happens when parametric models
are fitted when they are wrong. A similar analysis can be performed to
see what exactly the maximum likelihood procedure aims at in the
present problem with data of both known-class and unknown-class type.
If $n_k/sla N=/lambda_k$ and $N/sla(N+M)=a$, then maximising the full
likelihood (7.23) essentially means trying to fit the underlying
models $f_k(x)=f_k(x,/th_k)$ by minimising a linear combination of
Kullback-Leibler distances involving the </it true> densities $f_k$
and $f=/sg/pi_k^0 f_k$, viz.
$$a/sum_<k=1>^K/lambda_k I/lb f_k,f_k/lp.,/th_k/rp/rb+
(1-a)/ I/lbrace f,/sum_<k=1>^K/pi_k/,f_k/lp.,/th_k/rp/rbrace.$$
In particular the method is consistent in the idealised case where
$f_k(x)$ really </it is//> $f_k(x,/th_k^0)$ for suitable $/th_k^0$, $k=1/upto
K$.

More valuable is the information the above remark provides for cases
where the model is imperfect. A typically occurring situation is one
where the normal model is used for each class but where one or more of
the classes have two-peaked distributions, corresponding for example
to two principally different ways of writing a specified symbol. One
can find out what the updated parameter estimates really aim at by
minimising the above distance function. The situation is much more
complex than for the  case of fully classified data. In this latter
case the maximum likelihood estimates aim at the true mean vectors and
covariance matrices for the class distributions, even if the model is
incorrect, see 3.1.B. But in the present case the estimates go for
least false/sla most fitting values that are more complicated functions of
the underlying model. In particular, even if the parametric model happens to be
perfect for class 3, but incorrect for class 4, then parameter
estimates may go for the wrong values for class 3 too. Thus departures
from the assumed model may cause more unexpected problems than in the
traditional cases with fully observed data. This also points to the
virtue of avoiding a </it too> small value of $a=N/sla(N+M)$ above.

The particular cases that deal with normal models are likely to become
most popular among the updating methods. A possible danger with these
is occurrences of clearly non-normal behaviour in the classes, like
the bi-modality phenomenon noted above. There are ways of making the
methods more resistant to such departures from the model. One
possibility is to replace the simple normal model for each class with
a mixture of two normals, at least for some of the classes. The
updating methods developed in Section 7.2 are general enough to handle
such a case too, but become somewhat complicated since an EM
iteration loop then must be included within each ``outer'' EM
iteration step. A simpler solution that tends to be robust against
several non-normal facets is to use the multi-dimensional
$t$-distribution with a fixed number of degrees of freedom, say $6$,
instead of the normal; see Berger (1980a, page 395) for a definition.
A rapid Newton-Raphson iteration subroutine, starting at moment
estimators, provides the maximum likelihood estimates that are needed
for the class parameters within each EM step. These are of the robust
$M$-type.

In cases where class distributions contain obvious sub-clusters a
preliminary cluster analysis can be carried out to define proper
sub-classes. Cluster analysis can be useful in conjunction with
updating methods also in other respects.

The experience we have had time to gather so far promises us that the
potential of updating methods, both towards simplifying and reducing
the costly training stage and towards improved classification
performance, is vast. Further empirical and theoretical investigations
are needed before explicit, foolproof guidelines can be given
regarding model choice, the necessary size of initial training set,
whether the semi-parametric or nonparametric updating methods of
Section 7.3 can be expected to outperform parametric ones, etc., and
before we fully can assess the value of updating methods for pattern
recognition problems.


/chapter<8><Handling singular covariance matrices>
This chapter provides a solution to the problem that arises when
one tries to apply the ``best quadratic rule'' of discriminant
analysis, i.e.~the one presented in 3.2.B, in a case where one or more
of the estimated covariance matrices turn out to be singular,
i.e.~null-definite.

There are several </it ad hoc//> ways to solve this problem, by 
artificially making the covariance matrix positive definitive.
One could for example add a small value $/eps$ to the diagonal elements.
This will add $/eps$ to the eigenvalues, and creates positive 
definite-ness. This is a sub-optimal approach, however.
A singular covariance matrix reveals information about the underlying
probabilistic structure, which can be used to enhance discriminatory power.
It is felt that this information should be utilised in an optimal way.
/section<8.1><Exact theory><>
Assume that the covariance matrix $/sg_k$ for class $k$ has zero determinant.
The usual cause for such a phenomenon is 
that one of the $d$ components of a vector from class $k$ is constant.
More generally, $0$ is an eigenvalue, implying the existence of
a non-zero vector $u$ with $/sg_k u=0$, so that the variance 
$u'/sg_k u$ of $u'X$ is zero, i.e.~every vector in a training set
for class $k$ will fall in the $(d-1)$-dimensional linear subspace
$$A=/lb x/vert u'x=u_1/rx_1+/cdots+u_d/rx_d=a/rb,/eqno(8.1)$$
$a$ being the appropriate constant.

Let $P_k$ be the underlying probability measure for class $k$. Then $P_k(A)=1$,
in particular $P_k$ is not absolutely continuous w.r.t.~Lebesgue measure in
$/rr^d$. Assume for the moment that only one of the eigenvalues for
$/sg_k$ is zero, i.e.~$/sg_k$ has rank $d-1$, and assume that $P_k$
is normal, having in mind the proper ressurection of the classical procedure.
Then $P_k$ is absolutely continuous w.r.t.~$/lam_<d-1,A>$, Lebesgue measure
on the $(d-1)$-dimensional $A$.
Its density
$$/fkx=<dP_k/over d/lam_<d-1,A>>(/rx_1/upto /rx_d)$$
can be found as follows: Seek out a $j/leq d$ such that $u_j/not=0$ in
(8.1), entailing
$$/rx_j/equiv /lp/sum_<i:i/not=j>u_i/rx_i-a/rp/sla u_j$$
as an identity for points $x$ in $A$. The covariance matrix $/sg_<k,/ast>$
for $X_/ast=(/rX_1/upto /rX_<j-1>,/rX_<j+1>,/allowbreak/ldots,/rX_d)'$ is non-singular
since no linear relationship ties together the components of $X_/ast$.
Hence $X_/ast/sim N_<d-1>(/mu_<k,/ast>,/sg_<k,/ast>)$, say, implying that
$$/eqalign<f_k(/rx_1/upto /rx_d)&=N_<d-1>(/mu_<k,/ast>,/sg_<k,/ast>)
				  (/rx_1/upto/rx_<j-1>,/rx_<j+1>/upto /rx_d)/cr
				&=(2/pi)^<-(d-1)/sla 2>/vert/sg_<k,/ast>/vert
				  ^<-1/sla 2>/exp/lb-/halv(x_/ast-/mu_<k,/ast>)'/,
				  /sg_<k,/ast>^<-1>(x_/ast-/mu_<k,/ast>)/rb,/cr
>$$
for $x/in A$.

In order to apply the theory of Chapter 1 we need $K$ class densities 
w.r.t.~a common $/sigma$-finite measure. Even if $/sg_k$ happened to be the 
only singular covariance matrix, as above, we may still have outcomes
in the set $A$ for vectors from other classes, i.e.~$P_m(A)/gt0$ for some 
$m/not=k$. Hence the common measure must dominate both $/lam_<d-1,A>$
and the usual Lebesgue measure $/lam_d$ in $/rr^d$.

A suitable generalisation of the considerations above will involve 
lower-dimensional linear subsets
$$A_1/upto A_p/eqno(8.2)$$
with the property that $P_k(A_j)/gt0$ for some $k;/ j=1/upto p$.
We may take them to be disjoint, without loss of generality,
and assume further that these are all the lower-dimensional subsets with said property.
Let $A_j$ have dimensionality $/dim(j)$ and consider the measure 
$$/mu=/lam_<d,B>+/lam_</dim(1),A_1>+/cdots+/lam_</dim(p),A_p>/eqno(8.3)$$
on $/rr^d$, where $/lam_</dim(j),A_j>$ is the appropriate Lebesgue measure on
$A_j,/ j=1/upto p$, and $/lam_<d,B>$ is $d$-dimensional Lebesgue measure
on the set $B=/rr^d-/cup_<j=1>^p A_j$. (Of course $/lam_<d,B>$ is really
equivalent to $/lam_d$ on $/rr^d$, but it will be convenient to use its
restriction to $B$.)

Now consider the class distribution $P_k$, which is dominated by $/mu$ by 
assumptions. Let
$$/eqalignno<
p_<k,j>&=P_k(A_j),/ j=1/upto p,&(8.4)/cr
  1-p_k&=1-/sum_<j=1>^p p_<k,j>=P_k(B).&(8.5)/cr>$$
Let furthermore
$$/eqalignno<
X_<(k)>/mid/lb X_<(k)>/in A_j/rb&/sim g_<k,j>(x),/ j=1/upto p,&(8.6)/cr
  X_<(k)>/mid/lb X_<(k)>/in B/rb&/sim h_k(x),&(8.7)/cr>$$
where $X_<(k)>$ denotes a vector drawn from $P_k$, i.e.
$$/int_C g_<k,j>(x)/,d/lam_</dim(j),A_j>(x)=P_k(C/mid A_j)=P_k(C)/sla p_<k,j>$$
for $C/subset A_j$ and
$$/int_C h_k(x)/,d/lam_<d,B>(x)=P_k(C/mid B)=P_k(C)/sla(1-p_k)$$
for $C/subset B$. (We define $g_<k,j>=0$ when $p_<k,j>=0$ and
$h_k=0$ if $1-p_k=0$.) It follows that
$$/fkx=/cases< p_<k,j>/,g_<k,j>(x)&if $x/in A_j;/ j=1/upto p$,/cr
               (1-p_k)h_k(x)    &if $x/not/in/cup_<j=1>^p A_j$/cr>/eqno(8.8)$$
is the Radon--Nikodym derivative of $P_k$ w.r.t.~$/mu$ in (8.3).

Now the theory of Chapter 1 applies. (1.10) defines the </it optimal//> 
classification procedure with posterior probabilities $P(k/mid x)$ defined
by (1.9) and (8.8), provided these are known.

In order to apply this result in practice, going for the ML plug-in
rules with the parametric model that emerges when the components
$g_<k,j>$ and $h_k$ are assumed Gaussian, we should in principle know the 
lower-dimensional subsets $A_j$ </it a priori//>. The $A_j$'s that have 
$P_k(A_j)=1$ for som $k$ may be automatically detected, as these correspond
to eigenvectors for empirical covariance matrices with eigenvalues zero,
cf. (8.1). And if it is known in advance, or seen by inspection of
histograms, say, that the second component of vectors from class 7 assumes the
value $0$ for a fair position of a training set, then $/lb x/mid/rx_2=0/rb$
should be included as an $A_j$, etc.

Assume that sets $A_1/upto A_p$ have been arrived at, and that a training set
$X_1^<(k)>/upto/allowbreak X_<n_k>^<(k)>$ is available from class $k$.
Estimate $p_<k,j>$ by the simple observed frequency
$$/eqalignno<
/widehat p_<k,j>&=<1/over n_k>/sum_<i=1>^<n_k> I/lb X_i^<(k)>/in A_j/rb,
				/ j=1/upto p,&(8.9)/cr
    /widehat p_k&=/sum_<j=1>^p /widehat p_<k,j>.&(8.10)/cr>$$
(Obviously many of the $/widehat p_<k,j>$'s will be zero, which is good for 
discrimination.) Next use the $n_k/,/widehat p_<k,j>$ $X_i^<(k)>$'s that fell in
$A_j$ to arrive at a (singular) Gaussian description $N_d(/hm_<k,j>,/hs_<k,j>)$
, rank $(/hs_<k,j>)=/dim(j)$, and express $/widehat g_<k,j>(x)=/widehat g_<k,j>
(/rx_1/upto/rx_d)$ as a $/dim(j)$-dimensional non-singular Gaussian density
in $/dim(j)$ of the $d$ components, similarly to the arrangement in our
introductory example. $/widehat h_k(x)$ is just $N_d(/hm_k,/hs_k)(x)$
with $/hm_k,/hs_k$ obtained from the $n_k(1-/widehat p_k)$
observations that fell outside all $A_j$'s.

This defines a parametric estimate $/widehat/fkx$ of (8.8). The collection
$/widehat f_1(x)/upto/widehat f_K(x)$ now defines a classification rule of the type 
described in (3.4) -- (3.5).

In order to be effective in practice such a rule should not include too many 
lower-dimensional sets $A_j$; only the ``major'' ones having $P_k(A_j)$
``significantly positive'' should enter in (8.8). In the applications we 
have in mind only a small number of $A_j$'s are needed; the basic aim is 
to generalise the classical method making it capable of handling singular
covariance matrices in an optimal way.

An example encountered in a symbol recognition experiment may clarify matters.
The symbols were hand-written numerals $0,1/upto 9$. A feature extraction method
was chosen where the first component $/rX_1$ measured
the distance from the upper left hand corner of the smallest rectangular box
enclosing the symbol to the nearest end-node of the symbol. Experience from
training data revealed that $/rX_1$ was always zero for symbol ``7''; that
this happened in 13/
respectively; and that $/rX_1$ was never zero for symbols other than these.

Employing the general framework above we put $A=/lb(/rx_1/upto/rx_d)/mid
/rx_1=0/rb$, and have
$$f_k(/rx_1/upto/rx_d)=/cases<
	p_k/,N_<d-1>(/mu_<k,1>,/sg_<k,1>)(/rx_2/upto/rx_d)&if $x/in A$/cr
	(1-p_k)/,N_d(/mu_k,/sg_k)(/rx_1/upto/rx_d)&if $x/not/in A$/cr>$$
as class densities, where $p_k=P_k(A)$.

Assume that a candidate vector is observed for which $/rX_1=0$. Then the Bayes
rule picks the class having the highest $/pi_k/,f_k(/rX_1/upto/rX_d)$, 
i.e.~among
$$/displaylines<
0,0,/pi_2 0.13N_<d-1>(/mu_<2,1>,/sg_<2,1>)(/rX_2/upto/rX_d),/cr
/pi_3 0.21N_<d-1>(/mu_<3,1>,/sg_<3,1>)(/rX_2/upto/rX_d),0,0,/cr
0,/pi_7 N_<d-1>(/mu_<7,1>,/sg_<7,1>)(/rX_2/upto/rX_d),0,0./cr>$$

If an $X$ is observed with $/rX_1/gt0$, then it will never be assigned to 
class 7. Conversely, if 7 were the </it only//> class with 
$P_k/lb/rX_1=0/rb/gt0$ (even if this probability were lower than one), then
</it every//> $X$ with zero first component would be classified as 7.
/section<8.2><Concluding remarks><>
A sometimes satisfactory approximation to the solution advocated above
consists of pushing zero eigenvalues up to positivity. More precisely, 
if $/hs$ for a class has determinant zero, then $P/hs 
P'=D=/diag(/lam_1/upto/lam_d)$, say, where one or more of the 
$/lam_i$'s are zero. Add $/eps$ to these, producing $D_/eps$, and then 
let $/hs_/eps=P'D_/eps P$. $/eps$ could be some fraction of the 
smallest positive $/lam_i$. In effect a genuine spike in the 
probability distribution is approximated with a smoothed spike.

A disadvantage of this approximation solution, compared to the more 
complex one discussed above, is that interesting lower-dimensional 
hyperplanes, having positive probability for some of the class 
distributions, may go unnoticed.

Sometimes effective two-stage or multi-stage classifiers can be 
constructed where the first stage consists of an initial sorting into 
subgroups of classes using such lower-dimensional hyperplanes, or, more 
generally, discrete feature components. A related comment is that the 
solution outlined in 8.1 resembles those discussed in Chapter 9.


/chapter<9><Discrete feature vectors> 
Our exposition has devoted most attention to continuously distributed
feature vectors. Often valuable information about the objects to be
classified exist in the form of variables that can take on only a low
number of values, however. One example, in the symbol recognition context,
is the feature vector whose components $X_1/upto X_d$ record whether or
not the symbol intersects subsets $1/upto d$ of the smallest box
containing the symbol. If the box is divided into $5$ by $5$ sub-rectangles
in the natural way, for example, then $X=(X_1/upto X_<25>)$ is
a long vector of $0$'s and $1$'s that obviously carries discriminatory
information. Another example could be a component counting the number of
end-nodes of a (skeletonised) handwritten numeral, a variable with $0, 1, 2,
 3, 4$ as the possible outcomes. And of course discrete variables arise
naturally by grouping values of a continuous variable together.

Thus there is a need for theoretical and practical machinery for turning
discrete-valued information into classification procedures. Section 9.1
below considers some general methods for building classifiers when the
feature vectors have binary valued components, a situation of particular
importance. Then comments are offered on the more general polychotomous 
problem,
in Section 9.2. Finally the mixed variables case, where some components of
$X$ are discrete and some continuous, is treated, to a limited extent, in
Section 9.3.

Observe first that important parts of the theory developed in earlier
chapters apply also to discrete-valued vectors, so one needs not start
from scratch. Results in Chapters 1 and 2 about optimal procedures, in
terms of class distributions and </it a priori//> probabilities, are valid still;
in particular
$$C_0(x)=/cases<k&if $P(k/vert x)/ge P(l/vert x)$ for $l/not=k$,/cr
		D&if every $P(k/vert x)/le 1-c$/cr>/eqno(9.1)$$
is the classifier performing best w.r.t.~the loss function (1.3), where
$$P(k/vert x)=/pi_k f_k(x)/sla/sum_<l=1>^K/pi_l f_l(x)$$
as usual. Thus a natural approach is to model the class densities $f_k$ in
effective ways, and then use estimated versions in (9.1) as in Chapters 
3, 4, 5, and 7.
/section<9.1><Binary feature components><>
Suppose $X=(X_1/upto X_d)$ has components $X_i$ with possible values $0$
and $1$, i.e.~$X$ lies in $/Omega=/lb0,1/rb^d$, a sample space with $2^d$ 
elements. To carry out the natural approach alluded to above, the class
densities (or point probabilities)
$$f_k(x)=/Pr/lb X_1=x_1/upto X_d=x_d/ /vert/ X/ </rm is/ from/ class/ >k/rb
	/eqno(9.2)$$
must be estimated from training data (possibly supplemented with
unclassified $X$-vectors, as in Chapter 7). This is, at one level, easier
than estimating continuous class densities, since
$$/tf_k(x)=<1/over/nk>/sum_<j=1>^/nk I/lb/Xkj=x/rb,/ x=/lp x_1/upto x_d/rp 
/in/Omega/eqno(9.3)$$
provide natural estimates. These ``counting estimates'' are in fact
maximum likelihood estimators under the nonparametric model. This is only
fruitful when $/nk/sla2^d$ is of reasonable size, say above $5$. Unless
$d$ is modest, therefore, smoothing methods of some sort quickly become 
necessary.

We remark in passing that the theory of Chapter 4 suggests that a
``predictive alternative'' to (9.3), namely 
$$/lp/sum_<j=1>^/nk I/lb/Xkj=x/rb+1/rp/sla/lp/nk+2/rp,/eqno(9.4)$$
can be expected to work better.
/subsection<9.1.A><The Bahadur-Lazarsfeld expansion>
In many cases of practical importance the number of parameters in 
$f_k(x)$, i.e.~$2^d$, is too large to handle without any smoothing
or modelling. A very simple model is that of </it independence//>, where 
$$f_k(x)=/prod_<i=1>^d/lp1-p_<k,i>/rp^<1-x_i>p_<k,i>^<x_i>./eqno(9.5)$$
This model spesification carries only $d$ unknown parameters, viz.
$$p_<k,i>=/Pr/lb X_i=1/ /vert/ X</rm/ is/ from/ class/ >k/rb./eqno(9.6)$$
These parameters are effectively estimated by
$$/wh p_<k,i>=<1/over/nk>/sum_<j=1>^/nk I/lb X_<j,i>^<(k)>=1/rb,
/ i=1/upto d,/eqno(9.7)$$
or presumably somewhat better by $(/nk/wh p_<k,i>+1)/sla(/nk+2)$.

Of course the independence assumption cannot be expected to hold, but
(9.5) might be taken as a first-stage guess at $f_k(x)$, whereupon
corrections can be made suggested by data. The machinery of orthogonal
expansions, developed in Sections 5.3 and 5.4, provides a framework for
obtaining such corrections.

Let us drop ``$k$'' in the notation in what follows, and concentrate on
estimating an unknown $f(x)=f(x_1/upto x_d)$ on the basis of
$X_1/upto X_n$. The initial, rough description suggested above is
$$f_0(x)=/prod_<i=1>^d/lp1-p_<i0>/rp^<1-x_i>p_<i0>^<x_i>./eqno(9.8)$$
Here $p_<i0>$ could be the estimate $/wh p_<i0>=(/sum_<j=1>^n X_<j,i>+1
)/sla(n+2)$, or only a reasonable guess for $p_i=/Pr/lb X_i=1/rb$.

Now introduce 
$$y_i=<x_i-p_<i0>/over/sqrt<p_<i0>(1-p_<i0>)>>,/ i=1/upto d/eqno(9.9)$$
and the $2^d$ polynominals
$$/psi_<j_1/upto j_d>(x)=y_1^<j_1>/cdots y_d^<j_d>,/ j_1/upto j_d=0,/ 1.
/eqno(9.10)$$
These are orthonormal w.r.t.~$f_0$, i.e., writing $/uj$ for $(j_1/upto j_d)$,
$$/eqalign<
/sum_x/psi_/uj(x)/psi_</uj'>(x)f_0(x)&=E_<f_0>/psi_/uj(X)/psi_</uj'>(X)/cr
&=/prod_<i=1>^d E_<f_0>/lp<X_i-p_<i0>/over/sqrt<p_<i0>(1-p_<i0>)>>
/rp^<j_i+j_i'>/cr
&=I/lb j_1=<j_1>'/upto j_d=<j_d>'/rb./cr>$$
Of course this may be written as 
$/int/psi_</uj>/psi_</uj'>f_0 d/mu=I/lb/uj=/uj'/rb$, where $/mu$ is counting 
measure on $/Omega$, so that the general orthogonal expansion theory
of Section 5.4.A applies.

This theory leads one to consider expansions of $h=f/sla f_0$. Minimisation of
$$/sum_x/lb/sum_/uj c/lp j_1/upto j_d/rp/psi_/uj(x)-h(x)/rb^2 f_0(x)$$
w.r.t.~the coefficients, and an additional argument involving counting 
of parameters, show that
$$f(x)=f_0(x)/sum_<j_1/upto j_d> c(j_1/upto j_d)/ y_1^<j_1>/cdots
y_d^<j_d>,/eqno(9.11)$$
where
$$/eqalign<
c(j_1/upto j_d)&=/sum_x/psi_/uj(x)/ h(x) f_0(x)/cr
&= E_f Y_1^<j_1>/cdots Y_d^<j_d>./cr>$$
Natural estimators for these are
$$/hc(j_1/upto j_d)=<1/over n>/sum_<t=1>^n Y_<t,1>^<j_1>/cdots 
Y_<t,d>^<j_d>,/eqno(9.12)$$
where, of course, $Y_<t,i>=(X_<t,i>-p_<i0>)/sla/lb p_<i0>(1-p_<i0>)
/rb^<1/sla2>$.

There is a natural ordering of the $2^d$ polynominals in (9.10). One 
can arrange them and rename them in terms of their order 
$j_1+/cdots+j_d$:
$$/eqalign<
/psi_0(x)&=1;/cr
/psi_<1;i>(x)&=y_i;/cr
/psi_<2;i,j>(x)&=y_iy_j,/ i/lt j;/cr
/psi_<3;i,j,k>(x)&=y_iy_jy_k,/ i/lt j/lt k;/cr>$$
etc. The coefficients can be ordered accordingly:
$$/matrix<
/hfill       c_0&=&1,/hfill   &  /hfill /wh c_0&=&1;/cr
/hfill   c_<1;i>&=&E_fY_i=</textstyle p_i-p_<i0>/over/sqrt</textstyle p_<i0>(1-p_<i0>)>>,/hfill
			      &  /hfill /hc_<1;i>&=&</textstyle/wh p_i-p_<i0>/over
						/sqrt</textstyle p_<i0>(1-p_<i0>)>>;/cr
/hfill   c_<2;i,j>&=&E_fY_iY_j,/hfill
			      &  /hfill /hc_<2;i,j>&=&<1/over n>/sum_<t=1>^n
Y_<t,i>Y_<t,j>;/cr>/eqno(9.13)$$
etc. Here $/wh p_i=<1/over n>/sum_<t=1>^n X_<t,i>$ estimates $p_i$.

The $m$'th order expansion estimator of $f(x)$ is
$$/eqalignno<
/wh f_m(x)&=f_0(x)/Bigg/lbrace1+/sum_i/hc_<1;i>y_i+/sum_<i/lt j>/hc_<2;i,j>
y_iy_j/cr
&/qquad+/cdots+/sum_<i_1/lt/cdots/lt i_m>/hc_<m;/ i_1/upto i_m>
y_<i_1>/cdots y_<i_m>/Bigg/rbrace.&(9.14)/cr>$$ 
In a practical application values for $p_<i0>$ in (9.8) must be
chosen, as well as the order $m$.

The special case $p_<10>=/cdots=p_<d0>=<1/over2>$ leads to products of
$(2x_i-1)$ terms, and (9.11) then becomes what is known, for example in switching
theory, as the </it Rademacher-Walsh expansion//>. A disadvantage of
this choice is that a high order $m$ might be needed in (9.14) to capture 
the structure of the true $f(x)$.

A more popular choice is the data-based $p_<i0>=/wh p_i$ or
$/wh p_<i0>=/lp n/wh p_i+1/rp/sla/lp n+2/rp$, for which 
$f_0(x)=/wh f_0(x)$ becomes the natural estimator under independence,
and for which the first-order terms vanish. $/wh f_m$ in this form is
known as the </it Bahadur-Lazarsfeld expansion estimator//>.
Often the simple 
$$/wh f_2(x)=/wh f_0(x)/lb1+/sum_<i/lt j>/hc_<2;i,j>
<x_i-/wh p_i/over/sqrt</wh p_i/lp1-/wh p_i/rp>>
<x_j-/wh p_j/over/sqrt</wh p_j/lp1-/wh p_j/rp>>/rb,/eqno(9.15)$$
with $d+(<>^d_2)$ estimated parameters, or perhaps $/wh f_3$
with $d+(<>^d_2)+(<>^d_3)$ such, perform satisfactorily.
$/wh f_2$ corrects the independence model $/wh f_0$ by taking 
all pair-wise correlations into account.
/note<Remark 1.> There is a problem with (9.15) and its higher-order 
relatives, however, namely that the sampling variability for 
$/wh p_i$ is ignored in its construction. One gets for example
$c_<1;i>=E_f(X_i-/wh p_i)/sla /sqrt</wh p_i(1-/wh p_i)>$ as the optimal
first-order coefficients, and these are not equal to zero. Similarly
the stopping rule devised below rely on the numbers $p_<i0>$ to be given
independently of the data $X_1/upto X_n$, so that it does not, strictly
speaking, apply to the Bahadur-Lazarsfeld expansion given as above. This
is the problem studied and partially solved in Section 5.4.G, and a 
similar analysis can be carried out here, leading to a better
stopping rule for $m$. This analysis is not supplied here, however.
/subsection<9.1.B><Inclusion rules>
What order $m$ should be used in (9.14)? One rule for stopping$/sla$inclusion
can easily be devised based on the general theory of 5.4.A. Consider
the more general estimator
$$/wh f_J(x)=f_0(x)/sum_<(j_1/upto j_d)/in J>/hc
(j_1/upto j_d)/ y_1^<j_1>/cdots y_d^<j_d>,/eqno(9.16)$$
where $J$ is a subset of $/Omega=/lb0,1/rb^d$. Its integrated (summed)
squared error is
$$/eqalign</ise(J)&=/sum_x/lb/wh f_J(x)-f(x)/rb^2/sla f_0(x)/cr
&=/sum_x/lsb/ /sum_</uj/in J>/lb/hc/lp/uj/rp-c/lp/uj/rp/rb/psi_/uj(x)
-/sum_</uj/not/in J>c/lp/uj/rp/psi_/uj(x)/ /rsb^2/ f_0(x)/cr
&=/sum_</uj/in J>/lb/hc/lp/uj/rp-c/lp/uj/rp/rb^2+/sum_</uj/not/in J>
c(j)^2,/cr>$$
the expectation of which is
$$/eqalign<
/mise(J)&=/sum_</uj/in J><1/over n>/tau/lp/uj/rp^2+/sum_</uj/not/in J> c(j)^2/cr
&=/sum_</uj/in/Omega> c/lp/uj/rp^2-A(J),/cr>$$
i.e.~a large value of
$$A(J)=/sum_<j/in J>/lb c/lp/uj/rp^2-<1/over n>/tau(j)^2/rb/eqno(9.17)$$
is favourable. Here
$$/eqalignno<
/tau/lp/uj/rp^2&=</rm Var>_f/lp Y_1^<j_1>/cdots Y_d^<j_d>/rp/cr
&=/sum_x/lb y_1^<j_1>/cdots y_d^<j_d>-c/lp/uj/rp/rb^2/ f(x),&(9.18)/cr>$$
an unbiased estimator for which is
$$/wh/tau/lp/uj/rp^2=<1/over n-1>/sum_<t=1>^n/lb
Y_<t,1>^<j_1>/cdots Y_<t,d>^<j_d>-/wh c/lp/uj/rp/rb^2./eqno(9.19)$$
It follows that
$$/hA(J)=/sum_</uj/in J>/lb/hc/lp/uj/rp^2-<2/over n>/wh/tau/lp/uj/rp^2/rb
/eqno(9.20)$$
is unbiased for $A(J)$.

A natural procedure is to choose the set $J$, among those under 
consideration, that gives maximal $/hA(J)$. $/hf_m$ of (9.14)
corresponds to $J=J_m$ consisting of all $/uj$  whose order 
$j_1+/cdots+j_d$ is at most $m$. Thus $/hf_3$ is preferred to $/hf_2$,
for example, if $/hA(J_3)/gt/hA(J_2)$, i.e.~if
$$/sum_<i/lt j/lt k>/lb/lp/hc_<3;/ i,j,k>/rp^2-<2/over n>/lp 
/wh/tau_<3;/ i,j,k>/rp^2/rb/gt0.$$

Other inclusion rules can easily be invented for specific purposes,
see the discussion in Section 5.4. The classification task should not be 
forgotten: our proclaimed goal is to construct classifiers with optimal 
performance. Thus one could vary expansion lengths $m_1/upto m_K$ for the 
$K$ class distributions and estimate error rate, or minimal class separation,
for each combination, and finally select the best one.
/subsection<9.1.C><Chow expansions> The simplest, non-trivial
Bahadur-Lazarsfeld expansion estimator is $/hf_2$ of (9.15), involving
estimates of $d(d+1)/sla2$ parameters. Even this might be too heavy a 
burden on the data; compare the example mentioned in the introduction
where a not unreasonably chosen feature vector had $d=25$ components.
Variations could be considered where a fair portion of the second-order
interaction coefficients are set to zero on a priori grounds. One might
also choose to replace $/hc(j_1/upto j_d)$ in (9.16) with zero
if the hypotesis $c(j_1/upto j_d)=0$ is not rejected by 
$/sqrt<n>/hc(j_1/upto j_d)/sla /wh/tau(j_1/upto j_d)$ 
(this statistic is close to standard normality under said
hypothesis).

The approach of 9.1.A treated the $d$ components equally, perhaps 
even </it too//> fairly. Sometimes a natural ordering of the components
exists, and such knowledge can guide one in building effective models 
with few parameters. An example where the components of $X$ proceed in 
an orderly fashion is a Markov chain, a stochastic process in
which the distribution of $X_i$ given $(X_1/upto X_<i-1>)$
only depends upon $X_<i-1>$. In such a situation one would have
$$f/lp x_1/upto x_d/rp=/Pr/lb X_1=x_1/rb/prod_<i=2>^d/Pr/lb X_i=
x_i/ /vert/ X_<i-1>=x_<i-1>/rb,$$
requiring spesification of only $2d-1$ parameters.

Of course the Markov chain situation is somewhat simplistic, but 
it illustrates a principle for building parsimonious expansions. As an 
example, suppose it is possible a priori to find an ordering of
$X$ in which the distribution of $X_i$ given $(X_1/upto X_<i-1>)$
only depends upon at most two given ones among $X_1/upto X_<i-1>$, say
$X_<a(i)>$ and $X_<b(i)>$. Then
$$/eqalign<
f(x)&=/Pr/lb X_1=x_1/rb/ /Pr/lb X_2=x_2/ /vert/ X_1=x_1/rb/cr
&/qquad/prod_<i=3>^d/Pr/lb X_i=x_i/ /vert/ X_<a(i)>=x_<a(i)>, 
X_<b(i)>=x_<b(i)>/rb,/cr>/eqno(9.21)$$
a distribution determined by $4d-5$ parameters, which is dramatically 
less than Bahadur-Lazarfeld's $/hf_2$ for long vectors $(X_1/upto X_d)$.

Expansions of the above type are called </it Chow expansions//>.
In practice the ``dependence tree'' must be specified, for example,
in (9.21) the ``variables of importance'' to $X_i$, namely $X_<a(i)>$
and $X_<b(i)>$ must be pin-pointed. Methods aiming at finding such dependence 
structures are offered in Chow and Liu (1966).

As an example of the Chow expansion approach, consider a ``binary image''
consisting of a matrix $/lb X_<(i,j)>/rb$ of $0$'s and $1$'s
of a $N/times N$ grid (say). Then one could postulate that the distribution 
of $X_<(i,j)>$, given the rest of the picture, only depends upon its eight
immediate neighbours. This actually characterises a two-dimensional Markov
process, and one is led to so-called auto-binomial or auto-Bernoulli
models; see Besag (1974, 1986).

It is clear how models like (9.21) and the binary image model above
give one the possibility of estimating class distributions $f_1/upto f_K$
for even long feature vectors of binary components. This provides a basis
for construction of classification procedures, and a real challenge to the statistician.
/section<9.2><General discrete feature vectors><> Suppose $X$ has
components $X_1/upto X_d$ with $X_i$ taking on values in $/lb0,1/upto M_i-1/rb$.
The task of modelling
$$/fk/lp x_1/upto x_d/rp =/Pr/lb X_1=x_1/upto X_d=x_d/ /vert</rm/ class/ > k
/rb/eqno(9.22)$$
for $x$ in $/Omega=/prod_<i=1>^d/lb0,1/upto M_i-1/rb$ is even more
difficult than the special case $M_i/equiv 2$ considered in 9.1.
Several methods can be written down, but it is presumably the case that 
no single method can be expected to work well in every application. No
complete review will be attempted here; discussion is limited to
a short description of a couple of approaches and some references 
to literature. 

The simplest model assumption is again that of independence, for which
$$/fk(x)=f_<k,ind>(x)=/prod_<i=1>^d/Pr/lb X_i=x_i/ /vert</rm/ class/ >k/rb,
/eqno(9.23)$$
and each product term can be estimated easily and usually with 
sufficient precision. A natural idea, along the lines of the
previous section, is to expand $f_<k,ind>(x)$ using correction terms
corresponding to interactions between the components. This is
indeed possible. Bahadur (1961b) provides a natural generalisation
of the Bahadur-Lazarsfeld expansion (9.14).
Butler and Kronmal (1985) provide a general framework using discrete 
Fourier functions, for which fast computational schemes are available
(see Karpovsky, 1978).

The discussion of Chow expansions (Section 9.1.C) extends to the 
present situation. Again the crux is to construct a ``dependence tree'',
summarising the dependence structure among $X_1/upto X_d$. If these are 
related in some way to spatial location, for example, then close neighbours
could be modelled to be dependent and distant ones not.

There are ``kernel methods'' of smoothing sparse multinomial tables, and
that in several ways are similar to the continuous case kernel methods 
reviewed in Section 5.2. A standard reference is Aitchison and Aitken (1976).
Their method is somewhat improved in Brown and Rundell (1985), and 
Titterington (1985) places kernel methods for categorical variables in 
a larger context of general smoothing techniques.

Nearest neighbour methods can also be constructed, again parallelling 
the case of continuous variables. A problem here is deciding on a suitable 
distance function; see Hand (1981, Section 5.2) for further comments.

We mention finally the logistic regression type approach that avoids
specifying concrete models for the class distributions (it does so
only implicitly). This approach models the posterior probabilities
$P(k/vert x)=/pi_k/fk(x)/sla/sum_<l=1>^K/pi_l f_l(x)$ directly, for example
as
$$P(k/vert x)=/exp/lp/alpha_k+/beta_k'x/rp/sla/sum_<l=1>^K/exp
/lp/alpha_l+/beta_l'x/rp./eqno(9.24)$$
The parameters, which are unique when a constraint like $/alpha_K=0$,
$/beta_K=0$ is imposed on them, can be estimated from the data, see
e.g.~Haggstrom (1985). Also quadratic generalisations of (9.24)
can be explored, in which $/alpha_k+/beta_k'x$ is replaced by
$/alpha_k+/beta_k'x+x'/Gamma_kx,/ k=1/upto K$.
/section<9.3><Mixed discrete and continuous variables><>
It is only natural to exploit discrete and continuous measurements in
tandem, when both carry information related to class membership.
Experience has shown that just adding one or two yes$/sla$no
variables to a modestly successful continuous feature vector can lead to 
much better classification results---using appropriate methods for
mixed discrete and continuous variables.

Write $X=(A_1/upto A_d,/ Y_1/upto Y_c)'$ for a vector with
$d$ discrete components followed by $c$ continuous ones.
In principle one may treat $(A_1/upto A_d)'$ as a single variable $A$ 
with possible values $a=(a_1/upto a_d)$ in some finite space, say $/Omega$,
and we will partly do that in what follows, for ease of notation.

It is convenient to represent the class density $/fkx$ as follows:
$$/fkx=/fk(a,y)=g_k(a)h_k(y/vert a)./eqno(9.25)$$
We have studied several models and estimating procedures for $g_k(a)$
already in this chapter, and they are valid here too: when data vectors
$$/xkj=/lp A_j^<(k)'>,Y_j^<(k)'>/rp',/ j=1/upto n_k$$
have been observed, use $A_1^<(k)>/upto A_<n_k>^<(k)>$
to obtain $/wh g_k(a)$, say. It remains to model the distribution of
$Y$ given $A=a$, for each outcome $a$.

The simplest model puts down a normal distribution for $h_k(y/vert a)$. This model 
has already been considered in Sections 3.2.D, 4.2.E, and 10.2.C 
(see also (10.19)).

Formally, assume that
$$Y/mid(A=a)/sim N_c/lp/mu_a^<(k)>,/sg_a^<(k)>/rp,/ a/in/Omega/eqno(9.26)$$
when $X=(A, Y')'$ comes from class $k$. If enough $X$-vectors have been
observed with $A_j^<(k)>=a$, then $/mu_a^<(k)>$ and $/sg_a^<(k)>$ can be 
estimated in the usual way from the $Y_j^<(k)>$-vectors, see Sections
3.2.D and 4.2.E for details. Very often it becomes necessary to pool
the $/sg_a^<(k)>$ matrices, however, to reduce the number of parameters to 
be estimated. One should first pool across $a$-outcomes, getting
$/sg_a^<(k)>=/sg^<(k)>$ in the model above, and if necessary one could
model the $K$ $/sg^<(k)>$-matrices using fewer than the full set of
$K<d(d+1)/sla2>$ parameters, for example by putting 
$/sg^<(1)>=/cdots=/sg^<(K)>=/sg$. Details about maximum likelihood 
estimation, for example, are as in 3.2.D.

Usually the conditional mean vectors $/mu_a^<(k)>$ are more important than
$/sg_a^<(k)>$ for discrimination. Furthermore, $/mu_a^<(k)>$ usually
needs only few $Y_j^<(k)>$ vectors, whose accompanying $A_j^<(k)>$'s
were $a$, to be reasonably well estimated. But it may also happen that too 
few, or even none, were available, in cases where $A=(A_1/upto A_d)'$
can take on a large number of different values. In such cases even 
$/mu_a^<(k)>$ must be modelled or smoothed.

One possibility is to view $/mu_a^<(k)>$ as changing smoothly with 
$a/in/Omega$, perhaps like $/gamma_0+/sum_<i=1>^d/gamma_ia_i$. Parameter
estimates can be computed for $/gamma_0/upto/gamma_d$.

Of course there are other reasonable models around than (9.26), but 
this ``discrete times normal'' is the most important one, and has 
produced good results in several applications.
/section<9.4><Concluding remarks><>
Because of the persistent success of Fisher's linear discriminant
rule and its quadratic generalisation in cases where discrimination is
``easy'', these methods have been quite popular also with categorical
data. Much can be won by using methods particularly suited to discrete
variables, however.

Early papers on the subject of classification and probability 
distribution estimation for data with binary valued components are Bahadur
(1961a, b) and Chow (1962), both using the orthogonal expansion
approach. There are naturally several specific orthogonal expansion
methods, and only some of them have been mentioned here, along with
our own variations. A method of Martin and Bradley (1972) has
performed well in some published simulation studies, but is not
expected to be as good as the Bahadur-Lazarsfeld method of 9.1 when
several classes are involved. (The simulation studies mostly use $K=2$
classes.) Kernel type estimators for categorical variables were first
studied by Hills (1967) and Aitchison and Aitken (1976). A good review
of later advances using the kernel approach is available in Hand
(1982, Ch.~4), see also Titterington (1985) for smoothing procedures
in a larger context.

Machine preprocessing of a symbol can lead to its representation as a
large matrix of $1$'s and $0$'s, where a ``$1$'' corresponds to the
event that the symbol intersected the subrectangle in question. One
particular approach to statistical symbol recognition is to model this
large collection of variables directly, rather than deriving feature
components from it. One might try out Markov random field models, as
briefly mentioned in 9.2. This represents a real challenge:
parsimonious models must be chosen, and parameter estimation
procedures must be devised. There has been an intense, recent interest
in Markov random field models for image type data, see for example
Geman and Geman (1984) and Besag (1986), but the methods have been
more oriented towards optimal analysis of a single picture at a time
than to the quicker methods one needs in symbol recognition
applications that should be able to handle several pictures a second.
Ripley (1986, Section 4) reviews and extends the so-called Serra
calculus that can be used to quickly summarise binary images, and such
ideas might lead to classification methods too. (If the gridding of
the image is very fine, and the symbols are represented by their
skeletons as is often done, then one might as well model the full
continuous phenomenon as a random curve in the plane. Helgeland and
Hjort (1986) report preliminary success using Fourier description
techniques to model such curves.)

Another approach to classification that is well suited to discrete
data is that of Breiman, Friedman, Olshen, and Stone (1984). Their
CART methods (classification and regression trees) involve automatic
construction of (large) decision trees.

Most nonparametric density estimaton methods made up with continuous
data in mind have some analogue for categorical data. The $k-NN$
method mentioned in 9.2 is one such example. With a suitable choice of
distance function that method also works for the mixed variable case.

There are links to other chapters in this report. The important
updating procedures of Chapter 7 can be adapted to the present
discrete state of affairs, in particular Section 7.3 provides clues
for this. It becomes necessary, also for the purpose of updating, to
guard against outliers, which means that separate outlier detection
programs must be constructed for the type of data considered here.
Such a method has been derived for the discrete-times-normal model of
(9.25), (9.26), but is not described here.

There are non-trivial problems associated with estimating the
parameters of the various models considered. Only the simpler ones,
for illustration, were considered in 3.2.E and 4.2.C.

For exploration of feature choices and reduction of dimensionality
problems one needs separability measures for categorical and mixed
data. Some such are being put forward in Section 10.2. Chapter 10 also
gives ideas for what should constitute a ``standard output'' for discrete
and mixed-type data.

Error rates for classifiers should be routinely computed. The general
principles presented in Chapter 12 apply here, and in particular the
shortcuts for computing leave-one-out error rate estimates described
in 12.1.B and 12.1.C can be generalised to similar shortcuts for some
of the discrete-times-normal models.

Finally we mention that goodness-of-fit methods, both graphical and
formal, as discussed in Chapter 11, can be devised also for the models
considered here.


/chapter<10><Class descriptions and separability measures>
When a set of features has been chosen, natural and important questions are: 
how can the resulting class densities be efficiently described, and how 
well are the classes separated? The present and the next chapter aim at 
providing diagnostic tools that answer these questions.

The serious experimenter will for a given pattern recognition problem 
try out several feature extraction methods, and should ideally be able to 
tell, even if based on a modestly sized training set, what kind of 
parametric or nonparametric density estimation method that adequately and 
efficiently will capture the data structure; whether an extended training 
phase is necessary, and if so, of what size; and finally something about 
the error rates that will result on a future large set of new objects, 
i.e.~about the success of the chosen set of features.

Note that while estimation of class densities and estimation of 
separability between classes a priori concern non-related matters and 
therefore seem to invite separate and independent study, the final aim of 
constructing a classifier with good performance makes them related after 
all. If preliminary study shows that classes ``3'' and ``7'' are far from 
each other, and also far from all other classes, then a rough and 
non-sophisticated description of $f_3$ and $f_7$ suffices, and one may 
concentrate on estimating the others. In particular it can happen that a 
secondary stage of training, in the form of automatic unsupervised 
updating by methods described in Chapter 7, or of a factual second set of 
classified objects, is only necessary for some of the classes.

It is thus hoped that proper and experienced use of the methods presented 
in the following can provide the basis for the construction of the best 
classifier and at the same time keep the amount of costly training to 
a minimum.
/vfill/eject
/section<10.1><Basic descriptive measures for classes><>
Some statistics are proposed below for ``basic description'' of the 
distribution of feature vectors from a class. Let us write $/XtXn$ 
instead of $X_1^<(k)>/upto X_<n_k>^<(k)>$ for these vectors. An extreme 
point of view is that nothing but the collection of these $n$ 
$d$-dimensional vectors describes the distribution $f$, another is that 
the density estimate $/hf$ that is ultimately used, be it parametric or 
nonparametric, provides the best description. The first description is not 
advantageous because of its size and its lack of guidelines as to what to do
next. Here we look for sound summarising statistics that can guide us in
the choice of method to produce the final $/hf$, and that more or less
routinely can be computed and inspected for each choice of features.
/subsection<10.1.A><Robust estimation of mean and covariance matrix>
The mean $/mu=EX$ and the covariance matrix $/sg=E(X-/mu)(X-/mu)'$ are
essential. The usual estimates are 
$$/eqalign<
/hm&=/ndel/sum_<t=1>^n X_t,/cr
/hs&=<1/over n-1>/sum_<t=1>^n (X_t-/hm)(X_t-/hm)'.>/eqno(10.1)$$

It can sometimes be advantageous to use versions that are robust against the
possible dominance of a few extreme $X$'s, depending upon the further use.
Suppose for example that the distribution of $X$'s from the class in
question is very nearly normal with parameters $/mu_0$ and $/sg_0$ except for
some low percentage of persistent ``wilder ones'' that come from a more
spread-out distribution. The perhaps optimal way to go then would be
to model the underlying $f$ as an appropriate mixture, estimate
parameters, and use the resulting $/hf$ in the classification program.
But another and usually simpler solution is also satisfactory:
estimate $/mu_0$ and $/sg_0$ using robust versions of (10.1), say by trimming
or truncating extreme vectors.  Then use the $N_d(/hm_0,/hs_0)$ for $/hf$,
but in tandem with some outlier criterion, cf.~Chapter 6, and be prepared to 
accept some low percentage of future objects being classified as ``out'' .

It is useful to include such robust estimates $/hm_0, /hs_0$  in the 
``standard output''. If $/hm$ 
and $/hs$ differ much from $/hm_0$ and $/hs_0$, then ``wild values'' are
present, and this fact about the underlying distribution should be taken
into account in one way or another when the classifier is to be
constructed. We propose the following two-stage trimming procedure. Compute
first </it trimmed means//> for each component, avoiding the 3/
ones to the right and the left, and similarly calculate
</it upgraded, trimmed standard deviations//> (this is explained below).
Thus we possess first-stage estimates $/tm_<0,i>,/tss_<0,i>$, say,
and can go on to compute distances
$$/eqalign<
R_t^2&=(X_t-/tm_0)'/,/ts_0^<-1>(X_t-/tm_0)/cr
     &=/sum_<i=1>^d(X_<t,i>-/tm_<0,i>)^2/sla /tss_<0,i>^2/cr>$$
for the whole data set, $t=1/upto n$. Finally discard those of the
$n$ vectors that have disproportionally large $R_t^2$, and compute mean
$/hm_0$ and covariance matrix $/hs_0$ in the usual way for the remaining ones.

It is usually sufficient to determine by a simple ``inspection by
eye'' whether some $R_t^2$ is disproportionately large or not. The
numbers $R_t^2/sla d$ have expectation approximately equal to $1$ and
variance approximately equal to $2/sla d$, unless they are among the
extremists that we look for. A histogram of the $R_t^2/sla d$ values
will usually quickly reveal any outlandish vectors that are present. A
more formal rule, motivated by the rough approximation
$$/Pr/lb/max_<t/le n>R_t^2/le/lam/rb/doteq/,/Gamma_d(/lam)^n,$$
is this: discard $X_t$ whenever
$R_t^2/gt/lam_0=/Gamma_d^<-1>(.95^<1/sla n>)$.

$/lbrack$Let us define properly the first-stage estimators $/tm_0$, $/ts_0=
</rm diag>(/tm_<0,1>^2/upto/tm_<0,d>^2)$ used above. They are based on
coordinate-wise trimming of 3/
$a=</rm/ closest/ integer/ to/ >.03n$, and let
$$/tm_<0,i>=<1/over n-2a>/sum_<t=1>^n X_<t,i>/,I/lb X_<(a+1),i>/le
X_<t,i>/lt X_<(n-a),i>/rb$$
be the corresponding trimmed averages. Here $X_<(1),i>/le/cdots/le
X_<(n),i>$ are the ordered observations. Next consider
$$/wt/tau_<0,i>^2=<1/over n-2a-1>/sum_<i=1>^n/lp
X_<t,i>-/tm_<0,i>/rp^2 I/lb X_<(a+1),i>/le X_<t,i>/le
X_<(n-a),i>/rb,$$
i.e./ $/wh/tau_<0,i>$ is the standard deviation for the trimmed data
set. But $/wt/tau_<0,i>$ obviously underestimates the true standard
deviation for $X_<t,i>$, and has to be upgraded. To find out by how
much, consider the limiting values of $/tm_<0,i>$ and $/wt/tau_<0,i>$:
$$/eqalign<
/tm_<0,i>&/totop</rm a.s.>/mu_<0,i>=<1/over .94>
/int_<F_i^<-1>(.03)>^<F_i^<-1>(.97)> x/,f_i(x)/,dx,/cr
/wt/tau_<0,i>^2&/totop</rm a.s.>/tau_<0,i>^2=<1/over .94>
/int_<F_i^<-1>(.03)>^<F_i^<-1>(.97)> (x-/mu_<0,i>)^2/,
f_i(x)/,dx./cr>$$ 
If $f_i(x)$, the density for the $i$'th component, happens to be
normal $(/mu_i,/ssg_i^2)$, then $F_i^<-1>(.97)=/mu_i+c/ssg_i$ and
$F_i^<-1>(.3) =/mu_i-c/ssg_i$, where $c=/Phi^<-1>(.97)=1.8808$. It
follows that
$$/eqalign<
/tm_<0,i>&/totop</rm a.s.>/ <1/over
.94>/int_<-c>^c(/mu_i+/ssg_iy)/phi(y)/,dy=/mu_i/cr 
/wt/tau_<0,i>^2&/totop</rm a.s.>/ <1/over .94>/int_<-c>^c
y^2/,/phi(y)/,dy=.7277/,/ssg_i^2./cr>$$
Accordingly we use $/tss_<0,i>^2=1.3742/,/wt/tau_<0,i>^2$, or 
$$/tss_<0,i>=1.17225/,/wt/tau_<0,i>$$
as our upgraded, trimmed standard deviation estimate.$/rbrack$

More sophisticated methods for obtaining robust estimates $/tm_0$,
$/ts_0$ than those discussed above can be found in Huber (1981), for
example. The simpler procedure proposed here has however proved
effective in practice.

There are obvious dangers involved in performing a classical
discriminant analysis based on $N_d(/mu_k,/sg_k)$ descriptions with
standard estimators (10.1), without checking for outliers. A couple of
extreme vectors in the training set may inflate the $/sg$-estimate
unduly, and make later classification suboptimal. Also, other parts of
the ``standard output'', e.g.~the generalised standard deviation and
class separability measures, to be discussed below, would suffer from
such an inflated $/sg$-estimate. On these grounds we advocate robust
estimators in normal theory based discriminant analysis.
/subsection<10.1.B><Further ``standard output''>
Of course the one-dimensional </it histograms//> for the data set provide
valuable immediate information too, also about possible ``wild'' vectors.
These should also be included in the standard output.
Similarly scatterplots, for two components at a time, are useful.

A ``small'' $/sg$ matrix is good for classification, since it means less
interference or confusion with other classes. It is convenient to study 
its ``smallness'' via its $d$ eigenvalues. Thus compute
$$/hl_1/geq/ldots/geq/hl_d,$$
the eigenvalues of $/hs$. For inspection purposes we recommend displaying
the roots
$$/hl_1^<1/sla 2>/geq/cdots/geq/hl_d^<1/sla 2>/eqno(10.2)$$
instead, these being standard deviations for the $d$ principal components,
i.e.~the components of
$$Z_t=/hP(X_t-/hm),/ t=1/upto n,/eqno(10.3)$$
where
$$/hP/hs/hP'=</rm diag>/lb/hl_1/upto/hl_d/rb/eqno(10.4)$$
is the spectral decomposition. In particular, the root $/hl_i^<1/sla 2>$
is on the same scale as the $i$'th component of $X$.

To get a single measure for the spread of the distribution, consider the
</it concentration ellipsoid//>
$$/lb x:(x-/mu)'/sg^<-1>(x-/mu)/leq c/rb,$$
cf.~the nonparametric outlier test mentioned in 6.1.C. (The ellipsoid
above contains e.g.~90/
$c=/Gamma_d^<-1>(0.90)$, if $f$ is normal.) The volume of the ellipsoid is
proportional to $/vsgv^<1/sla 2>$, and provides such a measure of
spread. Again it is easier to judge such a value when it is normalised to
the original scale of $X$-component values, so we settle on
$/vsgv^<1/sla 2d>$, which can be termed the </it generalised
standard deviation//> for the $d$-dimensional distribution of $X$.

An estimate of $/vsgv^<1/sla 2d>$ is obtained by ``plugging in'',
i.e.
$$/vhsv^<1/sla 2d>=/lp/hl_1/cdots/hl_d/rp^<1/sla 2d>.$$
This estimator has however a tendency towards underestimating the real
$/vsgv^<1/sla 2d>$, for example, one can show that 
$$E/vhsv=<m-1/over m>/cdots <m-(d-1)/over m>/vsgv,$$
where $m=n-1$ is the ``degrees of freedom'' for
$$A=(n-1)/hs=/sum_<t=1>^n(X_t-/hm)(X_t-/hm)'./eqno(10.5)$$
One can use the unbiased estimator
$(/prod_<i=0>^<d-1><m/over m-i>)/vhsv$
for $/vsgv$ instead, and then take the $1/sla 2d$ root.

We shall however devise a slightly different ``adjusting for finite sample
size'' trick. We encounter from time to time different powers of the
individual standard deviation $(/ssg_<ii>)^<1/sla 2>$ and the generalised
standard deviation $/vsgv^<1/sla 2d>$, e.g.~$/ssg_<ii>,
(/ssg_<ii>)^2,1/sla(/ssg_<ii>)^<1/sla2>,/vsgv,
/vsgv^<1/sla 2>,/vsgv^<1/sla d>,$ and later on in this
chapter $/invert</sg_1>/sla/invert</sg_2>$ and $/invert</sg_1>^<1/sla d>
/sla/invert</sg_2>^<1/sla d>$ for two classes. It is therefore natural to
consider these powers and products on the logarithmic scale, and provide
an adjustment formula there. If $/log/vsgv^/ast $ is unbiased for
$/log/vsgv$, then $/log(/vsgv^/ast )^p$ is unbiased for $/log/vsgv^p$, for
each $p$!

The correction formula we derive now is based on exact theory under
normality, but is reasonable also without normality. Thus the starting
point is $/hs=A/sla m$ where $A$ of (10.5) is $</rm Wishart>_d(m,/sg),
m=n-1$. Thus $/invert<A>/sla/vsgv$ is distributed as the product of $d$
independent $/chi^2$-variables with degrees of freedom $m,m-1/upto m-(d-1)$,
cf.~Mardia, Kent and Bibby (1979, Ch.~3). And since
$$/eqalignno<
E/log/chi_/nu^2&=/int_0^/infty x^<</nu/over 2>-1>e^<-x/sla 2>/log/,x/,dx/sla
		2^</nu/over 2>/Gamma/lp</nu/over 2>/rp/cr
	     &=/log 2+/psi/lp</nu/over 2>/rp,&(10.6)/cr>$$
where $/psi(x)=</partial/over/partial x>/log/Gamma(x)=/Gamma'(x)/sla/Gamma(x)$,
it follows that
$$/eqalign<
E/log/vhsv&=E/log/lb/invert<A>/sla m^d/rb/cr
	  &=/log/vsgv+/sum_<i=0>^<d-1>/lb/log2+/psi/lp<m-i/over2>/rp/rb-
		d/log m/cr
	  &=/log/vsgv-B_d(m),/cr>$$
say, where
$$B_d(m)=d/log <m/over 2>-/sum_<i=0>^<d-1>/psi/lp<m-i/over 2>/rp./eqno(10.7)$$
Accordingly, $/log/vhsv+B_d(m)$ is unbiased for $/log/vsgv$, and
$$/vsgv^/ast =/vhsv e^<B_d(m)>/eqno(10.8)$$
becomes our corrected estimator of $/vsgv$. In particular we want
$$/lp/vsgv^/ast /rp^<1/sla2d>=/vhsv^<1/sla2d>/exp/lb<1/over2d>B_d(m)/rb
/eqno(10.9)$$
displayed as the estimate of the generalised standard deviation.

/medskip
/note<Remark 1.> The $/psi(x)=/Gamma'(x)/sla/Gamma(x)$ function is
available in standard mathematical software libraries (and in tables). An
approximation may also be used. One has
$$/psi(z)=/log z-<1/sla 2/over z>-<1/sla 12/over z^2>+<1/sla 120/over z^4>
-<1/sla252/over z^6>+/cdots,$$
see Abramowitz and Stegun (1964, Ch.~6), so that
$$/psi(<m/over2>)/doteq/log<m/over 2>-<1/over m>-<1/sla 3/over m^2>+
<2/sla 15/over m^4>./eqno(10.10)$$
Hence
$$/eqalignno<
B_d(m)&/doteq/log<m/over m-1>/cdots<m/over m-(d-1)>/cr
      &/quad+/sum_<i=0>^<d-1>/lb<1/over m-i>+<1/sla 3/over(m-i)^2>-
		<2/sla15/over(m-i)^4>/rb.&(10.11)/cr>$$
Notice that $B_d(m)$ approaches zero for large $m$.

/medskip
/note<Remark 2.> We have remarked that a small value of
$/vsgv^<1/sla2d>$ is good for classification, in particular a small
eigenvalue may lead us to study a corresponding principal component with
good discriminatory power. A small value of $/vsgv^<1/sla 2d>$ is also
good news for the training phase, because it becomes easier to ``learn'' a
distribution using a modest sample size when its concentration around the
mean is high. While it is useful to compare the $/invert</sg_k>^<1/sla2d>$
numbers for the different classes, the generalised standard deviation
$/vsgv^<1/sla2d>$ itself is not a good measure for whether the underlying
distribution is ``easy to learn'', however, since it depends upon the
particular scale used, say meter or centimeter. A better indicator is the
</it generalised coefficient of variation//>,
$$</rm CV>=/vsgv^<1/sla 2d>/sla/inpar/mu,/eqno(10.12)$$
$/inpar/mu$ being the length of $/mu$. (Another possible generalisation
of the one-dimensional $/ssg/sla/mu$ is $(/mu'/sg^<-1>/mu)^/mh$.)
CV can be estimated to guide one when planning the size of a possible
second stage training phase, but is not perfect for that purpose, because
of obvious difficulties for $/mu$ close to zero. Let us nevertheless
provide an estimator for CV, correcting the obvious
$/vhsv^<1/sla2d>/sla/inpar</hm>$ for biasedness on the log scale.
Its logarithm is $<1/over 2d>/log/vhsv-/halv/log/inpar/hm^2$, and
$$E/inpar/hm^2=E/sum_<i=1>^d</hm_i>^2=/sum_<i=1>^d(</mu_i>^2+<1/over
n>/ssg_<ii>)=/inpar/mu^2+<1/over n></rm Tr>(/sg).$$
Hence $E/log/inpar/hm^2/doteq/log/inpar/hm^2+<1/over/inpar/mu^2><1/over
n></rm Tr>(/sg)$, and $/log/inpar/hm^2-<1/over n></rm
Tr>(/hs)/sla/inpar/hm^2$ is better for $/log/inpar/mu^2$, i.e.
$/inpar/hm/exp/lb-<1/over2n></rm Tr>(/hs)/sla/inpar/hm^2/rb$ is better
for $/inpar/mu$. Combining this with (10.8) we reach
$$</rm CV>^/ast =</vhsv^<1/sla2d>/over/inpar/hm>/exp/lb<1/over2d>B_d(m)+
<1/over2n><</rm Tr>(/hs)/over/inpar/hm^2>/rb/eqno(10.13)$$
as the recommended estimate of CV.

/smallskip
We have so far discussed uses for the estimated $/sg$ matrix. Another
useful entity is the </it correlation matrix//> $R$ with elements $/rho_<ij>=
/ssg_<ij>/sla(/ssg_<ii>/ssg_<jj>)^<1/sla2>$, estimates for which are 
$/wh/rho_<ij>=/wh/ssg_<ij>/sla(/wh/ssg_<ii>/wh/ssg_<jj>)^<1/sla2>$.

A low correlation may (at best) mean approximate independence between the
coordinates in question, i.e. invites the possibility of using a
multiplicative structure in the density, for increased efficiency, cf.
remarks made in 5.4. A classifier based on $N_d(/mu_k,/sg_k)$ descriptions
may sometimes work better, if training sets are small, if off-diagonal
elements in $/sg_k$ are put to zero whenever the corresponding
$/wh/rho_<ij>$'s are below say $/halv$ in absolute value. (A more formal
rule can be based on the fact that if coordinates $i$ and $j$ of $X$
really are statistically independent, then 
$$/sqrt<n>/wh/rho_<ij>/totop<D>N(0,1),$$
regardless of the underlying distributions,
so that only $/rho_<ij>$'s with $/invert</wh/rho_<ij>>/gt1.96/sla/sqrt<n>$,
say, are ``significant''.)

If a correlation is high, on the other hand, one of the components 
``explains'' much of the other one, and if this is seen in all classes
one might remove one of the components in the feature vector.

We have discussed in this section various entities that could constitute
``standard output'' for each class after the (initial or final)
training phase: $n_k$, estimates of $/mu_k$, $/sg_k$,
$/invert</sg_k>^<1/sla 2d>$, $CV_k=/invert</sg_k>^<1/sla 2d>/sla/inpar</mu_k>$,
estimated root eigenvalues $(/hl_<k,1>)^<1/sla 2>/upto(/hl_<k,d>)^<1/sla 2>$,
and $</wh R>_k$.

It is also desirable to include diagnostics for each class that tell one 
whether some specific ``class description method'', like normality or the
multiplicative beta model with cosine expansions, cf.~5.4, is adequate. 
This topic is treated, to a limited extent, in Chapter 11.
/vfill/eject
/section<10.2><Error rates and separability measures><>
/subsection<10.2.A><Error rates> One should ideally be able to predict,
with some precision, the error rates that will result from a classifier
on the basis of measures for separability between classes, and that can be 
evaluated even before the classifier is put to work. While
this may be hoping for too much, such distance measures can be constructed,
and are useful for several reasons. One can sometimes sort out the ``safe'' 
cases, i.e. classes that are almost never confused with others, and on the
other hand, one can also spot the hard cases, where a pair of classes are to
close to hope for good discrimination at all. In such a case one might consider
momentarily pooling two close ones together, aiming only at classifying such
objects into the new extended class, with the present choice of features.
A second stage could utilise a new feature that aimed specifically at 
discriminating these two. Thus a table of distances (of some sort) are helpful
in devising hierarchical classification schemes. They are also essential 
when attacking the difficult problem of ``reducing dimensionality'',
i.e.~how to select a small number of feature components among a large number
such that error rates become as low as possible.

Let us concentrate on a particular pair of classes, say ``$1$'' and ``$2$'',
with class densities $f_1$ and $f_2$. We aim for the moment at good measures
for the distance between $f_1$ and $f_2$ as such, without regard to the
prior probabilities. Thus think of ``$1$'' and ``$2$'' as the only classes 
present, and take them to be equally likely. The optimal allocation rule
assigns label ``$2$'' to a new $X$ if it falls in
$$A=/lb x:f_2(x)/gt f_1(x)/rb,/eqno(10.14)$$
and to ``$1$'' otherwise. The two </it conditional//> error rates are
$$/eqalignno<
/eps(1/to2)&=/Pr/lb X/in A/vert X/sim f_1/rb=/int_A f_1,/cr
/eps(2/to1)&=/Pr/lb X/not/in A/vert X/sim f_2/rb=/int_</sim A>f_2,&(10.15)/cr>
$$
and the overall, </it unconditional//> error rate becomes
$$/eps=/halv/int_A f_1+/halv/int_</sim A>f_2.$$
Now observe that
$$/int/invert<f_1-f_2>=/int_A(f_2-f_1)-/int_</sim A>(f_2-f_1)=
2/int_A(f_2-f_1),$$
so that
$$/eps=/halv/int_A f_1+/halv-/halv/int_A f_2=/halv-<1/over4>/int
/invert<f_1-f_2>./eqno(10.16)$$

Thus the $L_1$ distance $/int/invert<f_1-f_2>$ is a very natural distance
measure in the present classification problem context. It is zero only if 
$f_1=f_2$, when no discrimination is possible, and is at its maximum $2$
when $f_1$ and $f_2$ are singular (supported on disjoint sets), and then
perfect discrimination is possible.  If $f_1$ and $f_2$ are Gaussian with
means $/mu_1$ and $/mu_2$, and common covariance matrix $/sg$, for example,
then one can show that
$$/int/invert<f_1-f_2>=4/Phi/lp</dl/over2>/rp-2=2/Phi/lp-</dl/over2>,
</dl/over2>/rp,$$
where
$$/dl=/lb(/mu_1-/mu_2)'/sg^<-1>(/mu_1-/mu_2)/rb^<1/sla 2>/eqno(10.17)$$
is the Mahalanobis distance. This corresponds by (10.16) to a classic formula
$$/eps=/Phi/lp-/halv/dl/rp/eqno(10.18)$$
relating, in this particular setting, error rate to distance.

Note that (10.16) is completely general, for example, the densities may be
point mass distributions for discrete variables, or mixed. In the simplest
mixed-variable model, for example, where $X=(A,Y)$ and
$$f_k(x)=f_k(a,y)=q_k(a)N_d/lp/mu^<(k)>_a,/sg/rp(y),$$
cf.~Section 9.3, one has
$$/eqalign<
/int/invert<f_1-f_2> dx&=/sum_a/int_</rr^d>/vert q_1(a)N_d/lp/mu_a^<(1)>,/sg/rp(y)/cr
		&/qquad/qquad/qquad-q_2(a)N_d/lp/mu_a^<(2)>,/sg/rp(y)/vert dy,
/cr>$$
which can be evaluated explicitly, giving error rate
$$/eqalignno<
/eps&=/halv/sum_a q_1(a) /Pr/lb N/gt/halv/dl_a+<1/over/dl_a>/log
	<q_1(a)/over q_2(a)>/rb/cr
    &/quad+/halv/sum_a q_2(a)/Pr/lb N/gt/halv/dl_a-<1/over/dl_a>/log
	<q_1(a)/over q_2(a)>/rb&(10.19)/cr>$$
where $/dl_a^2=(/mu_a^<(1)>-/mu_a^<(2)>)'/sg^<-1>(/mu_a^<(1)>-/mu_a^<(2)>)$
are the squared Mahalanobis distances in $Y$-space for the various 
$A$-outcomes. ($N$ denotes a standard normal variable.)

Note further that (10.16) makes a strong point in favour of the $L_1$ distance
as a ``universal'' one, in that it does not depend upon dimension or scale.
If the vector $X$ is transformed in any one-to-one manner to a $Y$ the $L_1$
distance is preserved. A $6$-dimensional $L_1$ distance observed in Germany
can be directly compared to a $3$-dimensional Norwegian one, since they relate
to probabilities of the same type of event, by (10.16).

It is fair to add here, however, that the fact that $/int/invert<f_1-f_2>$ can be 
explicitly evaluated only for the simplest models limits its wider 
applicability, and causes us to look for other distance measures, to cope with
normal distributions with unequal covariance matrices, for example.

It has been usual in the literature to display tables of estimated Mahalanobis
distances, using either the ``two-pooled'' estimated covariance matrix for
the two classes in question, or the ``grand-pooled'' one where a common 
$/sg$ is used for all pairs. It is somewhat easier to relate to and immediately
interpret distances, say in the range $2$ to $15$, than to tiny error rates.

The problem with such practice is that the discriminatory power of unequal 
covariance matrices is not seen, and more generally, that the underlying 
statistical model used to estimate $f_1$ and $f_2$, for example normality,
is never exactly correct.

It is still useful to produce model-based estimates of error rates, however.
Suppose the training set gives us $/hf_1$ and $/hf_2$, based on
some model and method, parametric or nonparametric. Accompanying
$$/eideal=/halv/int_A f_1+/halv/int_</sim A>f_2=/halv-<1/over4>
/int/invert<f_1-f_2>$$
is the model based
$$/emodel=/halv/int_</hA>/hf_1+/halv/int_</sim/hA>/hf_2 =
/halv-<1/over4>/int/invert</hf_1-/hf_2>,/eqno(10.20)$$
where
$$/hA=/lb x:/hf_2(x)/gt/hf_1(x)/rb/eqno(10.21)$$
is the deciding region actually used for candidate vectors to come. The
real error probability to be seen in future is
$$/eqalignno<
/etrue &=/halv/int_</hA>f_1+/halv/int_</sim/hA>f_2/cr
&=/halv-/halv/int_</hA>(f_2-f_1)/cr
&=/halv-/halv/int_</hA>/lp/hf_2-/hf_1+f_2-/hf_2-f_1+/hf_1/rp/cr
&=/emodel+/halv/int_</hA>/lp f_1-/hf_1/rp+/halv/int_</sim/hA>
/lp f_2-/hf_2/rp.&(10.22)/cr>$$
Observe that $/emodel$ is not necessarily a good estimate of 
$/eideal$, the two will in practice be close only when the method is 
nonparametric and the sample size is large. More important for the present 
discussion is that $/emodel$ often provides a reasonable estimate of $/etrue$, 
since
$$/int_/hA f_1/approx/int_/hA/wh f_1,/,/int_</sim/hA>f_2/approx
/int_</sim/hA>/wh f_2/eqno(10.23)$$
might hold, even with a coarse model leading to $/wh f_1$ and $/wh f_2$.
Asking for (10.23) to hold is certainly less than asking for $/wh f_i$
to be uniformly close to $f_i$, and often less than
$$/int_/hA/wh f_1/approx/int_A f_1,/,/int_</sim/hA>/wh f_2/approx
/int_</sim A>f_2.$$

We should also point out that the true error rate and its components 
$/int_</hA>f_1$, $/int_</sim/hA>f_2$ can be estimated too from the data, 
and without relying on the model that produced $/hf_1$ and $/hf_2$. The
nonparametric estimation of error rates is in fact the topic of Chapter 12.
While the ultimate assessment of the accuracy of a classifier should use
such methods, we intend here to provide tools to make a reasonable prediction
of its accuracy even before it is put to use, and maybe based on a preliminary
training set only. The nonparametric error rate estimates use cross-validation
or resampling techniques, and much work is needed before they can be found.
The intended use of the separability measures below includes that of quickly 
exploring, and perhaps later discarding, many choices of feature vectors.
/note<Example.> Let us illustrate the three different error rates by means of
a simple example. Observe first that as the training set's size grows
$/wh f_i$ converges to a certain function $f_i^0$ that is ``least
false'' under the model used, i.e.~closest to the true $f_i$ according
to the distance measure associated with the estimation procedure used,
e.g.~the Kullback-Leibler distance for maximum likelihood estimators,
cf.~Section 3.1. This means that $/hA$ of (10.21) eventually becomes
$$A^0=/lb x/colon f_2^0(x)/gt f_1^0(x)/rb,$$
and that 
$$/eqalign<
/emodel&/to/halv/int_<A^0>f_1^0+/halv/int_</sim A^0>f_2^0,/cr
/etrue&/to/halv/int_<A^0>f_1+/halv/int_</sim A^0>f_2./cr>$$

Now consider a situation where $f_1/sim N(0,1)$ and $f_2/sim
N(a,/tau^2)$, with $a/gt 0$ and $/tau/gt1$, but where the model
specifies
$f_1/sim N(/mu_1,/ssg^2)$, $f_2/sim N(/mu_2,/ssg^2)$. As the training
set grows $/hm_1$ converges to $0$ and $/hm_2$ to $a$, whereas the
common $/hssg$ converges to the least false value
$/ssg=(/halv+/halv/tau^2)^<1/sla 2>$. Hence the deciding region $/hA$
for future objects tends to $A^0=/lb x/colon x/gt/halv a/rb$, and 
$$/eqalign<
/emodel&/to/halv</rm Pr>/lb N(0,/ssg^2)/gt/halv a/rb +/halv</rm Pr>/lb
N(a,/ssg^2)/le/halv a/rb/cr
&=/Phi/lp-/halv<a/over/ssg>/rp,/cr
/etrue&/to/halv</rm Pr>/lb N(0,1)/gt/halv a/rb+/halv</rm
Pr>/lb N(a,/tau^2)/le/halv a/rb/cr
&=/halv/Phi/lp-/halv a/rp+/halv/Phi/lp-/halv<a/over/tau>/rp./cr>$$
Let us also compute $/eideal$: The optimal rule classifies $X$ to
``class 2'' if it falls in $$A=/lb x/colon/halv/lp1-<1/over/tau^2>/rp
x^2+<a/over/tau^2>x-/halv<a^2/over/tau^2>-/log/tau/gt0/rb,$$
which is seen after simplifying $/log/lb f_2(x)/sla f_1(x)/rb$.
One has
$$A=/lp-/infty,<-a-/tau D/over/tau^2-1>/rp/cup/lp<-a+/tau D/over/tau^2-1>,
/infty/rp$$
where $D=/lb a^2+2(/tau^2-1)/log/tau/rb^<1/sla2>$, and finds
$$/eqalign<
/eideal&=/halv-/halv/Phi/lp<a+/tau D/over/tau^2-1>/rp
-/halv/Phi/lp<-a+/tau D/over/tau^2-1>/rp/cr
&/quad+/halv/Phi/lp<-a/tau+D/over/tau^2-1>/rp+/halv/Phi
/lp<a/tau+D/over/tau^2-1>/rp./cr>$$

Figures (i)---(iv) display the three different error rates
discussed for the four cases $a=1, 3, 5, 7$, as a function of $/tau$.
$/emodel$ is at least reasonably close to $/etrue$ when
$1/le/tau/le3$,
which corresponds to cases where the model for $f_1$ and $f_2$ isn't
terribly wrong.
/pageinsert
/vbox to 20truecm<>/vfill
/note<Figures> (i)--(iv). </it Error rates (in percent) when discriminating
between $N(0,1)$ and $N(a,/tau^2)$, and the imperfect model specifies equal 
variance//>. The figures suggest that $/emodel$ can be used as a rough
estimate of $/etrue$ for usually occurring distances, even when the 
model is somewhat incorrect.
/endinsert
/pageinsert
/vbox to 20truecm<>/vfill
/endinsert
/subsection<10.2.B><Affinity, distance measures, and a generalised Mahalanobis 
distance>
The discussion above motivates studying $/emodel$ of (10.20), for 
various methods.
/note<Equal covariance matrices.> The classical normal model with common 
covariance matrix $/sg$, for example, gives
$$/emodel=/int_/hA/hf_1=/int_</sim/hA>/hf_2=/Phi/lp-/halv/hd/rp,
$$
where $/hd^2=(/hm_1-/hm_2)'/hs^<-1>(/hm_1-/hm_2)$ and
$$/hs=<n_1-1/over n_1+n_2-2>/hs_1+<n_2-1/over n_1+n_2-2>/hs_2$$
(say), $n_i$ being the size of the training set for class $i$.

It is good practice in this and similar cases to correct the plug-in version 
of $/hd^2$ above for sampling variability, making it less biased, before
computing $/emodel = /Phi(-/halv/hd)$, or before displaying $/hd$ in a table
of distances. While it is true that $/emodel$ as defined in (10.20)
happens to be $/Phi(-/halv/hd)$ with the uncorrected $/hd^2$, this will give an
incorrect (and too optimistic, as we shall se below) picture
of $/etrue$. This becomes clear when one considers the case $f_1$ close to
$f_2$. Then $/hd$ severely over-estimates the true $/dl$, and the two
right hand probabilities in (10.23) become too small, compared with the
true left hand probabilities, and that will be close to $/halv$. Thus we help
(10.23) and a fortiori $/emodel/approx/etrue$, if we use $/wh/emodel$
instead, with the distance measure corrected for bias. This remark also applies to distance measures for other models considered below.

In the present case $/hs=<1/over m>A$, where $A=(n_1-1)/hs_1+(n_2-1)/hs_2$
is $</rm Wishart>_d(/sg,m)$ and $m=n_1+n_2-2$. Since $EA^<-1>=/sg^<-1>/sla
/lp m-/lp d+1/rp/rp$ and $</rm VAR>(/hm_1-/hm_2)=(<1/over n_1>+<1/over n_2>)
/sg$ we obtain
$$/eqalign<
E/hd^2&=/sum_<i,j><m/over m-(d+1)>/ssg^<ij>/lb(/mu_<1,i>-/mu_<2,i>)
(/mu_<1,j>-/mu_<2,j>)+/lp<1/over n_1>+<1/over n_2>/rp/ssg_<ij>/rb/cr
&=<m/over m-(d+1)>/lb/dl^2+/lp<1/over n_1>+<1/over n_2>/rp
</rm Tr>(/sg^<-1>/sg)/rb,/cr>$$
so that
$$/bd^2=<m-(d+1)/over n>/,/hd^2-/lp<1/over n_1>+<1/over n_2>/rp d/eqno(10.24)$$
is recommended, along with $/wh/emodel=/Phi(-/halv/bd)$.
/note<Different covariance matrices.> Next consider the important case
with two normal distributions $(/mu_1,/sg_1)$ and $(/mu_2,/sg_2)$.
Here no closed form solution can be obtained for the error rate, and we must
look for non-obvious distance measures, hopefully generalising (10.17)
in a natural way.

First  we review a short piece of linear algebra. There is a standard 
decomposition of $/sg_1$ as $/sg_1=P'D_1P$, where $D_1$ is diagonal with 
the eigenvalues of $/sg_1$; in fact $P'$ has as its columns the normalised 
eigenvectors of $/sg_1$. Hence we can form the unique symmetric inverse
square root of $/sg_1$, namely $/sg_1^/mh=P'D_1^/mh P$, 
in particular $/sg_1^/mh/sg_1/sg_1^/mh=I$. Next consider
$/sg_1^/mh/sg_2/sg_1^/mh$. It is symmetric and non-negative 
definite, so there is a spectral decomposition
$$Q/sg_1^/mh/sg_2/sg_1^/mh Q'=D,$$
where $D$ is diagonal with the eigenvalues of $/sg_1^/mh/sg_2/sg_1^/mh$, 
and $Q'$ the corresponding orthonormal eigenvectors. In sum we have 
$$F/sg_1 F'=I/,,/, F/sg_2 F'=D/eqno(10.25)$$
with $F=Q/sg_1^/mh$. Note that $/invert<D>=/invert<F>^2/invert</sg_2>=
/invert</sg_2>/sla/invert</sg_1>$.

Now transform the feature vector $X$ to $Y=F(X-/mu_1)$. Then $Y$ is 
$N_d(0,I)$ if it is from ``$1$'' and $N_d(b,D)$ if it is from ``$2$'',
where $b=F(/mu_2-/mu_1)$. The $L_1$ distance and the error rate are 
unaffected by this transformation of data. Hence
$$/eqalign<
/int/invert<f_1-f_2>dx&=/int/vert N_d(0,I)(y)-N_d(b,D)(y)/vert/,dy/cr
&=/int(2/pi)^/mdh e^<-/halv y'y>/Bigl/vert 1-/invert<D>^/mh 
e^</halv y'y-/halv(y-b)'D^<-1>(y-b)>/Bigr/vert/,dy/cr
&=E/Big/vert 1-</invert</sg_1>^<1/sla 2>/over/invert</sg_2>^<1/sla 2>>/exp
	/halv/sum_<i=1>^d/lb Y_i^2-<1/over/lam_i>/lp Y_i-b_i/rp^2/rb/Big/vert,
/cr>$$
where the expectation is w.r.t. $Y$ being $N_d(0,I)$. Accordingly one might
consider evaluating (an approximation to) the error rate using numerical
integration or perhaps by simulating a large number of $Y$'s in the 
expectation formula above. This has turned out to involve unsuspected
numerical unstability problems, and we recommend computing the error rate, 
if strongly needed, using 
$$/eqalignno<
/eps&=/halv/Pr/lb N_d(b,D)(Y)/gt N_d(0,I)(Y)/vert Y/sim N_d(0,I)/rb/cr
&/quad+/halv/Pr/lb N_d(b,D)(Y)/lt N_d(0,I)(Y)/vert Y/sim N_d(b,D)/rb/cr
&=/halv/Pr/lb/sum_<i=1>^d/lp Y_i^2-<1/over/lam_i>(Y_i-b_i)^2/rp-/log
	/invert<D>/gt 0/,/vert Y/sim N_d(0,I)/rb/cr
&/quad+/halv/,/Pr/lb/sum_<i=1>^d/lp Y_i^2-(b_i+/lam_i^<1/sla 2>Y_i)^2/rp+/log
	/invert<D>/gt 0/,/vert Y/sim N_d(0,I)/rb,&(10.26)/cr>$$
and then obtaining these probabilities using simulation. One should use as
many as 10,000 replicates of $Y=(Y_1/upto Y_d)$ to secure a reasonable
accuracy to two decimal places.

Such a process, calculating for each pair $(/mu_k,/sg_k),(/mu_l,/sg_l)$ a
probability by simulating 10,000 or more outcomes, may appear frighteningly 
complicated and cumbersome, compared to a classical one like computing 
$/wh/dl$ or $/bar/dl$ from given parameter estimates and then $/eps=/Phi
(-/halv/dl)$. It is however a perfectly feasible task for the electronic 
computer of the 80ies and 90ies, and is not difficult to build into a
sophisticated discriminant analysis package system. One should bear in
mind that the classical theory, featuring linear methods and Mahalanobis
distances, was developed from 1936 to about 1960, and intended for the
mechanical calculator.

We will nevertheless proceed to obtain a reasonable generalisation of
the classical Mahalanobis distance for this case of unequal covariance
matrices.

For a pair of densities $f_1$ and $f_2$, introduce the </it affinity//>
between them by
$$/aff=/aff(f_1,f_2)=/int(f_1f_2)^<1/sla 2>./eqno(10.27)$$
It is related to a distance measure, namely
$$/lb/int/lp/sqrt<f_1>-/sqrt<f_1>/rp^2/rb^<1/sla 2>=(2-2/aff)^<1/sla
2>,$$
and provides valuable information about discrimination possibilities too,
though in a less direct manner than the $L_1$ distance. The main force of
the affinity measure is that it can be explicitly computed in a much
wider range of models than can $/int/invert<f_1-f_2>$, and that it more readily
is generalised to a measure for more than two classes, viz.
$/int(f_1/cdots f_K)^<1/sla K>$.

An affinity close to $1$ means that $f_1$ is close to $f_2$, and it is
difficult to discriminate. Affinity close to zero, on the other hand,
means that $f_1$ is small when $f_2$ is large, and vice versa, i.e.~a low
error rate should result.

A more precise result is 
$$/halv/lb1-/lp1-/aff^2/rp^<1/sla 2>/rb/le/eps/le/halv/aff./eqno(10.28)$$
Only the second half, which is the most useful, will be proven here. One
has /def/xx#1<f_1^<1/sla 2>#1f_2^<1/sla 2>>$/invert</xx<->>/lp/xx<+>/rp
/ge/invert</xx<->>/invert</xx<->>$, i.e. $/invert<f_1-f_2>/ge/lp/xx<->/rp^2$.
Consequently
$$/int/lb/invert<f_1-f_2>+2(f_1f_2)^<1/sla 2>/rb/ge2,$$
which indeed translates into $/halv-<1/over4>/int/invert<f_1-f_2>/le/halv
/int(f_1f_2)^<1/sla 2>$. Thus we have a useful upper bound on the error
rate, whenever we know the affinity.

It is not difficult to show that
$$/aff/lb N_d(/mu_1,/sg),N_d(/mu_2,/sg)/rb=e^<-<1/over 8>/dl^2>,/eqno(10.29)$$
with $/dl$ again as in (10.17). In the more general case, one obtains
$$/eqalign<
/aff&/lb N_d(/mu_1,/sg_1),N_d(/mu_2,/sg_2)/rb/cr
&/qquad=/int/lb N_d(0,I)(y)/,N_d(b,D)(y)/rb^<1/sla 2>dy,/cr>$$
using the transformation described in (10.25) and below. Clever algebra gives
$$/eqalign<
&/aff=/invert<D>^<-1/sla 4>/int(2/pi)^<-d/sla 2>/exp-<1/over4>/lb y'y+
(y-b)'D^<-1>(y-b)/rb dy/cr
&=/invert<D>^<-1/sla4>/int(2/pi)^<-d/sla2>/invert<T>^<-1/sla2>/exp-/halv
/Big/lbrace y'T^<-1>y-b'D^<-1>y/cr
&/qquad/qquad/qquad/qquad/qquad/qquad/qquad+/halv b'D^<-1>b/Big/rbrace dy/invert<T>^<1/sla2>/cr
&=/invert<D>^<-1/sla4>/invert<T>^<1/sla2>e^<-<1/over4>b'D^<-1>b>
Ee^</halv b'D^<-1>Y>,/cr>$$
writing $T=2(I+D^<-1>)^<-1>=2D(I+D)^<-1>$, and where the $Y$ is $N_d(0,T)$. 
Further simplification involves 
$$/eqalign<
/sg&=/halv(/sg_1+/sg_2),/cr
/dl^2&=(/mu_1-/mu_2)'/sg^<-1>(/mu_1-/mu_2)./cr>/eqno(10.30)$$
One gets
$$/aff=</invert</sg_1>^<1/sla4>/invert</sg_2>^<1/sla4>/over/invert</sg>^<1/sla2>
>e^<-<1/over8>/dl^2>=e^<-<1/over8>(/dl^2+/gamma^2)>,/eqno(10.31)$$
where
$$/gamma^2=4/log</invert</sg>/over/invert</sg_1>^<1/sla2>/invert</sg_2>^<1/sla2>
>/eqno(10.32)$$
is a measure for the discrepancy between the two covariance matrices.

It becomes natural then, in view of the generalisation from (10.29) to (10.31),
and the close connection between error rate and affinity, to define the 
</it generalised Mahalanobis distance//> as
$$/omega=(/dl^2+/gamma^2)^<1/sla2>,/eqno(10.33)$$
combining in a Pythagorean way a $/dl$-component for difference in mean
vectors and a $/gamma$-component for difference in covariance matrices.

Below we report on a limited study of the usefulness of the new distance 
$/omega$. Of course the construction and motivation above provide some 
justification for $/omega$ as a proper generalisation of the classical
$/dl$, but it is for example unclear how much direct information it offers
for the error rate. It turns out that for the few and rather regular cases 
considered,
$$/eps/doteq/Phi/lp-/halv/omega/rp/eqno(10.34)$$
holds as a reasonable approximation, thus lending extra support to the term
generalised Mahalanobis distance.

The first class of cases considered was the rather opposite one to the 
classical, namely that of equal means but different $/sg$ matrices,
and that for simplicity were taken to be proportional. In this case
an exact expression for the error rate can be found, using 
$$/log<N_d/lp0,/lam/sg/rp(x)/over N_d(0,/sg)(x)>=-/halv d/log/lam+/halv
/lp1-<1/over/lam>/rp x'/sg^<-1>x.$$
One arrives at 
$$/eps=/halv/Pr/lb/chi_d^2/sla d/gt</lam/log/lam/over/lam-1>/rb+
/halv/Pr/lb/chi_d^2/sla d/ge</log/lam/over/lam-1>/rb.$$
Now for given values $/omega=0.5,1.0/upto6.0$ the appropriate values for 
$/lam$ were found from
$$/omega/lb N_d(/mu,/lam/sg),N_d(/mu,/sg)/rb=/gamma
=/lb2d/log<(1+/lam)^2/over4/lam>/rb^<1/sla2>=/omega.$$
Then the exact error rate $/eps$ was computed and compared to
$/Phi(-/halv/omega)$. The results are displayed in Table 1 below for
dimension $d=1,3,5,10,1000$, and show very good agreement. A pleasant and perhaps 
unexpected fact is that $/eps=/Phi(-/halv/gamma)=/Phi(-/halv/omega)$
</it exactly//> in the limiting case where $d/rightarrow/infty$. Also included in the 
table is the general upper bound $/halv/aff=/halv e^<-<1/over8>/omega^2>$
from (10.28). It appears to be too high in this situation to be of any use.

/pageinsert
/note<table 1.></it/ Error rates (in percent) for two opposite 
situations//>. Situation ($i$) concerns discriminating between
$N_d(/mu_1,/sg)$ and 
$N_d(/mu_2,/sg)$, while ($ii$) concerns discriminating between $N_d(/mu,/sg)$
and $N_d(/mu,/lambda/sg)$. The error rates are given as a function of the new
distance $/omega$, for dimensions $d=1$, $3$, $5$, $10$.

/noindent/hbox to/hsize</hfil/vbox</halign<&/quad/hfil#/cr
/hfil $/omega$&$/eps(/omega)$&$d=1$&$d=3$&$d=5$&$d=10$&$d=1000$&
upper bound/cr
/noalign</smallskip>
0.5&40.1294&41.4463&40.6239&40.4297&40.2809&40.1254&48.4617/cr
1.0&30.8537&32.9313&31.6388&31.3299&31.0951&30.8574&44.1248/cr
1.5&22.6627&24.6405&23.4222&23.1219&22.8954&22.6660&37.7420/cr
2.0&15.8655&16.9981&16.3116&16.1314&15.9997&15.8660&30.3265/cr
2.5&10.5650&10.5885&10.5576&10.5481&10.5536&10.5649&22.8917/cr
3.0& 6.6807& 5.8640& 6.2675& 6.4024& 6.5312& 6.6789&16.2326/cr
3.5& 4.0059& 2.8631& 3.3675& 3.5702& 3.7674& 4.0032&10.8133/cr
4.0& 2.2750& 1.2286& 1.6182& 1.8098& 2.0113& 2.2720& 6.7668/cr
4.5& 1.2224& 0.4630& 0.6888& 0.8256& 0.9866& 1.2196& 3.9780/cr
5.0& 0.6210& 0.1532& 0.2580& 0.3359& 0.4412& 0.6186& 2.1968/cr
5.5& 0.2980& 0.0445& 0.0847& 0.1210& 0.1786& 0.2963& 1.1397/cr
6.0& 0.1350& 0.0114& 0.0244& 0.0384& 0.0649& 0.1339& 0.5554/cr>>/hfil>
/medskip
The rates for situation ($i$), $/eps(/omega)=/Phi(-/omega/sla 2)$, are exact. The
table shows that this is a good approximation also for situation ($ii$).
Also given in the table are upper bounds, affinity$/sla2$. These appear
to be too large to be of value.
/note<table 2.></it/ Error rates (in percent) for $N_d(/mu_1,/sg_1)$
versus $N_d(/mu_2,/sg_2)$ in a situation where difference in covariance matrices
equals difference in means (as measured by the new distance $/omega$)//>.
The error rates are given as a function of $/omega$, for dimensions $d=1$, $3$,
$5$, $10$, $1000$. The right hand column gives the approximation 
$/eps(/omega)=/Phi(-/omega/sla2)$. The table suggests that this approximation 
is fairly good.

/noindent/hbox to/hsize</hfil/vbox</halign<&/quad/hfil#/cr
$/omega$&$d=1$&$d=3$&$d=5$&$d=10$&approximation/cr
/noalign</smallskip>
0.5&41.530&40.45&39.95&40.45&40.1294/cr
1.0&32.610&31.09&31.80&30.40&30.8537/cr
1.5&24.045&22.88&22.95&24.65&22.6627/cr
2.0&16.540&16.01&15.00&15.90&15.8655/cr
2.5&10.385&10.80&10.25& 9.40&10.5650/cr
3.0& 5.810& 6.26& 6.10& 7.55& 6.6807/cr
3.5& 2.935& 3.15& 3.15& 3.25& 4.0059/cr
4.0& 1.165& 1.46& 2.25& 2.15& 2.2750/cr
4.5& 0.395& 0.81& 0.70& 1.35& 1.2224/cr
5.0& 0.185& 0.32& 0.25& 0.45& 0.6210/cr
5.5& 0.045& 0.14& 0.00& 0.20& 0.2980/cr
6.0& 0.000& 0.04& 0.00& 0.00& 0.1350/cr>>/hfil>
/endinsert
The second class of cases considered is the one where the difference in means 
and in $/sg$ matrices are of equal importance, as measured by $/dl$ and
$/gamma$. Thus for each value of $/omega$ a situation was found where 
$/dl=/omega/sla/sqrt<2>$ and $/gamma=/omega/sla/sqrt<2>$. We took for 
simplicity  the case of $N_d(0,I)$ versus $N_d(b,/lam I)$, with 
$b=(b_0/upto b_0)'$. Table 2 below displays approximations to the
exact error rate, using the simulation device of (10.26), and again they
are in good agreement with $/Phi(-/halv/omega)$.

On these grounds we recommend evaluating estimates of $/omega=(/dl^2+
/gamma^2)^<1/sla2>$ for each pair of classes under consideration 
(estimates better than the plug-in ones are developed in 10.3 below),
in the general process of feature evaluation. One can also set up
an accompanying table of tentative pairwise error rates, using
$/Phi(-/halv/omega)$.
/note<Remark 3.> Chernoff (1973)  provides another distance measure between
multivariate normal distribution with unequal covariance matrices, and in the
same spirit as our $/omega$, having also $/eps=/Phi(-/halv/omega)$ in
mind. His distance entails $/eps=/Phi(-/halv/omega)$ </it exactly//>,
but where $/eps$ is the error rate for the best </it linear//> rule,
whereas the present treatment has considered the more interesting 
</it optimal//> classification rule, which is the quadratic one.
/subsection<10.2.C><Other models> Suppose the merits of two different feature 
vector choices are to be compared, and that perhaps very different 
probability measures are used to describe the resulting class densities. It
is natural to turn to affinities then, because (10.27) can be explicitly 
evaluated for a wide range of models, and because of its invariance under 
general  transformations and its relation to error rates, cf.~(10.28).
One can also calculate </it distances//> through
$$/omega=/lsb-8/log/lb/aff(f_1,f_2)/rb/rsb^<1/sla2>,/eqno(10.35)$$
in order to have meaningful measures on the same scale as Mahalonobis distances.
But of course one cannot expect the $/eps/doteq/Phi(-/halv/omega)$ 
approximation to hold in any generality.

Imagine, for example, that method one gives $f_k$ and $f_l$ that are reasonably 
normal, resulting in distance
$$/omega_<k,l>=(/dl_<k,l>^2+/gamma_<k,l>^2)^<1/sla2>,$$
say, whereas method two gives discrete components and discrete probability
distributions $g_k$ and $g_l$. Then 
$$/aff(g_k,g_l)=/sum_x/lb g_k(x)g_l(x)/rb^<1/sla2>$$
is merely a finite sum, and $/omega_<k,l>$ above can be compared to 
$/lbrack-8/log/sum_x/lb g_k(x)g_l(x)/rb^<1/sla2>/rbrack^<1/sla2>$.

Of course these distances must be estimated in practice, and plug-in estimates 
should be  corrected for bias, as in 10.3 below.

Of perhaps particular importance in the symbol recognition context are feature 
vectors with both discrete and continuous components, cf.~Chapter 9.
Consider $X=(A,Y)$ where
$$/fkx=q_k(a)N_c/lp/mu_a^<(k)>,/sg_a^<(k)>/rp.$$
Then
$$/aff(f_k,f_l)=/sum_a/lb q_k(a)q_l(a)/rb^<1/sla2>e^<-<1/over8>/lb/dl^2_a(k,l)+
/gamma_a^2(k,l)/rb>,$$
where
$$/eqalign<
/dl_a^2(k,l)&=/lp/mu_a^<(k)>-/mu_a^<(l)>/rp'/lb/halv/lp/sg_a^<(k)>+/sg_a^<(l)>/rp
/rb^<-1>/lp/mu_a^<(k)>-/mu_a^<(l)>/rp,/cr
/gamma_a^2(k,l)&=4/log</invert</halv/lp/sg_a^<(k)>+/sg_a^<(l)>/rp>/over
/invert</sg_a^<(k)>>^<1/sla2>/invert</sg_a^<(l)>>^<1/sla2>>./cr>$$
/subsection<10.2.D><More than two classes> Assume there are $K$ classes present, 
with a priori probabilities $/pi_1/upto/pi_K$. Consider
$$A_<k,l>=/lb x:/pi_kf_k(x)/lt/pi_lf_l(x)/rb$$
and
$$A_k=/lb x:/pi_kf_k(x) </rm/ is/ not/ largest/ among/ >/pi_lf_l(x)
/rb=/bigcup_<l/not=k>A_<k,l>.$$
Then
$$/eps(k/to l)=P_k(A_<k,l>)$$
is the error probability of an $X$-vector from $k$ falsely being closer
to class $l$, in their own two-class competition. The error rates
discussed earlier in this chapter have been of the type 
$$/eps_<k,l>=</pi(k)/over/pi(k)+/pi(l)>/eps(k/to l)+
</pi(l)/over/pi(k)+/pi(l)>/eps(l/to k),$$
but with $/pi(k)=/pi(l)$, since we aimed at ``objective''
interclass distances.

Of interest is the overall error rate for objects of type $k$, i.e.
$/eps_k=P_k(A_k)$ (we rule out doubt and outlier possibilities now).
One has 
$$/eps_k/le/sum_<l/not=k>/eps(k/to l),$$
but it is otherwise difficult to obtain more precise results for 
$/eps_k$, in any generality.

The best way to estimate $/eps_k$, conditional error rates $P_k/lbrace/pi_lf_l
(x)=/mo<m/le K>/pi_m$/break$f_m(X)/rbrace$, total error rate
$/sum_k/pi_k/eps_k$, and doubt and 
outlier rates, is by some cross-validation scheme as in Chapter 12, but perhaps 
only when the features have been found sufficiently promising to warrant
full implementation and ``full learning'', i.e.~with estimates from a serious,
as opposed to a preliminary, training set. We reemphasise that the methods
of this chapter also have been intended for preliminary and exploratory 
analysis, and for planning of hierarchical classification procedures.
/section<10.3><Estimating the generalised Mahalanobis distance><> The
previous section discussed distances between populations, and in particular 
provided
$$/omega=(/dl^2+/gamma^2)^<1/sla2>$$
as such a measure to be estimated for each pair of classes, where
$$/eqalign<
/dl^2&=(/mu_1-/mu_2)'/sg^<-1>(/mu_1-/mu_2),/cr
/gamma^2&=4/log</invert</sg>/over/invert</sg_1>^<1/sla2>/invert</sg_2>^<1/sla2>>
,/cr>$$
and $/sg=/halv(/sg_1+/sg_2)$. While the $/omega$ distance was derived for and
is particularly relevant for multinormal densities, it is reasonable even
without this assumption.
/subsection<10.3.A><Adjustments for bias> The plug-in estimates are
$$/eqalignno<
/hd_0^2&=(/hm_1-/hm_2)'/hs^<-1>(/hm_1-/hm_2),/cr
/hg_0^2&=4/log</invert</hs>/over/invert</hs_1>^<1/sla2>/invert</hs_2>^<1/sla2>>
,&(10.36)/cr>$$
where $/hm_i=<1/over n_i>/sum_<t=1>^<n_i>X_t^<(i)>$, $/hs_i=<1/over n_i-1>
/sum_<t=1>^<n_i>(X_t^<(i)>-/hm_i)(X_t^<(i)>-/hm_i)'$, and
$/hs=/halv(/hs_1+/hs_2)$. These have a non-negligible bias for finite $n_1$ and 
$n_2$, and should be adjusted.

This turns out to be difficult to do in an exact manner, under normality, 
because of the complicated distribution for $/hs$, and we will be content
with an approximation, boldly taking 
$$/hs=/halv<A_1/over m_1>+/halv<A_2/over m_2>/doteq c</rm Wishart>_d(K,m)
/eqno(10.37)$$
for appropriately chosen constant $c$ and matrix $K$, where $m_1=n_1-1$,
$m_2=n_2-1$ and $A_i$ is $</rm Wishart>_d(/sg_i,m_i)$. ($/doteq$ here means
``approximately distributed as.'') Such approximations have a long tradition and
are  known to be good in the one-dimensional case, effectively coming close
to the distribution of a linear combination of independent chi-squares.

Taking expectations in (10.37) one gets $/halv/sg_1+/halv/sg_2=/sg=cmK$, i.e.
$$/hs /doteq c</rm Wishart>_d(<1/over cm>/sg,m)=<A/over m>,/eqno(10.38)$$
where $A$ is $</rm Wishart>_d(/sg,m)$, and it remains to specify a value
for $m$ that makes the approximation good. It is natural to get the second 
moments to agree as much as possible. More specifically, if $y$ is a 
$d$-vector then (10.37) and (10.38) imply
$$/halv<1/over m_1>y'A_1y+/halv<1/over m_2>y'A_2y /doteq <1/over m>y'Ay,$$
and $y'A_iy/sim(y'/sg_iy)/chi_<m_i>^2$, $y'Ay/sim(y'/sg y)/chi_m^2$.
Having equal variances means
$$<1/over m>=/lp<z_1/over z_1+z_2>/rp^2<1/over m_1>+/lp<z_2/over z_1+z_2>/rp^2
<1/over m_2>,$$
where $z_i=y'/sg_iy$. A good value for $m$ hinges on $r=<z_1/over z_2>$, which  
for example is equal to $/ssg_<1,ii>/sla/ssg_<2,ii>$ for one choice of $y$.
To settle on some reasonable average value one can consider $r$ on the
transformed scale after applying the simultaneous diagonalisation trick;
$F/sg_1F'=I$ and $F/sg_2F'=D=</rm diag>/lb/lam_1/upto/lam_d/rb$ gives 
$r=v'D^<-1>v/sla v'v$ for $v=D^<1/sla2>F^<-1>y$. In particular
$<1/over/lam_1>/upto<1/over/lam_d>$ are extreme choices, and a natural
average value is
$$/overline r=/lp<1/over/lam_1>/cdots<1/over/lam_d>/rp^<1/sla d>=/invert<D>^<-1/sla d>=
/lp</invert</sg_1>/over/invert</sg_2>>/rp^<1/sla d>./eqno(10.39)$$
(Another choice is $/overline r=<1/over 
d>(<1/over/lam_1>+/cdots+<1/over/lam_d>)/allowbreak=
<1/over d></rm Tr>(D^<-1>)/allowbreak=<1/over d></rm Tr>((F')^<-1>/sg_2^<-1>F^<-1>)$
$=<1/over d></rm Tr>(/sg_1/allowbreak/sg_2^<-1>)$.)

Now solve
$$<1/over m>=/lp</overline r/over /overline r+1>/rp^2<1/over m_1>+/lp<1/over /overline r+1>/rp^2
<1/over m_2>$$
for $m$, i.e.
$$m=<(/overline r+1)^2m_1m_2/over m_1+/overline r^2m_2>./eqno(10.40)$$
This will be our choice for $m$, only with $/overline r$ replaced with an estimate,
namely
$$/overline r=</invert</hs_1>^<1/sla d>/exp/lb<1/over d>B_d(m_1)/rb/over
/invert</hs_2>^<1/sla d>/exp/lb <1/over d>B_d(m_2)/rb>,/eqno(10.41)$$
following (10.7)--(10.9). Note that the nominator and denominator
have essentially been computed already as part of the ``standard output''
describing each class, cf.~10.1.

We are now in a position to compute approximate expectations of the naive
plug-in estimators (10.36).

Firstly, since $</rm VAR>(/hm_1-/hm_2)=<1/over n_1>/sg_1+<1/over n_2>/sg_2$,
one has
$$/eqalign<
E/hd_0^2&=E/sum_<i,j>(/hm_<1,i>-/hm_<2,i>)(/hm_<1,j>-/hm_<2,j>)/wh/ssg^<ij>/cr
&/doteq/sum_<i,j>/lb(/mu_<1,i>-/mu_<2,i>)(/mu_<1,j>-/mu_<2,j>)+
<1/over n_1>/ssg_<1,ij>+<1/over n_2>/ssg_<2,ij>/rb EmA^<-1>/cr
&=<m/over m-(d+1)>/lb/dl^2+/sum_<i,j>/lp<1/over n_1>/ssg_<1,ij>+<1/over n_2>
/ssg_<2,ij>/rp/ssg^<ij>/rb/cr
&/doteq<m/over m-(d+1)>/lb/dl^2+<1/over n_0>/sum_<i,j>/ssg_<1,ij>/ssg^<ij>+
<1/over n_0>/sum_<i,j>/ssg_<2,ij>/ssg^<ij>/rb/cr
&=<m/over m-(d+1)>/lb/dl^2+<2/over n_0></rm Tr>(/sg/sg^<-1>)/rb,/cr>$$
i.e.
$$</hd>^2=<m-(d+1)/over m>/,</hd>_0^2-<2d/over n_0>/eqno(10.42)$$
should be cured for bias. Another approximation was used here to obtain a 
reasonable simple formula, namely that of equating $<1/over n_1></rm Tr>
(/sg_1/sg^<-1>)+<1/over n_2></rm Tr>(/sg_2/sg^<-1>)$ to $<1/over n_0>
</rm Tr>(/sg_1/sg^<-1>)+<1/over n_0></rm Tr>(/sg_2/sg^<-1>)$.
This holds exactly for $n_0=(s+1)n_1n_2/sla(n_1+sn_2)$, where
$s=</rm Tr>(/sg_1/sg^<-1>)/sla$/break$</rm Tr>(/sg_2/sg^<-1>)$. Using the
diagonalisation transformation that led to (10.39) again,
$$/eqalign<
s&=</Tr(I+D)^<-1>/over/Tr D(I+D)^<-1>>=<<1/over d>/sum_1^d 1/sla(1+/lam_i)/over
<1/over d>/sum_1^d/lam_i/sla(1+/lam_i)>/cr
&/qquad/doteq</lb/prod_1^d 1/sla(1+/lam_i)/rb^<1/sla d>/over/lb/prod_1^d
/lam_i/sla(1+/lam_i)/rb^<1/sla d>>=/invert<D>^<-1/sla d>=
</invert</sg_1>^<1/sla d>/over/invert</sg_2>^<1/sla d>> = </overline r>,/cr>$$
and we may conveniently use $s=/overline r$ of (10.41) once more to get
$$n_0=<(</overline r>+1)n_1n_2/over n_1+</overline r>n_2>.$$
$n_0$ is always between $n_1$ and $n_2$, and is in particular equal to the
common value $n_1=n_2$ when these are equal, and in this balanced case the 
second  approximation in $E/hd_0^2$ is an identity.

Secondly, consider
$$/wh/gamma_0^2=4/log/invert</hs>-2/log/invert</hs_1>-2/log/invert</hs_2>,$$
where, under normality still,
$$/log/invert</hs_1>=/log/invert<<A_1/over m_1>>/,</buildrel D/over=>/,/log
/invert</sg_1>+/sum_<i=0>^<d-1>/log/chi_<m_1-i>^2-d/log m_1$$
and similarily for $/log/invert</hs_2>$ and $/log/invert</hs>$, using
$/hs /doteq A/sla m$ again. By results of 10.1, see e.g.~(10.6)--(10.8),
one obtains
$$/eqalign<
E/wh/gamma_0^2&/doteq/gamma^2-2d/log<m^2/over m_1m_2>+/sum_<i=0>^<d-1>/lb
4/psi/lp<m-i/over 2>/rp-2/psi/lp<m_1-1/over2>/rp-2/psi/lp<m_2-i/over2>/rp/rb/cr
&=/gamma^2-4B_d(m)+2B_d(m_1)+2B_d(m_2)./cr>$$
The bias-adjusted estimator becomes
$$/eqalignno<
/wh/gamma^2&=/wh/gamma_0^2+4B_d(m)-2B_d(m_1)-2B_d(m_2)/cr
&=/wh/gamma_0^2+2d/log<m^2/over m_1m_2>-/sum_<i=0>^<d-1>/lb4/psi
(<m-i/over2>)-2/psi(<m_1-i/over2>)-2/psi(<m_2-i/over2>)/rb.&(10.43)/cr>$$
Combining (10.41) and (10.43) we reach the desired distance estimator
$$/wh/omega=(/hd^2+/wh/gamma^2)^<1/sla2>./eqno(10.44)$$
/subsection<10.3.B><Amount of noise in the distance estimates> It is of interest
to assess the sampling variability of $/wh/omega$ and other
estimators appearing above. Below is a short study
of this problem, using limit distribution techniques that are useful also in
other circumstances. The results can be used to give approximate confidence
intervals for a $/dl$, a $/gamma$, or an $/omega$, one can put up a test for 
$/omega/le4$ against $/omega/gt4$, etc. The main objective is however
to get information about the ``level of noise'', so that one can be warned,
for example, against basing decisions on unreliable estimates.

We begin with a lemma on the limiting distribution of the inverse sample 
covariance matrix.
/Blemma Let $/hs=/sum_<t=1>^n(X_t-/hm)(X_t-/hm)'/sla(n-1)$ be
the sample covariance matrix based on an i.i.d.~sample from a
distribution in $/rr^d$ possessing finite fourth order moments, and let 
$EX=/mu$, $</rm VAR>(X)=/sg=(/ssg_<ij>)$. Then $/sqrt<n>(/hs-/sg)/totop<D>S$,
a stochastic symmetric matrix with elements $S_<ij>$ whose distribution is 
multinormal, and with 
$$</rm cov>(S_<ij>,S_<kl>)=E/lb(X_i-/mu_i)(X_j-/mu_j)-/ssg_<ij>/rb
/lb(X_k-/mu_k)(X_l-/mu_l)-/ssg_<kl>)/rb.$$
If $/sg^<-1>=(/ssg^<ij>)$ exists, then
$$/sqrt<n>(/hs^<-1>-/sg^<-1>)/totop<D>T=-/sg^<-1>S/sg^<-1>./eqno(10.45)$$

If in particular $X$ is normal, then 
$$</rm cov>(S_<ij>,S_<kl>)=/ssg_<ik>/ssg_<jl>+/ssg_<il>/ssg_<jk>/eqno(10.46)$$
and
$$</rm cov>(T_<ij>,T_<kl>)=/ssg^<ik>/ssg^<jl>+/ssg^<il>/ssg^<jk>./eqno(10.47)$$
/Elemma
/Bproof The first part is elementary and well known, and follows 
essentially from the central limit theorem, upon noting that $/sqrt<n>
(/hs-/hs_0)/totop<P>0$, where $/hs_0=/sum_<t=1>^n(X_t-/mu)(X_t-/mu)'
/sla(n-1)$. The normality assumption implies (10.46), see e.g.~Mardia, Kent and 
Bibby (1979, p.~92).

It is clear that $/sqrt<n>(/hs^<-1>-/sg^<-1>)$ has a limiting multinormal
distribution, since the inverse transformation is continuously 
differentiable. It remains to find it. One has
$$(I+A)^<-1>=I-A+A^2-A^3+/cdots$$
for matrices $A$ whose sup-norm $/max_<i,j>/invert<a_<ij>>$ is
less than one. Hence, writing $/hs=/sg+/eps_n$,
$$/eqalign<
/hs^<-1>&=/lp/sg/lp I+/sg^<-1>/eps_n/rp/rp^<-1>/cr
&=/lb I-/sg^<-1>/eps_n+/lp/sg^<-1>/eps_n/rp^2-/cdots/rb/sg^<-1>/cr
&=/sg^<-1>-/sg^<-1>/eps_n/sg^<-1>+O_p/lp<1/over n>/rp,/cr>$$
which implies (10.45), since $/sqrt<n>/eps_n/totop<D>S$.

Finally, under normality, 
$$/eqalign<
</rm cov>(T_<ij>,T_<kl>)&=E/sum_<u,v>/sum_<u',v'>/ssg^<iu>S_<uv>/ssg^<vj>
/ssg^<ku'>S_<u'v'>/ssg^<v'l>/cr
&=/sum_<u,v>/sum_<u'v'>/ssg^<iu>/ssg^<vj>/ssg^<ku'>/ssg^<v'l>
(/ssg_<uu'>/ssg_<vv'>+/ssg_<uv'>/ssg_<vu'>)/cr
&=/sum_<u,u'>/sum_<v,v'>(/ssg^<iu>/ssg_<uu'>/ssg^<u'k>)(/ssg^<jv>/ssg_<vv'>
/ssg^<v'l>)/cr
&/qquad+/sum_<u,v'>/sum_<u',v>(/ssg^<iu>/ssg_<uv'>/ssg^<v'l>)(/ssg^<jv>
/ssg_<vu'>/ssg^<u'k>)/cr
&=/ssg^<ik>/ssg^<jl>+/ssg^<il>/ssg^<jk>./qed/cr>$$
/medskip
We are now in a position to study the limit distributions of $/hd_0^2$ and 
$/wh/gamma_0^2$ of (10.36). Let us for simplicity consider the
</it balanced case//> where $n_1=n_2=n$. Then 
$$/lsb
	/matrix<	/sqrt<n>(/hm_1-/mu_1)/cr
			/sqrt<n>(/hs_1-/sg_1)/cr
			/sqrt<n>(/hm_2-/mu_1)/cr
			/sqrt<n>(/hs_2-/sg_2)
		>
/rsb
/totop<D>
/lsb
	/matrix<	M_1/cr
			S_1/cr
			M_2/cr
			S_2
		>
/rsb,/eqno(10.48)$$
where $(M_1,S_1)$ and $(M_2,S_2)$ are independent, $</rm cov>(M_<1,i>,
M_<1,j>)=/ssg_<1,ij>$, $S_1$, $S_2$ have covariance structure
as in the lemma above, and, finally, where
$$</rm cov>(M_<1,i>, S_<1,jk>)=E/lp X_i^<(1)>-/mu_<1,i>/rp/lb/lp X_j^<(1)>-
/mu_<1,j>/rp/lp X_k^<(1)>-/mu_<1,k>/rp-/ssg_<1,jk>/rb./eqno(10.49)$$
Now write
$$/hm_1=/mu_1+/dl_1,/,/hm_2=/mu_2+/dl_2,/,/hs_1=/sg_1+/eps_1,/,/hs_2=/sg_2+/eps_2,$$
where the error terms $/dl_1$, $/dl_2$, $/eps_1$, $/eps_2$ are of order
$O_p(n^<-1/sla2>)$, and consider $/hd_0^2$ first:
$$/eqalign<
/hd_0^2&=(/mu_1-/mu_2+/dl_1-/dl_2)'(/sg+/eps)^<-1>(/mu_1-/mu_2+/dl_1-/dl_2)/cr
&=(/mu_1-/mu_2)'(/sg+/eps)^<-1>(/mu_1-/mu_2)+2(/mu_1-/mu_2)'(/sg+/eps)^<-1>
(/dl_1-/dl_2)/cr
&/qquad+(/dl_1-/dl_2)'(/sg+/eps)^<-1>(/dl_1-/dl_2)/cr
&=(/mu_1-/mu_2)'(/sg^<-1>-/sg^<-1>/eps/sg^<-1>)(/mu_1-/mu_2)/cr
&/qquad+2(/mu_1-/mu_2)'/sg^<-1>(/dl_1-/dl_2)+O_p/lp/ndel/rp./cr>$$
Here $/eps=/halv(/eps_1+/eps_2)$. It follows already that 
$$/eqalignno<
/sqrt<n>(/hd_0^2-/dl^2)&=-(/mu_1-/mu_2)'/sg^<-1>(/sqrt<n>/eps)/sg^<-1>
(/mu_1-/mu_2)/cr
&/qquad+2(/mu_1-/mu_2)'/sg^<-1>(/sqrt<n>/dl_1-/sqrt<n>/dl_2)+O_p/lp<1/over
/sqrt<n>>/rp/cr
&/totop<D>-(/mu_1-/mu_2)'/sg^<-1>/lp/halv S_1+/halv S_2/rp/sg^<-1>(/mu_1-/mu_2)/cr
&/qquad+2(/mu_1-/mu_2)'/sg^<-1>(M_1-M_2).&(10.50)/cr>$$

Next consider $/hg_0^2$ of (10.36), which is slightly more complicated.
If $A=(a_<ij>)$ is a symmetric matrix, then one may show that
$$</partial/over/partial a_<ii>>/log/invert<A>=a^<ii>,</partial/over
/partial a_<ij>>/,/log/invert<A>=2a^<ij>,$$
if $A^<-1>=/lp a^<ij>/rp$ exists. Hence
$$/eqalign<
/hg^2_0&=4/log/invert</sg+/eps>-2/log/invert</sg_1+/eps_1>-2/log/invert</sg_2
+/eps_2>/cr
&=/gamma^2+/sum_i/lb4/ssg^<ii>/lp/halv/eps_<1,ii>+/halv/eps_<2,ii>/rp-
2/ssg^<ii>_1/eps_<1,ii>-2/ssg^<ii>_2/eps_<2,ii>/rb/cr
&/quad+/sum_<i/lt j>/lb8/ssg^<ij>/lp/halv/eps_<1,ij>+/halv
/eps_<2,ij>/rp-4/ssg_1^<ij>/eps_<1,ij>-4/ssg_2^<ij>/eps_<2,ij>/rb
+O_p/lp<1/over n>/rp,/cr>$$
so that
$$/eqalignno<
/sqrt<n>(/hg_0^2-/gamma^2)&/totop<D>/sum_i/lb/lp2/ssg^<ii>-2/ssg_1^<ii>/rp
S_<1,ii>+/lp2/ssg^<ii>-2/ssg_2^<ii>/rp S_<2,ii>/rb/cr
&/qquad+/sum_<i/lt j>/lb/lp4/ssg^<ij>-4/ssg_1^<ij>/rp S_<1,ij>
+/lp4/ssg^<ij>-4/ssg_2^<ij>/rp S_<2,ij>/rb.&(10.51)/cr>$$

The limiting variances can now in principle be read off (10.50) and
(10.51), and consistent estimates of these can be computed.
The expressions involved become fairly complex, however, and to gain insight
we proceed under the additional restriction that the underlying 
populations are </it normal//>, under which $M_1$, $S_1$, $M_2$, $S_2$ of
(10.48) are </it independent//>, and (10.46), (10.47) are valid.

Matters are technically simplified when we pass from $X^<(1)>$ and $X^<(2)>$
data to $F(X^<(1)>-/mu_1)$ and $F(X^<(2)>-/mu_1)$, where $F$ is the matrix
appearing in the simultaneous diagonalisation procedure, 
$$F/sg_1F'=I/,,/,F/sg_2F'=D=/diag(/lam_i).$$
The scheme is to evaluate answers of interest in terms of the simpler 
populations $N_d(0,I)$ and $N_d(b,D)$, where $b=F(/mu_2-/mu_1)$, and then 
transform back to $N_d(/mu_1,/sg_1)$ and $N_d(/mu_2$,$/sg_2)$.
It follows from (10.50) and (10.51) that 
$$/eqalignno<
/sqrt<n>/lp/hd_0^2-/delta^2/rp&/totop<D>/,<one>/,+/,<two>,/cr
/sqrt<n>/lp/hg_0^2-/gamma^2/rp&/totop<D>/,<three>,&(10.52)/cr
/sqrt<n>/lp/ho_0^2-/omega^2/rp&/totop<D>/,<one>/,+/,<two>/,+/,<three>,/cr>$$
say, where
$$/eqalign<
<one>/,&=-b'/lp/halv I+/halv D/rp^<-1>/lp/halv S_1^</ast >+/halv S_2^</ast >/rp
/lp/halv I+/halv D/rp^<-1>b,/cr
<two>/,&=-2b'/lp/halv I+/halv D/rp^<-1>(M_1^</ast >-M_2^</ast >),/cr
<three>/,&=2/sum_i/lb/lp<2/over 1+/lam_1>-1/rp S_<1,ii>^</ast >+
/lp<2/over1+/lam_i>-<1/over/lam_i>/rp S_<2,ii>^</ast >/rb,/cr>$$
and where $</rm VAR>(M_1^</ast >)=I$, $</rm VAR>(M_2^</ast >)=D$,
$$</rm Var>/,S_<1,ii>^</ast >=2,/,</rm Var>/,S_<1,ij>^</ast >=1,
/,</rm Var>/,S_<2,ii>^</ast >=2/lam_i^2,/,</rm Var>/,S_<2,ij>^</ast >=/lam_i/lam_j.$$
The $d(d+1)/sla 2$ different elements of each $S^</ast >$ matrix are mutually 
independent, cf.~(10.46).

Here judicious algebraic work gives 
$$/eqalign<
</rm Var>(two)&=4b'/lp/halv I+/halv D/rp^<-1>(I+D)
/lp/halv I+/halv D/rp^<-1>b/cr
&=8b'/lp/halv I+/halv D/rp^<-1>b/cr
&=8a'F'(F')^<-1>/sg^<-1>F^<-1>Fa=8/dl^2,/cr>$$
writing $a=/mu_1-/mu_2$;
$$/eqalign<
</rm Var>(one)&=</rm Var>/Bigg/lbrack/sum_ib_i^2/lp<2/over1+/lam_i>/rp^2
/lp/halv S_<1,ii>^</ast >+/halv S_<2,ii>^</ast >/rp/cr
&/qquad/qquad+2/sum_<i/lt j>b_ib_j<2/over1+/lam_i><2/over1+/lam_j>
/lp/halv S_<1,ij>^</ast >+/halv S_<2,ij>^</ast >/rp/Bigg/rbrack/cr
&=/sum_ib_i^4<4/over(1+/lam_i)^4>(2+2/lam_i^2)/cr
&/qquad/qquad+2/sum_<i/not=j>b_i^2b_j^2 4<1/over(1+/lam_i)^2>
<1/over(1+/lam_j)^2>(1+/lam_i/lam_j)/cr
&=8/sum_<i,j>b_j^2 b_j^2<1/over(1+/lam_i)^2><1/over(1+/lam_j)^2>
(1+/lam_j/lam_j)/cr
&=8/lb b'(I+D)^<-2>b/rb^2+8/lb b'(I+D)^<-1>D(I+D)^<-1>b/rb^2/cr
&=/halv(a'/sg^<-1>/sg_1/sg^<-1>a)^2+/halv(a'/sg^<-1>/sg_2/sg^<-1>a)^2,/cr>$$
since, for example,
$$/eqalign<
b'(I+D)^<-1>D(I+D)^<-1>b&=<1/over4>a'F'(F')^<-1>/sg^<-1>F^<-1>F/sg_2
F'(F')^<-1>/sg^<-1>F^<-1>Fa/cr
&=<1/over4>a'/sg^<-1>/sg_2/sg^<-1>a;/cr>$$
$$/eqalign<
</rm Var>(three)&=4/sum_i/lb/lp<1-/lam_i/over1+/lam_i>/rp^2 2
+/lp</lam_i-1/over/lam_i(1+/lam_i)>/rp^2 2/lam_i^2/rb/cr
&=16/sum_i/lp<1-/lam_i/over1+/lam_i>/rp^2/cr
&=16/Tr/lb(I-D)(I+D)^<-1>(I-D)(I+D)^<-1>/rb/cr
&=4/Tr/lb F(/sg_1-/sg_2)F'(F')^<-1>/sg^<-1>F^<-1>F(/sg_1-/sg_2)F'
(F')^<-1>/sg^<-1>F^<-1>/rb/cr
&=4/Tr/lb/sg^<-1>(/sg_1-/sg_2)/sg^<-1>(/sg_1-/sg_2)/rb./cr>$$

As a particular case it follows that when $/sg_1=/sg_2=/sg$ and
$/dl=/big/lbrace(/mu_1-/mu_2)'/sg^<-1>$/break$(/mu_1-/mu_2)/big/rbrace^<1/sla 2>$ is the
appropriate (classical Mahalanobis) distance, then
$$/sqrt<n>(/hd_0^2-/delta^2)/totop<D>N(0,/dl^4+8/dl^2),$$
or equivalently
$$/sqrt<n>(/hd_0-/dl)/totop<D>N(0,2+<1/over4>/dl^2)./eqno(10.53)$$
A confidence interval for the unknown $/dl$ with approximate level 95/%
is therefore
$$/hd_0/pm1.96/lp2+<1/over4>/hd_0^2/rp^<1/sla2>/sla/sqrt<n>,$$
and the estimated standard deviation $(2+<1/over4>/hd_0^2)^<1/sla2>/sla
/sqrt<n>$ provides information about the preciseness of $/hd_0$.

The paragraph above was phrased in terms of the direct plug-in estimator 
$/hd_0$ of (10.36). It is clear from (10.42), however, that the
bias-corrected and recommended version $/hd$ is within $O_p(<1/over n>)$
of $/hd_0$ as $n$ tends to infinity, i.e.~$/sqrt<n>(/hd-/dl)$ has the same
limit distribution as $/sqrt<n>(/hd_0-/dl)$, and the paragraph can be reread
with $/hd$ in lieu of $/hd_0$. Similarly
$$/sqrt<n>(/hg^2-/gamma^2)/totop<D>/,<three>,$$
whose variance was just found, since $/sqrt<n>(/hg^2-/hg_0^2)/totop<P>0$
using (10.43), and
$$/sqrt<n>(/ho^2-/omega^2)/totop<D>/,<one>/,+/,<two>/,+/,<three>.$$
Some further efforts produce $</rm cov>(<one>,/,<three>)$ and
$$/eqalignno<
</rm Var>(/,<one>/,+/,<two>/,+/,<three>)&=8/dl^2+/halv
(a'/sg^<-1>/sg_1/sg^<-1>a)^2/cr
&+/halv(a'/sg^<-1>/sg_2/sg^<-1>a)^2+4/Tr/lb/sg^<-1>(/sg_1-/sg_2)/sg^<-1>
(/sg_1-/sg_2)/rb/cr
&-2a'/sg^<-1>(/sg_1-/sg_2)/sg^<-1>(/sg_1-/sg_2)/sg^<-1>a.&(10.54)/cr>$$
The estimated standard deviation for the distance estimate $/ho$ of $/omega$
becomes $/wh V^<1/sla2>/sla(2/ho/sqrt<n>)$, where $/wh V$ is the plug-in
estimate of the long variance formula (10.54).
/section<10.4><Concluding remarks><>
The usual estimators (10.1) of mean vector and covariance matrix are
both sensitive to extreme data points, and because of their significance
in several classification procedures, not only classical normal-theory
based ones, robustified versions of them are recommended. The robust
estimators invented in 10.1 are very simplistic, but appear to do their
job satisfactorily. They are based on trimming away data points that
are ``obviously extreme''. Of course both ``genuine outliers'', coming from
objects outside the class studied, and ``unfortunate extremes'' that
really come from the class, can be trimmed away. This latter fact
points out that the robust $/wt/sg_0$ does not, strictly speaking, estimate
the true covariance matrix $/sg$ for the class. The classifier using
trimmed-robust estimates nevertheless usually provides the best practical
solution, when combined with an outlier criterion for future vectors,
as explained in Section~10.1.

The criterion given in 10.1 for what constitute ``obviously
extreme'' data points is only meant as a crude guideline, but again,
it is usually adequate. It can be made more accurate by including
robust estimates also for non-diagonal terms in the first-stage
covariance matrix estimator $/wt/sg_0$. This should be done in cases
of doubt.

Another approach to robust estimation of a crudely normal
density is to replace it with a $d$-dimensional $t$-distribution
(see e.g./ Berger, 1980a, p.~395), using an ad hoc-chosen value for the 
degrees of freedom, say $6$. The parameters of this distribution are then
estimated by the maximum likelihood method.

More sophisticated methods of robustly estimating means and covariance
matrices are given in Huber~(1981). A relevant comment here,
again, is that we are not specifically interested in $/mu$ and $/sg$ </it
per se>,
what matters more is an estimate of the full density.

Chapter~6 gave some methods for detecting outliers/sla extreme
observations among future samples. Some of these explicitly assumed that
the classes had normal density descriptions. It is then of obvious
importance that the parameters appearing there are robustly estimated.
A too large $/hs$-matrix for a class could mean that future vectors never
would qualify as outliers.

There are many published measures of distance between classes.
We have proposed a distance based on affinity and motivated by similarity
to Mahalanobis distances. These are intended also for quick
exploration of a chosen feature extraction method: Which classes are
well separated, and which are easily confused with others? How small
is the smallest distance? What matters in the end is error rates, and
another type of ``quick output'' could be a table of directly estimated
error rates based on model assumptions, computed by simulation. A third
alternative is perhaps the conceptually simplest one: Employ methods of
Chapter~12 to get a table of estimates for the true error rates. This
may be difficult for various reasons, however, and the sometimes simpler
methods proposed above often suffice.

Tables of inter-class distances are valuable in construction of
hierarchical classifiers and in reduction of dimensionality problems.


/chapter<11><Checking model assumptions>
/section<11.1><Introduction><>
The last chapter discussed separability measures between classes that
could be estimated using a (perhaps preliminary) training set, giving among 
other things the experimenter a chance to sort out ``safe'' classes (easily
separated from all other classes) and ``confusions'' (pairs of classes that
are too easily mixed). While there is no need to model the data structure
for a safe class in any sophisticated way, as even a ``box'' description may 
suffice for discrimination purposes, it becomes all the more important
to use a good density estimation strategy (parametric, semiparametric or 
nonparametric, and which one?) for classes in danger of being incorrectly 
classified.

Chapter 10 also provided basic summarising statistics for each class, and
that are of help in the modelling strategy. Below we go further and develop
tools that can check various assumptions that different strategies are based
on, e.g.~the normality assumption. These measures can also be computed
for each class once a set of features have been chosen, using a training set.

A perhaps ideal system could run through various sets of assumptions, 
calculating for each some measure of appropriateness, and then automatically
produce the estimated $/wh f$ using the most fitting model.
This program is not pursued to its end here, because the list of alternative 
models is not very extensive (yet), and because it is difficult to find 
satisfactory omnibus measures for adequacy. It is also important to retain
a minimum amount of human expert interaction within such a computer
expert system framework, even if it is nice to have a completely automatic
procedure as an alternative.

Let us also stress the point that even rough and ``philosophically incorrect''
model assumptions, like normality, can lead to very good classifiers. The
best quadratic rule happens to be optimal (in the class of </it all//> 
procedures) under $N_d(/mu_k,/sg_k)$ assumptions, but is in any case best 
among members in a fairly large class, and may obviously work well even if
the densities are non-normal. Thus checking the appropriateness of the best
quadratic rule is not quite the same as checking wether the normal density
appropriately fits each class. It would be taking an inappropriately high stand
if one refused to use the classifier that assigns to a candidate vector $x$ the
label $k$ that has highest $/pi_kN_d(/hm_k,/hs_k)(x)$, for the single reason
that the hypotesis of multivariate normality was rejected by some strong
test.

Not every aspect of an obtained density estimate $/wh f$ is important for 
discrimination. To reveal more precisely what matters, consider
once more the basic set-up involving $K$ classes with densities $f_k$ and 
prior probabilities $/pi_k$. The rule that minimises overall error rate
$$/eps=/sum_<k=1>^K/pi_k/,/eps(k),$$
among </it all//> rules, is the one that assigns $x$ to class $k$ whenever $x$
falls in 
$$A_k=/lb x/colon/pi_k/,f_k(x)=/max_l/pi_lf_l(x)/rb,/eqno(11.1)$$
see Chapter 1. $/eps(k)$ above is the conditional error rate for class $k$,
and becomes 
$$/eps(k)=1-P_k(A_k)=/int_</sim A_k>f_k/eqno(11.2)$$
for the optimal rule. The rule used in practice instead assigns $x$ to class
$k$ when $x$ falls in 
$$/wh A_k=/lb x/colon/pi_k/wh f_k(x)=/max_l/pi_l/wh f_l(x)/rb,/eqno(11.3)$$
and the true error rate for future objects is 
$$/etrue=/sum_<k=1>^K/pi_k/etrue(k),$$
where
$$/etrue(k)=/int_</sim/wh A_k>f_k./eqno(11.4)$$
Accordingly, the actually used rule will be close to the ideal one if only
$$/int_<A_k>f_k/approx/int_</wh A_k>f_k/eqno(11.5)$$
for each class. So $/wh f_k$ is not necessarily required to be uniformly
close to $f_k$, it matters only that $A_k/approx/wh A_k$. A good density
estimation strategy is one where the borders $/lb x/colon/pi_kf_k(x)=
/pi_l/,f_l(x)/rb$ are precisely estimated, or, rather, the part of the
borders that lie within probable land.

A method for checking the adequacy of a parametric family of distributions
could estimate the left and right hand sides of (11.5) and compare them.
While there are good methods available for estimating $/int_</wh A_k>f_k$,
cf.~Chapter 12, it is very difficult to estimate $/int_<A_k>f_k$ 
nonparametrically for anything less than large training sets, however.
One must therefore be content with reasonable checking procedures for
$f_k/approx/wh f_k$.
/section<11.2><Checking normality><>
It is useful to develop some diagnostic tools that tell one whether basic
aspects of the normal distribution assumption are reflected in the data.
Such tools should also indicate specific departures from normality,
if present, so that more adequate models can be guessed at. The value 
of tests for parametric assumptions may primarily lie in such suggestions
as to what to do next, compare the general remarks in the introductory
section, where it is warned that only border-defining consequences 
of a modelling strategy are important for the classification problem.
/subsection<11.2.A><Distribution of distance and direction>
The two perhaps most important aspects of the normality assumption are
its consequences for the distribution of </it distance//> from the centre
and the distribution of </it direction//> from the centre. Specifically,
if $X$ is drawn from a $N_d(/mu,/sg)$, then $(X-/mu)'/sg^<-1>(X-/mu)$ is
distributed as a chi-square with $d$ degrees of freedom, and the
transformed variable $Y=/sg^<-1/sla2>(X-/mu)$, which is $N_d(0,I)$,
is isotropically distributed in $/rr^d$, i.e.~all directions are equally
likely.

These aspects can be checked for. Introduce the ``squared radii''
$$R_t^2=(X_t-/hm)'/hs^<-1>(X_t-/hm),/ t=1/upto n./eqno(11.6)$$
These are nearly distributed as $/chi_d^2$, and are nearly independent.
More precisely, because $X_t-/hm$ is $N_d/lp0,<n-1/over n>/sg/rp$ and
$/hs=<A/over m>$ with $A/sim</rm Wishart>_d(/sg,m),/,m=n-1$, one has that
$(<n/over n-1>)^<1/sla2>(X_t-/hm)'(<A/over m>)^<-1>
(<n/over n-1>)^<1/sla2>(X_t-/hm)$ is $T^2(d,m)$, 
cf.~Mardia, Kent and Bibby (1979, Ch.~3), i.e., by the connection
to Fisher from Hotelling (op.cit., p.~74),
$$/eqalignno<
R_t^2&/sim<n-1/over n>T^2(d,m)=<n-1/over n><dm/over m-(d+1)>F(d,m-(d+1))/cr
&=<n-1/over n><n-1/over n-d>dF(d,n-d),&(11.7)/cr>$$
which is close to the $/chi_d^2$ for moderate and large $n$.

Graphical checks for the distribution of distances from centre can now be 
devised. The histogram of $R_t^2$ values can be compared to the density
implied by (11.7), or a $QQ$ plot, consisting of the $n$ points 
$(G_<n,d>^<-1>(<t/over n+1>),R_<(t)>^2)$, where
$R_<(1)>^2/lt/cdots/lt R_<(n)>^2$ and $G_<n,d>$ is the cumulative
distribution of (11.7), can be inspected. The $QQ$ plot should give
approximately the straight line $R_<(t)>^2=G_<n,d>(<t/over n+1>)$.
Again it is useful with some standard quantitative outfit, however. We
have found it convenient and useful to display the ten relative
frequencies $M_j/sla n,/ j=1/upto 10$, of squared radii falling in the ten
decile cells $(c_<j-1>,c_j)$, where $G_<n,d>(c_j)=j/sla10$ (and
$c_0=0,/,c_<10>=/infty$). These relative frequencies should all be
approximately equal to $1/over 10$, under the normality assumption.

Note here that the Karl Pearson goodness of fit statistic for these cell
frequencies, $10n/sum_<j=1>^<10>(<M_j/over n>-<1/over10>)^2$, has a 
different limiting distribution that the expected $/chi_9^2$,
since the normalised empirical process
$$/sqrt<n>/lsb<1/over n>/sum_<t=1>^n I/lb R_t^2/le u/rb-G_<n,d>(u)/rsb
/eqno(11.8)$$
has a different (and more complicated) limit process than the perhaps 
expected $W^0(/Gamma_d(u))$, $W^0$ being the Brownian bridge, cf.~11.2.B
below.

We develop a rigorous test for this distance-from-centre aspect of normality
in 11.2.B, using a split-sample version of the empirical process above,
but note here that studying the structure of the ten cell frequencies
$M_j/sla n$, compared to $<1/over10>/upto<1/over10>$, may be more
informative than the observed significance value.

Let us next supply a method of checking the distribution of directions from 
centre. Transform data to pseudo-observations 
$$Y_t=/hs^<-1/sla2>/lp X_t-/hm/rp,/eqno(11.9)$$
$t=1/upto n$, cf.~5.4.C. These have $<1/over n>/sum_<t=1>^n Y_t=0$
and $<1/over n-1>/sum_<t=1>^n Y_tY_t'=I$, i.e.~pretend to be like 
$N_d(0,I)$ variables. If they were, then the observations, when projected 
into the unit sphere surface, would be uniformly distributed on this
surface. Define
$$Z_t=Y_t/sla/inpar<Y_t>,$$
therefore, where $/inpar<Y_t>^2=(X_t-/hm)'/hs^<-1>(X_t-/hm)=R_t^2$.
Consider the matrix
$$M(n)=<1/over n>/sum_<t=1>^n Z_tZ_t'=<1/over n>/sum_<t=1>^nY_tY_t'/sla R_t^2
/eqno(11.10)$$
and its $d$ eigenvalues $/hm_1/ge/cdots/ge/hm_d$. These have sum one, and 
should all be about equal to $<1/over d>$ under the assumption that the 
$Z_t$'s are uniformly distributed on the unit surface, and in fact
$$<1/over2>nd(d+2)/sum_<i=1>^d/lp/hm_i-<1/over d>/rp^2/totop<D>
/chi^2_<d(d+1)/sla 2-1>.$$
Even taking the noise involved in the initial transformation (11.9) into
account it still holds that the $/hm_i$'s should be about equal, and 
a marked departure from this indicates a preference of some direction
of $Y$ over others.

There are many published tests for multivariate normality, some looking
for heavy tails and others looking at skewness and kurtosis, see e.g.~Mardia
(1974). We re-emphasise the point, however, that moderate deviations from
normality are not necessarily  dramatic in their consequences for the 
performance of the normal-theory based classifier.
/subsection<11.2.B><A rigorous test based on distances> The reason
why the empirical process (11.8) is difficult to study is that the squared
radii $R_t^2$ of (11.6) are dependent, via $/hm$ and $/hs$.
Durbin (1973) obtains limit process results for some empirical 
processes when parameters are estimated, but these can not be applied
here since the noise is in the transformation to $R_t^2$. Koziol (1982)
has obtained the limiting normal process for $/sqrt<n>/lsb<1/over n>/sum_1^n
I/lb R_t^2/le u/rb-/Gamma_d(u)/rsb$; it has covariance structure
$$/Gamma_d(u)/lb1-/Gamma_d(v)/rb-<2/over d>u/gamma_d(u)
v/gamma_d(v), u/le v.$$

Below we manage to find the limit distribution of a similar process, using a 
split-sample device.

Divide the sample at hand into two sets, say $X_1^<(0)>/upto X_<n_0>^<(0)>$
and $X_1/upto X_n$. The first part is used to get estimates $/hm_0$ and 
$/hs_0$, say. Define next, for the second part,
$$R_</ast t>^2=(X_t-/hm_0)'/hs_0^<-1>(X_t-/hm_0)./eqno(11.11)$$
These are independent, given the past observations,
with common distribution
$$G_</ast >(u)=Pr_</ast >/lb/lp X-/hm_0/rp'/hs_0^<-1>/lp X-/hm_0/rp/le u/rb,
/eqno(11.12)$$
the $/ast $ indicating probabilities in the conditional framework given
$X_1^<(0)>/upto X_<n_0>^<(0)>$. The empirical process we shall study is
$$Z_n(u)=/sqrt<n>/lsb<1/over n>/sum_<t=1>^n I/lb R_</ast t>^2/le u/rb-
/Gamma_d(u)/rsb./eqno(11.13)$$

Some remarks are in order here before we proceed. $n_0$ and $n$ above can be 
chosen by the experimenter; we suggest taking $n_0$ to be $1/over2$ or 
$3/over4$ of the whole sample. The division itself is random,
i.e.~$X_1^<(0)>/upto X_<n_0>^<(0)>$ are supposed to form a random
sample from the whole initial sample. This may appear somewhat artificial,
and indeed may be for some applications, but is not so unnatural in the 
pattern recognition world, where it is usual to have two-stage training sets,
and perhaps designed just for this type of use, with ``training set'' and
``test set''. In other cases an initial set is provided first,
and examination of this one, and perhaps of interclass distances etc.~as 
outlined in Chapter 10, dictates the size of the second phase set. It
should also be pointed out that this two-stage thinking is but an
example of the cross-validation practice and philosophy that are
essential for many other testing purposes too, cf.~Chapter 12. 
A final comment is that it is usual and sound practice to ``reverse sides''
afterwards, i.e.~train on the test set and test on the training set,
or to cycle the subsets in other ways, as in the ``leave a quarter out''
method that tests each of four equally sized subsets using sample estimates 
from the appropriate complementary $<3/over4>$ sample.

The proof of the following theorem happens to require several pages, but after
that a class of rigorous tests for the normality hypothesis is derived as a 
corollary.
/Btheorem<> Write
$$Z_n(u)=Z_<n,1>(u)+Z_<n,2>(u),$$
where
$$/eqalignno<
Z_<n,1>(u)&=/sqrt<n>/lsb<1/over n>/sum_<t=1>^n I/lb R_</ast t>^2/le u/rb-G_/ast (u)
/rsb,/cr
Z_<n,2>(u)&=/sqrt<n>/lb G_/ast (u)-/Gamma_d(u)/rb.&(11.14)/cr>$$
Then $(Z_<n,1>,Z_<n,2>)/totop<D>(Z_1,Z_2)$ in the function space
$D_0/lsb0,/infty/rsb/times D_0/lsb0,/infty/rsb$, where the elements of 
$D_0/lsb0,/infty/rsb$ are those right continuous functions on $/lsb0,/infty/rsb$
that vanish at zero and infinity and have left hand limits. Here
$$Z_1=W^0/lb/Gamma_d(.)/rb$$
where $W^0$ is the Brownian bridge,
$$Z_2=h(.)N$$
where $N$ is standard normally distributed and independent of $Z_1$, and
$$h(u)=/halv/lb/Gamma_d(u)-/Gamma_<d+2>(u)/rb/sqrt<c>(2d)^<1/sla2>.
/eqno(11.15)$$
c is the limit of the sample size quotient $n/sla n_0$. Finally
$$Z_n(.)/totop<D>Z_1+Z_2=W^0/lb/Gamma_d(.)/rb+h(.)N,/eqno(11.16)$$
a zero mean Gaussian process with covariance  function 
$$</rm cov>/lb Z(u),Z(v)/rb=/Gamma_d(u)/lb1-/Gamma_d(v)/rb
+h(u)/,h(v),/,u/le v./eqno(11.17)$$
/Etheorem
/Bproof Because $R_</ast 1>^2/upto R_</ast n>^2$ are i.i.d.~with distribution
$G_/ast $ in the conditional framework $Z_<n,1>$ is close to $W^0/lb G_/ast (.)/rb$ in
distribution, by the classical Doob-Donsker theorem for empirical processes,
see Billingsley (1968, Section 13). But
$$G_/ast (u)=G(u;/hm_0,/hs_0),$$
say, where
$$G(u;/tm,/ts)=Pr_</mu,/sg>/lb/lp X-/tm/rp'/ts^<-1>/lp X-/tm/rp/le u/rb,
/eqno(11.18)$$
and $/hm_0$ and $/hs_0$ converge to the true parameters $/mu$ and $/sg$ (almost
surely). Hence
$$G_/ast (u)/totop</rm a.s.>G(u;/mu,/sg)=/Gamma_d(u),$$
and it follows that $Z_<n,1>$ becomes close to $W^0/lb/Gamma_d(.)/rb$
in distribution. That
$$Z_<n,1>/totop<D>W^0/lb/Gamma_d(.)/rb/eqno(11.19)$$
can indeed be proved rigorously using random time change arguments
as in Billingsley (1968, Section 17).

It is now possible to construct tests based on $Z_<n,1>$ alone, for 
example,
$$/sqrt<n>/,</rm sup>_<u/ge0>/invert<<1/over n>/sum_<t=1>^n I/lb R_</ast t>^2/le 
u/rb-G_/ast (u)>/totop<D>/inpar<W^0>=</rm sup>_<0/le u/le1>/invert<W^0(u)>,$$
and
$$/sum_<j=1>^m<(M_j-nq_j)^2/over nq_j>/totop<D>/chi_<m-1>^2,$$
where $q_j$ is the $G_/ast $-probability of the cell $(c_<j-1>,c_j)$ and
$M_j$ is the number of observations in this cell; 
$0=c_0/lt/cdots/lt c_m=/infty$. These are awkward tests, however, since no
explicit expressions for the observed $G_/ast (u)$ is available. Since
$G_/ast $ is close to the well known $/Gamma_d$, more natural tests result if
$G_/ast $ is replaced by $/Gamma_d$ above. This requires the limit distribution
for the original $Z_n=Z_<n,1>+Z_<n,2>$, and we must study $Z_<n,2>$.

A first order Taylor expansion gives 
$$/eqalignno<
G_/ast (u)&/doteq/Gamma_d(u)+/sum_i</partial/over/partial/mu_i>G(u;/mu,/sg)
(/hm_<0,i>-/mu_i)/cr
&+/sum_i</partial/over/partial/ssg_<ii>>G(u;/mu,/sg)
(/wh/ssg_<0,ii>-/ssg_<ii>)/cr
&+/sum_<i/lt j></partial/over/partial/ssg_<ij>>G(u;/mu,/sg)
(/wh/ssg_<0,ij>-/ssg_<ij>),&(11.20)>$$
and shows, after properly tending to remainder terms, that 
$/sqrt<n>/lb G_/ast (u)-/Gamma_d(u)/rb$ has a limiting normal distribution.
It turns out to be technically convenient to evaluate these difficult
derivatives of (11.18) for the simpler case $/mu=0,/,/sg=I$ first,
and then generalise afterwards. 

Thus consider
$$H(u;/dl,/eps)=Pr_<0,I>/lb(Y-/dl)'(I+/eps)^<-1>(Y-/dl)/le u/rb$$
where $Y$ is $N_d(0,I)$, and $/dl=(/dl_i)$ and $/eps=(/eps_<ij>)$
are close to zero, $/eps_<ij>=/eps_<ji>$. The task is to find
$</partial/over/partial/dl_i>H(u;/dl,/eps)/mid_<(0,0)>$ and
$</partial/over/partial/eps_<ij>>H(u;/dl,/eps)/mid_<(0,0)>$.
This can be managed using the following trick. Write
$$/eqalign<
H(u;/dl,/eps)&=/int_<y/colon(y-/dl)'(I+/eps)^<-1>(y-/dl)/le u>g(y)/,dy/cr
&=/int_<z/colon z'z/le u>g/lp/dl+(I+/eps)^<1/sla2>z/rp
/invert<(I+/eps)^<1/sla2>>/,dz,/cr>$$
where $g(y)=(2/pi)^<-d/sla2>/exp(-/halv/sum_<i=1>^d y_i^2)=
/phi(y_1)/cdots/phi(y_d)$ is the $N_d(0,I)$ density, and where the 
transformation $z=(I+/eps)^<-1/sla2>(y-/dl)$ has been used.
Derivatives can now be evaluated by differentiating under the integral 
sign for each fixed $u$.

Firstly,
$$/eqalignno<
</partial/over/partial/dl_i>H(u;/dl,/eps)/mid_<(0,0)>&=
</partial/over/partial/dl_i>/mid_<(0,0)>/int_<z'z/le u>g(/dl+z)/,dz/cr
&=/int_<z'z/le u>-z_i g(z)/,dz=0.&(11.21)>$$
(This could also have been found using 
$$/eqalign<
H(u;/dl,0)&=Pr/lb/chi_d^2/lp/sum_<i=1>^d/dl_i^2/rp/le u/rb/cr
&/doteq/Gamma_d(u)+/halv/lb/Gamma_<d+2>(u)-/Gamma_d(u)/rb
/sum_<i=1>^d/dl_i^2.)/cr>$$
Secondly, using the fact that 
$$</partial/over/partial a_<ii>>/invert<A>=/invert<A>a^<ii>,/,
</partial/over/partial a_<ij>>/invert<A>=2/invert<A>a^<ij>$$
for an invertible symmetric matrix $A$, where $A^<-1>=(a^<ij>)$, it 
holds that
$$/invert<I+/eps>/doteq/invert<I>+/sum_i/eps_<ii>,$$
i.e.~$/invert<I+/eps>^<1/sla2>/doteq1+/halv/sum_i/eps_<ii>$ for symmetric 
matrices $/eps$ with elements close to zero. This implies
$$/eqalign<
H(u;0,/eps)&=/int_<z'z/le u>g/lp(I+/eps)^<1/sla2>z/rp
/invert<I+/eps>^<1/sla2>/,dz/cr
&/doteq/int_<z'z/le u>g/lp(I+/eps)^<1/sla2>z/rp/lb1+/halv Tr(/eps)/rb/,dz/cr
&=/lb1+/halv Tr(/eps)/rb/int_<z'z/le u>g(z)/exp/lp-/halv z'/eps z/rp/,dz/cr
&/doteq/lb1+/halv Tr(/eps)/rb/lb/Gamma_d(u)-/halv/int_<z'z/le u>
z'/eps zg(z)/,dz/rb,/cr>$$
and therefore
$$/eqalign<
H(u;0,/eps)-H(u;0,0)&/doteq/halv Tr(/eps)/Gamma_d(u)-/halv
/sum_<i,j>/eps_<ij>/int_<z'z/le u>z_iz_jg(z)/,dz/cr
&/qquad-/halv/sum_i/eps_<ii>/,/halv/sum_<i,j>/eps_<ij>/int_<z'z/le u>
z_iz_jg(z)/,dz/cr
&/doteq/halv/sum_i/eps_<ii>/Gamma_d(u)-/halv/sum_i/eps_<ii>/int_<z'z/le u>
z_i^2g(z)/,dz./cr>$$
All this goes to show that
$$/eqalignno<
</partial/over/partial/eps_<ii>>H(u;/dl,/eps)/mid_<(0,0)>
	&=/halv/Gamma_d(u)-/halv b,/cr
</partial/over/partial/eps_<ij>>H(u;/dl,/eps)/mid_<(0,0)>
	&=0,/,i/not=j,&(11.22)/cr>$$
where $b=/int z_i^2g(z)I/lb z'z/le u/rb/,dx$. But
$$/eqalign<
db&=/int/sum_iz_i^2g(z)I/lb/sum_iz_i^2/le u/rb/,dz/cr
&=E/chi_d^2I/lb/chi_d^2/le u/rb/cr
&=/int_0^ux/gamma_d(x)/,dx=d/Gamma_<d+2>(u),/cr>$$
writing $/gamma_d(x)=x^<d/sla2-1>e^<-x/sla2>/sla/lb2^<d/sla2>
/Gamma(d/sla2)/rb$ for the $/chi_d^2$ density, i.e.
$$</partial/over/partial/eps_<ii>>H(u;/dl,/eps)/mid_<(0,0)>=/halv/lb/Gamma_d(u)-
/Gamma_<d+2>(u)/rb.$$

These efforts make it possible for us to find the coefficients appearing in 
the Taylor expansion (11.20). Going back to the $(/mu,/sg)$ probability
(11.18), let $Y=F(X-/mu)$, where $F/sg F'=I$ (i.e.~$F=/sg^<-1/sla2>$).
For small $/dl$ and $/eps$, introduce $a=F/dl$ and $/phi=F/eps F'$. 
Then
$$/eqalign<
/lp X-(/mu+/dl)/rp'&(/sg+/eps)^<-1>/lp X-(/mu+/dl)/rp/cr
&=(Y-a)'(F')^<-1>/lb F^<-1>(F')^<-1>+F^<-1>/phi(F')^<-1>/rb^<-1>F^<-1>(Y-a)/cr
&=(Y-a)'(I+/phi)^<-1>(Y-a),/cr>$$
and
$$/eqalign<
G(u;/mu+/dl,/sg+/eps)&=Pr_</mu,/sg>/lb/lp X-(/mu+/dl)/rp'(/sg+/eps)^<-1>
/lp X-(/mu+/dl)/rp/le u/rb/cr
&=Pr_<0,I>/lb(Y-a)'(I+/phi)^<-1>(Y-a)/le u/rb/cr
&/doteq Pr_<0,I>/lb Y'Y/le u/rb+/sum_i/halv/lb/Gamma_d(u)-/Gamma_<d+2>(u)/rb
/phi_<ii>/cr
&=/Gamma_d(u)+/halv/lb/Gamma_d(u)-/Gamma_<d+2>(u)/rb/,Tr(/phi),/cr>$$
because of (11.18) and (11.19). The consequence is that
$$G(u;/hm_0,/hs_0)/doteq/Gamma_d(u)+/halv/lb/Gamma_d(u)-/Gamma_<d+2>(u)/rb/,
Tr/lb F/lp/hs_0-/sg/rp F'/rb,$$
i.e., finally,
$$/eqalign<
Z_<n,2>(u)&=/sqrt<n>/lb G_/ast (u)-/Gamma_d(u)/rb/cr
&/doteq/halv/lb/Gamma_d(u)-/Gamma_<d+2>(u)/rb/,Tr/lb/sg^<-1>/sqrt<n>
/lp/hs_0-/sg/rp/rb./cr>$$
If $n/sla n_0=c$, and $n$ and $n_0$ tend to infinity keeping their
relationship stable, then 
$$Z_<n,2>(u)/totop<D>/halv/lb/Gamma_d(u)-/Gamma_<d+2>(u)/rb/sqrt<c>/,
Tr/lp/sg^<-1>S_0/rp,$$
where $S_0$ is the limiting random matrix of $/sqrt<n_0>/lp/hs_0-/sg/rp$.
Under normality the elements of $S_0$ have $/cov(S_<0,ij>,S_<0,kl>)$
$=/ssg_<ik>/ssg_<jl>+/ssg_<il>/ssg_<jk>$, see the lemma in 10.3.B. Hence
$$/eqalign<
</rm Var>/lp/sg^<-1>S_0/rp&=/sum_<i,j>/sum_<k,l>/ssg^<ij>/ssg^<kl>
(/ssg_<ik>/ssg_<jl>+/ssg_<il>/ssg_<jk>)/cr
&=Tr(/sg^<-1>/sg/sg^<-1>/sg)+Tr(/sg^<-1>/sg/sg^<-1>/sg)=2d./cr>$$
$/lb Z_<n,2>/rb$ can also be seen to be a tight sequence of elements in the
space $D_0/lsb0,/infty/rsb$, so that we have actually proved
$$Z_<n,2>/totop<D>h(.)N(0,1),/eqno(11.23)$$
where $h$ is as in (11.15).

The proof of the theorem will be completed when it is shown that $Z_1$ and
$Z_2$ are independent, i.e.~that $Z_<n,1>$ and $Z_<n,2>$ become independent
in the limit. (That $X_n+Y_n/totop<D>X+Y$ in $D_0/lsb0,/infty/rsb$ when
$(X_n,Y_n)/totop<D>(X,Y)$ in the product space holds provided
$X$ and $Y$ have continuous sample paths, as here, since the function
$(x,y)/to x+y$ is continuous in the Skohorod topology on this subset.)

The arguments of the proof can be copied to show that 
$/alpha Z_<n,1>(u)+/beta Z_<n,2>(v)$ at any rate has a normal limit, 
so that at least $(Z_<n,1>,Z_<n,2>)/totop<D>$ some $(Z_1,Z_2)$ that is 
bi--Gaussian, and it is sufficient to prove $EZ_1(u)Z_2(v)=0$, for all
$u$ and $v$.

But it is easily seen that 
$$E_/ast Z_<n,1>(u)Z_<n,2>(v)=0,$$
using $/ast $ to denote the conditional framework again, and in which the whole 
$Z_<n,2>$ process is known. Hence $EZ_<n,1>(u)Z_<n,2>(v)=0$ for all $n$,
and that $EZ_1(u)Z_2(v)=0$ too follows rigorously if 
$/lb Z_<n,1>(u)Z_<n,2>(v)/rb$ is uniformly integrable, which
it is;
$$/eqalign<
EZ_<n,1>(u)^2Z_<n,2>(v)^2&=EG_/ast (u)/lb1-G_/ast (u)/rb Z_<n,2>(v)^2/cr
&/le En/lb G_/ast (v)-/Gamma_d(v)/rb^2,/cr>$$
and an exact version of the Taylor expansion (11.20) guarantees a bound
on the last expression. /Eproof
We are finally in a position to derive test statistics for normality
based on the squared radii (11.11). In fact a multitude of tests can be
constructed based on the limit process result (11.16), an example of the
``omnibus'' type (consistent in all alternatives) being a
Kolmogorov-Smirnov like
$$D_n=</rm sup>_<u/ge0>/invert<<1/over n>/sum_<t=1>^n I/lb R_</ast t>^2/le
u/rb-/Gamma_d(u)>,$$
for which $/sqrt<n>D_n/totop<d></rm sup>_<0/le v/le1>/invert<W^0(v)+h/lb
/Gamma_d^<-1>(v)/rb N>$.
The distribution of the latter is complicated but can be obtained by a
laborious conditioning argument. We shall be content here to contruct
$/chi^2$ type tests.

Let then $0=c_0/lt c_1/lt/cdots/lt c_m=/infty$ be a division of the half
line into $m$ cells $(c_<j-1>,c_j/rbrack$, let
$p_j=/Gamma_d(c_<j-1>,c_j/rbrack$ be the (approximate) probability that
a $R_/ast ^2$ falls in cell no.~$j$, and let $M_j/sla n$ be the observed
frequency in it. By the theorem
$$/eqalignno<
/sqrt<n>/lp<M_j/over n>-p_j/rp&/totop<D>K_j=(Z_1+Z_2)(c_j)-(Z_1+Z_2)
(c_<j-1>)/cr
&=W^0/lb/Gamma_d(c_j)/rb-W^0/lb/Gamma_d(c_<j-1>)/rb+/lb h(c_j)-
h(c_<j-1>)/rb N,&(11.24)/cr>$$
and the corresponding $m$-vector converges to $K=(K_1/upto K_m)'$.
(11.17) can be used to obtain 
$$/cov(K_j,K_l)=p_j(/dl_<j,l>-p_l)+w_jw_l,/eqno(11.25)$$
where
$$w_j=h(c_j)-h(c_<j-1>)=/halv/lsb p_j-/lb/Gamma_<d+2>(c_j)-/Gamma_<d+2>
(c_<j-1>)/rb/rsb/lp<n/over n_0>/rp^<1/sla2>/sqrt<2d>./eqno(11.26)$$

$K$ has a singular covariance matrix since $/sum_<j=1>^m K_j=0$. The
shortened $K_0=(K_1/upto K_<m-1>)'$ has 
$$T=A+w_0w_0',$$
say, where $A=/lp p_j(/dl_<jl>-p_l)/rp$ and $w_0=(w_1/upto w_<m-1>)'$.
By a matrix identity sometimes attributed to Bartlett,
$$T^<-1>=A^<-1>-A^<-1>w_0w_0'A^<-1>/sla/lp1+w_0'A^<-1>w_0/rp.$$
Here $A^<-1>=(<1/over p_j>/dl_<jl>+<1/over p_m>)$, and
$$/eqalign<
K_0'A^<-1>K_0=/sum_<j=1>^m K_j^2/sla
p_j&,/,w_0'A^<-1>w_0=/sum_<j=1>^mw_j^2/sla p_j,/cr
w_0'A^<-1>K_0&=/sum_<j=1>^mw_jK_j/sla p_j/cr>$$
are readily demonstrated. Hence
$$/eqalign<
K_0'(A+w_0w_0')^<-1>K_0&=K_0'A^<-1>K_0 - <K_0'A^<-1>w_0w_0'A^<-1>K_0/over
1+w_0'A^<-1>w_0>/cr
&=/sum_<j=1>^mK_j^2/sla p_j-/lp/sum_<j=1>^m<w_j/over p_j>K_j/rp^2/sla
/lp1+/sum_<j=1>^mw_j^2/sla p_j/rp^2,/cr>$$
and this quadratic form is $/chi_<m-1>^2$ distributed.
It is also the limit in distribution of the test statistic
$$P=/sum_<j=1>^m<(M_j-np_j)^2/over np_j>-/lb/sum_<j=1>^m<w_j/over p_j>
/sqrt<n>/lp<M_j/over n>-p_j/rp/rb^2/sla/lb1+/sum_<j=1>^mw_j^2/sla
p_j/rb./eqno(11.27)$$
The test that rejects the normality assumption if $P$ is greater than
the upper $/alpha$-point of the $/chi_<m-1>^2$ has (asymptotically)
level $/alpha$.

If for example $c_1/upto c_9$ are the deciles of the chi-square,
i.e.~$c_j=/Gamma_d^<-1>/lp<j/over10>/rp$, and $M_j$ is the observed
number of $(X_t-/hm_0)'/hs_0^<-1>(X_t-/hm_0)$ in
$(c_<j-1>,c_j/rbrack$, then
$$P=10n/sum_<j=1>^<10>/lp<M_j/over n>-<1/over10>/rp^2- <100n/over
1+10/sum_<j=1>^<10> w_j^2>/lb/sum_<j=1>^<10> w_j/lp<M_j/over
n>-<1/over10>/rp/rb^2/eqno(11.28)$$
is a ``corrected Pearson statistic'' for goodness of fit, and can be
referred to a $/chi_9^2$ to test normality.
/section<11.3><Checks for other density estimation schemes><>
If some parametric model other than the normal is being tried out,
parameter estimates should of course be included in the general class
description output (see Chapter 10), but also quantitative diagnostics
for the adequacy of the model used, as in the normal case studied above.
Such are not provided here, but are in principle easy to invent for a
given parametric model. 

We go on to study checking procedures tailored to some of the special
and less parametric methods described in Section 5.4.

Consider for example the multiplicative cosine expansion estimator of
5.4.A. There new variables $V_t=/wh P(X_t-/hm)$ were created, and the
components were fitted to a cosine expansion density. The crucial and
bold assumption underlying the method is that the transformed components
of $V_t$ become (nearly) independent. More precisely, if $X$ is an
upcoming vector from the same class it is assumed that the $d$
components of $V=/wh P(X-/hm)$ are independent.

To check this assumption and also the ``total fit'' of the model without
new vectors one can resort to the ``split sample trick'', see 11.2.B and
the discussion that follows equation (11.13). Divide the sample in two
and write it as $X_1^<(0)>/upto X_<n_0>^<(0)>$, $X_1/upto X_n$, with
$n_0$ and $n$ adding to the size of the full sample and not necessarily
of equal size. Obtain estimates $/hm_0,/hs_0,/wh P_0$ based on the first
set of $X$'s and consider ``new'' transformed variables
$$V_t=/wh P_0(X_t-/hm_0),/ t=1/upto n.$$
These have the same distribution and are independent, given the past
initial set. One should now draw histograms and scatterplots for the
coordinate components of the $V_t$'s. Tests for independence could be
run, as well as tests of fit to the individual cosine expansion density
estimates, involving coefficients $c_i(j)$ estimated from the first set
only, cf.~(5.70)---(5.74). In particular the correlation coefficients
could be studied, and a nice output could include the estimated smooth
density curve for component no.~$i$ along with the corresponding
histogram.

A pleasant point to make is that the orthogonal expansion procedures
partly take care of their own model-checking through the stopping rule.
The estimator $/wh g_m=g_0/sum_<j=0>^m/wh c(j)/psi_j$ estimates an
approximation $g_m=g_0/sum_<j=0>^m c(j)/psi_j$ to the density $g$, and
the procedure for choosing which coefficients to include can be seen as
a dynamic goodness-of-fit method where the hypotheses $c(j)=0$ are
scrutinised, using $/wh c(j)$ and a variance estimate.

Consider the third order corrected normal density estimate, for example.
It is equal to (5.88) provided $T$ of (5.89) is positive, and is simply
$N_d(/hm,/hs)(x)$ if $T$ is negative. Built into this procedure is a
test for normality, with certain deviations from it in mind, viz.~$H_0:$
``$/gamma_i=0$, $/dl_<i,j>=0$, $/eps_<i,j,l>=0$, for every $i$, $j$,
$l$''. $H_0$ is true under normality, and the $T$ diagnostic weighs
the evidence $/wh/gamma_i$, $/hd_<i,j>$, $/wh/eps_<i,j,l>$ against it.
/section<11.4><Concluding remarks><> There is a bewilderingly long list of 
published tests for multivariate normality. As usual there is no single 
best method; each one has its individual advantages and disadvantages. 
While it is a simple fact of life that the multivariate normal model can 
never be completely correct, it is also true that methods based on it can 
be extremely useful. This is a particularly pertinent comment in the 
discriminant analysis context, where even a rough description of a class 
distribution suffices, if the class in question is reasonably well 
separated from the others. So a rejected hypothesis of normality is not a 
decree stating that normal-based procedures cannot be used. On the other 
hand, if the data support the normal model for classes then there may be 
little gain in using more sophisticated procedures.

We have found it convenient and useful to split the normality question 
into two, one testing distance from centre and one testing direction from 
centre. Studying not only the finally calculated significance values but 
also other informative output from the procedures, say the ten decile 
ratios $M_j/sla n$ mentioned in Section 11.2.B, makes it possible to 
pinpoint what is wrong with normality in the case at hand, and often 
suggests what kind of model or method will work better.

The rigorous split-sample test for normality based on distances is new. 
Related procedures, most often advocated for tentative use (as opposed to 
rigorous significance testing), are discussed in Andrews, 
Gnanadesikan, and Warner (1973) and in Royston (1983).

One may of course develop goodness of fit tests, including graphical 
procedures, for many other models than the multinormal. Some general 
comments were given in 11.3. One can in particular put up such checking 
procedures for the important discrete-times-normal models studied in 
Chapter 9.


/chapter<12><Estimating error rates>
Consider once more the usual framework involving $K$ classes,
whose class distributions will be denoted by $F_1/upto F_K$ in the 
following, and with prior probabilities $/pi_1/upto/pi_K$. Let 
$/hC(x)=/hC(x;/,z)$ be some classification procedure under consideration,
based on the training set
$$z=/lb x_1^<(k)>/upto x_</nk>^<(k)>;/ k=1/upto K/rb/eqno(12.1)$$
that happened to be available. $/hC(x)$ assigns one of the values
$1/upto K, D, /out$ (say) to the feature vector $x$.

The present chapter is concerned with the basic problem of assessing
the performance of the classifier. The best measures of quality,
CPU-considerations excluded, are the various error and success rates
involved. We shall in particular take interest in 
$$/pmc_k=/Pr_<F_k>/lb/hC/lp X_0;/,z/rp/not= k/rb,/ k=1/upto
K./eqno(12.2)$$
$/pmc_k$ is the probability of misclassifying an upcoming object from
class $k$, and is also termed the (conditional) </it error rate for
class $k$//>. The unconditional or </it total error rate//> is
$$/pmc=/Pr_F/lb/hC/lp X_0;/,z/rp/not= </rm class>/lp X_0/rp/rb
=/sum_<k=1>^K/pi_k/pmc_k./eqno(12.3)$$ 
In (12.2) $X_0$ is drawn from $F_k$ and in (12.3) $X_0$ comes from the
mixture distribution $F=/sum_<k=1>^K/pi_kF_k$. Note that $/pmc_k$ and
$/pmc$ are defined conditionally on the actually observed training set
$z$.

The sections below give methods for estimating $/pmc_k$ and $/pmc$,
using cross-val/-i/-da/-tion and resampling procedures. We concentrate on
$/pmc_k$ and $/pmc$ in our discussion, but the reasoning and methods
apply equally well to the </it success rates//>
$$/eqalign<
/pcc_k&=/Pr_<F_k>/lb/hC/lp X_0;/,z/rp=k/rb=1-/pmc_k,/cr
/pcc&=/Pr_F/lb/hC/lp X_0;/,z/rp=</rm class>/lp X_0/rp/rb=
/sum_<k=1>^K/pi_k /pcc_k,/cr>$$
the </it doubt rate//>
$$/pd_k=/Pr_<F_k>/lb/hC/lp X_0;/,z/rp=D/rb,$$
the similarly defined </it outlier rate//>, and the </it partial error rates//> 
$$/pmc(k/to l)=/Pr_<F_k>/lb/hC/lp X_0;/,z/rp=l/rb.$$
/section<12.1><The leave-one-out method><>
/subsection<12.1.A><Apparent error rates and cross validation>
An obvious possibility is to run the classifier on the training set,
and count errors. Thus the </it apparent error rates//>
$$ /apmc_k=<1/over /nk>/sum_<j=1>^</nk>I/lb/hC/lp /xkj ;/,z/rp/not= k/rb/eqno(12.4)$$
can be computed. These tend to be too optimistic, i.e.~smaller than
the true error rates $ /pmc_k$, since the same data set $z$ has been
used both to construct and evaluate $/hC( x;/,z)$.
We would rather see $/hC(./,;/,z)$ at work on a set of new vectors $X_0$
from class $k$, and indeed, this is a very good method of
obtaining estimates of $ /pmc_k$, </it if//> such additional data are
available. 

Consider the usual situation, where no additional data are available.
We can simulate the ``new $X_0$ from class $k$'' situation by setting
part of the training data $z$ aside as a </it test set//>, constructing
the classifier using the complementary </it learning set//>, and then
run it on the test set.

The extreme version of this is the </it leave-one-out//> method,
introduced by Lachenbruch and Mickey (1968). Leave a particular
$/xkj $ out of the training set, reducing it to 
$$z_<(k,j)>=z-/lb /xkj /rb,$$
say. Construct the classifier $/hC$ anew, but this time based on
$z_<(k,j)>$, according to the same method that led one to the full-set
$/hC(./,;/,z)$. Next observe whether it correctly classifies $/xkj $
as coming from class $k$, or not. This produces the leave-one-out
estimates, of </it cross-validation//> estimates,
$$/hpmc_k^</cross>=<1/over /nk>/sum_<j=1>^</nk> I/lb/hC/lp /xkj ;/
z_<(k,j)> /rp/not=k/rb./eqno(12.5)$$

It seems to be quite a laborious task to arrive at (12.5) in practice,
since they involve a total of $N=/sum_<k=1>^K n_k$ different discriminant
rules. Often algebraic manipulations can provide shortcuts, however,
making it possible to compute $/hpmc_k^</cross>$ without actually
constructing every $/hC(./,;/,z_<(k,j)>)$.
/subsection<12.1.B><Leave-one-out estimates for the best quadratic
rule>
The remark above is nicely illustrated below for the best quadratic
rule $/hC_Q$, which is based on finding the maximal one among
$/pi_k/hf_k(k)=/pi_k N_d(/hm_k,/hs_k)(x)$, say with
$$/hm_k=<1/over /nk>/sum_<j=1>^</nk>/xkj ,/ /hs_k=a/nkp
/sum_<j=1>^</nk>/twice</lp /xkj -/hm_k/rp>'.$$
Here $a(/nk)$ is either $1/sla(/nk-1)$ or $1/sla n_k$.

$/hC_Q( /xkj ;/,z_<(k,j)>)$, needed in (12.5), is the class
index that gives maximum among
$$/pi_1/hf_1/lp x_j^<(k)>/rp/upto /pi_k/hf_<k,(j)>/lp
/xkj /rp/upto /pi_K/hf_K/lp /xkj /rp,$$
and these numbers are computed (at some subroutine level) using the
</it full//> classifier, except
$$/hf_<k,(j)>/lp /xkj /rp=N_d/lp/hm_<k,(j)>,/hs_<k,(j)>/rp/lp
/xkj /rp./eqno(12.6)$$ 
$/hf_<k,(j)>$ is $/hf_k$, but with estimates inserted for $/mu_k$ and
$/sg_k$ based on only $z-/lbrace /xkj /rbrace$.

But a shortcut is available for (12.6). One has
$$/eqalign<
/hm_<k,(j)>&=<1/over /nk-1>/sum_<i/not=j>x^<(k)>_i/cr
&=/hm_k-<1/over /nk-1>/lp /xkj -/hm_k/rp/cr>$$
and
$$/eqalign<
/hs_<k,(j)>&=a/lp /nk-1/rp/sum_<i/not=j>/twice</lp x_i^<(k)>-
/hm_<k,(j)>/rp>'/cr
&=<a/lp /nk-1/rp/over a/nkp>/lb/hs_k-</nk/over /nk-1>a/nkp
	/twice</lp /xkj -/hm_k/rp>'/rb./cr>$$
Introduce
$$/dl_<k,j>^2=/lp /xkj -/hm_k/rp'/hs_k^<-1>/lp /xkj -/hm_k/rp.
	/eqno(12.7) $$
Matrix algebra gives
$$/invert</hs_<k,(j)>>=/lp <a/lp /nk-1/rp/over a/nkp>/rp^d
	/invert</hs_k>/lb 1-</nk/over /nk-1> a/nkp/dl_<k,j>^2
	/rb$$ 
and
$$/hs_<k,(j)>^<-1>=<a/nkp/over a/lp n_k-1/rp>/lb/hs_k^<-1>+
	c_<k,j> </nk/over /nk-1> a/nkp/hs_k^<-1>
	/twice</lp /xkj -/hm_k /rp>'/hs_k^<-1>/rb,$$
where $c_<k,j>=/lbrace 1-</nk/over /nk-1> a(n_k)/dl_<k,j>^2/rbrace^<-1>$. It
follows that 
$$/eqalign<
/lp /xkj-/hm_<k,(j)>/rp'&/hs_<k,(j)>^<-1>/lp/xkj-/hm_<(k,j)>/rp/cr
&=/lp</nk/over/nk-1>/rp^2/lp/xkj-/hm_k/rp'/hs_<k,(j)>^<-1>/lp/xkj-
	/hm_<k,j>/rp /cr
&=/lp</nk/over/nk-1>/rp^2<a/nkp/over a/lp/nk-1/rp>/lb /dl_<k,j>^2+
	c_<k,j> </nk/over/nk-1> a/nkp/dl_<k,j>^4/rb/cr
&=/dl_<k,j>^2/lp1+r_<k,j>/rp,/cr>$$
say, writing
$$r_<k,j>=/lp</nk/over/nk-1>/rp^2<a(/nk)/over a(/nk-1)>/sla
/lb1-</nk/over/nk-1>/ a(/nk)/,/delta_<k,j>^2/rb-1./eqno(12.8)$$
By these efforts one arrives at
$$/eqalignno<
/hf_<k,(j)>/lp/xkj/rp&=/lp2/pi/rp^<-d/sla2>/lp<a(/nk-1)/over a(/nk)>
/rp^<-d/sla2>/invert</hs_k>^<-1/sla2>/cr
&/qquad/lb1-</nk/over/nk-1>/ a(/nk)/,/delta_<k,j>^2/rb^<-1/sla2>/exp
/lb-/halv/,/delta_<k,j>^2/lp1+r_<k,j>/rp/rb/cr
&=/hf_k/lp/xkj/rp/ B_<k,j>,&(12.9)/cr>$$
where
$$B_<k,j>=/lp<a(/nk)/over a(/nk-1)>/rp^<d/sla2>/lb1-</nk/over/nk-1>/ a(/nk)
/,/delta_<k,j>^2/rb^<-1/sla2>/exp/lp-/halv/,/delta_<k,j>^2 r_<k,j>/rp.
/eqno(12.10)$$

Suppose first that $a(/nk)=1/sla(/nk-1)$. Then
$r_<k,j>=/nk^2(/nk-2)/sla/lbrace(/nk-1)^3-/nk(/nk-1)/,/delta_<k,j>^2/rbrace-1$
and
$$/hf_<k,(j)>/lp/xkj/rp=/hf_k/lp/xkj/rp/lp</nk-2/over/nk-1>/rp^<d/sla 2>/lb
1-</nk/over(/nk-1)^2>/,/delta_<k,j>^2/rb^<-1/sla2>/exp/lp-/halv/,
/delta_<k,j>^2 r_<k,j>/rp./eqno(12.11)$$
If on the other hand $/hsk$ is used with multiplicative factor
$a/nkp=1/sla/nk$, then $r_<k,j>=/nk/sla(/nk-1-/delta_<k,j>^2)-1$ and
$$/hf_<k,(j)>/lp/xkj/rp=/hf_k/lp/xkj/rp/lp</nk-1/over/nk>/rp^<d/sla2>
/lb1-<1/over/nk-1>/delta_<k,j>^2/rb/exp/lp-/halv/delta_<k,j>^2
r_<k,j>/rp./eqno(12.12)$$

It is seen that $B_<k,j>$ is close to one for $/nk$ large, in which case
$$/wh f_<k,(j)>/lp x_j^<(k)>/rp/doteq/ /hf_k/lp x_j^<(k)>/rp,$$
so that the leave-one-out estimates (12.5) become close to the 
apparent error-rates (12.4).

The formulae above make it possible to compute the $/hpmc_k^</cross>$,
$k=1/upto K$, with little extra cost, at least in terms of computing
time. One does not need all the individual classifiers $/hC_Q(./,;/,Z_<(k,j)>)$;
it suffices to compute, for each $/xkj$,
$$/pi_1/hf_1/lp/xkj/rp/upto/pi_k/hf_k/lp/xkj/rp/upto/pi_K/hf_K/lp/xkj/rp$$
and find the maximiser (this gives a zero or one in the sum appearing 
in (12.4)); and then
$$/pi_1/hf_1/lp/xkj/rp/upto/pi_k/hf_k/lp/xkj/rp B_<k,j>/upto/pi_K/hf_K
/lp/xkj/rp$$
%
%
and find the maximiser (this gives a zero or one in the sum appearing in
(12.5)).
/subsection<12.1.C><Leave-one-out estimates for the best linear rule> Let us
also record how 
$$/hpmc_k^</cross>=<1/over/nk>/sum_<j=1>^/nk I/lb/hC_L/lp/xkj;/,z_<(k,j)>
/rp/not= k/rb/eqno(12.13)$$
can be computed for the best linear rule $/hC_L$, without 
actually constructing the individual classifiers.

Single out a particular $/xkj$ for attention. We must find out whether
$/hC_L(./,;/,z-/lbrace/xkj/rbrace)$ classifies $/xkj$ correctly or not, in as 
simple a fashion as possible.

$/hC_L(./,;/,z)$ works by evaluating
$$/hf_l(x)=</rm const.>/exp/lb-/halv/lp x-/hm_l/rp'/hs^<-1>
/lp x-/hm_l/rp/rb,$$
$l=1/upto K$, corresponding to the $N_d(/mu_l,/sg)$ assumption,
and then proceeds by normalising to $/pi_l/hf_l(x)/sla/sum_<l=1>^K/pi_l
/hf_l(x)$ and finding the maximiser and the maximum among these numbers.
Needed now are expressions for $/hf_<l,(k,j)>(/xkj)$, where
$/hf_<l,(k,j)>$ uses parameter estimates $/hm_<l,(k,j)>$ and $/hs_<(k,j)>$
that are based on $z$ without $/xkj$.

Let $S_k=/sum_<i=1>^k/twice<(x_i^<(k)>-/hm_k)>'$. The pooled estimate
of $/sg$ has the form $/hs=a(N)/sum_<k=1>^KS_k$, usually with 
$a(N)=1/sla(N-K)$ but sometimes with $a(N)=1/sla N$. After going through
matrix manipulations similar to those needed in the previous subsection 
one arrives at the following:
$$/eqalignno<
/lp/xkj-/hm_l/rp'&/hs_<(k,j)>^<-1>/lp/xkj-/hm_l/rp/cr
&=/lp/xkj-/hm_l/rp/hs^<-1>/lp/xkj-/hm_l/rp/lb1+r_<k,j>(l)/rb&(12.14)/cr>$$
for $l/not=k$, while
$$/eqalignno<
/lp/xkj-/hm_<k,(k,j)>/rp'&/hs_<(k,j)>^<-1>/lp/xkj-/hm_<k,(k,j)>/rp/cr
&=/lp/xkj-/hm_k/rp'/hs^<-1>/lp/xkj-/hm_k/rp/lb1+r_<k,j>(k)/rb.&(12.5)/cr>$$
Expressions for $r_<k,j>(1)/upto r_<k,j>(K)$ are given below;
they are complicated, but easy to program. It follows that 
$$/hf_<l,(k,j)>/lp/xkj/rp=</rm const./ >/hf_l/lp/xkj/rp B_l/lp/xkj/rp,
/eqno(12.16)$$
writing
$$B_l/lp/xkj/rp=/exp/lb-/halv r_<k,j>(l)/lp/xkj-/hm_l/rp'/hs^<-1>
/lp/xkj-/hm_l/rp/rb.$$
This means that the leave-one-out error rates (12.13) and the corresponding 
apparent error rates can be computed simultaneously, with only
one pass over the data.

It remains to provide formulae for the correction factors $r_<k,j>(l)$
appearing above. One can derive
$$1+r_<k,j>(l)=<a(N)/over a(N+1)>/lsb1+<</nk/over/nk-1>a(N)/over
1-</nk/over/nk-1>a(N)d_k'/hs^<-1>d_k>/ </lp d_k'/hs^<-1>d_l/rp^2/over
d_l'/hs^<-1>d_l>/rsb/eqno(12.17)$$
for $l/not=k$, and
$$/eqalignno<
1+r_<k,j>(k)&=<a(N)/over a(N-1)>/lp</nk/over/nk-1>/rp^2/lsb1+<</nk/over
/nk-1>a(N)/over1-</nk/over/nk-1>a(N)d_k'/hs^<-1>d_k> d_k'/hs^<-1>d_k/rsb/cr
&=<a(N)/over a(N-1)>/lp</nk/over/nk-1>/rp^2/sla/lsb1-</nk/over/nk-1>
a(N)d_k'/hs^<-1>d_k/rsb,&(12.18)/cr>$$
where $d_l=/xkj-/hm_l$, $d_k=/xkj-/hm_k$.

Observe that $r_<k,j>(l)$ is close to zero for $N=/sum_<k=1>^K/nk$ 
large, i.e., the $B_l(/xkj)$ factors in (12.16) are close to one,
so that cross validation rates in such a case will be very similar to 
apparent rates.
/subsection<12.1.D><Leave-one-out estimates for some/nl 
nonparametric classification procedures> Some nonparametric 
classification methods are so
intricate that no shortcuts are visible for the evaluation of the 
leave-one-out rates (12.5). In such cases one must resort to the 
laborious and CPU-time consuming task of actually constructing each
$/hC(./,;/,z_<(k,j)>)$,  or use simpler and less effective estimation
procedures, like the leave-a-quarter-out method: divide
$z$ into subsets $z_1$, $z_2$, $z_3$, $z_4$ of about equal size,
then run $/hC(./,;/,z_1/cup z_2/cup z_3)$ on $z_4$, and similarly for the
other three choices of quarter-sized test set. Average in the end to obtain 
error rate estimates.

It is, however, reasonably simple to carry out the leave-one-out
program for kernel type discriminant analysis  and
for the $k$-nearest neighbour method (see Sections 5.2 and 5.5).

Consider first kernel type density estimates of the form
$$/hf_k(x)=<1/over/nk>/sum_<j=1>^</nk>K_<h_k>/lp x-x_j/rp,$$
where $h_k=(h_<k,1>/upto h_<k,d>)$ is a vector of smoothing parameters for
the $d$ directions, and where $K_h(u)=K(<u_1/over h_1>/upto
<u_d/over h_d>)/sla(h_1/cdots h_d)$ for some basis kernel $K$.
The complete kernel classification method includes the difficult
step of choosing $h_k$ based on the training data $z$, and 
leave-one-outing quickly becomes tiresome if $h_k$ has to be 
found anew for each $z-/lbrace/xkj/rbrace$, say. But if one uses the same $h_k$
vectors, $k=1/upto K$, also for the reduced-by-one training sets, then
leave-one-out rates are easy to compute. This is seen from the fact that
$$/eqalign<
/hf_<k,(j)>/lp/xkj/rp&=<1/over/nk-1>/sum_<i/not=j>K_<h_k>/lp x_<(j)>^<(k)>
 -x_i^<(k)>/rp/cr
&=</nk/over/nk-1>/hf_k/lp/xkj/rp-<1/over/nk-1>K_<h_k>(0)./cr>$$

Finally consider the $k$-nearest-neighbour scheme, in which $/hC_<k-NN>(x;/,z)$
is the class being most frequent in the set $S_k(x;/,z)$ consisting of
the $k$ nearest neighbours, in $z$, to $x$. Then
$/hC_<k-NN>(/xkj;/,z-/lbrace/xkj/rbrace)$ is based on inspecting the set
$S_k(/xkj;/,z-/lbrace/xkj/rbrace)$, which happens to be the same set as
$S_<k+1>(/xkj;/,z)-/lbrace/xkj/rbrace$. Hence a </it k-NN> software package
can be used to evaluate
$$/hpmc_k^</cross>=<1/over/nk>/sum_<j=1>^</nk> I/lb/hC_<k-NN>/lp/xkj;
z-/lb/xkj/rb/rp/rb,$$
needing essentially the same amount of time as running the full-set 
$/hC_<k-NN>(./,;/,z)$ on the training set, i.e., obtaining cross validation
error rates is no more expensive than obtaining apparent error rates.
/note<Remark 1.> It is very useful to study not only the error rates
$/pmc_k$, but the full </it confusion matrix//>. If
$/hC(./,;/,z)$ is run on a second set $z'$, then
$$M(k,l)=/sum_<<>^</ / x'</rm/ in/ >z'>_<</rm from/ class/ >k>>I/lb/hC(x';/,z)=l/rb
/eqno(12.19)$$
is the number among those from class $k$ in $z'$ that end up 
being classified as $l$. The matrix of these numbers is the confusion
matrix, from which further information can be extracted. Note that
(12.19) is useful also for $l=/doubt$ and $l=/out$.

It is of course possible to run $/hC(./,;/,z)$ on the training set $z$ itself;
indeed this is often done in practice because it is easy and does not
require further programming. The result is however an
``apparent confusion matrix'' which tends to be too optimistic;
the off-diagonal elements usually become larger when $/hC(./,;/,z)$
is run on a second, independent set of the same size.

A more realistic picture of the classifier's ability to discriminate
future objects is provided by the </it cross-validated//> confusion matrix,
consisting of elements
$$M_</cross>(k,l)=/sum_<j=1>^</nk>I/lb/hC/lp/xkj;/,z_<(k,j)>/rp=l/rb./eqno(12.20)$$
The machinery needed to obtain (12.5) produces (12.20) too with little extra 
cost.
/note<Remark 2.> We have concentrated on estimating the $/pmc_k$ rates
above. Similar reasoning applies to other rates, however.

Of particular importance is the </it total//> error rate $/pmc$ of (12.3).
If the $/pi_k$'s are known, then
$$/eqalign<
/hpmc^</cross>&=/sum_<k=1>^K/pi_k/hpmc_k^</cross>/cr
&=/sum_<k=1>^K/pi_k<1/over/nk>/sum_<j=1>^</nk>I/lb/hC/lp/xkj;/,z_<(k,j)>/rp
/not=k/rb/cr>$$
is the canonical estimator. If the $/pi_k$'s are estimated themselves
from the training set, as $/nk/sla N$, then
$$/hpmc^</cross>=<1/over N>/sum_<x/in z>I/lb/hC/lp x;/,z-/lb x/rb/rp
/not=</rm class>(x)/rb$$
is the natural, overall leave-one-out estimate.
/section<12.2><Bootstrapping and other resampling methods><>
The cross validation technique produces estimators for $/pmc_k$ and $/pmc$ 
that are nearly unbiased, but that may have high variability.
Experimental evidence for this has been offered several times in the 
literature, see e.g.~Efron (1983). Thus there may be a scope for
improvement over cross validation.

This section discusses bootstrap estimation of the $/pmc_k$'s,
and a particular, related method called the $.632$ estimator.
/subsection<12.2.A><Bootstrap and jackknife estimators>
We know that the apparent error rate $/apmc_k$ usually
underestimates the true error rate $/pmc_k$, i.e.~the
random variable
$$/op_k=/pmc_k-/apmc_k,/eqno(12.21)$$
called the </it optimism//> for $/apmc_k$, is usually positive.
$$/eqalign<
/op_k=/op_k(z;/,F_k)&=/pmc_k(z;/,F_k)-/apmc_k(z)/cr
&=/Pr_<F_k>/lb/hC(X_0;/,z)/not=k/rb-/Pr_</wh F_k>/lb/hC(X_0;/,z)/not=
k/rb/cr>/eqno(12.22)$$
is the same equation, but spells out the dependence upon the 
actually observed training set and the distribution of future vectors from 
class $k$. $/wh F_k$ is the usual empirical distribution with
point mass $1/sla/nk$ in each $/xkj$.

An idea of Efron (1983) is as follows. Let
$$/som_k=/som_k(F_1/upto F_K)=E_<F_1/upto F_K>/op_k(Z;/,F_k)/eqno(12.23)$$
be the expectation of $/op_k$ over all training sets 
$Z=/lbrace X_1^<(k)>/upto X_</nk>^<(k)>;/ k=1/upto K/rbrace$.
($Z$ is a collection of random variables, $z$ is the actually
observed training data.) If $/som_k$ were known, then $/apmc_k+/som_k$ would
be a natural estimator. In practice $/som_k$ must be estimated,
however, and for each estimator $/wh /som_k$ there is a corresponding
estimator $/hpmc_k=/apmc_k+/wh /som_k$ of the true error rate.

One way of estimating $/som_k$ is by cross validation,
$$/wh /som_k^</cross>=<1/over/nk>/sum_<j=1>^</nk> I/lb/hC/lp/xkj;/,z_<(k,j)>/rp
/not=k/rb-/apmc_k,/eqno(12.24)$$
producing in its turn $/hpmc_k^</cross>=/apmc_k+/wh /som_k^</cross>$,
which is the leave-one-out estimator (12.5) again.

Another way of estimating the expected optimism (12.23) is by simply inserting
the natural nonparametric estimates $/wh F_1/upto/wh F_K$ for
$F_1/upto F_K$:
$$/wh /som_k^</boot>=/som_k/lp/wh F_1/upto/wh F_K/rp=E_</ast>/op_k
/lp Z^</ast>;/wh F_k/rp./eqno(12.25)$$
Here $Z^</ast>=/lbrace X_1^<(k)/ast>/upto X_</nk>^<(k)/ast>;/ k=1/upto K/rbrace$ is a
``pseudo training set'', drawn with $X_j^<(k)/ast>$ from $/wh F_k$,
and $E_</ast>$, $/Pr_</ast>$ mean operations w.r.t.~this distribution.
$/wh /som_k^</boot>$ is called the bootstrap estimator of $/som_k$, and
is in practice computed by simulating a large number of
$Z^/ast$ sets.
The </it bootstrap estimator//> of the true error rate becomes
$$/hpmc_k^</boot>=/apmc_k+/wh /som_k^</boot>./eqno(12.26)$$

The following rewriting is practical when it comes to the computation of
$/wh /som_k^</boot>$:
$$/eqalignno<
</rm op>_k/lp Z^/ast;/wh F_k/rp&=/pmc_k/lp Z^/ast;/wh F_k/rp-/apmc_k/lp Z^/ast/rp/cr
&=/Pr_</wh F_k>/lb/hC/lp X_0;/,Z^/ast/rp/not=k/rb-<1/over/nk>/sum_<j=1>^</nk>
I/lb/hC/lp X_j^<(k)/ast>;/,Z^/ast/rp/not=k/rb/cr
&=/sum_<j=1>^</nk>/lp<1/over/nk>-P_<k,j>^/ast/rp I/lb/hC/lp/xkj; Z^/ast/rp/not=k
/rb.&(12.27)/cr>$$ 
Here $P_<k,j>^/ast$ is the proportion among $X_1^<(k)/ast>/upto X_</nk>^<(k)/ast>$
that fall in $/xkj$. Perhaps $B=200$ simulated values of (12.27) can 
now be obtained (in the computer), their average is (a very good
approximation to) $/wh /som_k^</boot>$.
/note<Remark 3.> When we refer to a random training set $Z$, as in
(12.23) and later on, we mean a random set $/lbrace X_1^<(k)>/upto
X_</nk>^<(k)>;/ k=1/upto K/rbrace$ with $N=n_1+/cdots+n_K$ independent
variables,
and $/xkj$ drawn from $F_k$. In particular we consider the class sizes
$n_1/upto n_K$ as fixed. A slightly different situation, and that
also could be analysed, is one where $N$ $X$'s are drawn independently 
from the mixture distribution $/pi_1F_1+/cdots+/pi_KF_K$.

The above viewpoint is natural since we take aim at the conditional, true
error rates $/pmc_k$. The second approach is more appropriate when the 
class sizes $n_1/upto n_K$ also are random and when one focusses on the
</it total//> error rate $/pmc=/sum_<k=1>^K/pi_k/pmc_k$. This is the 
approach chosen by Efron (1983).

In adherence to the above also the bootstrap training samples $Z^/ast$ are made
up with $n_1$ from $/wh F_1$, $n_2$ from $/wh F_2$, and so on.
/Enote
Of course (12.26) is a computationally intensive estimator of $/pmc_k$, but
it is absolutely feasible in practice. It may nevertheless be of interest to 
derive approximations that require less computing. Such an approximation
is the </it jackknife estimator//>
$$/hpmc_k^</jack>=/apmc_k+/wh /som_k^</jack>,/eqno(12.28)$$
where
$$/eqalignno<
/wh /som_k^</jack>&=<1/over/nk>/sum_<j=1>^</nk>I/lb/hC/lp/xkj;/,z_<(k,j)>
/rp/not=k/rb/cr
&/quad-<1/over/nk^2>/sum_<j=1>^</nk>/sum_<l=1>^</nk>I/lb/hC/lp/xkj;
z_<(k,l)>/rp/not=k/rb.&(12.29)/cr>$$
That $/wh /som_k^</jack>$ approximates $/wh /som_k^</boot>$ follows from a
quadratic expansion of (12.27) and from properties of the multinomial 
distribution. It requires $K$-sample generalisations of some of the theory in Efron (1982, Ch.~6).

Experience of Efron (1982, 1983) and others indicates the following: (i)/ 
$/hpmc_k^</cross>$ is almost unbiased, but has a somewhat large
variance; (ii)/ $/hpmc^</jack>$, an approximation to the natural
nonparametric estimator $/hpmc^</boot>$, is similar to $/hpmc_k^</cross>$
in form and value; (iii)/ $/hpmc_k^</boot>$ may have a slight, downward bias, but is less variable than the others, and performs 
usually better (e.g.~in terms of mean squared error).
/subsection<12.2.B><The $.632$ estimator> There are various bootstrap
and resampling schemes around that can compete with the ``basic bootstrap''
discussed above. Efron (1983) discusses a ``randomised bootstrap'' and
a ``double bootstrap''. He discussed only the two-class problem
and the </it total//> error rate (12.3), and it is an interesting and
non-trivial exercise to generalise the randomised bootstrap and the double
bootstrap in an appropriate way, to cover $K$-class problems
and conditional error rates. We concentrate now, however, on providing
such a generalisation of a third bootstrap version of Efron's, namely
the so-called $.632$ estimator. This estimator performed best in
several sampling experiments.

$/pmc_k$ is the error rate for a future $X_0$ from $F_k$, while the apparent 
rate $/apmc_k$ may be seen as the error rate for an $X_0$ drawn from
$/wh F_k$, i.e.~from the ($k$'th part of the) training set, cf.~(12.22).
$/pmc_k$ is usually larger then $/apmc_k$ because an $X_0$ from $F_k$ will 
lie some distance away from $z^<(k)>=/lbrace x_1^<(k)>/upto x_</nk>^<(k)>/rbrace$.
To introduce a notion of this distance, define first a system of 
neighbourhoods $S_k(x,/Delta)$ around points $x$ in the feature space, by
$$/Pr_<F_k>/lb X_0/in S_k/lp x,/Delta/rp/rb=/Delta, 0/le/Delta/le1./eqno(12.30)$$
One could, for example, take
$$S_k(x,/Delta)=/lb y:(y-x)'/sg_k^<-1>(y-x)/le r_k(x,/Delta)^2/rb,/eqno(12.31)$$
where $r_k(x,/Delta)$ is determined from (12.30). Thus $r_k(x,0)=0$ and
$r_k(x,1)=/infty$. Now let
$$/delta_k(x_0,z)=/inf/lb/Delta:/ x_0/in/bigcup_<j=1>^</nk>S_k
/lp/xkj,/Delta/rp/rb,/eqno(12.32)$$
the distance from a point $x_0$ to the nearest point in $z^<(k)>$. 
We may view $/pmc_k$ as the probability of misclassifying an $X_0$ lying a
random distance $/delta_k(X_0,z)$ away from $z^<(k)>$.

Consider
$$R_k(/Delta)=/Pr_<F_k>/lb/hC/lp X_0; Z/rp/not=k/ /vert/ /delta_k
/lp X_0,Z/rp=/Delta/rb,/eqno(12.33)$$
the conditional misclassification probability for $X_0$ vectors lying a 
distance $/Delta$ away from $Z^<(k)>=/lbrace X_1^<(k)>/upto X_</nk>^<(k)>/rbrace$.
We aim at estimating $/som_k=E </rm op>_k(Z; F_k)$, and therefore have a random
training set $Z$ in (12.33). Observe that
$$/eqalignno<
R_k(0)&=/Pr_<F_k>/lb/hC(X_0;/,Z)/not=k/ /vert/ X_0</rm/ is/ in/ >Z^<(k)>/rb/cr
&=<1/over/nk>/sum_<j=1>^</nk>/Pr/lb/hC/lp/xkj;/,Z/rp/not=k/rb/cr
&=E/ /apmc_k.&(12.34)/cr>$$
Note also that
$$/eqalign<
E/ /pmc_k(Z;/,F_k)&=E/lb E/ /pmc_k(Z;/,F_k)/ /vert/ /delta_k(X_0,Z)/rb/cr
&=/int_0^1 R_k(/Delta)dD_k(/Delta),/cr>$$
so that
$$/som_k=/int_0^1/lb R_k(/Delta)-R_k(0)/rb/ dD_k(/Delta)./eqno(12.35)$$
Here $D_k(.)$ is the distribution of $/delta_k(X_0,Z)$.

The content of Remark 3 applies again, now to (12.33)--(12.35).

The ``bootstrap equivalent'' to (12.33) is
$$R_k^<(B)>(/Delta)=</rm Prob>_k/lb/hC/lp X_0^/ast;/,Z^/ast/rp/not=/ k/vert/ /delta_k
/lp X_0^/ast,Z^/ast/rp=/Delta/rb./eqno(12.36)$$
Here $X_0^/ast$ is drawn from $/wh F_k$ and independently of $Z^/ast$, and the 
probability is w.r.t.~</it both//> the random $Z$ and the following 
$X_0^/ast$, $Z^/ast$. $R_k^<(B)>(/Delta)$ will be in close agreement with 
$R_k(/Delta)$, see Efron (1983) for some evidence.

Cases where $/delta_k^/ast =/delta_k(X_0^/ast,Z^/ast)/gt 0$ are qualitatively different 
from those where $/delta_k^/ast =0$. One has
$$/eqalignno<
</rm Prob>_k/lb/delta_k^/ast/gt0/rb&=</rm Prob>_k/lb X_0^/ast </rm/ is/ new/ to/ >X_1^<(k)/ast>/upto
X_</nk>^<(k)/ast>/rb/cr
&=E_k/ </rm Prob>_k/lb X_j^<(k)/ast>/not=X_0^/ast,/ j=1/upto/nk/ /vert/ Z, X_0^/ast/rb/cr
&=/lp1-<1/over/nk>/rp^</nk>/cr
&/doteq e^<-1>=.368.&(12.37)/cr>$$
Now consider, for the given training set $z$,
$$/heps_k=/heps_k(z)=/Pr_</ast,k>/lb/hC/lp X_0^/ast;/,Z^/ast/rp/not=k/ /vert/ /delta_k^/ast 
/gt0/rb,/eqno(12.38)$$
which in practice can be computed by simulation, for example as
$$/sum_<b=1>^B I/lb/hC/lp X_0^</ast b>;/,Z^</ast b>/rp/not=k/rb/sla/sum_<b=1>^B
I/lb/dl_k^</ast b>/gt0/rb/eqno(12.39)$$
for a large number $B$ of outcomes $Z^</ast b>$, $X_0^</ast b>$.
A remarkable fact, to be commented on later, is that
$$</rm Prob>_k/lb/dl_k^/ast/gt/Delta/ /vert/ /dl_k^/ast/gt0/rb/doteq</rm Prob>_k
/lb/dl_k/gt.632/Delta/rb/eqno(12.40)$$
holds, to a good approximation, at least for small values of $/Delta$.
Here $/dl_k=/dl_k(X_0,/,Z)$.

The observations above combine to the following motivation for an estimator
of $/som_k$, the bias of $/apmc_k$:
$$/eqalignno<
E/lb/heps_k(Z)-/apmc_k(Z)/rb&=/int_</Delta/gt0> R_k^<(B)>(/Delta)dD_k^<(B)>(/Delta)-R_k(0)/cr
&/doteq/int_</Delta/gt0>/lb R_k(/Delta)-R_k(0)/rb dD_k^<(B)>(/Delta)/cr
&/doteq/int_</Delta/gt0>/lb R_k(/Delta)-R_k(0)/rb dD_k(.632/Delta)/cr
&=/int_</Delta/gt0>/lb R_k/lp</Delta/over/ .632>/rp-R_k(0)/rb dD_k(/Delta)/cr
&/doteq/int_</Delta/gt0><1/over.632>/lb R_k(/Delta)-R_k(0)/rb dD_k(/Delta)/cr
&=/som_k/sla.632.&(12.41)/cr>$$
This suggests
$$/wh /som_k^<(.632)>=.632/lp/heps_k-/apmc_k/rp$$
as an estimator of $/som_k$, and, finally,
$$/eqalignno<
/hpmc_k^<(.632)>&=/apmc_k+/wh /som_k^<(.632)>/cr
&=.368/apmc_k+.632/heps_k&(12.42)/cr>$$
as an estimator of $/pmc_k$.

Several comments are in order here. $D_k^<(B)>(.)$ is the distribution
of $/dl_k^/ast$, conditional on $/dl_k^/ast$ being positive, cf.~(12.40).
The final approximation involved in (12.41) assumes that the function
$R_k(/Delta)-R_k(0)$ is reasonably linear for small $/Delta$, i.e.,
the values that count for the $D_k$ distribution of $/dl_k=/dl_k(X_0,Z)$.
Under these circumstances, assuming also that (12.36) is close to (12.33)
and that (12.40) holds, will $/hpmc_k^<(.632)>$ be reasonably unbiased
for $/pmc_k$.

It remains only to provide evidence for (12.40). From the system of 
neighbourhoods (12.30), define
$$T_k(x_0,/Delta)=/lb x:/,x_0/in S_k(x,/Delta)/rb.$$
Then, from the definition (12.41),
$$/eqalign<
</rm Prob>_k/lb/dl_k(X_0,Z)/gt/Delta/rb&=</rm Prob>_k/lb/Xkj/not/in
T_k(X_0,/Delta),/ j=1/upto/nk/rb/cr
&=/int/lsb1-F_k/lb T_k(x_0,/Delta)/rb/rsb^</nk> d/,F_k(x_0)./cr>$$
$F_k/lbrace T_k(x_0,/Delta)/rbrace=/Pr_<F_k>/lbrace x_0/in S_k(X,/Delta)/rbrace$
will be close to $/Delta$, for small $/Delta$, by the defining property
(12.30). Hence
$$</rm Prob>_k/lb/dl_k/gt/Delta/rb/doteq(1-/Delta)^</nk>./eqno(12.43)$$
$/lbrack$Consider the specific system (12.31), for example. Then
$S_k(x,/Delta)$ is an ellipsoid with volume approximately equal to
$/Delta/sla/fk(x)$, which means that the radius $r_k(x,/Delta)$ appearing
there is proportional to $/lbrace/Delta/sla/fk(x)/rbrace^<1/sla d>$,
$d$ being the feature vector dimension. Accordingly,
$$/eqalign<
/Pr_<F_k>/lb X/in T_k(x_0,/Delta)/rb&=/Pr_<F_k>/lb(X-x_0)'/sg_k^<-1>
(X-x_0)/le r_k(X,/Delta)^2/rb/cr
&/doteq/Pr_<F_k>/lb(X-x_0)'/sg_k^<-1>(X-x_0)/le r_k(x_0,/Delta)^2/rb
=/Delta,/cr>$$
at least for small $/Delta$.$/rbrack$

Next consider
$$/eqalign<
</rm Prob>_k/lb/dl_k^/ast/gt/Delta/ /vert/ /dl_k^/ast/gt0/rb&=
</rm Prob>_k/lb X_j^<(k)/ast>/not/in T_k/lp X_0^/ast,/Delta/rp,/ j=1/upto/nk
/ /vert/ /dl_k^/ast/gt0/rb/cr
&=/int E/lsb1-F_k/lb T_k(x_0,/Delta)/rb/rsb^</nk^/ast>d/,F_k(x_0)/cr
&/doteq E(1-/Delta)^</nk^/ast>,/cr>$$
where $/nk^/ast$ is the number of data points $/xkj$, $j=1/upto/nk$ that are
visited by a bootstrap sample $X_1^<(k)/ast>/upto X_</nk>^<(k)/ast>$.
$$/nk^/ast =/sum_<j=1>^</nk> I/lb/xkj</rm/ is/ visited/ >/rb$$
has expectation $/nk/lbrace1-(1-<1/over/nk>)^</nk>/rbrace/doteq
.632/nk$, and is close to being distributed as a binomial $(/nk,/,.632)$.
Hence
$$/eqalignno<
</rm Prob>_k/lb/dl_k^/ast/gt/Delta/ /vert/ /dl_k^/ast/gt0/rb&/doteq
E(1-/Delta)^<</rm Bin>(/nk,/,.632)>/cr
&=(1-/,.632/Delta)^</nk>.&(12.44)/cr>$$
This verifies (12.40).
/note<Remark 4.> $/dl_k=/dl_k(X_0,Z)$, the distance from a new $X_0$ to
the training set, has (12.43) as its approximate distribution, in
particular
$$E/dl_k(X_0,Z)/doteq<1/over/nk+1>,$$
on the scale from zero to one set out by (12.30) and (12.32).
Also, by (12.35),
$$/som_k/doteq/int_0^1/lb R_k(/Delta)-R_k(0)/rb/,/nk(1-/Delta)^</nk-1>/,d/Delta,$$
and, since $R_k(/Delta)$ must be reasonably linear in the short range of
$/Delta$-values that count,
$$/som_k/doteq/int_0^1 c_kd/,D_k(/Delta)=<c_k/over/nk+1>./eqno(12.45)$$
Here $c_k$ is the slope of $R_k(/Delta)$ just past $/Delta=0$.
This provides some insight into the de-bias factor $/som_k$, and can also
be utilised to construct alternative estimates, for example using
$R_k^<(B)>(/Delta)/doteq R_k(/Delta)$.
/section<12.3><Concluding remarks><> There is an enormous amount of published
papers on the topic of error rate estimation. Toussaint (1974)
contains a bibliography with 188 references, see also Kanal (1974).
Hand (1986a) presents an updated review of the topic.

The leave-one-out or the sometimes simpler leave-a-quarter-out methods
usually perform quite well. Only in vital situations where more care
is called for (when the ``final classification procedure'' is to be
assessed, chosen after much initial browsing, for example), or in
unstable situations with many classes and/sla or small training sets,
are the more sophisticated and computer intensive methods of 12.2
called for. In many cases even the apparent error rates will do, say
in an initial exploring phase when a lot of different feature
extraction methods are being tried out, preferably in combination with
estimation of interclass distances, see Chapter 10.

That there are shortcuts available to easier evaluation of cross
validation error estimates for the best linear and the best quadratic
rules seems to be part of the statistical folklore, but we haven't
found good references; hence Sections 12.1.B and 12.1.C.

Another approach to estimating error rates, often presented in
not-so-new textbooks, relies on trusting the parametric models used
to estimate class distributions. If the normality assumptions hold, for
example, then the error rate is only a (complex) function of the
parameters, and the corrections to the plug-in estimators, based
usually on asymptotic expansions, can be made. These methods are 
sensitive to the correctness of the assumed parametric models, of course,
and are generally not recommended.

The main reference for bootstrap and other resampling based
competitors to cross validation is Efron (1983). Our contribution here
is to generalise some of his methods to the $K$-class situation and to
conditional error rates.

There is a more efficient way  of computing the $.632$ estimator than
via (12.39), see Efron's Section 8. Also supplied by Efron (op.~cit.)
is an interesting connection from the $.632$ estimator to the
halfsample cross validation estimator.

Chapter 7 considered the surprising potential of unclassified feature
vectors, which sometimes can be easily available, to updating and even
refining class distribution estimates. Unclassified vectors can also
be utilised to estimate error rates for classification procedures. The
basic idea is presented in Fukunaga and Kessell (1973), and consists
of representing error rate as an integral of the $x$-conditional error
rate weighted by the mixture distribution $f(x)$. Incoming,
unclassified vectors can then be used to estimate $f$, say by one of
the nonparametric methods of Chapter 5. This leads to average
conditional error rate estimates. Recent references are Hand (1986a, 1986b).

We have concentrated attention on performance for a given classification 
rule. It is also interesting to estimate the </it Bayes//> or </it ideal 
error rate//> that would have resulted if class densities and prior 
probabilities were known exactly, cf.~Chapter 1 and the discussion 
involving $/etrue$ and $/eideal$ in Section 10.2.A. The </it 
ideal success probability//> is
$$/eqalign<
/pcc_</rm ideal>&=/sum_<k=1>^K/pi_k/,/pcc_<k,</rm ideal>>/cr
&=/sum_<k=1>^K/pi_k/int I/lb/pi_k/fkx</rm/ is/ largest>/rb/fkx/,dx/cr
&=/int/lb/max_<k/le K>/pi_k/fkx/rb dx/cr
&=/int/lb/max_<k/le K>P(k/mid x)/rb/,f(x)/,dx,/cr>$$
and is a benchmark against which classifiers can be measured.
$/pcc_</rm ideal>$ can be estimated in various ways, and unclassified vectors 
can again be of use. There are also lower and upper bounds available for 
this probability.


/vfill/supereject
/begingroup
/nonfrenchspacing
/header<References>
/parindent=0pt
/parskip=/medskipamount
/def/auth</rm>
/def/art</rm>
/def/from</rm>
/def/publ</it>
/def/book#1</rm In: /it #1>
/def/vol#1<</bf #1>>
/def/pages#1-#2</rm, #1--#2>
/def/page#1</rm, p./ #1>

/auth Abramowitz, M. and Stegun, I.A. (1964).
/publ Handbook of Mathematical Functions.
/from U.S. Government Printing Office, Washington, D.C.

/auth Aitchison, J. and Aitken, C.G.G. (1976).
/art  Multivariate binary discrimination by the kernel method.
/publ Biometrika /vol<63>/pages 413-420.

/auth Aitchison, J., Habbema, J.D.F. and Kay, J.W. (1977).
/art  A critical comparison of two methods of statistical discrimination.
/publ Applied Statistics /vol<26>/pages 15-25.

/auth Akaike, H. (1974).
/art  A new look at statistical model identification.
/publ IEEE Transactions on Automatic Control /vol<AC-19>/pages 716-723.

/auth Anderberg, M.R. (1973).
/publ Cluster Analysis for Applications.
/from Academic Press, New York.

/auth Anderson, J.A. (1969).
/art  Discrimination between $k$ populations with constraints on the
      probabilities of misclassification.
/publ Journal of the Royal Statistical Society Series B /vol<31>/pages 123-139.

/auth Anderson, T.W. (1958).
/publ An Introduction to Multivariate Analysis.
/from Wiley, New York.

/auth Andrews, D.F., Gnanadesikan, R., and Warner, J.L. (1973).
/art  Methods for assessing multivariate normality.
/book Multivariate Analysis-III,
/from ed. P.R. Krishnaiah, Academic Press, New York/pages 95-116.

/auth Bahadur, R.R. (1961a).
/art  A representation of the joint distribution of responses to $n$
      dichotomous items.
/book Studies in Item Analysis and Prediction,
/from ed. H. Solomon, Stanford University Press, California/pages 158-167.

/auth Bahadur, R.R. (1961b).
/art  On classification based on responses to $n$ dichotomous items.
/book Studies in Item Analysis and Prediction,
/from ed. H. Solomon, Stanford University Press, California/pages 169-176.

/auth Berger, J. (1980a).
/publ Statistical Decision Theory: Foundations, Concepts, and Methods.
/from Springer, Berlin.

/auth Berger, J. (1980b).
/art  Improving on inadmissible estimators in continuous exponential
      families with applications to simultaneous estimation of
      gamma scale parameters.
/publ Annals of Statistics /vol<8>/pages 545-571.

/auth Besag, J. (1974).
/art  Spatial interaction and the statistical analysis of lattice
      systems (with discussion).
/publ Journal of the Royal Statistical Society Series B /vol<36>/pages 192-236.

/auth Besag, J. (1986).
/art  On the statistical analysis of dirty pictures (with discussion).
/publ Journal of the Royal Statistical Society Series B /vol<48>/pages 259-302.

/auth Billingsley, P. (1968).
/publ Convergence of Probability Measures.
/from Wiley, New York.

/auth Bosq, D. (1970).
/art  Contribution a le th/'eorie de l'estimation fonctionelle.
/publ Publications de l'Institut de Statistique de l'Universit/'e de Paris /vol<19>/pages 1-177.

/auth Box, G.E.P., and Tiao, G.C. (1973).
/publ Bayesian Inference in Statistical Analysis.
/from Addison-Wesley, Menlo Park.

/auth Breiman, L., Meisel, W., and Purcell, E. (1977).
/art  Variable kernel estimates of multivariate densities.
/publ Technometrics /vol<19>/pages 135-144.

/auth Breiman, L., Friedman, J.H., Olshen, R.A., and Stone, C.J. (1984).
/publ Classification and Regression Trees.
/from Wadsworth, Belmont, California.

/auth Brown, P.J. and Rundell, P.W.K. (1985).
/art  Kernel estimates for categorical data.
/publ Technometrics /vol<28>/pages 293-299.

/auth Br}ten, K., Holb{k-Hanssen, E. and Taxt, T. (1986a).
/publ Statistical Symbol Recognition: Development of a System.
/from Report No. 777, Norwegian Computing Centre, Oslo.

/auth Br}ten, K., Holb{k-Hanssen, E. and Taxt, T. (1986b).
/art  A general software system for supervised statistical classification 
      of symbols.
/publ Proceedings of the Eighth International Conference on Pattern 
      Recognition, Paris, October 27-31 1986.

/auth Butler, W.J. and Kronmal, R.A. (1985).
/art  Discrimination with polychotomous predictor variables using orthogonal functions.
/publ Journal of the American Statistical Association /vol<80>/pages 443-448.

/auth Chernoff, H. (1973).
/art  Some measures for discriminating between normal multivariate 
      distributions with unequal covariance matrices.
/book Multivariate Analysis-III,
/from ed. P.R. Krishnaia/-h, Academic Press, New York/pages 337-344.

/auth Chow, C.K. (1962).
/art  A recognition method using neighbor dependence.
/publ IRE Transactions on Electronic Computing /vol<EC-11>/pages 683-690.

/auth Chow, C.K. and Liu, C.N. (1966).
/art  An approach to structure adaptation in pattern recognition.
/publ IEEE Transactions on Systems, Science, and Cybernetics /vol<SCC-2>/pages 73-80.

/auth Crain, B.R. (1976).
/art  Matrix density estimation.
/publ Communications in Statistics---Theory /& Methods /vol<5>/pages 89-96.

/auth Dempster, A.P., Laird, N.M., and Rubin, D.B. (1977).
/art  Maximum likelihod from incomplete data via the EM algorithm (with discussion).
/publ Journal of the Royal Statistical Society Series B /vol<39>/pages 1-38.

/auth Duda, R.J. and Hart, P.E. (1973).
/publ Pattern Classification and Scene Analysis.
/from Wiley, New York.

/auth Duin, R.P.W. (1976).
/art  On the choice of smoothing $j$ parameters for Parzen estimators of probability density functions.
/publ IEEE Transactions on Computers /vol<C-25>/pages 1175-1179.

/auth Durbin, J. (1973).
/art  Weak convergence of the sample distribution function when parameters are estimated.
/publ Annals of Statistics /vol<1>/page 279-290.

/auth Efron, B. (1982).
/publ The Jackknife, the Bootstrap, and Other Resampling Plans.
/from SIAM NSF-CBMS Monograph /#38.

/auth Efron, B. (1983).
/art  Estimating the error rate of a prediction rule: improvement on cross-validation.
/publ Journal of the American Statistical Association /vol<78>/pages 316-331.

/auth Eisenberger, I. (1964).
/art  Genesis of bimodal distributions.
/publ Technometrics /vol<6>/pages 357-363.

/auth Fisher, R.A. (1936).
/art  The use of multiple measurements in taxonomic problems.
/publ Annals of Eugenics /vol<7>/pages 179-188.

/auth Friedman, J.H., Stuetzle, W., and Schroeder, A. (1984).
/art  Projection pursuit density estimation.
/publ Journal of the American Statistical Association /vol<78>/pages 599-608.

/auth Friedman, J.H. (1986).
/art  Classification and multiple regression through projection pursuit.
/publ Journal of the American Statistical Association
/from (to appear).

/auth Fukunaga, K. and Kessell, D.L. (1973).
/art  Nonparametric Bayes error estimation using unclassified samples.
/publ IEEE Transactions on Information Theory /vol<IT-19>/pages 434-440.

/auth Fukunaga, K. and Narenda, P.M. (1975).
/art  A branch and bound algorithm for computing $k$-nearest neighbours.
/publ IEEE Transactions on Computers /vol<C-24>/pages 750-753.

/auth Gates, G.W. (1972).
/art  The reduced nearest neighbour rule.
/publ IEEE Transactions on Information Theory /vol<IT-18>/pages 431-432.

/auth Geisser, S. (1982).
/art  Bayesian discrimination.
/book Handbook of Statistics,
/from vol. /vol<2>, eds. P.R. Krishnaiah and L.N. Kanal, North-Holland/pages 101-120.

/auth Geman, S. and Geman, D. (1984).
/art  Stochastic relaxation, Gibbs distributions, and the Bayesian restoration of images.
/publ IEEE Transactions on Pattern Analysis and Machine Intelligence
      /vol<PAMI-6>/pages 721-741.

/auth Glick, N.D. (1974).
/art  Consistency condition for probability estimators and integrals
      of density estimates.
/publ Utilitas Matematica /vol<6>/pages 61-74.

/auth Greblicki, W. (1978).
/art  Asymptotically optimal pattern recognition procedures with
      density estimates.
/publ IEEE Transactions on Information Theory /vol<IT-24>/pages 250-251.

/auth Haggstrom, G.W. (1983).
/art  Logistic regression and discriminant analysis by ordinary least squares.
/publ Journal of Business /& Economic Statistics /vol<1>/page 229-238.

/auth Hand, D.J. (1981).
/publ Discrimination and Classification.
/from Wiley, Chichester.

/auth Hand, D.J. (1982).
/publ Kernel Discriminant Analysis.
/from Research Studies Press, Chichester.

/auth Hand, D.J. (1986a).
/art  Recent advances in error rate estimation.
/publ Pattern Recognition Letters /vol<4>/pages 335-346.

/auth Hand, D.J. (1986b).
/art  An optimal error rate estimator based on average conditional
      error rate: Asymptotic results.
/publ Pattern Recognition Letters /vol<4>/pages 347-350.

/auth Hand, D.J. and Batchelor, B.G. (1978).
/art  An edited condensed nearest neighbor rule.
/publ Information Science /vol<14>/pages 171-180.

/auth Harris, J.F., Kittler, J., Llewellyn, B. and Preston, G. (1982).
/art  A modular system for interpreting binary pixel representatitons of
      line-structured data on maps.
/publ Cartographia /vol<19>, No./ 2.

/auth Hart, P.E. (1968).
/art  The condensed nearest neighbour rule.
/publ IEEE Transactions on Information Theory /vol<IT-14>/pages 515-516.

/auth Hathaway, R.J. (1968).
/art  A constrained formulation of maximum-likelihood estimation
      for normal mixture distributions.
/publ Annals of Statistics /vol<13>/pages 795-800.

/auth Helgeland, J. and Hjort, N.L. (1986).
/art  Optimal recognition of closed Gaussian curves in the plane (preliminary report).
/from Presented at The First World Congress of the Bernoulli Society,
      Tashkent, USSR.

/auth Henry, D.H. (1983).
/art  Multiplicative models in projection pursuit.
/from Technical Report Orion 25, Department of Statistics, Stanford University.

/auth Hermans, J. and Habbema, J.D.F. (1976).
/art  Manual for the ALLOC discriminant analysis programs.
/from Department of Medical Statistics, University of Leiden, Netherlands.

/auth Hills, M. (1967).
/art  Discrimination and allocation with discrete data.
/publ Applied Statistics /vol<16>/pages 237-250.

/auth Hjort, N.L. (1986).
/art  On frequency polygons and average shifted histograms in higher
      dimensions.
/from Technical Report No. 22-LCS, Department of Statistics, Stanford University.

/auth Huber, P.J. (1981).
/publ Robust Statistics.
/from Wiley, New York.

/auth Huber, P.J. (1985).
/art  Projection pursuit (with discussion).
/publ Annals of Statistics /vol<13>/pages 435-525.

/auth Jeffreys, H. (1961).
/publ Theory of Probability.
/from Third edition, Clarendon Press, Oxford.

/auth Jones, M.C. and Lotwick, H.W. (1984).
/art  A remark on ``Kernel density estimation using the fast Fourier
      transform''.
/publ Applied Statistics /vol<33>/page 120.

/auth Kamgar-Parsi, B. and Kanal, L.N. (1985).
/art  An improved branch and bound algorithm for computing $k$-nearest 
      neighbours.
/publ Pattern Recognition Letters /vol<3>/pages 7-12.

/auth Kanal, L. (1974).
/art  Patterns in pattern recognition; 1968--1974.
/publ IEEE Transactions on Information Theory /vol<IT-20>/pages 472-479.

/auth Karpovsky, M.G. (1978).
/publ Finite Orthogonal Series in Design of Digital Devices.
/from Halsted Press, New York.

/auth Kimura, F., Takashina, K., Tsuruoka, S., and Miyake, Y. (1987).
/art  Modified quadratic discriminant functions and the application to 
      Chinese character recognition.
/publ IEEE Transactions on Pattern Analysis and Machine Intelligence 
      /vol<PAMI-9>/pages 149-153.

/auth Koziol, J. (1982).
/art  A class of invariant procedures for assessing multivariate 
      normality. Technical report, Department of Mathematics /& Medicine, 
      UC San Diego, California.

/auth Krishnaiah, P.R. and Kanal, L.N. (editors) (1982).
/publ Handbook of Statistics 2: Classification, Pattern Recognition, and 
      Reduction of Dimensionality.
/from North-Holland, Amsterdam.

/auth Kronmal, R. and Tarter, M. (1968).
/art  Estimation of error rates in discriminant analysis.
/publ Technometrics /vol<10>/pages 1-11.

/auth Lachenbruch, P.A. and Mickey, M.R. (1968).
/art  Estimation of error rates in discriminant analysis.
/publ Technometrics /vol<10>/pages 1-11.

/auth Lehmann, E.L. (1983).
/publ Theory of Point Estimation.
/from Wiley, New York.

/auth Llewellyn, B., Preston, G., Kittler, J. and Harris, J.F. (1982).
/publ Pattern Recognition Methods for Symbols in Cartography.
/from Report 61/sla 82, Oxford University.

/auth Mardia, K.V. (1974).
/art  Applications of some measures of multivariate skewness and 
      kurtosis in testing normality and robustness studies.
/publ Sankhya Series B /vol<36>/pages 115-128.

/auth Mardia, K.V., Kent, J.T. and Bibby, J.M. (1979).
/publ Multivariate Analysis.
/from Academic Press, San Francisco.

/auth Mahalanobis, P.C. (1936).
/art  On the generalized distance in statistics.
/publ Proceedings from the National Institute of Science,
/from India, /vol<12>/pages 49-55.

/auth Martin, D.C. and Bradley, R.A. (1972).
/art  Probability models estimation and classification for multivariate 
      dichotomous populations.
/publ Biometrics /vol<28>/pages 203-222.

/auth Moran, M.A. and Murphy, B.J. (1979).
/art  A closer look at two alternative methods of statistical 
      discrimination.
/publ Applied Statistics /vol<28>/pages 223-232.

/auth Prakasa Rao, B.L.S. (1983).
/publ Nonparametric Functional Estimation.
/from Academic Press, London.

/auth Ranneby, B. (1984).
/art  The maximum spacing method: An estimation method related to the 
      maximum likelihood method.
/publ Scandinavian Journal of Statistics /vol<11>/pages 93-112.

/auth Rao, C.R. (1948).
/art  The utilization of multiple measurements in problems of biological classification.
/publ Journal of the Royal Statistical Society Series B /vol<10>/pages 159-193.

/auth Ripley, B.D. (1986).
/art  Statistics, images, and pattern recognition (with discussion).
/publ Canadian Journal of Statistics /vol<14>/pages 83-111.

/auth Royston, J.P. (1983).
/art  Some techniques for assessing multivariate normality based on the 
      Shapiro-Wilk ``W''.
/publ Applied Statistics /vol<32>/pages 121-133.

/auth Rudemo, M. (1982).
/art  Empirical choice of histograms and kernel density estimators.
/publ Scandinavian Journal of Statistics /vol<9>/pages 65-78.

/auth Ruiz, E.V. (1986).
/art  An algorithm for finding nearest neighbours in (approximately) 
      constant average time.
/publ Pattern Recognition Letters /vol<4>/pages 145-158.

/auth Schwarz, G. (1978).
/art  Estimating the dimension of a model.
/publ Annals of Statistics /vol<6>/pages 461-464.

/auth Scott, D.W. (1985a).
/art  Frequency polygons: Theory and applications.
/publ Journal of the American Statistical Association
      /vol<80>/pages 348-354.

/auth Scott, D.W. (1985b).
/art  Average shifted histograms: Effective nonparametric density 
      estimation in several dimensions.
/publ Annals of Statistics /vol<13>/pages 1024-1040.

/auth Scott, D.W., Thapia, R.A., and Johnson, J.R. (1980).
/art  Nonparametric probability density estimation by discrete maximum 
      penalized-likelihood criteria.
/publ Annals of Statistics /vol<8>/pages 820-832.

/auth Silverman, B.W. (1976).
/art  On a Gaussian process related to multivariate probability density 
      estimation.
/publ Mathematical Proceedings from the Cambridge Philosophical Society 
      /vol<80>/pages 135-144.

/auth Silverman, B.W. (1978).
/art  Choosing window width when estimating a density.
/publ Biometrika /vol<65>/pages 1-11.

/auth Silverman, B.W. (1982).
/art  Kernel density estimation using the fast Fourier transform.
/publ Applied Statistics /vol<31>/pages 93-99.

/auth Sj/"ostr/"om, M. and Wold, S. (1980).
/art  SIMCA: A pattern recognition method based on principal component 
      models.
/book Pattern Recognition in Practice,
/from E.S. Gelsema and L.N. Kanal (editors), North-Holland/pages 351-359.

/auth Titterington, D.M. (1985).
/art  Common structure of smoothing techniques in statistics.
/publ International Statistical Review /vol<53>/pages 141-170.

/auth Tittertington, D.M., Smith, A.F.M. and Makov, U.E. (1985).
/publ Statistical Analysis of Finite Mixture Distributions.
/from Wiley, Chichester.

/auth Toussaint, G.T. (1974).
/art  Bibliography on estimation of misclassification.
/publ IEEE Transactions on Information Theory /vol<IT-20>/pages 472-479.

/auth Van Ryzin, J. (1966).
/art  Bayes risk consistency of classification procedures using density 
      estimation.
/publ Sankhya Series A /vol<28>/pages 261-479.

/auth Wold, S. (1976).
/art  Pattern recognition by means of disjoint principal components 
      models.
/publ Pattern Recognition /vol<8>/pages 127-133.

/auth Wu, C.F.J. (1983).
/art  On the convergence properties of the EM algorithm.
/publ Annals of Statistics /vol<11>/pages 95-103.
/endgroup

/bye